\documentclass[11pt]{article}
\usepackage{amsmath,amsthm,amsfonts,amssymb,amsopn,amscd}     
\usepackage{epsfig}
\usepackage{graphicx,tikz}
\usetikzlibrary{arrows}
\usepackage{setspace}
\usepackage{stmaryrd}
\usepackage{framed}

\newcommand{\tikzmatht}[2][0.5]
{\vcenter{\hbox{\begin{tikzpicture}[scale=#1]\footnotesize #2
				 \end{tikzpicture}}}
}
\usetikzlibrary{arrows,positioning} 
\tikzset{
    >=stealth',
    arrow/.style={
           ->,
           thick,
           shorten <=2pt,
           shorten >=2pt,}
}

\definecolor{shadecolor}{gray}{0.9}
\definecolor{tintedcolor}{gray}{0.8}

\newenvironment{tinted}{%
  \MakeFramed {\FrameRestore}}%
 {\endMakeFramed}%

\renewenvironment{shaded}{%
  \MakeFramed {\FrameRestore}}%
 {\endMakeFramed}%

\newenvironment{graybox}%
       {\fboxsep=12pt\relax
        \begin{shaded}%
        \list{}{\leftmargin=12pt\rightmargin=\leftmargin\leftmargin=12pt\topsep=2pt\relax}%
        \expandafter\item
}%
       {\endlist\end{shaded}}%

\newenvironment{tintedbox}%
       {\fboxsep=12pt\relax
        \begin{tinted}%
        \list{}{\leftmargin=12pt\rightmargin=\leftmargin\leftmargin=12pt\topsep=2pt\relax}%
        \expandafter\item
}%
       {\endlist\end{tinted}}%

\DeclareMathOperator{\Tr}{Tr}
\DeclareMathOperator{\LTr}{LTr}
\DeclareMathOperator{\RTr}{RTr}

\def\qed{{\unskip\nobreak\hfil\penalty50
\hskip2em\hbox{}\nobreak\hfil$\square$
\parfillskip=0pt \finalhyphendemerits=0\par}\medskip}

\def\Ad{{\mathrm{\rm Ad}}}

\def\dim{{\mathrm{dim}}}
\def\End{{\mathrm{End}}}
\def\Hom{{\mathrm{Hom}}}
\def\Mat{{\mathrm{Mat}}}
\def\id{{\rm id}}
\def\eins{{\mathbf{1}}}

\newtheorem{theorem}{Theorem}[section]
\newtheorem{lemma}[theorem]{Lemma}

\newtheorem{corollary}[theorem]{Corollary}

\newtheorem{proposition}[theorem]{Proposition}

\theoremstyle{definition} 
\newtheorem{definition}[theorem]{Definition}
\newtheorem{example}[theorem]{Example}
\newtheorem{remark}[theorem]{Remark}

\def\RR{{\mathbb R}}
\def\CC{{\mathbb C}}
\def\NN{{\mathbb N}}

\def\ZZ{{\mathbb Z}}

\def\A{{\mathcal A}}
\def\B{{\mathcal B}}
\def\C{{\mathcal C}} 
\def\DHR{\C^{\rm DHR}}

\def\H{{\mathcal H}}

\def\sig{\sigma}
\def\LRA{\Leftrightarrow} 
\def\RA{\Rightarrow}
\def\Psit{\Psi_{\tau\otimes\tau}}
\def\Psis{\Psi_{\sig\otimes\sig}}
\def\Ad{\mathrm{Ad}}
\def\wh{\widehat}
\def\eps{\varepsilon} 
\def\inv{^{-1}} 
\def\loc{_{\rm loc}}
\def\can{_{\rm can}}
\def\ol{\overline}

\def\wt{\widetilde}
\def\opp{^{\rm opp}}
\def\mm{{\mathbf{m}}}
\def\mn{{\mathbf{n}}}
\def\BA{{\mathbf{A}}}
\def\BB{{\mathbf{B}}}
\def\BC{{\mathbf{C}}}
\def\BR{{\mathbf{R}}}
\def\scirc{{\, \circ\, }}
\def\ncirc{{}}
\def\RA{\Rightarrow}
\def\ql{_{\rm ql}}

\newcommand{\be}{\begin{equation}} 
\newcommand{\ee}{\end{equation}}
\newcommand{\bea}{\begin{eqnarray}} 
\newcommand{\eea}{\end{eqnarray}}
\newcommand{\ba}{\begin{array}} 
\newcommand{\ea}{\end{array}}
\newcommand{\eref}[1]{Eq.\ \!(\ref{#1})}
\newcommand{\nref}[1]{(\ref{#1})}
\newcommand{\sref}[1]{Sect.~\ref{#1}}
\newcommand{\dref}[1]{Def.~\ref{#1}}
\newcommand{\cref}[1]{Cor.~\ref{#1}}
\newcommand{\rref}[1]{Remark~\ref{#1}}
\newcommand{\xref}[1]{Example~\ref{#1}}
\newcommand{\lref}[1]{Lemma~\ref{#1}}
\newcommand{\pref}[1]{Prop.~\ref{#1}}
\newcommand{\tref}[1]{Thm.~\ref{#1}}
\newcommand{\stapel}[2]{\stackrel{\hbox{\small #1}}{#2}}

\title{\vskip-10mm \LARGE Tensor categories and endomorphisms \\ of von
  Neumann algebras \\[2mm] \Large (with applications to Quantum Field Theory)}

\author{{\sc Marcel Bischoff}\\ \scriptsize Institut f\"ur Theoretische Physik, Universit\"at G\"ottingen,
\\[-1.5mm] \scriptsize Friedrich-Hund-Platz 1, D-37077 G\"ottingen,
Germany \\[-1.5mm] \scriptsize as of Sept.\ 2014: Vanderbilt University, 
Department of Mathematics \\[-1.5mm] \scriptsize
1326 Stevenson Center, 
Nashville, TN 37240-0001, USA
\\[-1.5mm] \scriptsize {\tt marcel.bischoff@vanderbilt.edu}
\\[2mm]
{\sc Yasuyuki Kawahigashi} \\ \scriptsize Department of Mathematical Sciences, The
University of Tokyo, \\[-1.5mm] \scriptsize Komaba, Tokyo 153-8914, Japan; \\ \scriptsize Kavli IPMU
(WPI), The University of Tokyo \\[-1.5mm]
\scriptsize 5-1-5 Kashiwanoha, Kashiwa, 277-8583, Japan \\[-1.5mm] \scriptsize {\tt yasuyuki@ms.u-tokyo.ac.jp}
\\[2mm]
{\sc Roberto Longo} \\
\scriptsize Dipartimento di Matematica,
Universit\`a di Roma ``Tor Vergata'',\\[-1.5mm]
\scriptsize Via della Ricerca Scientifica, 1, I-00133 Roma, Italy \\[-1.5mm] \scriptsize {\tt longo@mat.uniroma2.it}
\\[2mm]
{\sc Karl-Henning Rehren} \\
\scriptsize Institut f\"ur Theoretische Physik, Universit\"at G\"ottingen,
\\[-1.5mm] \scriptsize Friedrich-Hund-Platz 1, D-37077 G\"ottingen, Germany
\\[-1.5mm] \scriptsize {\tt rehren@theorie.physik.uni-goettingen.de}}

\begin{document}

\maketitle

\begin{abstract}
Q-systems describe ``extensions'' of an infinite von Neumann
factor $N$, i.e., finite-index unital inclusions of $N$ into another
von Neumann algebra $M$. They are (special cases of) Frobenius
algebras in the C* tensor category of endomorphisms of $N$. We review 
the relation between Q-systems, their modules and bimodules as 
structures in a category on one side, and homomorphisms 
between von Neumann algebras on the other side. We then elaborate 
basic operations with Q-systems (various decompositions in the 
general case, and the centre, the full centre, and the braided product 
in braided categories), and illuminate their meaning in the von
Neumann algebra setting. The main applications are in local quantum 
field theory, where Q-systems in the subcategory of DHR endomorphisms of 
a local algebra encode extensions $\A(O)\subset\B(O)$ of local
nets. These applications, notably in conformal quantum field theories 
with boundaries, are briefly exposed, and are discussed in more
detail in two separate papers \cite{BKL,BKLR}.  
\end{abstract}

\newpage 

\tableofcontents

\vfill

\noindent
{\footnotesize Supported by the Grants-in-Aid for Scientific
Research, JSPS. Supported by the German Research
Foundation (Deutsche Forschungsgemeinschaft (DFG)) through the
Institutional Strategy of the University of G\"ottingen.  
The hospitality and support of the Erwin Schr\"odinger International
Institute for Mathematical Physics, Vienna, is gratefully
acknowledged. R.L. thanks the Alexander von Humboldt Foundation and
the European Research Council for support.}

\newpage

\section{Introduction}
\setcounter{equation}{0}
Q-systems have first appeared in \cite{L94} as a device to characterize
finite-index subfactors $N\subset M$ of infinite (type $III$) von
Neumann algebras, generalizing the Jones theory of type $II$
subfactors \cite{J,Ko,P}. A Q-system is a triple 
$$\BA=(\theta,w,x),$$ 
where $\theta$ is a unital endomorphism of $N$ and $w\in\Hom(\id_N,\theta)$, 
$x\in\Hom(\theta,\theta^2)$ are a pair of intertwiners whose algebraic 
relations guarantee that $\theta$ is the dual canonical endomorphism 
(\sref{s:subf}) associated with a subfactor $N\subset M$. 

Notice that the data of the Q-system pertain only to $N$, so the
Q-system actually characterizes $M$ as an ``extension'' of $N$.  
In fact, the larger algebra $M$ along with the embedding of $N$ into
$M$ can be explicitly reconstructed (up to isomorphism) from the
data. One issue in this work is a generalization to Q-systems for
extensions $N\subset M$ where $M$ may have a finite centre, i.e., $M$
is a direct sum of infinite factors. 

Subfactors are, apart from their obvious mathematical interest, also
of physical interest since they describe, e.g., the embedding of a 
physical quantum sub-system in larger system. In this context, it is
essential that the algebras are C* or von Neumann algebras, since 
quantum observables are always (selfadjoint) elements of such
algebras, and in relativistic QFT, local observables generate factors
of type $I\!I\!I$. 

From a category point of view, a Q-system is the same as a special C*
Frobenius algebra in a (strict, simple) C* tensor category. 
In the case at hand, the category would be (a subcategory of) the
category $\End_0(N)$ of endomorphisms of $N$ with finite dimension. 
This is actually the most general situation, since every (rigid, countable) 
abstract C* tensor category can be realized as a full subcategory of
$\End_0(N)$ \cite{Y03}. 

In a more general setting (notably without assuming the C* structure
which is naturally present in the case of $\End_0(N)$) abstract tensor
categories and Frobenius algebras have been extensively studied by
many mathematicians \cite{D91,JS,ML,BK,D02}, 
and interesting ``derived'' structures have been discovered and
classified, notably when the underlying tensor category is braided, or
even modular \cite{M00,KO,O,M03,M04,ENO,DMNO}. 

A connection to physics of this more general setting is provided by
\cite{TFT,TFT1} where a formulation of two-dimensional 
({\em Euclidean}) conformal 
quantum field theory on Riemannian surfaces is developped in terms of a 
three-dimensional ``topological quantum field theory'' which is a 
cobordism theory between pairs of Riemannian surfaces. The authors observed,
among a wealth of other results, that the modules and bimodules of the
representation category of the underlying chiral theory play a
prominent role in the classification of one-dimensional boundaries
between Riemannian surfaces.  

From the von Neumann algebra point of view, an important class of
braided tensor subcategories of $\End_0(N)$ naturally arises in the
algebraic formulation of {\em relativistic} Quantum Field Theory (QFT). 
Namely, a distinguished class of positive-energy 
representations of local QFT can be described in terms of
endomorphisms of the C* algebra $\A$ of quasi-local observables. 
These DHR endomorphisms
are the objects of a braided C* tensor category \cite{DHR,FRS}.  
By restricting attention to a von Neumann algebra $N=\A(O)$ of local
observables, one obtains a braided tensor subcategory of $\End_0(N)$.  
In this context, Q-systems describe finite-index extensions 
$\A\subset\B$ of quantum field theories, and $\B$ is local if and 
only if the Q-system is commutative w.r.t.\ the braiding. 

Our main motivation for the present work was the study of boundary
conditions in relativistic conformal QFT in two spacetime dimensions,
as discussed in detail in the compagnon papers \cite{BKL,BKLR}. Boundaries in 
relativistic quantum field theories \cite{LR04,LR09,CKL13}, with
observables that are Hilbert space operators subject to
the principle of locality (or rather causality), have been 
analyzed much less than in the Euclidean setting. Very little is known
about an apriori relation between Euclidean and Lorentzian
boundaries. Yet, our treatment of boundaries in relativistic 
two-dimensional conformal QFT shows that precisely the same
mathematical structures, namely the chiral representation category, 
its Q-systems and their modules and
bimodules, control the boundary conditions in both situations. We
address in particular the case of ``hard'' boundaries in \cite{BKL}
and ``transparent'' or ``phase boundaries'' (defects) in \cite{BKLR}.
In this work, we shall concentrate on the underlying mathematical
theory, with only scattered remarks about the relevance in QFT. 
A brief exposition of these physical applications will be given in 
\sref{s:QFT}. 

While large portions of the category side of this work are 
reformulations from \cite{FFRS06,KR08}, our original contribution is 
the elaboration of the relation between the abstract category notions
and the von Neumann algebra setting and subfactor theory. A prominent 
issue is our proof of \tref{t:centchar} (a characterization of the 
central projections of an extension $N\subset M$, which is given by
the braided product of two full-centre Q-systems in a modular
category). This theorem is implicitly present, but widely scattered 
in the work of \cite{FFRS04,FFRS06,FFRS07,TFT1,TFT,KR08}. Our proof 
is much more streamlined, because it benefits from substantial
simplifications in the C* setting, where one can exploit positivity
arguments in crucial steps.  

This theorem is relevant for phase boundaries in relativistic 
two-dimensional conformal QFT because it classifies the boundary
conditions in terms of chiral data \cite{BKLR}, very much the same as
in the Euclidean setting \cite{TFT1}. 

Other original contributions in this paper concern Q-systems for
extensions $N\subset M$ when $M$ is not a factor, a situation that
naturally occurs in several applications, as well as 
the characterization of various types of decompositions of
Q-systems (\sref{s:Qcentral}--\ref{s:interm}) in terms of algebraic 
properties of projections in $\Hom(\theta,\theta)$.

In \sref{s:homo}, we review the basic notions concerning endomorphisms
and homomorphisms of infinite von Neumann algebras, with special
emphasis on the notions of conjugates and dimension. 

\sref{s:frob} is devoted to the category structure, and to
the correspondences between Q-systems and algebra extensions, and
between bimodules between Q-systems and homomorphisms between the
corresponding extensions. 

\sref{s:Qsystem} is the main part of the paper. We introduce various
operations with Q-systems (decompositions, braided products, centres
and full centre), and investigate their meaning in the setting of von
Neumann algebras. 

\sref{s:QFT} contains an exposition of the appearance of braided and
modular C* tensor categories in the DHR theory of superselection
sectors in Algebraic QFT, and reviews the relevance of Q-systems for
issues like extensions and boundary conditions.

\section{Homomorphisms of von Neumann algebras}
\setcounter{equation}{0}
\label{s:homo}

Let $N$ and $M$ be two von Neumann algebras, and $\alpha$, $\beta$ a
pair of homomorphisms $:N\to M$. (Without further mentioning, the notion 
``homomorphism'' will include the * and unit-preserving properties
$\alpha(n^*)=\alpha(n)^*$ and $\alpha(\eins_N)=\eins_M$.)  
An operator $t\in M$ such that
$$t\cdot\alpha(n)=\beta(n)\cdot t\qquad\mathrm{for\;all}\quad n\in N$$
is called an {\bf intertwiner}, writing $t:\alpha\to\beta$ or $t\in
\Hom(\alpha,\beta)$. Clearly, if $t\in\Hom(\alpha,\beta)$, then
$t^*\in\Hom(\beta,\alpha)$; $\Hom(\alpha,\beta)$ is a complex vector
space, and $\Hom(\alpha,\alpha)$ is a C*-algebra.   

A homomorphism $\alpha: N \to M$ is composed with a homomorphism
$\beta:M\to L$, such that $\beta\circ\alpha:N\to L$. 

Likewise, for any three homomorphisms $\alpha,\beta,\gamma: N \to M$ and
intertwiners $t\in \Hom(\alpha,\beta)$ and $s\in \Hom(\beta,\gamma)$,
the product in $M$ gives an intertwiner $s\cdot t\in \Hom(\alpha,\gamma)$. 

These structures turn the endomorphisms of a von Neumann algebra $N$ into
a strict tensor category $\End(N)$, and the homomorphisms between von
Neumann algebras $N,M,\dots$ into a strict tensor 2-category, where
the concatenation of morphisms is the product of intertwiners:
$s\scirc t:=s\cdot t$, the monoidal product of objects is the composition of 
endomorphisms: $\beta\times\alpha:=\beta\circ\alpha$, and the monoidal 
product of morphisms $t_i:\alpha_i\to\beta_i$ is the product 
$$t_1\times t_2 = t_1\cdot\alpha_1(t_2) = \beta_1(t_2)\cdot t_1:\qquad
\tikzmatht{
  \fill[black!10] (.8,-1.5) rectangle (1.8,2.5); 
  \fill[black!15] (-.8,-1.5) rectangle (.8,2.5); 
  \fill[black!20] (-1.8,-1.5) rectangle (-.8,2.5); 
    \draw[thick] (-.8,2.5)--(-.8,-1.5);           
    \draw[thick] (.8,2.5)--(.8,-1.5); 
  \fill[white](-1.3,0) rectangle (-0.3,1); 
    \draw (-1.3,0) rectangle node{$t_1$} (-0.3,1); 
  \fill[white](0.3,0) rectangle (1.3,1); 
    \draw (0.3,0) rectangle node{$t_2$} (1.3,1); 
    \node at (-.7,2) [left] {$\beta_1$};
    \node at (-.7,-1.1) [left] {$\alpha_1$};
    \node at (.8,2) [right] {$\beta_2$};
    \node at (.8,-1.1) [right] {$\alpha_2$};
} = 
\tikzmatht{
  \fill[black!10] (.8,-1.5) rectangle (1.8,2.5); 
  \fill[black!15] (-.8,-1.5) rectangle (.8,2.5); 
  \fill[black!20] (-1.8,-1.5) rectangle (-.8,2.5); 
    \draw[thick] (-.8,2.5)--(-.8,-1.5); 
    \draw[thick] (.8,2.5)--(.8,-1.5); 
  \fill[white](-1.3,.6) rectangle (-0.3,1.6); 
    \draw (-1.3,.6) rectangle node{$t_1$} (-0.3,1.6); 
  \fill[white](0.3,-.6) rectangle (1.3,.4); 
    \draw (0.3,-.6) rectangle node{$t_2$} (1.3,.4); 
    \node at (-.7,2) [left] {$\beta_1$};
    \node at (-.7,-1.1) [left] {$\alpha_1$};
    \node at (.8,2) [right] {$\beta_2$};
    \node at (.8,-1.1) [right] {$\alpha_2$};
} = 
\tikzmatht{
  \fill[black!10] (.8,-1.5) rectangle (1.8,2.5); 
  \fill[black!15] (-.8,-1.5) rectangle (.8,2.5); 
  \fill[black!20] (-1.8,-1.5) rectangle (-.8,2.5); 
    \draw[thick] (-.8,2.5)--(-.8,-1.5); 
    \draw[thick] (.8,2.5)--(.8,-1.5);
  \fill[white](-1.3,-.6) rectangle (-0.3,.4);
   \draw (-1.3,-.6) rectangle node{$t_1$} (-0.3,.4);
  \fill[white](0.3,.6) rectangle (1.3,1.6); 
    \draw (0.3,.6) rectangle node{$t_2$} (1.3,1.6);
    \node at (-.7,2) [left] {$\beta_1$};
    \node at (-.7,-1.1) [left] {$\alpha_1$};
    \node at (.8,2) [right] {$\beta_2$};
    \node at (.8,-1.1) [right] {$\alpha_2$};
}
$$
(This graphical notation, directly appealing to the underlying tensor
category point of view, will render the structure of many algebraic
computations more transparent. Its basic rules are self-explaining from
this example: Different shades indicate different von Neumann algebras, 
and we usually reserve the lightest shade for $N$, lines are
homomorphisms, boxes and similar symbols to appear later are
intertwiners, the monoidal product is horizontal juxtaposition, 
and the concatenation product is read from the bottom to the top. 
The operator adjoint is represented by up-down reflection.)  

Notice that {\em as operators}, $t\times 1_\alpha = t$ is the same operator 
in a different intertwiner space, whereas $1_\alpha\times t=\alpha(t)$. 
To enhance readability, we shall occasionally suppress
the concatenation symbol and write simply $s\scirc t$ as the operator
product $st$. 

Because all intertwiner spaces $\Hom(\alpha,\beta)$ are linear subspaces of
the target von Neumann algebra, they inherit its weak and norm topologies.
In particular, $\End(N)$ is a C* tensor category, and the
self-intertwiners $\Hom(\alpha,\alpha)$ form a C* algebra. Important
consequences are that $t^*\scirc t \equiv t^* t$ is a positive operator in
$\Hom(\beta,\beta)$, and that $t^*\scirc t=0$ implies $t=0$. 

\subsection{Endomorphisms of infinite factors}
\label{s:endo}
\setcounter{equation}{0}
A von Neumann algebra $N$ is a {\bf factor} iff its centre $N'\cap
N\equiv\Hom(\id_N,\id_N) =\CC\cdot \eins_N$. Since $\id_N$ is the monoidal
unit in the tensor category, this is the same as saying that the
category $\End(N)$ is simple. 

These elementary facts can be supplemented by further structure. 
If $u:\alpha\to\beta$ is unitary, $\alpha$ and $\beta$ are said to be
{\bf unitarily equivalent}. The unitary equivalence class of $\alpha$ is
called the {\bf sector} $[\alpha]$. An endomorphism $\alpha$ is
{\bf irreducible} iff 
$\Hom(\alpha,\alpha)=\CC\cdot \eins_N$. 

In an {\bf infinite} ($\LRA$ purely infinite, type $I\!I\!I$) von Neumann 
factor acting on a separable Hilbert space (which we shall henceforth
assume throughout), every projection $e\neq 0$ can be written as
$e=ss^*$ where $s^*s=1$, and one can always choose decompositions of
the unit $1=\sum_i s_is_i{}^*$ such that $s_i{}^*s_j = \delta_{ij}$. 
The algebra generated by bounded quantum mechanical observables (= the
algebra $\B(\H)$ of all bounded operators) does not share this
property; instead, the local algebras of quantum field theory are
generically infinite von Neumann factors.   

Thanks to this property, one can define  

(i) an inclusion relation for endomorphisms: $\beta\prec\alpha$ iff there is
$s:\beta\to\alpha$ with $s^*\ncirc s=1_{\beta}$. 

(ii) subobjects: if $e:\alpha\to\alpha$ is a projection, then 
there is a sub-endomorphism $\alpha_s$ defined by the choice of $s$
such that $s\ncirc s^*=e$, $s^*\ncirc s=1$, and putting 
$$\alpha_s(\cdot) = s^*\alpha(\cdot)s:\qquad
\tikzmatht{
  \fill[black!10] (-1.2,-2) rectangle (1.2,2);
    \draw[thick] (0,-2)--(0,2);  
    \node at (-.5,-1.6) {$\alpha_s$}; \node at (-.5,1.6) {$\alpha_s$};
  \fill[white] (0,1.2)--(-0.3,.7)--
          (0.3,.7)--(0,1.2); 
    \draw[thick] (0,1.2) [right] node{$s^*$}--(-0.3,.7)--
          (0.3,.7)--(0,1.2); 
  \fill[white] (0,-1.2)--(-0.3,-.7)--
          (0.3,-.7)--(0,-1.2); 
    \draw[thick] (0,-1.2) [right] node{$s$}--(-0.3,-.7)--
          (0.3,-.7)--(0,-1.2); 
    \node at (0,0) [right] {$\alpha$};
}.
$$
We refer to $\alpha_s\prec\alpha$ as the {\bf range} of $e$. We shall 
sometimes write $\alpha_e$ instead, in order to emphasize that the
unitary equivalence class of $\alpha_s$ does not depend on the choice of $s$. 
(Categories where subobjects exist are also called ``Karoubian'', thus
$\End(N)$ is Karoubian if $N$ is an infinite factor.)

(iii) direct sums of endomorphisms: 
$$\alpha(\cdot):=\sum_i s_i\alpha_i(\cdot)s_i{}^*:\qquad\sum_i
\tikzmatht{
  \fill[black!10] (-1.3,-2) rectangle (1.3,2);
    \draw[thick] (0,-2)--(0,2);  
    \node at (-.5,-1.7) {$\alpha$}; \node at (-.5,1.7) {$\alpha$};
  \fill[white] (0,.7)--(-0.3,1.2)--(0.3,1.2)--(0,.7); 
    \draw[thick] (0,.7)--(-0.3,1.2)--(0.3,1.2) [right]
          node{$s_i$}--(0,.7); 
  \fill[white] (0,-.7)--(-0.3,-1.2)--(0.3,-1.2)--(0,-.7); 
    \draw[thick] (0,-.7)--(-0.3,-1.2)--(0.3,-1.2) [right]
          node{$s_i^*$}--(0,-.7); 
    \node at (0,0) [right] {$\alpha_i$};
}
$$
is an endomorphism, $\alpha_i\prec\alpha$. Suppressing the dependence
on the isometries $s_i$, we write sloppily $\alpha\simeq\bigoplus_i\alpha_i$. 
Since the choice of the isometries $s_i$ is irrelevant for the
unitary equivalence class (sector) $[\alpha]$, the direct sum should
be understood as a direct sum of sectors. We emphasize this by writing also 
$$[\alpha]=\bigoplus\nolimits_i[\alpha_i].$$ 

\subsection{Homomorphisms and subfactors}
\label{s:subf}
\setcounter{equation}{0}
All notions of the preceding presentation can be transferred to homomorphisms 
$\varphi:N\to M$ where both $N$ and $M$ are infinite factors. 
Notice that intertwiners $t\in\Hom(\varphi_1,\varphi_2)$ are elements of $M$. 

Admitting several factors, one obtains a 2-category,
whose objects are the factors, the 1-morphisms are the
homomorphisms, and the 2-morphisms are their intertwiners. 

\medskip

If $N\subset M$ is a {\bf subfactor} (i.e., both $N$ and $M$ are
factors), then the identical map $\iota:N\to M$, $n\mapsto n$, 
is a nontrivial homomorphism, that describes the embedding of $N$ into
$M$. 

One can define \cite[Chap.\ 3]{LRo} a {\bf dimension} function on the
homomorphisms $N\to M$ when both $N$ and $M$ are infinite
factors, which is additive under direct sums and multiplicative under
composition. It is defined through the notion of {\bf conjugates}:
$\alpha:N\to M$ and $\ol\alpha:M\to N$ are said to be conjugates of
each other whenever there is a pair of intertwiners $N\ni w:\id_N\to
\ol\alpha \alpha$ and $M\ni \ol w:\id_M\to \alpha\ol\alpha$ satisfying
the {\bf conjugacy relations}  
\bea\label{conj}
\notag
(w^*\times 1_{\ol\alpha})\scirc(1_{\ol\alpha}\times \ol w) =
1_{\ol\alpha}: \qquad
\tikzmatht{
  \fill[black!18] (-1.7,-1.8) rectangle (1.7,1.8);
  \fill[black!10] (1,1.8)--(1,0) arc (360:180:.5) 
          arc(0:180:.5)--(-1,-1.8)--(-1.7,-1.8)--(-1.7,1.8);  
    \draw[thick] (1,1.8)--(1,0) arc(360:180:.5)--(0,0)
          arc(0:180:.5)--(-1,-1.8);  
    \node at (-0.5,0.4) [above] {$w^*$}; 
    \node at (0.5,-0.4) [below] {$\ol w$}; 
    \node at (.9,1.5) [right] {$\ol\alpha$};
}= 
\tikzmatht{
  \fill[black!18] (-.8,-1.8) rectangle (.8,1.8);
  \fill[black!10] (-.8,-1.8) rectangle (0,1.8);
    \draw[thick] (0,1.8)--(0,-1.8); 
    \node at (-.1,1.5) [right] {$\ol\alpha$};
}, \\ 
(1_\alpha\times w^*) \scirc (\ol w\times 1_\alpha) = 1_\alpha: \qquad
\tikzmatht{
  \fill[black!18] (-1.7,-1.8) rectangle (1.7,1.8);
  \fill[black!10] (-1,1.8)--(-1,0) arc (180:360:.5)
          arc(180:0:.5)--(1,-1.8)--(1.7,-1.8)--(1.7,1.8); 
    \draw[thick] (-1,1.8)--(-1,0) arc(180:360:.5)--(0,0)
          arc(180:0:.5)--(1,-1.8);  
    \node at (0.5,0.4) [above] {$w^*$}; \node at (-0.5,-0.4)
          [below] {$\ol w$}; 
    \node at (-1,1.5) [right] {$\alpha$};

}=
\tikzmatht{
  \fill[black!18] (-.8,-1.8) rectangle (.8,1.8);
  \fill[black!10] (0,-1.8) rectangle (.8,1.8);
    \draw[thick] (0,1.8)--(0,-1.8); 
    \node at (0,1.5) [right] {$\alpha$};
}.
\eea

Being self-intertwiners of $\id_N$, resp.\ $\id_M$, $w^*\ncirc w=d\cdot
\eins_N$ and $\ol w^* \ncirc\ol w=d'\cdot \eins_M$ are positive
scalars, and $w,\ol w$ can be normalized such that $d=d'$. 
The {\bf dimension} $\dim(\alpha)=\dim(\ol\alpha)$ is defined to be 
\be\label{inf}
\dim(\alpha)=\dim(\ol\alpha) := \inf_{(w,\ol w)} d
\ee
where the infimum is taken over all solutions $(w,\ol w)$ of the 
conjugacy relations \eref{conj} with $d=d'$. 
A solution saturating the infimum is 
called {\bf standard solution} or {\bf standard pair}. If $\alpha$ and 
$\beta$ are irreducible, every solution with $d=d'$ is standard, 
because $\dim\Hom(\id,\alpha\ol\alpha)=\dim\Hom(\ol\alpha\alpha)$=1. 
In the general case, standard solutions always exists,
and are unique up to unitary equivalence \cite{KL,LRo}. 

(Here is a simple explicit proof: For $[\alpha]=\bigoplus_i n_i[\alpha_i]$ and
$[\ol\alpha]=\bigoplus_i \ol n_i[\ol\alpha_i]$ with $\alpha_i,\ol\alpha_i$ 
irreducible, one may choose standard pairs $(w_i,\ol w_i)$ for 
$\alpha_i,\ol\alpha_i$ and orthonormal bases
$s^i_a\in\Hom(\alpha_i,\alpha)$, $\ol s^i_b\in\Hom(\ol\alpha_i,\ol\alpha)$.
Then the most general element of $\Hom(\id,\ol\alpha,\alpha)$ is of
the form  
$w=\sum_i\sum_{ab} c^i_{ab}\ol\alpha(s^i_a)\ol s^i_b w_i$, and similarly 
$\ol w=\sum_i\sum_{ab} c'^i_{ab} \alpha(\ol s^i_b)s^i_a \ol w_i$. These
solve the conjugacy relations iff the coefficient matrices satisfy 
$c'^i = (c^i)\inv{}^*$ (in particular, the multiplicities $\ol
n_i=n_i$ must be the same), and one has $d=\sum_i\dim(\alpha_i)\Tr(c^i)^*c^i$,
$d'=\sum_i\dim(\ol\alpha_i)\Tr(c^i)\inv{}^*(c^i)\inv$. The variational problem 
$d[c]d'[c]$ $\stackrel!=\min$ with $d=d'$ is solved by any family of unitary
matrices $c^i$.)

The conjugate of an endomorphism is unique up to unitary equivalence. 
Endomorphisms which do not have conjugates can be assigned the
dimension $\infty$.  

The dimension is always $\geq 1$, and a homomorphism $\alpha$ 
is an isomorphism iff $\dim(\alpha)=1$. In this case, 
$\alpha\inv$ is a conjugate of $\alpha$. More generally, the
dimension is the square root of the (minimal) index \cite{L89,L94}:
$$\dim(\alpha)^2=[M:\alpha(N)].$$

In particular, for a subfactor $N\subset M$, $\dim(\iota)$ is the
square root of the index $[M:N]$ \cite{J}. In this case,
$\iota\ol\iota\in\End(M)$ is called the {\bf canonical endomorphism}, 
and $\ol\iota\iota\in\End(N)$ the {\bf dual canonical endomorphism}.

\medskip

\begin{lemma} \cite{LRo} \label{l:addmult} 
{\rm (i)} Let $(w_1,\ol w_1)$ and $(w_2,\ol w_2)$ be standard pairs
for $(\alpha_1,\ol\alpha_1)$ and for $(\alpha_2,\ol\alpha_2)$,
respectively. Then 
$$w=\ol\alpha_1(w_2)w_1,\quad \ol w=\alpha_2(\ol w_1)\ol w_2$$ 
is a standard pair for $(\alpha_2\alpha_1,\ol\alpha_1\ol\alpha_2)$. \\
{\rm (ii)} Let $(w_i,\ol w_i)$ be standard pairs for 
$(\alpha_i,\ol\alpha_i)$, and $[\alpha]=\bigoplus_i[\alpha_i]$, 
$[\ol\alpha]=\bigoplus_i[\ol\alpha_i]$. Choose orthonormal isometries 
$s_i\in\Hom(\alpha_i,\alpha)$ and $\ol
s_i\in\Hom(\ol\alpha_i,\ol\alpha)$. Then 
$$w = \sum_i(\ol s_i\times s_i)\scirc w_i,\quad \ol w =
\sum_i(s_i\times \ol s_i)\scirc \ol w_i$$
is a standard pair for $(\alpha,\ol\alpha)$. 
\end{lemma}

\begin{corollary} \cite{LRo} \label{c:addmult}
The conjugate respects direct sums, and the
dimension is additive and multiplicative: 
$$\dim(\alpha_2\circ\alpha_1) = \dim(\alpha_2)\cdot\dim(\alpha_1),
\quad \dim(\alpha) = \sum_i\dim(\alpha_i)\quad\hbox{if}\quad
[\alpha]=\bigoplus_i[\alpha_i].$$ 
\end{corollary}

It should, of course, be emphasized that all the notions of direct sums,
subobjects, conjugates and dimension respect the unitary equivalence
relation.  

\medskip

\begin{definition}\label{d:traces}
The left and right {\bf traces} are the faithful positive maps 
\bea\label{traces}
\LTr_\alpha: \Hom(\alpha\beta,\alpha\beta')\to \Hom(\beta,\beta'),
\qquad\qquad\qquad \notag \\[-2mm] 
t\mapsto (w^*\times 1_{\beta})\scirc(1_{\ol\alpha}\times
t)\scirc (w\times 1_{\beta}) \qquad 
\tikzmatht{
  \fill[black!10] (-1.5,-1.5) rectangle (1.5,1.5);
  \fill[white] (-.5,-.5) rectangle (1,.5); 
    \draw[thick] (-.5,-.5) rectangle node{$t$} (1,.5); 
    \draw[thick] (0,.5)--(0,.8) arc(0:180:.5)--(-1,-.8)
          arc(180:360:.5)--(0,-.5) (.5,.5)--(.5,1.5) (.5,-.5)--(.5,-1.5); 
}\\
\RTr_\alpha: \Hom(\beta\alpha,\beta'\alpha)\to \Hom(\beta,\beta'),
\qquad\qquad\qquad \notag \\[-2mm] t\mapsto (1_{\beta}\times \ol w^*)\scirc(t\times
1_{\ol\alpha})\scirc (1_{\beta}\times \ol w) \qquad 
\tikzmatht{
  \fill[black!10] (-1.5,-1.5) rectangle (1.5,1.5);
  \fill[white] (-1,-.5) rectangle (.5,.5); 
    \draw[thick] (-1,-.5) rectangle node{$t$} (.5,.5); 
    \draw[thick] (0,.5)--(0,.8) arc(180:0:.5)--(1,-.8)
          arc(360:180:.5)--(0,-.5) (-.5,.5)--(-.5,1.5) (-.5,-.5)--(-.5,-1.5); 
}
\eea
with any standard solution $(w,\ol w)$ for the conjugate homomorphisms
$(\alpha,\ol\alpha)$.
\end{definition}

\begin{proposition} \label{p:traces} \cite[Lemma 3.7]{LRo}
Let $N$ and $M$ be infinite factors, and let the traces $\LTr_\alpha$ and
$\RTr_\alpha$ be defined w.r.t.\ a standard solution $(w,\ol w)$ of
the conjugacy relations for $\alpha:N\to M$ and $\ol\alpha:M\to
N$. 
 
The traces do not depend on the choice of the conjugate and of the
standard solution, and satisfy the trace property 
\bea \label{traceprop} 
\LTr_\alpha (s\times 1_{\beta'})\scirc t = \LTr_{\alpha'} t \scirc
(s\times 1_{\beta})\quad 
\tikzmatht{
  \fill[black!10] (-1.5,-1.5) rectangle (1.5,1.5);
  \fill[white] (-.5,-.8) rectangle (1,0); 
    \draw[thick] (-.5,-.8) rectangle node {$t$} (1,0); 
  \fill[white] (-.3,.3) rectangle (.3,.8); 
    \draw[thick] (-.3,.3) rectangle node {$s$} (.3,.8); 
    \draw[thick] (0,0)--(0,.3) (0,.8) arc(0:180:.5)--(-1,-.8)
    arc(180:360:.5)--(0,-.8) (.5,0)--(.5,1.5) (.5,-.8)--(.5,-1.5); 
} =
\tikzmatht{
  \fill[black!10] (-1.5,-1.5) rectangle (1.5,1.5);
  \fill[white] (-.5,.8) rectangle (1,0); 
    \draw[thick] (-.5,.8) rectangle node {$t$} (1,0); 
  \fill[white] (-.3,-.3) rectangle (.3,-.8); 
    \draw[thick] (-.3,-.3) rectangle node {$s$} (.3,-.8); 
    \draw[thick] (0,0)--(0,-.3) (0,-.8)
          arc(360:180:.5)--(-1,.8)
          arc(180:0:.5)--(0,.8) (.5,0)--(.5,-1.5) (.5,.8)--(.5,1.5); 
}, \nonumber \\
\RTr_\alpha (1_{\beta'}\times s)\scirc t = \RTr_{\alpha'} t \scirc
(1_{\beta}\times s) \quad 
\tikzmatht{
  \fill[black!10] (-1.5,-1.5) rectangle (1.5,1.5);
  \fill[white] (.5,-.8) rectangle (-1,0); 
    \draw[thick] (.5,-.8) rectangle node {$t$} (-1,0); 
  \fill[white] (.3,.3) rectangle (-.3,.8); 
    \draw[thick] (.3,.3) rectangle node {$s$} (-.3,.8); 
    \draw[thick] (0,0)--(0,.3) (0,.8) arc(180:0:.5)--(1,-.8)
    arc(360:180:.5)--(0,-.8) (-.5,0)--(-.5,1.5) (-.5,-.8)--(-.5,-1.5);
} =
\tikzmatht{
  \fill[black!10] (-1.5,-1.5) rectangle (1.5,1.5);
  \fill[white] (.5,.8) rectangle (-1,0); 
    \draw[thick] (.5,.8) rectangle node {$t$} (-1,0); 
  \fill[white] (.3,-.3) rectangle (-.3,-.8); 
    \draw[thick] (.3,-.3) rectangle node {$s$} (-.3,-.8); 
    \draw[thick] (0,0)--(0,-.3) (0,-.8) arc(180:360:.5)--(1,.8)
          arc(0:180:.5)--(0,.8) (-.5,0)--(-.5,-1.5) (-.5,.8)--(-.5,1.5); 
}.
\eea
for $s\in\Hom(\alpha',\alpha)$ and $t\in \Hom(\alpha\beta,\alpha'\beta')$ 
resp.\ $t\in \Hom(\beta\alpha,\beta'\alpha')$. For $\beta=\beta'=\id$, 
both traces coincide and are denoted 
$\Tr_\alpha:\Hom(\alpha,\alpha)\to\CC$. In particular, 
\be\label{Tr=dim}
\Tr_\alpha 1_\alpha = \dim(\alpha).
\ee
\end{proposition}

The latter property can in fact be adopted as an alternative
definition for standardness, since one also has 

\begin{proposition} \label{p:leftright} \cite[Lemma 3.9]{LRo}
Let $N$ and $M$ be infinite factors, and let the traces $\LTr_\alpha$ and
$\RTr_\alpha$ be defined as in \dref{d:traces} w.r.t.\ {\em any} (i.e., not 
necessarily standard) solution $(w,\ol w)$ of the conjugacy relations for 
$\alpha:N\to M$ and $\ol\alpha:M\to N$. Then $\LTr_\alpha$ and $\RTr_\alpha$ 
coincide if and only if $(w,\ol w)$ is standard. 
\end{proposition}

If $(w,\ol w)$ is not standard, the maps $\LTr_\alpha$ and $\RTr_\alpha$ on
$\Hom(\alpha,\alpha)\to\CC$ may happen to be traces, without being
equal. E.g., for reducible $\alpha$ every $n\in\Hom(\alpha,\alpha)$ 
gives rise to a deformation $w':= (1_{\ol\alpha}\times
n)\circ w$, $\ol w' := (n^*{}\inv\times 1_{\ol\alpha})\circ \ol w$ of
a standard pair $(w,\ol w)$, which still solves the conjugacy
relations. Then $\LTr'_\alpha$ and $\RTr'_\alpha$ defined with
$(w',\ol w')$ are traces if and only if $n^*\ncirc n$ is central in
$\Hom(\alpha,\alpha)$, while $(w',\ol w')$ is standard iff $n^*\ncirc
n=1_\alpha$. One has the following characterization \cite[Lemma 2.3]{LRo}: 

\begin{proposition} \label{p:dotstar}
Let $(w,\ol w)$ and $(w',\ol w')$ be solutions of the conjugacy
relations for $\alpha,\ol\alpha$ and for $\alpha',\ol\alpha'$, not
necessarily standard. Define $\LTr_\alpha$ as in \dref{d:traces}
w.r.t.\ these pairs. The following are equivalent: \\
{\rm (i)} For $t\in\Hom(\alpha,\alpha')$ and
$s\in\Hom(\alpha',\alpha)$, one has $\LTr_\alpha(st) =
\LTr_{\alpha'}(ts)$. \\
{\rm (ii)} For $t\in\Hom(\alpha,\alpha')$, one has 
$$\tikzmatht{
  \fill[black!10] (-1.5,-1.8) rectangle (1.5,.7);
    \draw[thick] (0,-.3) arc(180:0:.4); 
  \fill[white] (-.3,-.3) rectangle (.3,-.8); 
    \draw[thick] (-.3,-.3) rectangle node {$t$} (.3,-.8); 
    \draw[thick] (0,-.8) arc(360:180:.4)--(-.8,.7)
    (.8,-.3)--(.8,-1.8); 
    \node at (1.1,.3) {$\ol w'{}^*$};
    \node at (-1,-1.4) {$w$};
} =
\tikzmatht{
  \fill[black!10] (-1.5,-.7) rectangle (1.5,1.8);
    \draw[thick] (0,.3) arc(180:360:.4); 
  \fill[white] (-.3,.3) rectangle (.3,.8); 
    \draw[thick] (-.3,.3) rectangle node {$t$} (.3,.8); 
    \draw[thick] (0,.8) arc(0:180:.4)--(-.8,-.7)
    (.8,.3)--(.8,1.8); 
    \node at (1.1,-.3) {$\ol w$};
    \node at (-1,1.4) {$w'{}^*$};
}
\in \Hom(\ol\alpha,\ol\alpha').$$
The same is true, replacing $\LTr$ by $\RTr$ in {\rm (i)}, or
replacing $t$ by $s\in\Hom(\alpha',\alpha)$ in {\rm (ii)}. 

In particular, {\rm (ii)} holds if $(w,\ol w)$ and $(w',\ol w')$ are standard.
\end{proposition}

{\em Proof:} ``(i) $\RA$ (ii)'' is the statement of \cite[Lemma
2.3(c)]{LRo}, although the authors actually prove also the
converse. The proof proceeds by noting that 
$$\LTr_\alpha(st) = 
\tikzmatht{
  \fill[black!10] (-1,-1.5) rectangle (1.3,1.5);
    \draw[thick] (.5,.8) arc(0:180:.5)--(-.5,-.8)
    arc(180:360:.5)--(.5,.8);
  \fill[white] (.2,.2) rectangle (.8,.7);
    \draw[thick] (.2,.7) rectangle (.8,.2);
  \fill[white] (.2,-.2) rectangle (.8,-.7);
    \draw[thick] (.2,-.7) rectangle (.8,-.2);
    \node at (1.1,.4) {$s$};
    \node at (1.1,-.4) {$t$};
}
= \tikzmatht{
  \fill[black!10] (-2,-1.5) rectangle (1.6,1.5);
    \draw[thick] (-1.6,-.5) arc(180:360:.4)--(-.8,.2)
    arc(180:0:.4)--(0,-.9) arc(180:360:.4)--(.8,.2)
    arc(0:180:1.2)--(-1.6,-.5); 
  \fill[white] (.5,-.3) rectangle (1.1,-.8);
    \draw[thick] (.5,-.8) rectangle (1.1,-.3);
  \fill[white] (-.5,.1) rectangle (-1.1,-.4);
    \draw[thick] (-.5,-.4) rectangle (-1.1,.1);
    \node at (1.4,-.5) {$s$};
    \node at (-.2,-.2) {$t$};
},\qquad \RTr_{\alpha'}(ts) = 
\tikzmatht{
  \fill[black!10] (-1,-1.5) rectangle (1.3,1.5);
    \draw[thick] (.5,.8) arc(0:180:.5)--(-.5,-.8)
    arc(180:360:.5)--(.5,.8);
  \fill[white] (.2,.2) rectangle (.8,.7);
    \draw[thick] (.2,.7) rectangle (.8,.2);
  \fill[white] (.2,-.2) rectangle (.8,-.7);
    \draw[thick] (.2,-.7) rectangle (.8,-.2);
    \node at (1.1,.4) {$t$};
    \node at (1.1,-.4) {$s$};
}
= \tikzmatht{
  \fill[black!10] (-2,-1.5) rectangle (1.6,1.5);
    \draw[thick] (-1.6,.5) arc(180:0:.4)--(-.8,-.2)
    arc(180:360:.4)--(0,.9) arc(180:0:.4)--(.8,-.2)
    arc(360:180:1.2)--(-1.6,.5); 
  \fill[white] (.5,.3) rectangle (1.1,.8);
    \draw[thick] (.5,.8) rectangle (1.1,.3);
  \fill[white] (-.5,-.1) rectangle (-1.1,.4);
    \draw[thick] (-.5,.4) rectangle (-1.1,-.1);
    \node at (1.4,.5) {$s$};
    \node at (-.2,.2) {$t$};
}.$$
Now, (ii) trivially implies equality of the two expressions, hence
(i). Conversely, (i) implies (ii) because $(1_{\ol\alpha'}\times s)\scirc w'$ is an arbitrary element
of $\Hom(\id,\ol\alpha\alpha')$. 

The variants of the statement follow by obvious modifications. 

Finally, if $(w,\ol w)$ and $(w',\ol w')$ are standard, then \pref{p:traces}
implies (i), hence (ii). 
\qed

\medskip

For a single infinite von Neumann factor $N$, $\End_0(N)$ is the
full subcategory of $\End(N)$, whose objects are the 
endomorphisms of finite dimension. This is a ``rigid'' category since
left and right duals exist for all objects (namely, the conjugate). 

All intertwiner spaces 
$\Hom(\alpha,\beta)$ in $\End_0(N)$ are finite-dimensional, 
and $\Hom(\alpha,\alpha)$ are isomorphic with a direct
sum of matrix algebras $\bigoplus_\lambda \Mat_\CC(n_\lambda)$, where
$\lambda$ are the equivalence classes of irreducible sub-endomorphisms
of $\alpha$ and $n_\lambda$ their multiplicities in $\alpha$. 

Whenever $\alpha$ has finite dimension (and hence a conjugate 
$\ol\alpha$ exists), one can use a standard solution $(w,\ol w)$ to
define linear bijections (left and right {\bf Frobenius conjugations}) 
between the spaces $\Hom(\gamma_2,\alpha\gamma_1)$ and
$\Hom(\ol\alpha\gamma_2,\gamma_1)$, and between
$\Hom(\gamma_2,\gamma_1\alpha)$ and
$\Hom(\gamma_2\ol\alpha,\gamma_1)$,
$$\tikzmatht{
  \fill[black!10] (-1.5,-1.5) rectangle (1.5,1.5);
  \fill[white] (-.6,-.4) rectangle (.6,.4);
    \draw[thick] (-.6,-.4) rectangle (.6,.4);
    \draw[thick] (-.4,.4)--(-.4,1.5) (.4,.4)--(.4,1.5) (0,-.4)--(0,-1.5);
    \node at (-.8,1) {$\alpha$};
    \node at (.8,1) {$\gamma_1$}; 
    \node at (.4,-1) {$\gamma_2$};
} \mapsto 
\tikzmatht{
  \fill[black!10] (-1.5,-1.5) rectangle (1.5,1.5);
  \fill[white] (-.3,-.4) rectangle (.9,.4);
    \draw[thick] (-.3,-.4) rectangle (.9,.4);
    \draw[thick] (-.1,.4)--(-.1,.6) arc(0:180:.4)--(-.9,-1.5) 
(.7,.4)--(.7,1.5) (.3,-.4)--(.3,-1.5);
    \node at (-.5,-1) {$\ol\alpha$};
    \node at (1.1,1) {$\gamma_1$}; 
    \node at (.8,-1) {$\gamma_2$};
    \node at (-1,1.2) {$w^*$};
}, \qquad
\tikzmatht{
\fill[black!10] (-1.5,-1.5) rectangle (1.5,1.5);
\fill[white] (-.6,-.4) rectangle (.6,.4);
    \draw[thick] (-.6,-.4) rectangle (.6,.4);
    \draw[thick] (-.4,-.4)--(-.4,-1.5) (.4,-.4)--(.4,-1.5) (0,.4)--(0,1.5);
    \node at (-.8,-1) {$\ol\alpha$};
    \node at (.8,-1) {$\gamma_2$}; 
    \node at (.4,1) {$\gamma_1$};
} \mapsto 
\tikzmatht{
  \fill[black!10] (-1.5,-1.5) rectangle (1.5,1.5);
  \fill[white] (-.3,-.4) rectangle (.9,.4);
    \draw[thick] (-.3,-.4) rectangle (.9,.4);
    \draw[thick] (-.1,-.4)--(-.1,-.6) arc(360:180:.4)--(-.9,1.5) 
(.7,-.4)--(.7,-1.5) (.3,.4)--(.3,1.5);
    \node at (-.5,1) {$\alpha$};
    \node at (1.1,-1) {$\gamma_2$}; 
    \node at (.7,1) {$\gamma_1$};
    \node at (-1.1,-1.2) {$\ol w^*$};
}, 
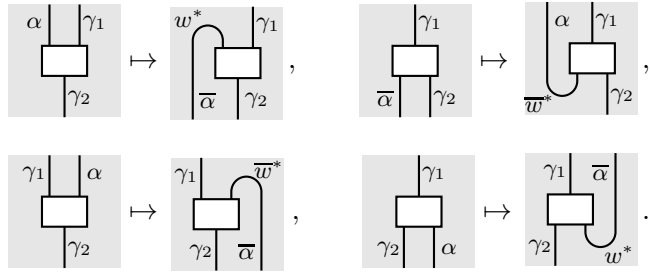
$$
$$
\tikzmatht{
  \fill[black!10] (-1.5,-1.5) rectangle (1.5,1.5);
  \fill[white] (-.6,-.4) rectangle (.6,.4);
    \draw[thick] (-.6,-.4) rectangle (.6,.4);
    \draw[thick] (-.4,.4)--(-.4,1.5) (.4,.4)--(.4,1.5) (0,-.4)--(0,-1.5);
    \node at (-.8,1) {$\gamma_1$};
    \node at (.8,1) {$\alpha$}; 
    \node at (.4,-1) {$\gamma_2$};
} \mapsto
\tikzmatht{
  \fill[black!10] (-1.5,-1.5) rectangle (1.5,1.5);
  \fill[white] (-.9,-.4) rectangle (.3,.4);
    \draw[thick] (-.9,-.4) rectangle (.3,.4);
    \draw[thick] (-.7,.4)--(-.7,1.5) (.1,.4)--(.1,.6)
    arc(180:0:.4)--(.9,-1.5) (-.3,-.4)--(-.3,-1.5); 
    \node at (-1.1,1) {$\gamma_1$};
    \node at (.5,-1) {$\ol\alpha$}; 
    \node at (-.7,-1) {$\gamma_2$};
    \node at (1.1,1.2) {$\ol w^*$};
}, \qquad
\tikzmatht{
  \fill[black!10] (-1.5,-1.5) rectangle (1.5,1.5);
  \fill[white] (-.6,-.4) rectangle (.6,.4);
    \draw[thick] (-.6,-.4) rectangle (.6,.4);
    \draw[thick] (-.4,-.4)--(-.4,-1.5) (.4,-.4)--(.4,-1.5) (0,.4)--(0,1.5);
    \node at (-.8,-1) {$\gamma_2$};
    \node at (.8,-1) {$\alpha$}; 
    \node at (.4,1) {$\gamma_1$};
} \mapsto 
\tikzmatht{
  \fill[black!10] (-1.5,-1.5) rectangle (1.5,1.5);
  \fill[white] (-.9,-.4) rectangle (.3,.4);
    \draw[thick] (-.9,-.4) rectangle (.3,.4);
    \draw[thick] (-.7,-.4)--(-.7,-1.5) (.1,-.4)--(.1,-.6)
    arc(180:360:.4)--(.9,1.5) (-.3,.4)--(-.3,1.5); 
    \node at (-1.1,-1) {$\gamma_2$};
    \node at (.5,1) {$\ol\alpha$}; 
    \node at (-.7,1) {$\gamma_1$};
    \node at (1,-1.2) {$w^*$};
}.
$$
These maps along with the ensuing equalities of the dimensions of the
intertwiner spaces,  
$$\dim\Hom(\gamma_2,\alpha\gamma_1)=\dim\Hom(\ol\alpha\gamma_2,\gamma_1), 
\qquad 
\dim\Hom(\gamma_2,\gamma_1\alpha) = \dim\Hom(\gamma_2\ol\alpha,\gamma_1),$$
are usually referred to as  {\bf Frobenius reciprocities}. 

\subsection{Non-factorial extensions} 
\label{s:nonfact}
\setcounter{equation}{0}
We want to extend our setup to $N$ being a factor, while $M$ is
admitted to be a properly infinite von Neumann algebra with finite
centre. For a related analysis, see \cite{FI} and \cite{BDH}. 

$M$ is a direct sum of finitely many infinite factors 
$$M=\bigoplus_i M_i.$$
The units of $M_i$ are the minimal central projections
$e_i$ of $M$. A homomorphism $\varphi:N\to M$ can then be written as
$$\varphi(n)=\bigoplus_i \varphi_i(n).$$
Unlike the direct sum of sectors involving isometric intertwiners,
cf.\ \sref{s:endo}, this is the true direct sum of homomorphisms
$\varphi_i:N\to M_i$, which is a homomorphism $N\to\bigoplus_i M_i$. 

Notice that the central projections $e_i\in M$ are self-intertwiners
of $\varphi$, but $e_i$ can {\em not} be split as $ss^*$ with 
isometries $s\in M$. Therefore, the direct sum of sectors
$[\varphi_i]$ as in \sref{s:endo} is not defined. 

\begin{proposition}\label{p:nonfact}
If all $\varphi_i:N\to M_i$ have conjugates $\ol\varphi_i$, then
a conjugate homomorphism $\ol\varphi:M\to N$ of $\varphi$ can be defined as 
$$\ol\varphi(m) = \sum_i s_i\ol\varphi_i(m_i)s_i^*$$
where $m=\bigoplus_i m_i$, $m_i\in M_i$, and $s_i$ are isometries 
in $N$ satisfying $s_i^*s_j=\delta_{ij}$ and $\sum_i s_is_i^*=\eins_N$. 
The dimension of $\varphi$ is 
\be\label{dim-nonf}
\dim(\varphi)=\big(\sum_i\dim(\varphi_i)^2\big)^{\frac12}.
\ee
\end{proposition}

The dimension $\dim(\varphi)$ is defined by the same infimum as
\eref{inf}, taken over all solutions $(w,\ol w)$ of the conjugacy
relations such that $w^*\ncirc w = d\cdot \eins_N$, $\ol w^*\ncirc \ol w 
=d\cdot \eins_M$. Notice that it is no longer additive, as in the factor case.

\medskip

{\em Proof:} One easily sees that the solutions of the conjugacy relations 
are parameterized by 
$$w=\sum_i\lambda_i\cdot s_iw_i,\qquad 
\ol w=\bigoplus_i\ol\lambda_i\inv\cdot \varphi_i(s_i)\ol w_i,$$ 
with parameters $\lambda_i\in\CC$. Here, $(w_i,\ol w_i)$ are solutions
for $(\varphi_i,\ol\varphi_i)$ satisfing $w_i^*\ncirc w_i=d_i\cdot
\eins_N$ and $\ol w_i^*\ncirc\ol w_i=d_i\cdot \eins_{M_i}$.  
Imposing $w^*\ncirc w=d\cdot \eins_N$ and $\ol w^*\ncirc \ol w=d\cdot \eins_M$ 
fixes the numerical coefficients by $\vert\lambda_i\vert^2=d/d_i$ and 
$d^2=\sum_i d_i^2$. This quantity is minimized if all $d_i$ are minimal, 
i.e., all $(w_i,\ol w_i)$ are standard, and $d_i=\dim(\varphi_i)$. 
This completes the proof. 
\qed
 
\begin{remark} \label{r:BDH} 
For standard pairs $(w,\ol w)$ of multiples of isometries satisfying
the minimality condition, the tracial properties (\pref{p:traces}, 
\pref{p:leftright}, and \pref{p:dotstar}) fail in general, when $M$ 
(or $N$) is not a factor.  
The authors of \cite{BDH} propose a different ``normalization condition'' 
(Eq.\ (4.3) in \cite{BDH}) for solutions to the conjugacy relations,
with $w^*\ncirc w\in N$ and $\ol w^*\ncirc \ol w\in M$ central but in
general not multiples of $\eins$. In the case of $N$ and $M$ both being
factors, their condition amounts to the equality of the left and right
traces, hence is equivalent to standardness by \pref{p:leftright}, but it
distinguishes different normalizations otherwise. In the case at hand,
it would rather fix $\vert\lambda_i\vert^2=1$, so that $\ol w^*\ncirc\ol w$ 
is no longer a multiple of an isometry. 
\end{remark}

\section{Frobenius algebras and modules}
\setcounter{equation}{0}
\label{s:frob}

We collect here some relevant results about the (simple, strict,
Karoubian) C* tensor category $\End_0(N)$ for an infinite von 
Neumann factor $N$. In fact, every full subcategory of $\End_0(N)$ can
be canonically completed so as to become a simple strict Karoubian 
C* tensor category with direct sums 
$$\C\subset\End_0(N).$$
This completion is precisely given by the constructions exposed 
in \sref{s:endo}. Without further specification, throughout this paper
$\C\subset\End_0(N)$ will denote a subcategory with the stated
properties. 

In the motivating application to QFT, as exposed in \sref{s:DHRendo},
$N$ will be the von Neumann algebra $\A(O)$ of observables localized
in some region $O$ of spacetime, which is known to be 
an infinite factor under very general assumptions. The assignment 
$O\mapsto\A(O)$ is called the local net of observables, and a 
distinguished class of positive-energy representations can be described 
by DHR endomorphisms \cite{DHR} of this net, which form a C* tensor 
category $\DHR(\A)$ (strict, simple, with subobjects, direct sums and 
conjugates). The DHR endomorphisms localized in $O$, when restricted 
to $\A(O)$, are in fact endomorphisms of $\A(O)$, and they have the 
same intertwiners as endomorphisms of the net and as elements of 
$\End(\A(O))$ \cite{GL96}. Therefore, they are the objects of a C*
tensor category $\DHR(\A)\vert_O$, which is a full
subcategory both of $\DHR(\A)$ and of $\End(N)$, $N=\A(O)$. 

In other words, if $\rho$ is localized in $O$, then one may safely
drop the distinction between $\rho\in\DHR(\A)$ and 
$\rho\in\DHR(\A)\vert_O\subset \End(N)$. 

Since $\dim(\rho)$ was defined in terms of intertwiners, one may
assign the same dimension to $\rho$ as a DHR endomorphism, and the
same properties (additivity and multiplicativity) remain valid. This  
definition coincides \cite{L89} with the ``statistical dimension'' 
originally defined in terms of the statistics operators \cite{DHR,FRS}.

It is physically most important that $\DHR(\A)$ is in fact 
a {\em braided} category, and in certain cases even {\em modular}. 
However, in our exposition, a braiding of the category $\C$ is not 
required before \sref{s:braid}, and the braided category is 
not required to be modular before \sref{s:modular}. 

\medskip

If the category $\C\subset\End_0(N)$ has only finitely many
equivalence classes of irreducible objects (sectors), then it is
called {\bf rational}. In this case, the  
structures discussed below admit only finitely many realizations, 
with complete classification available in many models. The {\bf global
dimension} of a rational tensor category is 
\be\label{globaldim}
\dim(\C):=\sum_{[\rho]\,\mathrm{irr}} \dim(\rho)^2,\ee
where the sum extends over the irreducible sectors of $\C$. 

\begin{graybox}
\begin{example} \vskip-5mm \label{x:Icat} (The Ising tensor category) 

In order to illustrate the ``rigidity'' of a C* tensor category (and
as a reference for further examples), we introduce the {\bf Ising
category}, which is one of two tensor categories with three
self-conjugate equivalence classes $[\id]$, $[\tau]$, $[\sig]$ of
irreducible objects with ``fusion rules'' $[\tau^2]=[\id]$,
$[\tau\circ\sig]= [\sig\circ\tau]=[\sig]$, $[\sig^2]=[\id]\oplus[\tau]$. 
It arises in QFT, e.g., as the category of DHR endomorphisms 
\sref{s:QFT} of the chiral Ising model.  

  The tensor category is specified by a choice
  of a representative in each class, an isometric
  intertwiner in each intertwiner space according to
  the fusion rules, and the action of the representative endomorphisms
  on the intertwiners. For all unitarily equivalent endomorphisms, the
  intertwiners are canonically related. (To specify a category in this manner, is sometimes refereed to as the ``Cuntz
  algebra approach''.) 

Because $\tau\circ\sig$ is unitarily equivalent to $\sig$, one can choose
$\tau$ in its equivalence class such that $\tau\circ\sig=\sig$. Because
$\tau^2$ is unitarily equivalent to the identity $\id$ and
$\tau^2\circ\sig=\sig$, it follows from irreducibility of $\sig$ that
$\tau^2=\id$. Therefore, $\Hom(\sig,\tau\sig)=\Hom(\id,\tau^2)=\CC\cdot\eins$. 
The remaining nontrivial intertwiner spaces are spanned by a pair of
orthogonal isometries $r\in\Hom(\id,\sig^2)$ and
$t\in\Hom(\tau,\sig^2)$, satisfying $rr^*+tt^*=1$, and 
$u\in\Hom(\sig,\sig\tau) =\Hom(\sig,\sig\tau)\times 1_\sig=
\Hom(\sig^2,\sig^2)$. Because $u^2\in\Hom(\sig,\sig\tau^2)=
\Hom(\sig,\sig)=\CC\cdot\eins$, one may choose $u=rr^*-tt^*$. 

Because $\tau(r)\in\Hom(\tau,\sig^2)$, one may choose $t=\tau(r)$, 
thus fixing the action of $\tau$: 
$$\tau(r)=t,\qquad \tau(t)=r,\qquad \tau(u)=-u.$$
$\sig(r)\in\Hom(\sig,\sig^3)$ and $\sig(\tau)\in\Hom(\sig\tau,\sig^3)$
are linear combinations of $r$ and $t$, resp.\ $ru$ and $tu$,
invariant under the action of $\tau$. Imposing $\sig^2(a)=rar^*+
t\tau(a)t^*$ for $a=r$ and $a=t$ suffices to fix all
coefficients up to an overall sign. For the Ising category one has  
$$\sig(r)=2^{-\frac12}(r+t),\qquad \sig(t)=2^{-\frac12}(r-t)u.$$
(The category of DHR endomorphisms of the $su(2)$ current algebra
at level 2 is specified by the opposite sign: $\sig(r)=-2^{-\frac12}(r+t)$,
$\sig(t)=-2^{-\frac12}(r-t)u$.)  

The dimensions are $\dim(\tau)=1$, $\dim(\sig)=\sqrt2$, and the global
dimension is $\dim(\C)=4$.
\end{example}
\end{graybox}

\subsection{C* Frobenius algebras}
\label{s:CFA}
\setcounter{equation}{0}
A Frobenius algebra $\BA=(\theta,w,x,\widehat w,\widehat x)$ in a 
C* tensor category (satisfying the unit, counit, associativity,
coassociativity and Frobenius relations \cite{FFRS06}) is called {\bf C* Frobenius
algebra} if the dual morphisms are given by the adjoint operators:
$\widehat w=w^*$ and $\widehat x = x^*$. By the latter property, the
unit and counit relations become equivalent, and so do the
associativity and coassociativity relations. 

More precisely, $\theta$ is an object of the C* category, and
$w\in\Hom(\id,\theta)$ and $x\in\Hom(\theta,\theta^2)$ are morphisms 
satisfying the relations 
\bea
\hbox{\bf unit property:} && (w^*\times 1_\theta)\scirc x
=(1_\theta\times w^*)\scirc x =1_\theta 
\notag \\ \label{unit} &&
\tikzmatht{
  \fill[black!10] (-1.2,-1.5) rectangle (1.2,1.2);
    \draw[thick] (.8,1.2)--(.8,.2) arc(360:180:.8)--(-.8,.6);
    \draw[thick] (0,-1.5)--(0,-.6);
  \fill[white] (-.8,.6) circle(.2);
    \draw[thick] (-.8,.6) circle(.2);
  \fill[black] (0,-.6) circle(.12);
    \node at (-.1,.6) {$w^*$};
    \node at (.4,-.9) {$x$};
} =
\tikzmatht{
  \fill[black!10] (-1.2,-1.5) rectangle (1.2,1.2);
    \draw[thick] (-.8,1.2)--(-.8,.2) arc(180:360:.8)--(.8,.6);
    \draw[thick] (0,-1.5)--(0,-.6);
  \fill[white] (.8,.6) circle(.2);
    \draw[thick] (.8,.6) circle(.2);
  \fill[black] (0,-.6) circle(.12);
} =
\tikzmatht{
  \fill[black!10] (-1,-1.5) rectangle (1,1.2);
    \draw[thick] (0,-1.5)--(0,1.2);
}
\\
\hbox{\bf associativity:} && (x\times 1_\theta)\scirc x
=(1_\theta\times x)\scirc x
\notag \\ \label{asso} &&
\tikzmatht{
  \fill[black!10] (-1.5,-1) rectangle (1.5,2);
    \draw[thick] (-1.2,2)--(-1.2,1.5) arc(180:360:.6)--(0,2);
    \draw[thick] (-.6,.9) arc(180:360:.9)--(1.2,2);
    \draw[thick] (.3,0)--(.3,-1);
  \fill[black] (-.6,.9) circle(.12);
  \fill[black] (.3,0) circle(.12);
    \node at (-.8,.5) {$x$};
    \node at (0,-.4) {$x$};
}=
\tikzmatht{
  \fill[black!10] (-1.5,-1) rectangle (1.5,2);
    \draw[thick] (1.2,2)--(1.2,1.5) arc(360:180:.6)--(0,2);
    \draw[thick] (.6,.9) arc(360:180:.9)--(-1.2,2);
    \draw[thick] (-.3,0)--(-.3,-1);
  \fill[black] (.6,.9) circle(.12);
  \fill[black] (-.3,0) circle(.12);
}=:
\tikzmatht{
  \fill[black!10] (-1.5,-1) rectangle (1.5,2);
    \draw[thick] (-1.2,2)--(-1.2,1.6) arc(180:360:1.2)--(1.2,2);
    \draw[thick] (0,-1)--(0,2); 
  \fill[black] (0,.4) circle(.12);
    \node at (.7,0) {$x^{(2)}$};
};
\\  
\hbox{\bf Frobenius property:} && (1_\theta\times x^*)\scirc (x\times 1_\theta) = x\ncirc x^*  =
(x^*\times 1_\theta)\scirc(1_\theta\times x)
\notag \\ \label{frob} &&
\tikzmatht{
  \fill[black!10] (-1.5,-1.5) rectangle (1.5,1.5);
    \draw[thick] (-1.2,1.5)--(-1.2,0) arc(180:360:.6)
    arc(180:0:.6)--(1.2,-1.5); 
    \draw[thick] (-.6,-.6)--(-.6,-1.5) (.6,.6)--(.6,1.5);
  \fill[black] (-.6,-.6) circle (.12); 
  \fill[black] (.6,.6) circle (.12);
}=
\tikzmatht{
  \fill[black!10] (-1.5,-1.5) rectangle (1.5,1.5);
    \draw[thick] (-.9,1.5)--(-.9,1.2) arc(180:360:.9)--(.9,1.5);
    \draw[thick] (-.9,-1.5)--(-.9,-1.2) arc(180:0:.9)--(.9,-1.5);
    \draw[thick] (0,-.3)--(0,.3);
  \fill[black] (0,-.3) circle (.12); \fill[black] (0,.3) circle (.12);
}=
\tikzmatht{
  \fill[black!10] (-1.5,-1.5) rectangle (1.5,1.5);
    \draw[thick] (1.2,1.5)--(1.2,0) arc(360:180:.6) arc(0:180:.6)--(-1.2,-1.5);
    \draw[thick] (.6,-.6)--(.6,-1.5) (-.6,.6)--(-.6,1.5);
  \fill[black] (-.6,.6) circle (.12); \fill[black] (.6,-.6) circle (.12);
}. 
\eea

In view of \eref{asso}, we also write $x^{(2)}$ for $(x\times
1_\theta)\scirc x =(1_\theta\times x)\scirc x$. 

Clearly, $w^* \ncirc w \in\Hom(\id,\id)=\CC$ is a multiple of $\eins$. 

\begin{definition} \label{d:standard}
If in addition, also $x^* x \in\Hom(\theta,\theta)$ is a multiple 
of $1_\theta$, the C* Frobenius algebra is called {\bf special}. If moreover, 
\be \label{standard}
w^*\ncirc w = d_{\BA} \cdot 1_\id
\qquad\hbox{and}\qquad 
x^*\ncirc x = d_{\BA} \cdot 1_\theta
\ee
with $d_{\BA}=\sqrt{\dim(\theta)}$, we call the C* Frobenius algebra
{\bf standard}. The number $d_{\BA}\geq 1$ is called the dimension of $\BA$.
\end{definition}

If $\alpha:N \to M$ and $\ol\alpha:M\to N$ are conjugate homomorphisms
between two factors, and 
$(w\in\Hom(\id_N,\ol\alpha\alpha),\ol w\in\Hom(\id_M,\alpha\ol\alpha))$ 
is a solution of the conjugacy relations, then 
$$(\theta=\ol\alpha\alpha,w,x=\ol\alpha(\ol w))$$ 
is a C* Frobenius algebra. It is automatically special because
$w^*\ncirc w\in\Hom(\id_N,\id_N)$ and $\ol w^*\ncirc \ol w
\in\Hom(\id_M,\id_M)$ are multiples of $\eins$. $(\theta,w,x)$ is standard
if and only if the pair $(w,\ol w)$ is standard. Therefore, 
standardness can not always be enforced by a scalar rescaling of a special C*
Frobenius algebra.  

Our aim is to prove \tref{t:Q=ext} which states that every standard 
C* Frobenius algebra is in fact of this type. 

Let us first comment on the independence of the above axioms.

\begin{lemma} \label{l:w>special} \cite{TFT1}
A Frobenius algebra is special, i.e., $x^*\ncirc x =\lambda\cdot
1_\theta$, if and only if $x^*\ncirc x\scirc w = \lambda\cdot w$ is a 
multiple of $w$. In particular, every Frobenius algebra with 
$\Hom(\id,\theta)$ one-dimensional is special.
\end{lemma}

{\em Proof:} $x^*\ncirc x =\lambda\cdot
1_\theta$ trivially implies $x^*\ncirc x\scirc w = \lambda\cdot w$. For
the converse conclusion:
$$
\tikzmatht{
  \fill[black!10] (-1.2,-1.5) rectangle (1.2,1.5);
    \draw[thick] (0,0) circle(.5); 
    \draw[thick] (0,-1.5)--(0,-.5) (0,1.5)--(0,.5); 
  \fill[black] (0,.5) circle(.12);
  \fill[black] (0,-.5) circle(.12);
} \stapel{\nref{unit}}=
\tikzmatht{
  \fill[black!10] (-1.3,-1.5) rectangle (1.3,1.5);
    \draw[thick] (0,-.6) circle(.5); 
    \draw[thick] (0,-1.5)--(0,-1.1) (0,-.1)--(0,.5) (-.5,1.5)--
(-.5,1) arc(180:360:.5)--(.5,1.2); 
  \fill[white] (.5,1.2) circle(.2); 
    \draw[thick] (.5,1.2) circle(.2); 
  \fill[black] (0,-1.1) circle(.12);
  \fill[black] (0,-.1) circle(.12);
  \fill[black] (0,.5) circle(.12);
} \stapel{\nref{frob}}=
\tikzmatht{
  \fill[black!10] (-1.5,-1.5) rectangle (1.1,1.5);
    \draw[thick] (-1,1.5)--(-1,0) arc(180:360:.4)arc(180:0:.4)--
(.6,-.5) arc(360:180:.6)--(-.6,-.4);
    \draw[thick] (0,-1.5)--(0,-1.1) (.2,1)--(.2,.4); 
  \fill[white] (.2,1) circle(.2); 
    \draw[thick] (.2,1) circle(.2);
  \fill[black] (0,-1.1) circle(.12);
  \fill[black] (-.6,-.4) circle(.12);
  \fill[black] (.2,.4) circle(.12);
} \stapel{\nref{asso}}=
\tikzmatht{
  \fill[black!10] (-1.5,-1.5) rectangle (1.1,1.5);
    \draw[thick] (-1,1.5)--(-1,-.4) arc(180:360:.6);
    \draw[thick] (.2,0) circle(.4);
    \draw[thick] (-.4,-1.5)--(-.4,-1) (.2,1)--(.2,.4);
  \fill[white] (.2,1) circle(.2); 
    \draw[thick] (.2,1) circle(.2); 
  \fill[black] (.2,-.4) circle(.12);
  \fill[black] (.2,.4) circle(.12);
  \fill[black] (-.4,-1) circle(.12);
} =\lambda \cdot
\tikzmatht{
  \fill[black!10] (-1.5,-1.5) rectangle (1.1,1.5);
    \draw[thick] (-1,1.5)--(-1,-.4) arc(180:360:.6)--(.2,.4);
    \draw[thick] (-.4,-1.5)--(-.4,-1);
  \fill[white] (.2,.4) circle(.2); 
    \draw[thick] (.2,.4) circle(.2); 
  \fill[black] (-.4,-1) circle(.12);
} \stapel{\nref{unit}}=\lambda \cdot 
\tikzmatht{
  \fill[black!10] (-.8,-1.5) rectangle (.8,1.5);
    \draw[thick] (0,1.5)--(0,-1.5);
}.$$
Standardness, however, is not automatic, as explained before. 

\begin{definition} \label{d:hom0}
For a Frobenius algebra $(\theta,w,x)$ in a simple C* tensor category,
$\Hom_0(\theta,\theta)$ is the subspace $\Hom(\theta,\theta)$ of
elements satisfying 
\be
\label{1to3} (1_\theta\times t)\scirc x = x\scirc t = (t\times 1_\theta)\scirc x: \qquad
\tikzmatht{
  \fill[black!10] (-1.3,-1.3) rectangle (1.3,1.3);
    \draw[thick] (-.8,1.3)--(-.8,.5) arc(180:360:.8)--(.8,1.3);
    \draw[thick] (0,-1.3)--(0,-.3); 
  \fill[white] (1.1,.4) rectangle (.5,.8);
    \draw[thick] (1.1,.4) rectangle (.5,.8);
  \fill[black] (0,-.3) circle (.12);
}=
\tikzmatht{
  \fill[black!10] (-1.3,-1.3) rectangle (1.3,1.3);
    \draw[thick] (-.8,1.3)--(-.8,1) arc(180:360:.8)--(.8,1.3);
    \draw[thick] (0,-1.3)--(0,.2); 
  \fill[white] (-.3,-.8) rectangle (.3,-.4);
    \draw[thick] (-.3,-.8) rectangle (.3,-.4);
  \fill[black] (0,.2) circle (.12);
}=
\tikzmatht{
  \fill[black!10] (-1.3,-1.3) rectangle (1.3,1.3);
    \draw[thick] (-.8,1.3)--(-.8,.5) arc(180:360:.8)--(.8,1.3);
    \draw[thick] (0,-1.3)--(0,-.3); 
  \fill[white] (-1.1,.4) rectangle (-.5,.8);
    \draw[thick] (-1.1,.4) rectangle (-.5,.8);
  \fill[black] (0,-.3) circle (.12);
}.
\ee
\end{definition}

We shall later identify this space with the self-morphisms of the
Frobenius algebra as a bimodule of itself (\sref{s:bim}), and exhibit
the importance of this space for the centre of the von Neumann algebra
extension $N\subset M$ associated with a Frobenius algebra
(\sref{s:Q}). 

\begin{lemma} \label{l:special} If $(\theta,w,x)$ is a Frobenius
  algebra in a simple C* tensor category, then 
$n:=x^*\ncirc x$ is a strictly positive
element of $\Hom_0(\theta,\theta)$. 
\end{lemma}

{\em Proof:} If $w^*\ncirc w=d\cdot 1_\id$, then $d\inv
\cdot w\ncirc w^*$ is a projection, hence $1_\theta\geq d\inv\cdot
w\ncirc w^*$, hence 
\be\label{strict}
x^*\ncirc x\geq d\inv\cdot x^*\scirc (1\times w\ncirc w^*)\scirc x
\stapel{\nref{unit}}= d\inv\cdot 1_\theta\ee  
is strictly positive. \eref{1to3} follows by associativity and the
Frobenius property: 
$$
\tikzmatht{
  \fill[black!10] (-1.5,-1.3) rectangle (1.5,1.5);
    \draw[thick] (1,.5) arc(0:360:.5);
    \draw[thick] (.5,1.5)--(.5,1) (.5,0) arc(360:180:.75)--(-1,1.5); 
    \draw[thick] (-.25,-1.3)--(-.25,-.75);
  \fill[black] (.5,1) circle(.12);
  \fill[black] (.5,0) circle(.12);
  \fill[black] (-.25,-.75) circle(.12);
} \stapel{\nref{asso}}=
\tikzmatht{
  \fill[black!10] (-1.5,-1.3) rectangle (1.5,1.5);
    \draw[thick] (-1,1.5)--(-1,.5) arc(180:360:.5) arc(180:0:.5)--(1,0)
arc(360:180:.75);
    \draw[thick] (.5,1.5)--(.5,1);
    \draw[thick] (.25,-1.3)--(.25,-.75);
  \fill[black] (.5,1) circle(.12);
  \fill[black] (-.5,0) circle(.12);
\fill[black] (.25,-.75) circle(.12);
} \stapel{\nref{frob}}=
\tikzmatht{
  \fill[black!10] (-1.3,-1.3) rectangle (1.3,1.5);
    \draw[thick] (.5,-.25) arc(0:360:.5);
    \draw[thick] (-.75,1.5) arc(180:360:.75); 
    \draw[thick] (0,-1.3)--(0,-.75) (0,.25)--(0,.75);
  \fill[black] (0,.75) circle(.12);
  \fill[black] (0,.25) circle(.12);
  \fill[black] (0,-.75) circle(.12);
} \stapel{\nref{frob}}=
\tikzmatht{
  \fill[black!10] (-1.5,-1.3) rectangle (1.5,1.5);
    \draw[thick] (1,1.5)--(1,.5) arc(360:180:.5) arc(0:180:.5)--(-1,0)
    arc(180:360:.75);
    \draw[thick] (-.5,1.5)--(-.5,1);
    \draw[thick] (-.25,-1.3)--(-.25,-.75);
  \fill[black] (-.5,1) circle(.12);
  \fill[black] (.5,0) circle(.12);
  \fill[black] (-.25,-.75) circle(.12);
} \stapel{\nref{asso}}=
\tikzmatht{
  \fill[black!10] (-1.5,-1.3) rectangle (1.5,1.5);
    \draw[thick] (0,.5) arc(0:360:.5);
    \draw[thick] (-.5,1.5)--(-.5,1) (-.5,0) arc(180:360:.75)--(1,1.5); 
    \draw[thick] (.25,-1.3)--(.25,-.75);
  \fill[black] (-.5,1) circle(.12);
  \fill[black] (-.5,0) circle(.12);
  \fill[black] (.25,-.75) circle(.12);
} . $$
\qed

\begin{corollary} \label{c:special}
The equivalent Frobenius algebra $(\theta,\wh w,\wh x)$ with 
$$\wh w := n^{\frac12}\scirc w,\qquad \wh x := (n^{-\frac12}\times
n^{-\frac12})\scirc x \scirc n^{\frac12}$$
is special. 
\end{corollary}

{\em Proof:} Replacing $w,x$ by $\wh w,\wh x$ clearly preserves the
unit property, associativity, and Frobenius property, and $\wh
w^*\ncirc \wh w \in \Hom(\id,\id)=\CC\cdot 1$. Specialness follows
because along with $n\in \Hom_0(\theta,\theta)$, also 
$n\inv \in \Hom_0(\theta,\theta)$, hence:
$$\wh x^*\ncirc \wh x = n^{\frac12}\scirc x^*\scirc (n\inv \times n\inv)\scirc
x \scirc n^{\frac12} \stapel{\nref{1to3}}= n^{-\frac12}\scirc x^*\ncirc x \scirc
n^{-\frac12} = n^{-\frac12}\scirc n \scirc n^{-\frac12} = 1_\theta.$$
\qed

In C* tensor categories, also the Frobenius property is not 
independent from the other relations. We shall now prove that \eref{frob} follows from
\eref{unit} and \eref{asso} along with the special property $x^*\ncirc
x =\lambda\cdot 1$. Notice that specialness is a relation in 
$\Hom(\theta,\theta)$, and is thus ``simpler'' than the Frobenius 
relation \eref{frob} in $\Hom(\theta^2,\theta^2)$. 
 
\begin{lemma} \label{l:UAS>F} 
\cite{LRo} In a C* tensor category, the Frobenius property is a
consequence of unit property, associativity, and specialness.
\end{lemma}

{\em Proof:}
Let $X:= (1_\theta\times x^*)\scirc (x\times 1_\theta) - x\ncirc x^* \in
\Hom(\theta^2,\theta^2)$. 
Then, if $x^*\ncirc x = d\cdot 1$, one has 
$$X^*X= (x^*\times 1_\theta)\scirc (1_\theta\times xx^*)\scirc (x\times
1_\theta) - d\cdot x\ncirc x^* \in
\Hom(\theta^2,\theta^2),$$
where for the two mixed terms the associativity relation has been
used. We define the map $\delta:\Hom(\theta^2,\theta^2) \to 
\Hom(\theta^2,\theta^2)$ given by 
$$\delta(T)=(x^*\times 1_\theta)\scirc (1_\theta\times T)\scirc (x\times
1_\theta): \quad 
\tikzmatht{
  \fill[black!10] (-1.5,-1.5) rectangle (1.5,1.5);
  \fill[white] (-.5,-.4) rectangle (1,.4); 
    \draw[thick] (-.5,-.4) rectangle node {$T$} (1,.4); 
    \draw[thick] (0,.4)--(0,.6) arc(0:180:.5)--(-1,-.6)
          arc(180:360:.5)--(0,-.4) (-.5,1.1)--(-.5,1.5) (-.5,-1.1)--(-.5,-1.5); 
    \draw[thick] (.5,.4)--(.5,1.5) (.5,-.4)--(.5,-1.5); 
  \fill[black] (-.5,1.1) circle(.12); \fill[black] (-.5,-1.1) circle(.12);
 }$$
$\delta$ is positive and faithful: Namely if $T=Y^*Y$ is positive,
then $\delta(T)=Z^*Z$ with $Z=(1_\theta\times Y)\scirc (x\times
1_\theta)$, hence $\delta(Y^*Y)$ is positive; and $\delta(Y^*Y)=0$
implies $Z=0$ from which it follows that $Y=0$ by the unit property
\eref{unit}. We apply $\delta$ to $T=X^*X$. Again, using the
associativity relation, one finds $\delta(X^*X)=0$. Hence $X=0$.  
\qed

\medskip

We make a little digression to report also the following observation: 
A triple satisfying only the unit property and associativity can be
``deformed'' in such a way that it is in addition special. Then the
Frobenius property also follows by \lref{l:UAS>F}, hence the deformed
triple is a C* Frobenius algebra. There enters in the proof, however,
a certain ``regularity condition'' which we do not quite know how to
control.  

The admitted deformations by any invertible element
$n\in\Hom(\theta,\theta)$ are defined via 
$$w\mapsto n^*{}{\inv} \scirc w,\quad x\mapsto (n\times n)\scirc x \scirc
n\inv,$$
obviously preserving \eref{unit} and \eref{asso}. The deformed triple
is standard if
$$x^*\scirc (n^*\ncirc n \times n^*\ncirc n)\scirc x = n^*\ncirc n.$$
We want to solve this eqution by iterating the following recursion: 
$$m_{k+1}:= x^*\scirc (m_k \times m_k)\scirc x,$$
starting with $m_0=1$, i.e., $m_1=x^*\ncirc x$. Clearly, each $m_k$ is a
positive element of $\Hom(\theta,\theta)$. It is even strictly
positive, because
$\big(\theta,w_{k+1}:=m_k^{-\frac12}\scirc w_k,x_{k+1}:=(m_k^{\frac12}\times
m_k^{\frac12})\scirc v_k\scirc m_k^{-\frac12}\big)$ 
is a sequence of triples satisfying \eref{unit} and \eref{asso}, and
$x_k^*x_k=m_k$ is strictly positive by \eref{strict}. The question is,
of course, whether $(m_k)_k$ converges. 

Now $\Hom(\theta,\theta)$ equipped with the product $m_1\ast
m_2=x^*\scirc (m_1\times m_2)\scirc x$ is an algebra. The algebra has the unit
$w\ncirc w^*$, and is associative by \eref{asso}. It is finite-dimensional,
because $\Hom(\theta,\theta)$ is finite-dimensional. Hence it is
isomorphic to some matrix algebra. W.r.t.\ this product,
$m_0=1_\theta$, $m_1=1_\theta\ast 1_\theta$, and 
$m_k= 1_\theta^{\ast 2^k}$. Because $m_k$ are strictly positive, they
cannot be zero, hence $1_\theta$ is not nilpotent w.r.t.\ the
$\ast$-product. Hence it has some largest eigenvalue, and hence some
multiple $\mu_0$ of $1_\theta$ has a largest eigenvalue 1, so that
$\mu_0^{\ast 2^k}$ converges to an idempotent $m$ w.r.t.\ the
$\ast$-product. This element therefore solves $x^*\scirc(m\times
m)\scirc x =
m$. If $m$ is strictly positive, then deforming the original triple
$(\theta,w,x)$ with $n=m^{\frac12}$, would give rise to a special
triple, which then satisfies the Frobenius property by \lref{l:UAS>F}. 
However, we only know that $m$ is positive as a limit of strictly 
positive elements of $\Hom(\theta,\theta)$. The ``regularity condition''
mentioned above is the absence of a kernel of the limit. (Actually, in
order to solve the equation, one may start from any initial element
$\mu_0$ (not necessarily a multiple of $1_\theta$), but in the most
general case, one will have even less control over the invertibility of
the limit.) 

\medskip

After this digression, we return to the main line of the section.

\subsection{Q-systems and extensions} 
\label{s:Q}
\setcounter{equation}{0}
\begin{definition} \label{d:Q} A {\bf Q-system} is a standard Frobenius
algebra $\BA=(\theta,w,x)$ in a simple strict C* tensor category
$\C$. Its {\em dimension} is $d_\BA=\sqrt{\dim(\theta)}$. 
\end{definition}

Even in the irreducible case, where the canonical endomorphism
$\theta$ fixes the intertwiner $w\in \Hom(\id,\theta)$ up to a complex
phase, there may be finitely many inequivalent $x\in
\Hom(\theta,\theta^2)$ \cite{IK}. 

From now on, we reserve the graphical representation
$$w= \tikzmatht{
  \fill[black!10] (-1,-1) rectangle (1,1);
    \draw[thick] (0,1)--(0,-.2);
  \fill[white] (0,-.2) circle(.2);
    \draw[thick] (0,-.2) circle(.2);
},\qquad x=
\tikzmatht{
  \fill[black!10] (-1.2,-1) rectangle (1.2,1);
    \draw[thick] (.8,1)--(.8,.5) arc(360:180:.8)--(-.8,1);
    \draw[thick] (0,-1)--(0,-.3);
  \fill[black] (0,-.3) circle(.12);
}, \qquad 
r := x\scirc w = 
\tikzmatht{
  \fill[black!10] (-1.2,-1) rectangle (1.2,1);
    \draw[thick] (.8,1)--(.8,.5) arc(360:180:.8)--(-.8,1);
}$$
for the intertwiners associated with a Q-system, i.e., $w$ and $x$ satisfy
\eref{unit}--\eref{frob} and \eref{standard}, and $(r,r)$ 
satisfies \eref{conj}. We shall freely
use these properties in the sequel.

For the irreducible case, and $\C=\End_0(N)$, this definition first
appeared in \cite{L94} as a characterization of 
subfactors $N\subset M$. In this section, we review and generalize
this work to the reducible case. The correspondence between Q-systems
and {\bf extensions} of a factor (= inclusions into a (possibly
non-factorial) von Neumann algebra) is the main reason for the study
of Q-systems. In quantum field theory, Q-systems in $\C=\DHR(\A)$
correspond to extensions $\A\subset\B$ of a given QFT. Non-factorial
extensions naturally arise, e.g., in the ``universal construction'' of
boundary conditions discussed in \cite{BKLR}, cf.\ \sref{s:transb}. 

An immediate consequence of standardness is the following:

\begin{corollary} \label{c:rrstandard}
Let $\BA=(\theta,w,x)$ be a Q-system, $r=x\scirc w$. Then $(r,\ol r=r)$ 
is a standard pair for $(\theta,\ol\theta=\theta)$. The left and right
Frobenius conjugations $\Hom(\theta,\theta^2)\to \Hom(\theta^2,\theta)$, 
$y\mapsto (r^*\times 1_\theta)\scirc(1_\theta\times y)$ and 
$y\mapsto (1_\theta\times r^*)\scirc(y\times 1_\theta)$ take $x$ to $x^*$. 
\end{corollary}

{\em Proof:} The conjugacy relations \eref{conj} follow by applying
the definition  
$\tikzmatht{
  \fill[black!10] (-1,-1) rectangle (1,1); 
    \draw[thick] (-.5,1)--(-.5,.4) arc(180:360:.5)--(.5,1); 
    \node at (0,-.5) {$r$};
} = 
\tikzmatht{
  \fill[black!10] (-1,-1) rectangle (1,1); 
    \draw[thick] (-.5,1)--(-.5,.6) arc(180:360:.5)--
          (.5,1) (0,-.6)--(0,.1);
  \fill[white] (0,-.6) circle(.2);
    \draw[thick] (0,-.6) circle(.2);
  \fill[black] (0,.1) circle(.12);
}$
in several ways to \eref{frob}, and \eref{unit}. $(r,\ol r=r)$ is
a standard pair because $r^*r=w^*x^*xw=d_\BA w^*w=d_\BA^2=\dim(\theta)$.  
\qed

\begin{remark} \label{r:rrspecial}
If $\BA=(\theta,w,x)$ is only special, $w^*\ncirc w=d_w\cdot 1$,
$x^*\ncirc x=d_x\cdot 1_\theta$, then $(r,r)$ still solves the conjugacy 
relations by the Frobenius and unit properties. Therefore, 
$r^*\ncirc r\geq \dim(\theta)$ by the definition of the dimension as an 
infimum. Hence, $d_wd_x\geq \dim(\theta)$ with equality if and 
only if $\BA$ is standard.  
\end{remark}

\medskip

Let $N\subset M$ be an infinite subfactor of finite index, and
$\iota:N\to M$ the embedding homomorphism. This gives rise to a
Q-system in the C* tensor category $\End_0(N)$ as follows. Because the
index $[M:N]$ is finite, the dimension $\dim(\iota)$ is finite, hence
there is a conjugate homomorphism $\ol\iota:M\to N$. Let 
$$w\in \Hom(\id_N,\ol\iota\iota)\subset N,\quad v\in
\Hom(\id_M,\iota\ol\iota)\subset M$$ be a standard solution of 
the conjugacy relations \eref{conj}. Then the triple
\be\label{Qsystem}
\BA=(\theta,w,x),\qquad \theta:=\ol\iota\iota\in \End_0(N), \quad w\in
N, \quad x:=\ol\iota(v)\in N
\ee
is a Q-system in $\End_0(N)$ of dimension $d_\BA=\dim(\iota)$. 
Graphically ``resolving'' $\theta=\ol\iota\circ\iota$, the
intertwiners $w$ and $x=\ol\iota(v)$ are displayed as 
$$\tikzmatht{
  \fill[black!10] (-1,-1) rectangle (1,1.5); 
    \draw[thick] (0,1.5)--(0,0);
  \fill[white] (0,0) circle(.2);
    \draw[thick] (0,0) circle(.2);
    \node at (0,-.6) {$w$}; \node at (.3,1) {$\theta$};
} \equiv 
\tikzmatht{
  \fill[black!10] (-1.3,-1) rectangle (1.3,1.5); 
  \fill[black!18] (-.3,1.5)--(-.3,0) arc(180:360:.3)--(.3,1.5);
    \draw[thick] (-.3,1.5)--(-.3,0) arc(180:360:.3)--(.3,1.5);
    \node at (-.6,1) {$\ol\iota$};
    \node at (.6,1) {$\iota$};
} ,\qquad
\tikzmatht{
  \fill[black!10] (-1.3,-1) rectangle (1.3,1.5); 
    \draw[thick] (-.8,1.5)--(-.8,.8) arc(180:360:.8)--
          (.8,1.5);
    \draw[thick] (0,0)--(0,-1); 
  \fill[black] (0,0) circle(.12);
    \node at (0,.6) {$x$};
} \equiv
\tikzmatht{
  \fill[black!10] (-1.3,-1) rectangle (1.3,1.5); 
  \fill[black!18] (-.8,1.5)--(-.8,.5) arc(180:235:.7)
  arc(55:0:.7)--(-.2,-1)--(1,-1)--(1,1.5);  
    \draw[thick] (-.8,1.5)--(-.8,.5) arc(180:235:.7)
    arc(55:0:.7)--(-.2,-1);
  \fill[black!10] (.8,1.5)--(.8,.5) arc(0:-55:.7)
  arc(125:180:.7)--(.2,-1)--(1.3,-1)--(1.3,1.5);  
    \draw[thick] (.8,1.5)--(.8,.5) arc(0:-55:.7)
    arc(125:180:.7)--(.2,-1);
  \fill[black!10] (-.4,1.5)--(-.4,.6) arc(180:360:.4)--(.4,1.5);  
    \draw[thick] (-.4,1.5)--(-.4,.6) arc(180:360:.4)--(.4,1.5);  
    \node at (0,.8) {$v$};
},$$
so that the unit, associativity and Frobenius properties are trivially
satisfied: 
$$\tikzmatht{
  \fill[black!10] (-1.3,-1) rectangle (1.3,1.5); 
  \fill[black!18] (.8,1.5) .. controls(.2,0) .. (.2,-1)--(-.2,-1)
  .. controls(-.2,0) .. (-.7,1) arc(210:45:.2) .. controls(0,.3)
  .. (.4,1.5);  
    \draw[thick] (.8,1.5) .. controls(.2,0) .. (.2,-1); 
    \draw[thick] (-.2,-1) .. controls(-.2,0) .. (-.7,1)
          arc(210:45:.2) .. controls(0,.3) .. (.4,1.5); 
} = 
\tikzmatht{
  \fill[black!10] (-1,-1) rectangle (1,1.5); 
  \fill[black!18] (-.2,1.5)--(-.2,-1)--(.2,-1)--(.2,1.5); 
    \draw[thick] (-.2,1.5)--(-.2,-1) (.2,-1)--(.2,1.5); 
} = 
\tikzmatht{
  \fill[black!10] (-1.3,-1) rectangle (1.3,1.5); 
  \fill[black!18] (-.8,1.5) .. controls(-.2,0) .. (-.2,-1)--(.2,-1)
  .. controls(.2,0) .. (.7,1) arc(-30:135:.2) .. controls(0,.3)
  .. (-.4,1.5); 
    \draw[thick] (-.8,1.5) .. controls(-.2,0) .. (-.2,-1); 
    \draw[thick] (.2,-1) .. controls(.2,0) .. (.7,1)
          arc(-30:135:.2) .. controls(0,.3) .. (-.4,1.5); 
},
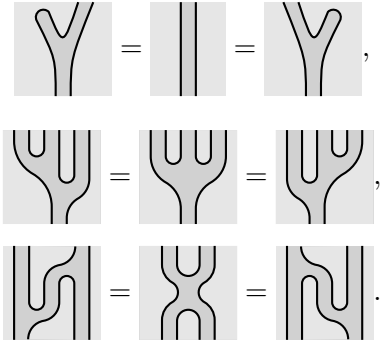
$$
$$
\tikzmatht{
  \fill[black!18] (-1.3,-1) rectangle (1.3,1.5); 
  \fill[black!10] (-1,1.5)--(-1,1.2) arc(180:240:1)
  arc(60:0:1)--(0,-1)--(-1.3,-1)--(-1.3,1.5); 
    \draw[thick] (-1,1.5)--(-1,1.2) arc(180:240:1)
    arc(60:0:1)--(0,-1); 
  \fill[black!10] (1,1.5)--(1,.2) arc(360:300:.6)
  arc(120:180:.6)--(.4,-1)--(1.3,-1)--(1.3,1.5); 
    \draw[thick] (1,1.5)--(1,.2) arc(360:300:.6)
    arc(120:180:.6)--(.4,-1); 
  \fill[black!10] (-.6,1.5)--(-.6,1) arc(180:360:.2)--(-.2,1.5); 
    \draw[thick] (-.6,1.5)--(-.6,1) arc(180:360:.2)--(-.2,1.5); 
  \fill[black!10] (.6,1.5)--(.6,.3) arc(360:180:.2)--(.2,1.5); 
    \draw[thick] (.6,1.5)--(.6,.3) arc(360:180:.2)--(.2,1.5); 
}=\tikzmatht{
  \fill[black!18] (-1.3,-1) rectangle (1.3,1.5); 
  \fill[black!10] (1,1.5)--(1,.8) arc(360:300:.8)
  arc(120:180:.8)--(.2,-1)--(1.3,-1)--(1.3,1.5); 
    \draw[thick] (1,1.5)--(1,.8) arc(360:300:.8)
    arc(120:180:.8)--(.2,-1); 
  \fill[black!10] (-1,1.5)--(-1,.8) arc(180:240:.8)
  arc(60:0:.8)--(-.2,-1)--(-1.3,-1)--(-1.3,1.5); 
    \draw[thick] (-1,1.5)--(-1,.8) arc(180:240:.8)
    arc(60:0:.8)--(-.2,-1); 
  \fill[black!10] (-.6,1.5)--(-.6,.8) arc(180:360:.2)--(-.2,1.5); 
    \draw[thick] (-.6,1.5)--(-.6,.8) arc(180:360:.2)--(-.2,1.5); 
  \fill[black!10] (.6,1.5)--(.6,.8) arc(360:180:.2)--(.2,1.5); 
    \draw[thick] (.6,1.5)--(.6,.8) arc(360:180:.2)--(.2,1.5); 
}=\tikzmatht{
  \fill[black!18] (-1.3,-1) rectangle (1.3,1.5); 
  \fill[black!10] (1,1.5)--(1,1.2) arc(360:300:1)
  arc(120:180:1)--(0,-1)--(1.3,-1)--(1.3,1.5); 
    \draw[thick] (1,1.5)--(1,1.2) arc(360:300:1)
    arc(120:180:1)--(0,-1); 
  \fill[black!10] (-1,1.5)--(-1,.2) arc(180:240:.6)
  arc(60:0:.6)--(-.4,-1)--(-1.3,-1)--(-1.3,1.5); 
    \draw[thick] (-1,1.5)--(-1,.2) arc(180:240:.6) arc(60:0:.6)--(-.4,-1);
  \fill[black!10] (.6,1.5)--(.6,1) arc(360:180:.2)--(.2,1.5); 
    \draw[thick] (.6,1.5)--(.6,1) arc(360:180:.2)--(.2,1.5); 
  \fill[black!10] (-.6,1.5)--(-.6,.3) arc(180:360:.2)--(-.2,1.5); 
    \draw[thick] (-.6,1.5)--(-.6,.3) arc(180:360:.2)--(-.2,1.5); 
},$$
$$
\tikzmatht{
  \fill[black!10] (-1.3,-1) rectangle (1.3,1.5); 
  \fill[black!18] (-1,1.5)--(-1,-1)--(1,-1)--(1,1.5); 
    \draw[thick] (-1,1.5)--(-1,-1) (1,1.5)--(1,-1); 
  \fill[black!10] (-.6,1.5)--(-.6,0) arc(180:360:.2)--(-.2,.5)
          arc(180:100:.5) arc(-80:0:.5); 
    \draw[thick] (-.6,1.5)--(-.6,0) arc(180:360:.2)--(-.2,.5)
          arc(180:100:.5) arc(-80:0:.5); 
  \fill[black!10] (.6,-1)--(.6,.5) arc(0:180:.2)--(.2,0)
          arc(0:-80:.5) arc(100:180:.5);
    \draw[thick] (.6,-1)--(.6,.5) arc(0:180:.2)--(.2,0)
          arc(0:-80:.5) arc(100:180:.5);
}=\tikzmatht{
  \fill[black!18] (-1.3,-1) rectangle (1.3,1.5);
  \fill[black!10] (.7,1.5)--(.7,1) arc(360:300:.55)
          arc(120:240:.3) arc(60:0:.55)--(.7,-1)--(1.3,-1)--(1.3,1.5);
    \draw[thick] (.7,1.5)--(.7,1) arc(360:300:.55)
          arc(120:240:.3) arc(60:0:.55)--(.7,-1);
  \fill[black!10] (-.7,1.5)--(-.7,1) arc(180:240:.55)
          arc(60:-60:.3) arc(120:180:.55)--(-.7,-1)--
          (-1.3,-1)--(-1.3,1.5); 
    \draw[thick] (-.7,1.5)--(-.7,1) arc(180:240:.55)
          arc(60:-60:.3) arc(120:180:.55)--(-.7,-1);
  \fill[black!10] (.3,1.5)--(.3,1) arc(360:180:.3)--(-.3,1.5); 
    \draw[thick] (.3,1.5)--(.3,1) arc(360:180:.3)--(-.3,1.5); 
  \fill[black!10] (.3,-1)--(.3,.-.5) arc(0:180:.3)--(-.3,-1); 
    \draw[thick] (.3,-1)--(.3,.-.5) arc(0:180:.3)--(-.3,-1); 
}=\tikzmatht{
  \fill[black!10] (-1.3,-1) rectangle (1.3,1.5); 
  \fill[black!18] (-1,1.5)--(-1,-1)--(1,-1)--(1,1.5); 
    \draw[thick] (-1,1.5)--(-1,-1) (1,1.5)--(1,-1); 
  \fill[black!10] (.6,1.5)--(.6,0) arc(360:180:.2)--(.2,.5)
          arc(0:80:.5) arc(260:180:.5); 
    \draw[thick] (.6,1.5)--(.6,0) arc(360:180:.2)--(.2,.5)
          arc(0:80:.5) arc(260:180:.5); 
  \fill[black!10] (-.6,-1)--(-.6,.5) arc(180:0:.2)--(-.2,0)
          arc(180:260:.5) arc(80:0:.5); 
    \draw[thick] (-.6,-1)--(-.6,.5) arc(180:0:.2)--(-.2,0)
          arc(180:260:.5) arc(80:0:.5); 
}.$$

Notice that the projections $d_\BA\inv\cdot ww^*$ and $d_\BA\inv\cdot xx^*$
have the same properties as the Jones projections in the type $I\!I$
case \cite{J}, satisfying the Temperley-Lieb algebra and starting the 
``Jones tunnel''. The Jones ``planar algebra'' \cite{J98}
  associated with a subfactor is the 2-category with two objects
  $N$ and $M$, whose 1-morphisms are sub-homomor\-phisms
  of alternating products of $\iota$ and $\ol\iota$, namely
  $\rho\prec(\ol\iota\iota)^n\in\End_0(N)$ for any $n\in\NN$,
  $\varphi\prec\iota(\ol\iota\iota)^n\in \Hom(N,M)$, etc., and
  whose 2-morphisms are their intertwiners.  

If $M=\bigoplus_i M_i$ is not a factor, and $\iota(n)=\bigoplus_i
\iota_i(n)$ as in \sref{s:nonfact}, then the Q-system defined by
\eref{Qsystem} can be computed with \pref{p:nonfact}:
\be\label{Qsum}
\theta(n)=\sum_i s_i\theta_i(n)s_i^*,\quad
w=\sum_i\sqrt{\frac{d}{d_i}}\cdot s_i\scirc w_i, 
\quad x= \sum_i\sqrt{\frac{d_i}{d}}\cdot (s_i\times s_i)\scirc
x_i\scirc s_i^* 
\ee
where $d=\sqrt{\dim(\theta)}=\sqrt{\sum_i\dim(\theta_i)}$, in
compliance with \eref{dim-nonf}. 

The projections $p_i=s_is_i^*\in\Hom(\theta,\theta)$ are elements of
$\Hom_0(\theta,\theta)$, cf.\ \dref{d:hom0}, i.e., they satisfy \eref{1to3}.

\medskip

The central result of this section is the converse to the construction
of a Q-system from an inclusion map $\iota:N\to M$: namely, the larger
von Neumann algebra $M$ can be reconstructed from $N$ and the Q-system.

\begin{tintedbox}
\begin{theorem} \vskip-5mm \label{t:Q=ext} \cite{L94} Let $N$ be an infinite
factor, and $\BA=(\theta,w,x)$ be a Q-system in $\End_0(N)$. Then there
is a von Neumann algebra $M$ and a homomorphism $\iota:N\to M$ with 
conjugate $\ol\iota:M\to N$ such that $\theta=\ol\iota\iota$, and 
a standard solution $(w,v)$ of the conjugacy relations \eref{conj}
such that $x=\ol\iota(v)$. The dimension $d_\BA$ equals the dimension
$\dim(\iota)=\sqrt{\dim(\theta)}$.  
\end{theorem}
\end{tintedbox}

{\em Proof:} 
The algebra $M$ is reconstructed from $N$ and the Q-system by
adjoining to $N$ one new element, called $v$, whose algebraic relations are
the same as those of the operator $v\equiv v\times
1_\iota\in\Hom(\iota,\iota\ol\iota\iota) = \Hom(\iota,\theta\iota)$ if
we {\em knew} that the Q-system comes from a conjugate solution
$(w,v)$ as before. Namely, $v$ satisfies the commutation relations: 
$$v\iota(n) = \iota\theta(n)v $$
with the elements $n\in N$ (where $\iota$ is the embedding of
$N$ into the larger algebra $M$), i.e., $v\in\Hom(\iota,\iota\theta)$,
its square is 
$$v^2 := \iota(x)v:\qquad 
\tikzmatht{
  \fill[black!10] (-1,-.5) rectangle (1.5,2); 
  \fill[black!18] (-1,-.5)--(.5,-.5)--(.5,0)--(0,.5)--(0,1)--
(-.5,1.5)--(-.5,2)--(-1,2)--(-1,-.5);  
  \fill[white] (0,1)--(.5,1.5)--(-.5,1.5)--(0,1); 
    \draw[thick] (0,1)--(.5,1.5)--(-.5,1.5)--(0,1); 
  \fill[white] (.5,0)--(1,.5)--(0,.5)--(.5,0); 
    \draw[thick] (.5,0)--(1,.5)--(0,.5)--(.5,0); 
    \draw[thick] (-.5,1.5)--(-.5,2);
    \draw[thick] (0,.5)--(0,1);
    \draw[thick] (.5,1.5)--(.5,2);
    \draw[thick] (1,.5)--(1,2);
    \draw[thick] (.5,-.5)--(.5,0);
    \node at (-.5,1) {$v$}; 
    \node at (0,0) {$v$}; 
    \node at (.7,-.2) {$\iota$}; 
} = 
\tikzmatht{
  \fill[black!10] (-1,-.5) rectangle (1.5,2); 
  \fill[black!18] (-1,-.5)--(0,-.5)--(0,0)--(-.5,.5)--(-.5,2)
--(-1,2)--(-1,-.5); 
    \draw[thick] (0,2)--(0,1.5) arc(180:360:.5)--(1,2); 
  \fill[white] (0,0)--(.5,.5)--(-.5,.5)--(0,0); 
    \draw[thick] (0,0)--(.5,.5)--(-.5,.5)--(0,0); 
    \draw[thick] (.5,.5)--(.5,1);
    \draw[thick] (-.5,.5)--(-.5,2);
    \draw[thick] (0,-.5)--(0,0);
  \fill[black] (.5,1) circle(.12);
    \node at (-.7,1.5) {$\iota$}; 
    \node at (.3,1.6) {$\theta$}; 
    \node at (1.2,1.6) {$\theta$};
},$$
and its adjoint is
$$v^*:=\iota(w^*x^*)v.$$ 
It follows from these relations that every element of $M$ can be written
in the form $\iota(n)v$ for some $n\in N$. The product thus defined is
associative by virtue of \eref{asso}, and it has a unit 
$\eins_M=\iota(w^*)v$ by virtue of \eref{unit}. The definition of $v^*$
implies the adjoint of a general element of $M$, namely
$(\iota(n)v)^*=v^*\iota(n^*)= \iota(w^*x^*\theta(n^*))v$. This turns
$M$ into a *-algebra, because the Frobenius property \eref{frob} ensures
that the adjoint is an anti-multiplicative involution. 

We have now constructed $M$ as a *-algebra. To see that it is in fact
a von Neumann algebra, one has to induce the weak topology from $N$ to
$M$ with the help of the faithful conditional expectation 
$\mu:M\to N$ given by
$$\mu(m)=d_\BA\inv\cdot
w^*\ol\iota(m)w,\qquad  \mu(vv^*)=d_\BA\inv\cdot \eins_N.
$$
Here,  
$$\ol\iota(\iota(n)v):= \theta(n)x$$
defines a conjugate homomorphism $\ol\iota:M\to N$, with $(w,v)$ as 
a standard solution of the conjugacy relations. $M$ is already weakly 
closed with respect to the induced topology because it is finitely
generated from $N$.  
\qed

\begin{remark} \label{r:charged} 
It may be convenient to consider, rather than the single generator
$v$ of the extension, 
the system of generators $\psi_\rho=\iota(w_\rho^*)v$  
(``charged intertwiners''), where $\rho\prec \theta$ is an irreducible
sub-endomorphism, and $w_\rho\in\Hom(\rho,\theta)$. By definition,
$\psi_\rho\in\Hom(\iota,\iota\circ\rho)$ which is equivalent to the
commutation relations (suppressing the embedding map $\iota$)
\be\label{charged-int}
\psi_\rho n = \rho(n)\psi_\rho \qquad (n\in N).
\ee
Every element of $M$ has an expansion $\sum n_\rho \psi_\rho$
into a basis of charged intertwiners with coefficients in $N$. The
Q-system controls the product and adjoint of charged intertwiners. 
\end{remark}

In the sequel, we shall always use Q-systems to characterize 
extensions $N\subset M$ of a given factor $N$. In particular, 
all properties of the embedding are encoded in the Q-system, see also
\sref{s:Qsystem}.

\begin{lemma} \label{l:irred}
For $\iota:N\to M$, the following are equivalent: \\
{\rm (i)} The extension is irreducible: $\iota(N)'\cap M=\CC\cdot \eins_M$; 
\\
{\rm (ii)} $\iota:N\to M$ is irreducible: $\Hom(\iota,\iota)=\CC\cdot \eins_M$; \\
{\rm (iii)} $\dim\Hom(\id_N,\bar\iota\iota)=1$. \\
Accordingly, we call a Q-system {\bf irreducible} iff
$\dim\Hom(\id_N,\theta)=1$. 
\end{lemma}

\begin{graybox}
\begin{example} \vskip-5mm \label{x:IQ} (Q-systems of the Ising category) 

The Ising category (cf.\ \xref{x:Icat}) has two irreducible
Q-systems: $(\id,1,1)$ with $M=N$, and $(\theta=\sigma^2,w=2^{\frac14} 
r,x=2^{\frac14}\sig(r)=2^{-\frac14}(r+t))$. In the latter case, the 
extension is $M=\iota(N)\vee\psi$, where $\psi=2^{\frac14}\iota(t^*)v$ 
satisfies the relations $\psi \iota(n)=\iota(\tau(n))\psi$, 
$\psi^*=\psi$, $\psi^2=1$. $M$ has an automorphism 
(fixing $N$, = gauge transformation) $\alpha:\psi\mapsto-\psi$. 
The conjugate $\ol\iota$ in the latter case takes
$\iota(n)$ to $\theta(n)=\sigma^2(n)$ and $\psi$ to
$\sigma^2(t^*)(r+t)=rt^*+tr^*$.
\end{example}
\end{graybox}

For an irreducible Q-system, $M$ is automatically a factor, because
$M'\cap M\subset \iota(N)'\cap M$. However, when $\Hom(\id_N,\theta)$
is more than one-dimensional, then $M$ may have a nontrivial centre,
as characterized by (ii) of the following Lemma.

\begin{definition}\label{d:simple}
We call the Q-system {\bf simple}%
\footnote{The term {\bf factorial} might be more appropriate in this
  context. ``Simple'', however, is more in line with standard category
  terminology, cf.\ \cref{c:simple=factorial}.},  
if the von Neumann algebra $M$ in \tref{t:Q=ext} is a factor. 
\end{definition}
We shall see the equivalence of this definition with the usual one
in \cref{c:simple=factorial} below. 

In the sequel, we give various characterizations of the relative
commutant $N'\cap M$ and of the centre of $M$.

\begin{lemma} \label{l:relc+cent} 
{\rm (i)} The relative commutant $N'\cap M$
is given by the elements $\iota(q)v$, $q\in\Hom(\theta,\id_N)$. \\
{\rm (ii)} $\iota(q)v$ is idempotent iff $(q\times q)\scirc x = q$: 
$\tikzmatht{
  \fill[black!10] (-1.3,-.8) rectangle (1.3,1.2);
    \draw[thick] (-.6,.8)--(-.6,.5) arc(180:360:.6)--(.6,.8);
  \fill[white] (-.6,1)--(-.8,.7)--(-.4,.7)--(-.6,1);
    \draw[thick] (-.6,1)--(-.8,.7)--(-.4,.7)--(-.6,1);
  \fill[white] (.6,1)--(.8,.7)--(.4,.7)--(.6,1);
    \draw[thick] (.6,1)--(.8,.7)--(.4,.7)--(.6,1);
    \draw[thick] (0,-.8)--(0,-.1);
  \fill[black] (0,-.1) circle (.12);
} = 
\tikzmatht{
  \fill[black!10] (-1.3,-.8) rectangle (1.3,1.2);
  \fill[white] (0,.8)--(-.2,.5)--(.2,.5)--(0,.8);
    \draw[thick] (0,.8)--(-.2,.5)--(.2,.5)--(0,.8);
    \draw[thick] (0,-.8)--(0,.5);
    \node at (.5,.7) {$q$};
}$, and it is selfadjoint iff $q^*=(1_\theta\times q)\scirc x\scirc w$: 
$\tikzmatht{
  \fill[black!10] (-1.3,-.8) rectangle (1.3,1.2);
    \draw[thick] (-.6,1.2)--(-.6,.2) arc(180:360:.6)--(.6,.8);
  \fill[white] (.6,.8)--(.8,.5)--(.4,.5)--(.6,.8);
    \draw[thick] (.6,.8)--(.8,.5)--(.4,.5)--(.6,.8);
\node at (1,.7) {$q$};
} = 
\tikzmatht{
  \fill[black!10] (-1.3,-.8) rectangle (1.3,1.2);
  \fill[white] (0,-.4)--(-.2,-.1)--(.2,-.1)--(0,-.4);
    \draw[thick] (0,-.4)--(-.2,-.1)--(.2,-.1)--(0,-.4);
    \draw[thick] (0,1.2)--(0,-.1);
    \node at (.5,0) {$q$};
}$.\\
{\rm (iii)} The centre of $M$ is given by the elements $\iota(q)v$, where
$q$ belongs to the subspace of $\Hom(\theta,\id_N)$ of elements
satisfying 
\be \label{hom-c}
(q\times 1_\theta)\scirc x=(1_\theta\times q)\scirc x: \qquad 
\tikzmatht{
  \fill[black!10] (-1.3,-1) rectangle (1.3,1.5);
    \draw[thick] (.6,1.5)--(.6,.2) arc(360:180:.6)--(-.6,.8);
  \fill[white] (-.6,1)--(-.8,.7)--(-.4,.7)--(-.6,1);
    \draw[thick] (-.6,1)--(-.8,.7)--(-.4,.7)--(-.6,1);
    \draw[thick] (0,-1)--(0,-.4); 
  \fill[black] (0,-.4) circle (.12);
    \node at (-1,1.1) {$q$};
} = 
\tikzmatht{
  \fill[black!10] (-1.3,-1) rectangle (1.3,1.5);
    \draw[thick] (-.6,1.5)--(-.6,.2) arc(180:360:.6)--(.6,.8);
  \fill[white] (.6,1)--(.8,.7)--(.4,.7)--(.6,1);
    \draw[thick] (.6,1)--(.8,.7)--(.4,.7)--(.6,1);
    \draw[thick] (0,-1)--(0,-.4);
  \fill[black] (0,-.4) circle (.12);
    \node at (1,1.1) {$q$};
}
\ee
In particular, the central projections are given by $\iota(q)v$ where
$q\in\Hom(\theta,\id_N)$ satisfies all the relations in {\rm (ii)} and 
{\rm (iii)}. 
\end{lemma}

{\em Proof:} We use the uniqueness of the representation $m=\iota(n)v$
for all three statements. Thus we write $c=\iota(q)v$ and characterize
the properties of $c$ in terms of $q$: \\ 
(i) For $c\in \iota(N)'\cap M$, the 
commutation relation $c\iota(n)=\iota(n)c$ reads $\iota(q\theta(n))v=
\iota(nq)v$. This is equivalent to $q\theta(n) = nq$. \\ 
(ii) 
Immediate from $(\iota(q)v)^2=\iota(q\theta(q)x)v$ and $(\iota(q)v)^* =
\iota(w^*x^*\theta(q^*))v$. \\ 
(iii) The commutation relation $cv=vc$ for $c\in M'\cap M$ reads 
$\iota(qx)v = \iota(\theta(q)x)v$, hence  $qx = \theta(q)x$. 
\qed

\begin{lemma} \label{l:charact}
{\rm (i)} The linear maps $\Hom(\theta,\id_N)\to\Hom(\theta,\theta)$, 
$\tikzmatht{
  \fill[black!10] (-.7,-.8) rectangle (.9,1.2);
  \fill[white] (0,.8)--(-.2,.5)--(.2,.5)--(0,.8);
    \draw[thick] (0,.8)--(-.2,.5)--(.2,.5)--(0,.8);
    \draw[thick] (0,-.8)--(0,.5);
    \node at (.5,.7) {$q$};
} \mapsto 
\tikzmatht{
  \fill[black!10] (-1,-.8) rectangle (1.3,1.2);
    \draw[thick] (-.6,1.2)--(-.6,.2) arc(180:360:.6)--(.6,.5);
  \fill[white] (.6,.8)--(.8,.5)--(.4,.5)--(.6,.8);
    \draw[thick] (.6,.8)--(.8,.5)--(.4,.5)--(.6,.8);
    \draw[thick] (0,-.8)--(0,-.4);
  \fill[black] (0,-.4) circle (.12);
    \node at (1,.8) {$q$};
}$, and 
$\Hom(\theta,\theta)\to\Hom(\theta,\id_N)$, 
$\tikzmatht{
  \fill[black!10] (-.8,-.8) rectangle (.8,1.2);
    \draw[thick] (0,-.8)--(0,1.2);
  \fill[white] (-.3,-.1) rectangle (.3,.5);
    \draw[thick] (-.3,-.1) rectangle node{$t$} (.3,.5);
} \mapsto 
\tikzmatht{
  \fill[black!10] (-.8,-.8) rectangle (.8,1.2);
    \draw[thick] (0,-.8)--(0,.5);
  \fill[white] (-.3,-.4) rectangle (.3,.2);
    \draw[thick] (-.3,-.4) rectangle node{$t$} (.3,.2);
  \fill[white] (0,.7) circle(.2);
    \draw[thick] (0,.7) circle(.2);
}
$,
define a bijection between 
$\Hom(\theta,\id_N)$ and the subspace of $\Hom(\theta,\theta)$ of
elements satisfying the first of \eref{1to3}:\\ 
\be\label{char1} 
\tikzmatht{
  \fill[black!10] (-1.3,-1.3) rectangle (1.3,1.3);
    \draw[thick] (-.8,1.3)--(-.8,.3) arc(180:360:.8)--(.8,1.3);
    \draw[thick] (0,-1.3)--(0,-.5); 
  \fill[white] (1.1,.2) rectangle (.5,.6);
    \draw[thick] (1.1,.2) rectangle (.5,.6);
  \fill[black] (0,-.5) circle (.12);
}=
\tikzmatht{
  \fill[black!10] (-1.3,-1.3) rectangle (1.3,1.3);
    \draw[thick] (-.8,1.3)--(-.8,.8) arc(180:360:.8)--(.8,1.3);
    \draw[thick] (0,-1.3)--(0,0); 
  \fill[white] (-.3,-.8) rectangle (.3,-.4);
    \draw[thick] (-.3,-.8) rectangle (.3,-.4);
  \fill[black] (0,0) circle (.12);
}.
\ee
\noindent {\rm (ii)} $q\in\Hom(\theta,\id_N)$ satisfies \eref{hom-c} iff
$t\in\Hom(\theta,\theta)$ (its image under the bijection  
in {\rm (i)}) satisfies also
\be \label{char2} 
\tikzmatht{
  \fill[black!10] (-1.3,-1.3) rectangle (1.3,1.3);
    \draw[thick] (-.8,1.3)--(-.8,.3) arc(180:360:.8)--(.8,1.3);
    \draw[thick] (0,-1.3)--(0,-.5); 
  \fill[white] (-1.1,.2) rectangle (-.5,.6);
    \draw[thick] (-1.1,.2) rectangle (-.5,.6);
  \fill[black] (0,-.5) circle (.12);
}=
\tikzmatht{
  \fill[black!10] (-1.3,-1.3) rectangle (1.3,1.3);
    \draw[thick] (-.8,1.3)--(-.8,.3) arc(180:360:.8)--(.8,1.3);
    \draw[thick] (0,-1.3)--(0,-.5); 
  \fill[white] (1.1,.2) rectangle (.5,.6);
    \draw[thick] (1.1,.2) rectangle (.5,.6);
  \fill[black] (0,-.5) circle (.12);
},\ee
i.e., iff $t\in \Hom_0(\theta,\theta)$.
\end{lemma}

{\em Proof:} (i) $\theta(q)x$ satisfies \eref{char1}: By associativity
$\tikzmatht{
  \fill[black!10] (-1.2,-1) rectangle (1.4,1.5);
    \draw[thick] (1,1)--(1,.8) arc(360:180:.6)--
          (-.2,1.5);
    \draw[thick] (.4,.2) arc(360:180:.7)--(-1,1.5);
  \fill[black] (-.3,-.5) circle(.12);
  \fill[black] (.4,.2) circle(.12);
    \draw[thick] (-.3,-1)--(-.3,-.5);
  \fill[white] (1,1.3)--(1.2,1)--(.8,1)--(1,1.3);
    \draw[thick] (1,1.3)--(1.2,1)--(.8,1)--(1,1.3);
} = 
\tikzmatht{
  \fill[black!10] (-1.2,-1) rectangle (1.4,1.5);
    \draw[thick] (-1,1.5)--(-1,.8) arc(180:360:.6)--
          (.2,1.5);
    \draw[thick] (-.4,.2) arc(180:360:.7)--(1,.9);
  \fill[black] (.3,-.5) circle(.12);
  \fill[black] (-.4,.2) circle(.12);
    \draw[thick] (.3,-1)--(.3,-.5);
  \fill[white] (1,.9)--(1.2,.6)--(.8,.6)--(1,.9);
    \draw[thick] (1,.9)--(1.2,.6)--(.8,.6)--(1,.9);
}$. 

The two maps invert each other: 
$$\tikzmatht{
  \fill[black!10] (-1,-1) rectangle (1,1.2); 
    \draw[thick] (0,-1)--(0,.5);
  \fill[white] (0,.8)--(.2,.5)--(-.2,.5)--(0,.8);
    \draw[thick] (0,.8)--(.2,.5)--(-.2,.5)--(0,.8);         
} \mapsto           
\tikzmatht{
  \fill[black!10] (-1.2,-1) rectangle (1.2,1.2);
    \draw[thick] (0,-1)--(0,-.2);
    \draw[thick] (-.5,1.2)--(-.5,.3) arc(180:360:.5)--(.5,.5);
  \fill[black] (0,-.2) circle(.12);
  \fill[white] (.5,.8)--(.7,.5)--(.3,.5)--(.5,.8);
    \draw[thick] (.5,.8)--(.7,.5)--(.3,.5)--(.5,.8);
} \mapsto           
\tikzmatht{
  \fill[black!10] (-1.2,-1) rectangle (1.2,1.2);
    \draw[thick] (0,-1)--(0,-.2);
    \draw[thick] (-.5,.7)--(-.5,.3) arc(180:360:.5)--(.5,.5);
  \fill[black] (0,-.2) circle(.12);
  \fill[white] (.5,.8)--(.7,.5)--(.3,.5)--(.5,.8);
    \draw[thick] (.5,.8)--(.7,.5)--(.3,.5)--(.5,.8);
  \fill[white] (-.5,.7) circle(.2);
    \draw[thick] (-.5,.7) circle(.2);
} = 
\tikzmatht{
  \fill[black!10] (-1,-1) rectangle (1,1.2);
    \draw[thick] (0,-1)--(0,.5);
  \fill[white] (0,.8)--(.2,.5)--(-.2,.5)--(0,.8);
    \draw[thick] (0,.8)--(.2,.5)--(-.2,.5)--(0,.8);         
}, $$
and
$$
\tikzmatht{
  \fill[black!10] (-1,-1.2) rectangle (1,1.2);
    \draw[thick] (0,1.2)--(0,-1.2);
  \fill[white] (-.3,-.2) rectangle (.3,.2);
    \draw[thick] (-.3,-.2) rectangle (.3,.2);
} \mapsto 
\tikzmatht{
  \fill[black!10] (-1.2,-1.2) rectangle (1.4,1.2);
    \draw[thick] (0,.7)--(0,-1.2);
  \fill[white] (-.3,-.6) rectangle (.3,-.2);
    \draw[thick] (-.3,-.6) rectangle (.3,-.2);
  \fill[white] (0,.7) circle(.2);
    \draw[thick] (0,.7) circle(.2);
} \mapsto 
\tikzmatht{
  \fill[black!10] (-1.2,-1.2) rectangle (1.4,1.2);
    \draw[thick] (0,-1.2)--(0,-.7);
    \draw[thick] (-.5,1.2)--(-.5,-.2) arc(180:360:.5)--(.5,.7);
  \fill[black] (0,-.7) circle(.12);
  \fill[white] (.2,-.2) rectangle (.8,.2);
    \draw[thick] (.2,-.2) rectangle (.8,.2);
  \fill[white] (.5,.7) circle(.2);
    \draw[thick] (.5,.7) circle(.2);
} 
\stapel{\nref{char1}}= 
\tikzmatht{
  \fill[black!10] (-1.2,-1.2) rectangle (1.4,1.2); 
    \draw[thick] (0,-.2)--(0,-1.2);
    \draw[thick] (-.5,1.2)--(-.5,.3) arc(180:360:.5)--(.5,.7);
  \fill[black] (0,-.2) circle(.12);
  \fill[white] (-.3,-.9) rectangle (.3,-.5);
    \draw[thick] (-.3,-.9) rectangle (.3,-.5);
  \fill[white] (.5,.7) circle(.2);
    \draw[thick] (.5,.7) circle(.2);
} = 
\tikzmatht{
  \fill[black!10] (-1,-1.2) rectangle (1,1.2);
    \draw[thick] (0,1.2)--(0,-1.2);
  \fill[white] (-.3,-.2) rectangle (.3,.2);
    \draw[thick] (-.3,-.2) rectangle (.3,.2);} .$$
(ii) ``If'': 
$\tikzmatht{
  \fill[black!10] (-1.2,-1.2) rectangle (1.2,1.2);
    \draw[thick] (0,-.4)--(0,-1.2) (.7,1.2)--(.7,.3) arc(360:180:.7)--(-.7,.5);
  \fill[white] (-.7,.7)--(-.9,.4)--(-.5,.4)--(-.7,.7);
    \draw[thick] (-.7,.7)--(-.9,.4)--(-.5,.4)--(-.7,.7);
  \fill[black] (0,-.4) circle(.12); 
}
:=\tikzmatht{
  \fill[black!10] (-1.2,-1.2) rectangle (1.2,1.2);
    \draw[thick] (0,-.6)--(0,-1.2) (.7,1.2)--(.7,.1)
          arc(360:180:.7) (-.7,.3)--(-.7,.6);
  \fill[white] (-.7,.8) circle(.2);
    \draw[thick] (-.7,.8) circle(.2);
  \fill[white] (-.4,.3) rectangle (-1,-.1);
    \draw[thick] (-.4,.3) rectangle (-1,-.1);
  \fill[black] (0,-.6) circle(.12); 
}
\stapel{\nref{char2}}= 
\tikzmatht{
  \fill[black!10] (-1.2,-1.2) rectangle (1.4,1.2);
    \draw[thick] (0,-.4)--(0,-1.2) (.7,1.2)--(.7,.3) arc(360:180:.7);
  \fill[white] (-.7,.5) circle(.2);
    \draw[thick] (-.7,.5) circle(.2);
  \fill[white] (.4,.6) rectangle (1,.2);
    \draw[thick] (.4,.6) rectangle (1,.2);
  \fill[black] (0,-.4) circle(.12);          
}
=
\tikzmatht{
  \fill[black!10] (-1,-1.2) rectangle (1,1.2);
    \draw[thick] (0,1.2)--(0,-1.2);
  \fill[white] (-.3,-.2) rectangle (.3,.2);
    \draw[thick] (-.3,-.2) rectangle (.3,.2);
}
=:
\tikzmatht{
  \fill[black!10] (-1.2,-1.2) rectangle (1.2,1.2);
    \draw[thick] (0,-.4)--(0,-1.2) (-.7,1.2)--(-.7,.3) arc(180:360:.7)--(.7,.5);
  \fill[white] (.7,.7)--(.9,.4)--(.5,.4)--(.7,.7);
    \draw[thick] (.7,.7)--(.9,.4)--(.5,.4)--(.7,.7);
  \fill[black] (0,-.4) circle(.12);           
}
$ by the unit property. 

``Only if'': $
\tikzmatht{
  \fill[black!10] (-1.2,-1.2) rectangle (1.4,1.2);
    \draw[thick] (0,-.4)--(0,-1.2) (.7,1.2)--(.7,.3)
          arc(360:180:.7)--(-.7,1.2);
  \fill[white] (-.4,.6) rectangle (-1,.2);
    \draw[thick] (-.4,.6) rectangle (-1,.2);
  \fill[black] (0,-.4) circle(.12);           
}:=
\tikzmatht{
  \fill[black!10] (-1.2,-1) rectangle (1.2,1.4);
    \draw[thick] (-1,1.4)--(-1,.8) arc(180:360:.6)--
          (.2,1);
    \draw[thick] (-.4,.2) arc(180:360:.7)--(1,1.4);          
  \fill[black] (.3,-.5) circle(.12);
  \fill[black] (-.4,.2) circle(.12);
    \draw[thick] (.3,-1)--(.3,-.5);
  \fill[white] (.2,1.1)--(0,.8)--(.4,.8)--(.2,1.1);
    \draw[thick] (.2,1.1)--(0,.8)--(.4,.8)--(.2,1.1);
} = 
\tikzmatht{
  \fill[black!10] (-1.2,-1) rectangle (1.2,1.4);
    \draw[thick] (1,1.4)--(1,.8) arc(360:180:.6)--(-.2,1);
    \draw[thick] (.4,.2) arc(360:180:.7)--(-1,1.4);
  \fill[black] (-.3,-.5) circle(.12);
  \fill[black] (.4,.2) circle(.12);
    \draw[thick] (-.3,-1)--(-.3,-.5);
  \fill[white] (-.2,1.1)--(0,.8)--(-.4,.8)--(-.2,1.1);
    \draw[thick] (-.2,1.1)--(0,.8)--(-.4,.8)--(-.2,1.1);
} \stapel{\nref{hom-c}}= 
\tikzmatht{
  \fill[black!10] (-1.2,-1) rectangle (1.4,1.4);
    \draw[thick] (1,1)--(1,.8) arc(360:180:.6)--(-.2,1.4);
    \draw[thick] (.4,.2) arc(360:180:.7)--(-1,1.4);
  \fill[black] (-.3,-.5) circle(.12);
  \fill[black] (.4,.2) circle(.12);
    \draw[thick] (-.3,-1)--(-.3,-.5);
  \fill[white] (1,1.1)--(1.2,.8)--(.8,.8)--(1,1.1);
    \draw[thick] (1,1.1)--(1.2,.8)--(.8,.8)--(1,1.1);
}=:\tikzmatht{
  \fill[black!10] (-1,-1.2) rectangle (1.2,1.2);
    \draw[thick] (0,-.4)--(0,-1.2) (.7,1.2)--(.7,.3)
          arc(360:180:.7)--(-.7,1.2);
  \fill[white] (.4,.6) rectangle (1,.2);
    \draw[thick] (.4,.6) rectangle (1,.2);
  \fill[black] (0,-.4) circle(.12);           
}$ by associativity. 
\qed

Thus, the relative commutant $N'\cap M$ and the centre $M'\cap M$ are
equivalently characterized by certain elements of $\Hom(\theta,\id_N)$ or of
$\Hom(\theta,\theta)$. In particular, the space
$\Hom_0(\theta,\theta)$, \dref{d:hom0}, is one way to characterize the centre of $M$. We shall come back to this in 
\sref{s:Qcentral} and \sref{s:Qirred}.

\begin{remark} \label{r:factor} 
The standardness property of the Q-system is not used in the
construction of the algebra $M$ in the proof of \tref{t:Q=ext}, 
and neither the (weaker) specialness property that $x^*\ncirc x$ 
is a multiple of $1_\theta$. These properties are only required to ensure 
that $v^*\ncirc v$ is a multiple of $\eins_M$, namely $x^* x=\ol\iota(v^*v)$.
Because $M'\cap M=\Hom(\iota\ol\iota,\iota\ol\iota)$, $v^*v$ is always
central in $M$, hence specialness is automatically satisfied if $M$ 
is a factor. 
\end{remark}

\subsection{The canonical Q-system} 
\label{s:canonical}
\setcounter{equation}{0}
Let $j:N\to j(N)$ be an antilinear isomorphism of factors. 
E.g., $j:n\mapsto n^*$ is an antilinear isomorphism of $N$ with
$j(N)=N\opp$ (the algebra with the opposite product), or a Tomita
conjugation $j=\Ad_J$ is an antilinear isomorphism of $N$ with
$j(N)=N'$. For $\C\subset \End_0(N)$, let $j(\C)$ the category with objects
$\rho^j \equiv j\circ\rho\circ j\inv$ ($\rho\in\C$) and with intertwiners
$j(t)$.  

We denote by $\C_1\boxtimes\C_2$ (the Deligne product) the completion of
the tensor product of categories $\C_1\otimes \C_2$ by direct sums. 

\begin{tintedbox}
\begin{proposition} \vskip-5mm \label{p:canonical} \cite{LR95} 
If $\C$ has only finitely many inequivalent irreducible objects
$\rho$, then there is a canonical irreducible Q-system $\BR$ 
in $\C\boxtimes j(\C)$ with 
$$[\Theta\can] = \bigoplus_\rho[\rho]\otimes[\rho^j],$$
where the sum runs over the equivalence classes of irreducible objects
of $\C$. Its dimension is given by $d_\BR^2=\dim(\Theta\can)=\dim(\C)$.  
Choosing isometries $T_\rho\in\Hom(\rho\otimes\rho^j,\Theta\can)$, the
Q-system is given by   
$$
W = d_\BR^{\frac12}\cdot T_\id,\qquad X=d_\BR^{-\frac12}\sum_{\rho,\sig,\tau}
\Big(\frac{d_\rho d_\sig}{d_\tau}\Big)^{\frac12}\cdot(T_\rho\times T_\sig)
\scirc \Big(\sum_a t_a\otimes j(t_{a})\Big) \scirc T_\tau^*,
$$
where the first sum extends over representatives of all sectors, and
the inner sums over $a$ extend over orthonormal bases of isometries
$t_a\in\Hom(\tau,\rho\sig)$.
\end{proposition}
\end{tintedbox}

Because of the anti-linearity of $j$, the sums over $a$ do not depend
on the choice of orthonormal bases $t_a$. A different choice of
$T_\rho$ gives a unitarily equivalent Q-system. 

One easily proves (cf.\ \pref{p:dotstar})

\begin{lemma} \cite{LRo} Choosing, for every $\rho\in\C$, a conjugate
  $\ol\rho\in\C$ and a standard pair $(w,\ol w)$, the assignment 
$$\rho\mapsto\ol\rho,\quad t\mapsto j\Big(\tikzmatht{
  \fill[black!10] (-1.5,-1.8) rectangle (1.5,.7);
    \draw[thick] (0,-.2) arc(180:0:.4); 
  \fill[white] (-.4,-.2) rectangle (.4,-.9); 
    \draw[thick] (-.4,-.2) rectangle node {$t^*$} (.4,-.9); 
    \draw[thick] (0,-.9) arc(360:180:.4)--(-.8,.7)
    (.8,-.2)--(.8,-1.8); 
    \node at (1.1,.3) {$\ol w_\rho^*$};
    \node at (-1,-1.4) {$w_\sig$};
}\Big) = j\Big(
\tikzmatht{
  \fill[black!10] (-1.5,-.7) rectangle (1.5,1.8);
    \draw[thick] (0,.2) arc(180:360:.4); 
  \fill[white] (-.4,.2) rectangle (.4,.9); 
    \draw[thick] (-.4,.2) rectangle node {$t^*$} (.4,.9); 
    \draw[thick] (0,.9) arc(0:180:.4)--(-.8,-.7)
    (.8,.2)--(.8,1.8); 
    \node at (1.1,-.3) {$\ol w_\sig$};
    \node at (-1,1.4) {$w_\rho^*$};
}\Big)
$$  
taking $\Hom(\rho,\sig)$ into $\Hom(\ol\rho^j,\ol\sig^j)$ is a
linear isomorphism between the C* tensor categories $\C\opp$ and $j(\C)$  
(the category equipped with the opposite monoidal product).
\end{lemma}

\begin{corollary} \label{c:opp}
The opposite tensor category $\C\opp$ can be realized
as $j(\C)\subset \End_0(N\opp)$ or $\End_0(N')$. Under this
isomorphism, the canonical Q-system in $\C\boxtimes j(\C)$ becomes a
Q-system in $\C\boxtimes \C\opp$ with $[\Theta]=\bigoplus
[\rho]\otimes[\ol\rho]$.
\end{corollary}

This is the way it is defined in the abstract setting (e.g., \cite[II,
Prop.\ 4.1]{M03}).

\subsection{Modules of Q-systems}
\label{s:modules}
\setcounter{equation}{0}
A {\bf module} ($\equiv$ left module) of a Q-system $\BA=(\theta,w,x)$
is a pair $\mm=(\beta,m)$, where $\beta$ is an object of the
underlying category and $m\in\Hom(\beta,\theta\beta)$,%
\footnote{More precisely, $(\beta,m^*)$ is a module and
  $(\beta,m)$ is a co-module. We do not make the distinction because
  the dualization is canonically given by the operator adjoint.} 
satisfying the relations 
\bea
\hbox{\bf unit property:} && 
(w^* \times 1_\beta)\scirc m = 1_\beta
\notag \\ \label{m-unit} &&
\tikzmatht{
  \fill[black!10] (-1,-1) rectangle (1,1); 
    \draw[thick] (.5,1)--(.5,-1); 
    \draw[thick] (-.5,.6)--(-.5,0) arc(180:270:.5)--(.33,-.5);
  \fill[white] (-.5,.6) circle(.2);
    \draw[thick] (-.5,.6) circle(.2);          
    \draw[thick] (.5,-.33) arc(90:270:.17); 
} = 
\tikzmatht{
  \fill[black!10] (0,-1) rectangle (1,1); 
    \draw[thick] (.5,1)--(.5,-1); 
} ,
\\ \hbox{\bf representation property:} && 
(1_\theta\times m)\scirc m = (x \times 1_\beta)\scirc m 
\notag \\ \label{repn} &&
\tikzmatht{
  \fill[black!10] (-1,-1) rectangle (1,1); 
    \draw[thick] (.7,1)--(.7,-1); 
    \draw[thick] (-.6,1)--(-.6,.3) arc(180:270:.8)--(.53,-.5);
    \draw[thick] (0,1)--(0,.5) arc(180:270:.4)--(.53,.1);
    \draw[thick] (.7,-.33) arc(90:270:.17); 
    \draw[thick] (.7,.27) arc(90:270:.17); 
} =
\tikzmatht{
  \fill[black!10] (-1,-1) rectangle (1,1); 
    \draw[thick] (.5,1)--(.5,-1); 
    \draw[thick] (-.5,.2)--(-.5,0) arc(180:270:.5)--
          (.33,-.5);
    \draw[thick] (-.8,1)--(-.8,.5) arc(180:360:.3)--
          (-.2,1);
  \fill[black] (-.5,.2) circle(.12);
    \draw[thick] (.5,-.33) arc(90:270:.17); 
}.
\eea

A module of a Q-system is called a {\bf standard module} if
$m^*\ncirc m$ is a multiple of $1_\beta$. (This property is automatic 
if $(\beta,m)$ is irreducible as a module, and in particular if
$\beta$ is irreducible as an endomorphism.)

A Q-system $\BA$ is also a standard $\BA$-module
$(\beta=\theta,m=x)$. 

Two modules $(\beta,m)$ and $(\beta',m')$ are {\bf equivalent}, when there
is an invertible $n\in\Hom(\beta,\beta')$ such that $m'\scirc n =
(1_\theta\times n)\scirc m$. They are unitarily equivalent if there is
a unitary such $n$.

\begin{lemma} \label{l:normmod}
{\rm (i)} If a module $\mm=(\beta,m)$ is standard, then (with
$d_\BA=$ the dimension of the Q-system) 
\be\label{modnorm}
m^*\ncirc m=
\tikzmatht{
  \fill[black!10] (-1,-1.3) rectangle (1,1.3); 
    \draw[thick] (.5,1.3)--(.5,-1.3);
    \draw[thick] (-.4,-.2)--(-.4,.2) arc(180:90:.5)--(.33,.7);
    \draw[thick] (-.4,-.2) arc(180:270:.4)--(.33,-.7);
    \draw[thick] (.5,.87) arc(90:270:.17);
    \draw[thick] (.5,-.53) arc(90:270:.17); 
} 
= d_\BA\cdot 1_\beta.
\ee
\noindent {\rm (ii)} Every module is equivalent to a standard module, i.e., there is
an invertible element $n$ of $\Hom(\beta,\beta)$ such that 
$(\beta,(1_\theta\times n)\scirc m \scirc n\inv)$ is a standard module.  
\end{lemma}

{\em Proof:} (i) follows from the representation property \eref{repn}
and $x^*\ncirc x=d_\BA\cdot 1_\theta$, which imply
$$m^*\scirc(1_\theta\times m^*\ncirc m)\scirc m =
\tikzmatht{
  \fill[black!10] (-1,-1.3) rectangle (1,1.3); 
    \draw[thick] (.5,1.3)--(.5,-1.3);
    \draw[thick] (-.7,-.5)--(-.7,.5) arc(180:90:.5)--(.33,1);
    \draw[thick] (-.4,-.2)--(-.4,.2) arc(180:90:.4)--(.33,.6);
    \draw[thick] (-.4,-.2) arc(180:270:.4)--(.33,-.6);
    \draw[thick] (-.7,-.5) arc(180:270:.5)--(.33,-1);
    \draw[thick] (.5,1.17) arc(90:270:.17); 
    \draw[thick] (.5,.77) arc(90:270:.17); 
    \draw[thick] (.5,-.43) arc(90:270:.17); 
    \draw[thick] (.5,-.83) arc(90:270:.17); 
}=
\tikzmatht{
  \fill[black!10] (-1,-1.3) rectangle (1,1.3); 
\draw[thick] (.5,1.3)--(.5,-1.3);
    \draw[thick] (-.4,.5) arc(180:90:.5)--(.33,1);
    \draw[thick] (-.4,-.5) arc(180:270:.5)--(.33,-1);
    \draw[thick] (-.4,0) circle(.5);
    \draw[thick] (.5,1.17) arc(90:270:.17);
    \draw[thick] (.5,-.83) arc(90:270:.17); 
  \fill[black] (-.4,.5) circle(.12);
  \fill[black] (-.4,-.5) circle(.12);
}=d_\BA\cdot 
\tikzmatht{
  \fill[black!10] (-1,-1.3) rectangle (1,1.3); 
    \draw[thick] (.5,1.3)--(.5,-1.3);
    \draw[thick] (-.4,-.2)--(-.4,.2) arc(180:90:.5)--(.33,.7);
    \draw[thick] (-.4,-.2) arc(180:270:.4)--(.33,-.7);
    \draw[thick] (.5,.87) arc(90:270:.17);
    \draw[thick] (.5,-.53) arc(90:270:.17); 
}
= d_\BA\cdot m^*\ncirc m.$$
For (ii), first we notice that $m^*\ncirc m$ is an invertible positive element
of $\Hom(\beta,\beta)$, because $e=d_\BA\inv\cdot w\ncirc w^*$ is a projection
in $\Hom(\theta,\theta)$, hence by the unit property,
$$m^*\ncirc m\geq m^*\scirc (e\times 1_\beta)\scirc m = d_\BA\inv \cdot 1_\beta.$$
Let $n\in \Hom(\beta,\beta)$ be the square root of $m^*\ncirc m$. Then by the 
representation property,
$$n\inv m^*\scirc (1_\theta\times n^2)\scirc mn\inv = n\inv m^*\scirc
(x^*\ncirc x\times 1_\beta)\scirc mn\inv = d_\BA\cdot 1_\beta.$$
\qed

\begin{lemma} \label{l:conjmod}
If $(\beta,m)$ is a standard module, then in addition to
the unit and representation relations, the relation  
\bea\label{standard-mod} (x^*\times 1_\beta)\scirc(1_\theta\times
m) = m\ncirc m^* = (1_\theta\times m^*)\scirc(x^*\times 1_\beta): \;\; 
\tikzmatht{
  \fill[black!10] (-1,-1.2) rectangle (1,1.2); 
    \draw[thick] (.7,1.2)--(.7,-1.2); 
    \draw[thick] (-.7,-1.2)--(-.7,0) arc(180:0:.4)
          arc(180:270:.4)--(.53,-.4) (-.3,.4)--(-.3,1.2);
  \fill[black] (-.3,.4) circle(.12);
    \draw[thick] (.7,-.23) arc(90:270:.17); 
} =
\tikzmatht{
  \fill[black!10] (-.8,-1.2) rectangle (.8,1.2); 
    \draw[thick] (.5,1.2)--(.5,-1.2); 
    \draw[thick] (-.5,1.2)--(-.5,.7) arc(180:270:.5)--(.33,.2);
    \draw[thick] (-.5,-1.2)--(-.5,-.7) arc(180:90:.5)--(.33,-.2);
    \draw[thick] (.5,.37) arc(90:270:.17); 
    \draw[thick] (.5,-.03) arc(90:270:.17); 
} =
\tikzmatht{
  \fill[black!10] (-1,-1.2) rectangle (1,1.2); 
    \draw[thick] (.7,1.2)--(.7,-1.2); 
    \draw[thick] (-.7,1.2)--(-.7,0) arc(180:360:.4)
          arc(180:90:.4)--(.53,.4) (-.3,-.4)--(-.3,-1.2);
  \fill[black] (-.3,-.4) circle(.12);
    \draw[thick] (.7,.57) arc(90:270:.17); 
} 
\eea
holds. This implies 
\be\label{mod-conj}
m^*=(r^*\times 1_\beta)\scirc (1_\theta\times m): \qquad 
\tikzmatht{
  \fill[black!10] (-1,-1) rectangle (1,1); 
    \draw[thick] (.5,1)--(.5,-1); 
    \draw[thick] (-.5,-1)--(-.5,0) arc(180:90:.5)--(.33,.5);
    \draw[thick] (.5,.67) arc(90:270:.17); 
} = 
\tikzmatht{
  \fill[black!10] (-1,-1) rectangle (1,1); 
    \draw[thick] (.7,1)--(.7,-1); 
    \draw[thick] (-.7,-1)--(-.7,0) arc(180:0:.4)
          arc(180:270:.4)--(.53,-.4);
    \draw[thick] (.7,-.23) arc(90:270:.17); 
},
\ee
and consequently
$$E:=d_\BA\inv\cdot (x^*\times 1_\beta)\scirc(1_\theta\times m) =
d_\BA\inv\cdot 
\tikzmatht{
  \fill[black!10] (-1,-1.5) rectangle (1.1,1.5); 
    \draw[thick] (.77,1.5)--(.77,-1.5); 
    \draw[thick] (.6,.2) arc(270:180:.4) arc(0:180:.4)
        --(-.6,-.6) arc(180:360:.4) arc(180:90:.4)--(.6,-.2);
    \draw[thick] (-.2,1)--(-.2,1.5) (-.2,-1)--(-.2,-1.5);
  \fill[black] (-.2,1) circle(.12);
  \fill[black] (-.2,-1) circle(.12);
    \draw[thick] (.77,.37) arc(90:270:.17); 
    \draw[thick] (.77,-.03) arc(90:270:.17); 
}
$$ 
is a self-adjoint idempotent, i.e., a projection in
$\Hom(\theta\beta,\theta\beta)$. 
\end{lemma}

{\em Proof:} The proof is very much the same as the proof of the
Frobenius property in \lref{l:UAS>F}, with $X$ replaced by 
$X':= (1_\theta\times m^*)\scirc (x\times 1_\beta) - m\ncirc m^* \in
\Hom(\theta\beta,\theta\beta)$, the associativity property of $x$
replaced by the representation property of $m$, and $\delta$ replaced
by the faithful positive map
$\delta':\Hom(\theta\beta,\theta\beta)\to \Hom(\theta\beta,\theta\beta)$
$$\delta'(T)=(x^*\times 1_\beta)\scirc (1_\theta\times T)\scirc (x\times
1_\beta): \qquad 
\tikzmatht{
  \fill[black!10] (-1.5,-1.5) rectangle (1.5,1.5);
  \fill[white] (-.5,-.4) rectangle (1,.4); 
    \draw[thick] (-.5,-.4) rectangle node {$T$} (1,.4); 
    \draw[thick] (0,.4)--(0,.6) arc(0:180:.5)--(-1,-.6)
    arc(180:360:.5)--(0,-.4) (-.5,1.1)--(-.5,1.5) (-.5,-1.1)--(-.5,-1.5); 
    \draw[thick] (.5,.4)--(.5,1.5) (.5,-.4)--(.5,-1.5); 
  \fill[black] (-.5,1.1) circle(.12); 
  \fill[black] (-.5,-1.1) circle(.12);
    \node at (.8,1) {$\beta$}; 
    \node at (.8,-1) {$\beta$};  
}
$$
The equation for $m^*$ then follows by left composition with $w^*\times
1_\beta$, and the statement about $E$ follows because $E=d_\BA\inv\cdot
m\ncirc m^*$ and $m^*\ncirc m=d_\BA\cdot 1_\beta$.
\qed

From now on, we reserve the graphical representation
$$m=
\tikzmatht{
\fill[black!10] (-1,-1.5) rectangle (1,1.2); 
    \draw[thick] (.5,1.2)--(.5,-1.5); 
    \draw[thick] (-.5,1.2)--(-.5,0) arc(180:270:.5)--(.33,-.5);
    \draw[thick] (.5,-.33) arc(90:270:.17); 
    \node at (-.2,.8){$\theta$};
    \node at (.8,.8) {$\beta$}; 
}
$$
for the intertwiner associated with a standard module $\mm = (\beta,m)$, i.e., $m$ satisfies
\eref{m-unit}, \eref{repn} and \eref{modnorm}, hence also \eref{mod-conj}. 
We shall freely use these properties in the sequel.

If $\BA$ is a Q-system in $\C=\End_0(N)$, corresponding to an 
extension $\iota:N\to M$, then every homomorphism $\varphi:N\to M$ of 
finite dimension gives rise to a standard module 
\bea \label{zeta}
(\beta,m) \equiv (\ol\iota\varphi,1_{\ol\iota}\times v\times
1_\varphi): \qquad
\tikzmatht{
\fill[black!10] (-1,-1.5) rectangle (1,1.2); 
    \draw[thick] (.5,1.2)--(.5,-1.5); 
    \draw[thick] (-.5,1.2)--(-.5,0) arc(180:270:.5)--(.33,-.5);
    \draw[thick] (.5,-.33) arc(90:270:.17); 
    \node at (-.2,.8){$\theta$};
} \equiv 
\tikzmatht{
  \fill[black!10] (-1.2,-1.2) rectangle (1.5,1.5); 
  \fill[black!18] (-.8,1.5)--(-.8,.5) arc(180:250:.9)
  arc(70:0:.9)--(.4,-1.2)--(.8,-1.2)--(.8,1.5);
    \draw[thick] (-.8,1.5)--(-.8,.5) arc(180:250:.9)
    arc(70:0:.9)--(.4,-1.2) (.8,-1.2)--(.8,1.5);
  \fill[black!10] (-.4,1.5)--(-.4,.4) arc(180:360:.4)--(.4,1.5);
    \draw[thick] (-.4,1.5)--(-.4,.4) arc(180:360:.4)--(.4,1.5);
    \node at (0,.3){$v$};
    \node at (1.2,-.3){$\varphi$};
    \node at (-.1,-.8){$\ol\iota$}; 
    \node at (-.2,1.2){$\iota$};
}\eea
of $\BA$. Notice that, as an operator in $N$, $m=\ol\iota(v)=x$. 
If $\C\subset\End_0(N)$ as specified in the beginning of the section, then the
same is true provided $\ol\iota\varphi$ belongs to $\C$. This
restriction on $\varphi$ is equivalent to the condition that 
$\varphi\prec\iota\rho$ with some $\rho\in\C$. 

The converse is also true: namely, we prove now that every 
standard module is of this form:

\begin{tintedbox}
\begin{proposition} \vskip-5mm \label{p:modulphi}
Every standard module $\mm=(\beta,m)$ of a simple Q-system 
$\BA=(\theta,w,x)$ in $\End_0(N)$ is unitarily equivalent to a 
standard module of the form $(\ol\iota\varphi,x)$ as in \eref{zeta},
where $\varphi$ is a homomorphism $\varphi:N\to M$.
\end{proposition}
\end{tintedbox}

(The same result was derived by \cite[Lemma 3.1]{EP03} by an exhaustion 
argument, using the known number of modules in the case of a 
{\em braided} category; our proof is more constructive, and does
not refer to a braiding.)

\medskip

{\em Proof:}  
Writing as before $\theta=\ol\iota\iota$, $m$ defines by left Frobenius 
conjugation an element 
$$e=d_\BA\inv\cdot (v^*\times 1_{\iota\beta})\scirc (1_\iota\times m)$$ 
of $\Hom(\iota\beta,\iota\beta)\subset M$. Then $1_{\ol\iota}\times
e$ equals, by \eref{standard-mod}, the projection $E=d_\BA\inv\cdot mm^*$ in \lref{l:conjmod}, hence $e$
is also a projection. Let $\varphi\prec \iota\beta$ be the sub-homomorphism 
$:N\to M$ corresponding to this projection, and 
$s\in \Hom(\varphi,\iota\beta)$ an isometry such that $e=ss^*$. 
By left Frobenius conjugation, $\tilde s :=(1_{\ol\iota}\times s^*)\scirc 
(w\times 1_\beta)\in
\Hom(\beta,\ol\iota\varphi)$. We claim that the range projection of 
$\tilde s$ equals $1_{\ol\iota\varphi}$.

Indeed, by inverting the definition of $e$, we have that 
$$m=d_\BA\cdot (1_{\ol\iota}\times s\ncirc s^*)\scirc (w\times 1_\beta),$$ 
hence 
$$\tilde s= d_\BA\inv\cdot (1_{\ol\iota}\times s^*)\scirc m.$$
Now, we use again \eref{standard-mod}: $m\ncirc m^* = d_\BA\cdot 1_{\ol\iota}\times e = d_\BA\cdot
1_{\ol\iota}\times ss^*$ to conclude
$$\tilde s\tilde s^* = d_\BA^{-2}\cdot (1_{\ol\iota}\times s^*)\scirc
m\ncirc m^*\scirc (1_{\ol\iota}\times s^*) = d_\BA\inv\cdot 
(1_{\ol\iota}\times s^*ss^*s) =d_\BA\inv\cdot 1_{\ol\iota\varphi}.$$ 
Thus, while $\varphi\prec \iota\beta$ by construction, we also have 
$\iota\beta\prec\varphi$, hence $\beta$ is equivalent to 
$\ol\iota\varphi$. It follows that $u:=\sqrt{d_\BA}\cdot \tilde s$ is 
a unitary $u\in\Hom(\beta,\ol\iota\varphi)$. Then, inserting 
$s=(1_\iota\times \tilde s^*)\scirc(w\times 1_\varphi)$ into 
$m=d_\BA\cdot (1_{\ol\iota}\times s\ncirc s^*)\scirc(w\times 1_\beta)$, one
arrives at 
$$m=(1_\theta\times u^*)\scirc (1_{\ol\iota}\times v\times 1_\varphi)\scirc u.$$
This proves the claim.
\qed

The homomorphism $\varphi$ corresponding to a module $\mm$ can be
explicitly computed: namely, $\varphi(n)\in M$ can be written as
$\varphi(n)=\iota(k)v$ for some $k\in 
N$. Applying $\ol\iota$ implies $\beta(n)=\theta(k)x$. Multiplying
$w^*$ from the right, implies $w^*\beta(n)=w^*\theta(k)x =
kw^*x=k$. Hence $\varphi:N\to M$ is given by
$$\varphi(n)=\iota(w^*\beta(n))v\in M.$$

Considering $\BA$ as a standard $\BA$-module
$(\beta=\theta,m=x)$, the corresponding homomorphism is
$\varphi=\iota:N\to M$. 

\medskip

The modules of a Q-system $(\theta,w,x)$ are the objects of 
the {\bf module category}. A {\bf morphism} between two modules $(\beta,m)$ 
and $(\beta',m')$ is an element $t\in \Hom(\beta,\beta')$ satisfying 
\bea\label{morph}
(1_\theta\times t)\scirc m = m' \scirc t:\qquad 
\tikzmatht{
  \fill[black!10] (-.8,-1.2) rectangle (1.6,1.2); 
    \draw[thick] (.5,1.2)--(.5,-1.2); 
    \draw[thick] (-.5,1.2)--(-.5,0) arc(180:270:.5)--(.33,-.5);
    \draw[thick] (.5,-.33) arc(90:270:.17); 
  \fill[white] (0,0) rectangle (1,.7);
    \draw[thick] (0,0) rectangle node{$t$}(1,.7);
    \node at (1.2,-.5) {$m$};
}=
\tikzmatht{
  \fill[black!10] (-.8,-1.2) rectangle (1.6,1.2); 
    \draw[thick] (.5,1.2)--(.5,-1.2); 
    \draw[thick] (-.5,1.2)--(-.5,1) arc(180:270:.5)--(.33,.5);
    \draw[thick] (.5,.67) arc(90:270:.17); 
  \fill[white] (0,0) rectangle (1,-.7);
    \draw[thick] (0,0) rectangle node{$t$}(1,-.7);
    \node at (1.2,.7) {$m'$};
}
\eea

It is obvious from the definition that the modules are closed under
right tensoring with $\rho\in\C$, namely $\mm\times 1_\rho\equiv
(\beta\circ\rho,\mm\times 1_\rho)$ is again a module, and the
corresponding homomorphism is $\varphi\circ\rho$. Moreover, the right
tensoring is compatible with the morphisms.
The category thus defined is therefore a right module category in the
sense of \cite[Def.\ 6]{O}. 

Clearly, every $s\in\Hom(\varphi,\varphi')$ defines a morphism 
$t=1_{\ol\iota}\times s$ between the associated standard modules. 
The converse is also true: 

\begin{tintedbox}
\begin{proposition} \vskip-5mm \label{p:modulint}
Every morphism $t$ between two standard modules
$(\ol\iota\varphi,m=1_{\ol\iota}\times v\times 1_{\varphi})$
and
$(\ol\iota\varphi',m'=1_{\ol\iota}\times v\times 1_{\varphi'})$
is of the form $t= 1_{\ol\iota}\times s$ where
$s\in\Hom(\varphi,\varphi')$.   
\end{proposition}
\end{tintedbox}

{\em Proof:} $s=d_\BA\inv\cdot\LTr_{\ol\iota}(t)$ does the job: 
$$\tikzmatht{
  \fill[black!10] (-1.2,-1.2) rectangle (1.5,1.5); 
  \fill[black!18] (.4,1.5)--(.4,.2) arc(360:180:.4)--
          (-.4,1.5)--(-.8,1.5)--(-.8,.1) arc(180:270:.6)
          arc(90:0:.6)--(.4,-1.2)--(.8,-1.2)--(.8,1.5);
    \draw[thick] (.4,1.5)--(.4,.2) arc(360:180:.4)--
          (-.4,1.5);
    \draw[thick] (-.8,1.5)--(-.8,.1) arc(180:270:.6)
          arc(90:0:.6)--(.4,-1.2); 
    \draw[thick] (.8,1.5)--(.8,-1.2);
  \fill[white] (0,.3) rectangle (1,1);
    \draw[thick] (0,.3) rectangle node{$t$} (1,1);
}
=\tikzmatht{
  \fill[black!10] (-1.2,-1.2) rectangle (1.5,1.5); 
  \fill[black!18] (.4,1.5)--(.4,1) arc(360:180:.4)--
          (-.4,1.5)--(-.8,1.5)--(-.8,.9) arc(180:270:.6)
          arc(90:0:.6)--(.4,-1.2)--(.8,-1.2)--(.8,1.5);
    \draw[thick] (.4,1.5)--(.4,1) arc(360:180:.4)--
          (-.4,1.5);
    \draw[thick] (-.8,1.5)--(-.8,.9) arc(180:270:.6)
          arc(90:0:.6)--(.4,-1.2); 
    \draw[thick] (.8,1.5)--(.8,-1.2);
  \fill[white] (0,-.2) rectangle (1,-.9);
    \draw[thick] (0,-.2) rectangle node{$t$}(1,-.9);
} 
\quad \Rightarrow \quad
\tikzmatht{
  \fill[black!10] (-1.2,-1.2) rectangle (1.5,1.5); 
  \fill[black!18] (-.8,1.5)--(-.8,-.1) arc(180:270:.6)
          arc(90:10:.6)--(.4,-1.2)--(.8,-1.2)--(.8,1.5);
    \draw[thick] (-.8,1.5)--(-.8,-.1) arc(180:270:.6)
          arc(90:10:.6)--(.4,-1.2); 
    \draw[thick] (.8,1.5)--(.8,-1.2);
  \fill[black!10] (-.4,1)--(-.4,0) arc(180:360:.4)--
          (.4,1) arc(0:180:.4);
    \draw[thick] (-.4,1)--(-.4,0) arc(180:360:.4)--
          (.4,1) arc(0:180:.4);
  \fill[white] (0,.1) rectangle (1,.8);
    \draw[thick] (0,.1) rectangle node{$t$}(1,.8);
} = 
\tikzmatht{
  \fill[black!10] (-1.2,-1.2) rectangle (1.5,1.5); 
  \fill[black!18] (-.8,1.5)--(-.8,.9) arc(180:270:.6)
          arc(90:0:.6)--(.4,-1.2)--(.8,-1.2)--(.8,1.5); 
    \draw[thick] (-.8,1.5)--(-.8,.9) arc(180:270:.6)
          arc(90:0:.6)--(.4,-1.2); 
    \draw[thick] (.8,1.5)--(.8,-1.2);
  \fill[black!10] (-.4,1) arc(180:360:.4) arc(0:180:.4);
    \draw[thick] (-.4,1) arc(180:360:.4) arc(0:180:.4);
  \fill[white] (0,-.2) rectangle (1,-.9);
    \draw[thick] (0,-.2) rectangle node{$t$}(1,-.9);
} 
= d_\BA\cdot 
\tikzmatht{
  \fill[black!10] (-.3,-1.2) rectangle (1.5,1.5); 
  \fill[black!18] (.4,1.5)--(.4,-1.2)--(.8,-1.2)--(.8,1.5);
    \draw[thick] (.4,1.5)--(.4,-1.2); 
    \draw[thick] (.8,1.5)--(.8,-1.2);
  \fill[white] (0,-.2) rectangle (1,.5);
    \draw[thick] (0,-.2) rectangle node{$t$}(1,.5);
}.
$$
\qed

We recognize that the argument in the proof of \lref{l:charact} is
just an instance of this general fact, namely \eref{char1} just states
that $t\in\Hom(\theta,\theta)$ is a morphism  
between $\BA=(\theta,x)$ as a $\BA$-module and itself, hence 
$t=1_{\ol\iota}\times s=\ol\iota(s)$ with $s=\iota(q)v\in \Hom(\iota,\iota)$. 

\begin{tintedbox}
\begin{corollary} \vskip-5mm \label{c:modulcat}
The module category of a simple Q-system in $\C\subset\End_0(N)$ is
equivalent to the full subcategory of $\Hom(N,M)$ whose objects are
the homomorphisms $\varphi\prec\iota\rho:N\to M$, $\rho\in\C$. 
\end{corollary}
\end{tintedbox}

In particular: 

\begin{corollary} \label{c:moduldec}
Let $\mm=(\beta,m)$ be a reducible module. The space of self-morphisms
of $\mm$ is a finite-dimensional C* algebra. If $p_i$ are minimal
projections in this algebra, and $p_i=t_it_i^*$ with isometries $t_i$,
then $\mm\simeq\bigoplus_i\mm_i$ with $\mm_i=(\beta_i,m_i)$, where
$\beta_i = t_i^*\beta(\cdot)t_i$ and $m_i= (1_\theta \times t_i^*)
\scirc m\scirc t_i$, i.e.,  
$$\beta = \sum_i t_i\beta_i(\cdot)t_i^*,\quad m= \sum_i (1_\theta
\times t_i)\scirc m_i\scirc t_i^*.$$
\end{corollary}

\begin{graybox}
\begin{example} \vskip-5mm \label{x:Imodules} (Modules in the Ising category)

The irreducible modules of the trivial Q-system are $(\rho,1)$ with
$\rho=\id,\sig,\tau$. The corresponding homomorphisms $:N\to M=N$ are $\varphi =
\rho$. 

The modules $(\beta,x=2^{-\frac14}(r+t))$ of the nontrivial Q-system
given in \xref{x:IQ} ($\theta=\sig^2\simeq\id\oplus\tau$) arising from $\varphi\prec\iota\rho$ are: 

{\rm (i)} $\rho=\id$: module $(\sig^2,x)$, homomorphism $\varphi=\iota$. 

{\rm (ii)} $\rho=\tau$: module $(\sig^2\tau,x)$,
homomorphism $\varphi=\iota\circ\tau$.

{\rm (iii)} $\rho=\sig$: The module ($\beta=\theta\sig=\sig^3,x)$ is
reducible: $\simeq(\beta_1=\sig,x)\oplus(\beta_2=\sig\tau,x)$ 
(with morphisms $\sig(r)\in\Hom(\beta_1,\beta)$ and 
$\sig(t)\in\Hom(\beta_2,\beta)$, respectively). For the 
submodule $(\sig,x)$, one computes $\varphi_1: n \mapsto r^*\sig(n)(r+t\psi)$, 
in particular, $r\mapsto 2^{-\frac12}(r+t\psi), t\mapsto
2^{-\frac12}(r-t\psi), u\mapsto \psi$.
For the submodule $(\sig\tau,x)$, $\varphi_2: n \mapsto 
r^*\sig\tau(n)(r+t\psi)$, in particular $r\mapsto 
2^{-\frac12}(r-t\psi), t\mapsto 2^{-\frac12}(r+t\psi), u\mapsto -\psi$. These
homomorphisms are surjective, hence isomorphisms, and $\varphi_2= 
\varphi_1\circ\tau=\alpha\circ\varphi_1$ ($\alpha =$ gauge
transformation $\psi\to-\psi$).
\end{example}
\end{graybox}

\subsection{Induced Q-systems and Morita equivalence}
\label{s:indQ}
\setcounter{equation}{0}
Let $\BA=(\theta,w,x)$ be a Q-system, defining an extension
$\iota:N\to M$. 

If $\mm=(\beta,m)$ is a standard module of $\BA$, and $\varphi:N\to M$
the corresponding homomorphism, we choose a conjugate homomorphism
$\ol\varphi:M\to N$ and a solution of the conjugacy relations
$w_\varphi$, $v_\varphi$. Then 
$$\BA_\varphi =
(\theta_\varphi,w_\varphi,x_\varphi)\quad\hbox{with}\quad
\theta_\varphi=\ol\varphi\varphi,\quad x_\varphi =
\ol\varphi(v_\varphi)$$
is a Q-system. We call $\BA_\varphi$ the 
{\bf Q-system induced by $\mm$}.

Notice that, by definition,
$\ol\varphi\varphi=\theta_\varphi=\ol\iota_\varphi\iota_\varphi$; but
the corresponding extension $\iota_\varphi:N\to M_\varphi$ should not
be confused with the homomorphism $\varphi:N\to M$, because
$\iota_\varphi(n)=n\in N\subset M_\varphi$, while $\varphi(n)\neq n\in
N\subset M$. 

\begin{lemma} \label{l:induced}
If a Q-system $\BA_2=(\theta_2,w_2,x_2)$ is induced by a
standard module $(\beta_1,m_1)$ of $\BA_1=(\theta_1,w_1,x_1)$,
then $\BA_1$ is induced by a standard module of $\BA_2$. 
\end{lemma}

{\em Proof:} 
By \pref{p:modulphi}, $(\beta_1,m_1)$ is of the form $\beta_1 =
\ol\iota_1\varphi_1$ and $m_1=x_1$, where $\varphi_1\prec
\iota_1\rho$ for some $\rho\in\C$. By definition of the induced
Q-system, $\iota_2 = \varphi_1$, and hence $\iota_1\prec
\iota_2\ol\rho$. Therefore, $\BA_1$ is induced by the module
$(\beta_2=\ol\iota_2\varphi_2,m_2=x_2)$ of $\BA_2$, where
$\varphi_2=\iota_1$. 
\qed

\begin{definition} \label{d:morita}
\cite[Def.\ 10]{O}
Two Q-systems in $\C$ are {\bf Morita equivalent} if their module
categories are equivalent, i.e., there exists an invertible functor between
the two module categories that commutes with the right tensoring by
$\rho\in\C$.  
\end{definition}

\begin{tintedbox}
\begin{proposition} \vskip-5mm \label{p:morita}
Two Q-systems are Morita equivalent if and only if one of
them is induced by a standard module of the other one (which implies
also the converse). 
\end{proposition}
\end{tintedbox}

{\em Proof:} If $\BA_2=(\theta_2,w_2,x_2)$ is induced by a standard module
$(\beta_1,m_1)$ of $\BA_1=(\theta_1,w_1,x_1)$, then $\iota_2 = \varphi_1\prec
\iota_1\rho$ for some $\rho\in\C$. Then the sub-homomorphisms
$\varphi$ of $\iota_2\rho'$ for some $\rho'\in\C$ are the same as the
sub-homomorphisms of $\iota_1\rho''$ for some $\rho''\in\C$. Then, by
\pref{p:modulphi} and \pref{p:modulint}, mapping the standard modules
$(\ol\iota_1\varphi,x_1)$ of $\BA_1$ to $(\ol\iota_2\varphi,x_2)$, and
morphisms $t_1=1_{\ol\iota_1}\times s$ to $t_2=1_{\ol\iota_2}\times s$, defines
a bijective functor that commutes with the right tensoring by $\rho\in\C$.  

Conversely, if $\BA_1$ and $\BA_2$ are Morita equivalent, then there
is a module $\mm_1$ of $\BA_1$ mapped by the bijective functor $F$ to $\BA_2$
as a module of itself. By \pref{p:modulphi}, $\mm_1= (\ol\iota_1\varphi,x_1)$ 
(up to equivalence) with $\varphi\prec\iota\rho$ for some $\rho\in\C$. 
We have to show that the Q-system $\BA_\varphi$ induced by $\varphi$
is equivalent to $\BA_2$. We first show that
$\theta_\varphi=\ol\iota_\varphi\iota_\varphi =\ol\varphi\varphi$
equals $\theta_2$ (up to unitary equivalence). 

For every $\sig\in\C$, one has $\Hom(\sig,\ol\varphi\varphi) \sim
\Hom(\varphi\sig,\varphi)$ by Frobenius reciprocity. Because $\varphi$
corresponds to $\iota_2$ under $F$, and $F$ commutes with right
tensoring by $\sig\in\C$, we further have
$\Hom(\varphi\sig,\varphi) \sim 
\Hom(\iota_2\sig,\iota_2) \sim \Hom(\sig,\theta_2)$, from which the
claim follows.  

Since the construction of the induced Q-system is invariant under the
isomorphism of module categories $F$, it follows that the Q-system 
induced by $\varphi$ from $\BA_1$ coincides with the Q-system induced
by $\iota_2$ from $\BA_2$, which is of course $\BA_2$.  
\qed

Thus, the Q-systems $\BA_\varphi$ induced from a Q-system $\BA$
precisely give the Morita equivalence class of $\BA$. However,
inequivalent $\varphi$ may induce equivalent Q-systems $\BA_\varphi$:
e.g., if $\BA=(\id,1,1)$ is the trivial Q-system, then all invertible
$\varphi$, hence $\ol\varphi\varphi=\id$, induce the trivial Q-system.

\subsection{Bimodules}
\label{s:bim}
\setcounter{equation}{0}
The identification \sref{s:modules} between standard modules (=
left modules) of a Q-system $\BA$ in $\C\subset\End_0(N)$ and
homomorphisms $N\to M$ of the associated pair of algebras works
exactly the same for standard right modules
$\mm=(\beta,m\in\Hom(\beta,\beta\theta))$ (satisfying the analogous
relations with the reversed tensor product). The correspondence is
then that every standard right module is of the form
$$(\beta = \varphi\iota,m=\varphi(v)),$$ 
where $\varphi:M\to N$ is a sub-homomorphism of $\beta\ol\iota$. 

In particular, a Q-system $\BA$ is also a standard right $\BA$-module
$(\beta=\theta,m=x)$, and the corresponding homomorphism is
$\varphi=\ol\iota:M\to N$. 

\medskip

By obvious generalizations of the arguments, one also treats bimodules. 
An $\BA_1$-$\BA_2$-{\bf bimodule} between two Q-systems is a triple
$\mm=(\beta,m_1\in\Hom(\beta,\theta_1\beta), m_2\in
\Hom(\beta,\beta\theta_2))$ such that $(\beta,m_1)$ is a left
$\BA_1$-module and $(\beta,m_2)$ is a right $\BA_2$-module, and the left
and right actions commute:
$$(1_{\theta_1}\times m_2)\scirc m_1 = (m_1\times 1_{\theta_2})\scirc
m_2: \qquad
\tikzmatht{
  \fill[black!10] (-1.4,-1.5) rectangle (1.4,1.2); 
    \draw[thick] (0,1.2)--(0,-1.5); 
    \draw[thick] (-1,1.2)--(-1,0) arc(180:270:.5)--(-.17,-.5);
    \draw[thick] (0,-.33) arc(90:270:.17); 
    \draw[thick] (1,1.2)--(1,.8) arc(0:-90:.5)--(.17,.3);
    \draw[thick] (0,.47) arc(90:-90:.17); 
    \node at (-.6,-1){$m_1$}; 
    \node at (.6,-.1){$m_2$};
} =
\tikzmatht{
  \fill[black!10] (-1.4,-1.5) rectangle (1.4,1.2); 
    \draw[thick] (0,1.2)--(0,-1.5); 
    \draw[thick] (-1,1.2)--(-1,.2) arc(180:270:.5)--(-.17,-.3);
    \draw[thick] (1,1.2)--(1,.2) arc(0:-90:.5)--(.17,-.3);
    \draw[thick] (0,-.3) circle(.17);
    \node at (-.5,-.7){$m$}; 
    \node at (.4,-1){$\beta$};
    \node at (-.6,.7){$\theta_1$}; 
    \node at (.6,.7){$\theta_2$};
} = 
\tikzmatht{
  \fill[black!10] (-1.4,-1.5) rectangle (1.4,1.2); 
    \draw[thick] (0,1.2)--(0,-1.5); 
    \draw[thick] (-1,1.2)--(-1,.8) arc(180:270:.5)--(-.17,.3);
    \draw[thick] (0,.47) arc(90:270:.17); 
    \draw[thick] (1,1.2)--(1,0) arc(0:-90:.5)--(.17,-.5);
    \draw[thick] (0,-.33) arc(90:-90:.17); 
    \node at (-.6,-.1){$m_1$}; 
    \node at (.6,-1){$m_2$};
}
$$
Equivalently, one may characterize the bimodule as a pair
$\mm=(\beta\in\C,m\in\Hom(\beta,\theta_1\beta\theta_2))$ satisfying 
\be \label{bimrep}
\tikzmatht{
  \fill[black!10] (-1.4,-1.5) rectangle (1.4,1.2); 
    \draw[thick] (0,1.2)--(0,-1.5); 
    \draw[thick] (-1,.8)--(-1,.2) arc(180:270:.5)--(-.17,-.3);
    \draw[thick] (1,.8)--(1,.2) arc(0:-90:.5)--(.17,-.3);
  \fill[white] (-1,.8) circle(.2);
    \draw[thick] (-1,.8) circle(.2);
  \fill[white] (1,.8) circle(.2);
    \draw[thick] (1,.8) circle(.2);
    \draw[thick] (0,-.3) circle(.17);
    \node at (-.5,-.7){$m$}; 
    \node at (.4,-1){$\beta$};
} = 
\tikzmatht{
  \fill[black!10] (-1,-1.5) rectangle (1,1.2);
    \draw[thick] (0,-1.5)--(0,1.2); 
    \node at (.4,-1) {$\beta$}; 
},\qquad 
\tikzmatht{
  \fill[black!10] (-1.8,-1.5) rectangle (1.8,1.2); 
    \draw[thick] (0,1.2)--(0,-1.5); 
    \draw[thick] (-.8,1.2)--(-.8,.5) arc(180:270:.5)--(-.17,0);
    \draw[thick] (.8,1.2)--(.8,.5) arc(0:-90:.5)--(.17,0);
    \draw[thick] (0,0) circle(.17);
    \draw[thick] (-1.5,1.2)--(-1.5,.2) arc(180:270:1)--(-.17,-.8);
    \draw[thick] (1.5,1.2)--(1.5,.2) arc(0:-90:1)--(.17,-.8);
    \draw[thick] (0,-.8) circle(.17);
    \node at (-.5,-1.2){$m$}; 
} = 
\tikzmatht{
  \fill[black!10] (-1.8,-1.5) rectangle (1.8,1.2); 
    \draw[thick] (0,1.2)--(0,-1.5); 
    \draw[thick] (-1,.1)--(-1,-.1) arc(180:270:.7)--(-.17,-.8);
    \draw[thick] (1,.1)--(1,-.1) arc(0:-90:.7)--(.17,-.8);
    \draw[thick] (0,-.8) circle(.17);
\draw[thick] (-1.5,1.2)--(-1.5,.6) arc(180:360:.5)--(-.5,1.2);
\draw[thick] (1.5,1.2)--(1.5,.6) arc(360:180:.5)--(.5,1.2);
  \fill[black] (-1,.1) circle(.12); 
  \fill[black] (1,.1) circle(.12);
    \node at (-.5,-1.2){$m$}; 
    \node at (-1.4,-.2){$x_1$}; 
    \node at (1.45,-.2){$x_2$};
}
\ee
Then $(\beta,m_1):=(1_{\theta_1}\times 1_\beta\times w_2^*)\scirc m$ is a
left $\BA_1$-module, $(\beta,m_2):=(w_1^*\times 1_\beta\times
1_{\theta_2})\scirc m$ is a right $\BA_2$-module, 
and their actions commute. 

A bimodule of a Q-system is called a {\bf standard bimodule} if
$m^*\ncirc m$ is a multiple of $1_\beta$. 

A Q-system $\BA$ is also a standard $\BA$-$\BA$-bimodule
$\BA=(\beta=\theta,m=x^{(2)})$.  

One proves the analogs of \lref{l:normmod},
\lref{l:conjmod}, \pref{p:modulphi}, 
\pref{p:modulint} and \cref{c:modulcat}
in more or less exactly the same way (replacing the left trace in 
\pref{p:modulint} by the right trace for the right module action):  

\begin{tintedbox}
\begin{proposition} \vskip-5mm \label{p:standbim}
{\rm (i)} Every bimodule is equivalent to a standard bimodule. The
normalization  
of a standard bimodule is $m^*\ncirc m = d_{\BA_1}d_{\BA_2}\cdot 1_\beta$. 
The adjoint of a bimodule is obtained by Frobenius reciprocity. \\
{\rm (ii)} \cite{EP03} Every standard bimodule is unitarily equivalent to
a bimodule of the form $\beta = \bar\iota_1\varphi\iota_2$ where 
$\varphi:M_2\to M_1$ is a sub-homomorphism of $\iota_1\rho\bar\iota_2$ 
for some $\rho\in\C\subset\End_0(N)$ (e.g., $\rho=\beta$), and 
$$b=1_{\bar\iota_1}\times v_1\times 1_\varphi \times v_2\times
1_{\iota_2}: \qquad 
\tikzmatht{
  \fill[black!10] (-1.8,-1.5) rectangle (1.8,1.2); 
    \draw[thick] (0,1.2)--(0,-1.5); 
    \draw[thick] (-1,1.2)--(-1,.2) arc(180:270:.5)--(-.17,-.3);
    \draw[thick] (1,1.2)--(1,.2) arc(0:-90:.5)--(.17,-.3);
    \draw[thick] (0,-.3) circle(.17);
    \node at (.5,-.7){$m$};
    \node at (-1.4,.7){$\theta_1$}; 
    \node at (1.4,.7){$\theta_2$};
} \equiv 
\tikzmatht{
  \fill[black!10] (-1.8,-1.2) rectangle (1.8,1.5); 
  \fill[black!18] (-1.6,1.5)--(-1.6,.5) arc(180:250:.9)
          arc(70:0:.9)--(-.4,-1.2)--(0,-1.2)--(0,1.5); 
    \draw[thick] (-1.6,1.5)--(-1.6,.5) arc(180:250:.9)
          arc(70:0:.9)--(-.4,-1.2); 
  \fill[black!10] (-1.2,1.5)--(-1.2,.6) arc(180:360:.4)--
          (-.4,1.5);
    \draw[thick] (-1.2,1.5)--(-1.2,.6) arc(180:360:.4)--
          (-.4,1.5);
  \fill[black!25] (1.6,1.5)--(1.6,.5) arc(360:290:.9)
          arc(110:180:.9)--(.4,-1.2)--(0,-1.2)--(0,1.5);
    \draw[thick] (1.6,1.5)--(1.6,.5) arc(360:290:.9)
          arc(110:180:.9)--(.4,-1.2); 
  \fill[black!10] (1.2,1.5)--(1.2,.6) arc(360:180:.4)--
          (.4,1.5);
    \draw[thick] (1.2,1.5)--(1.2,.6) arc(360:180:.4)--
          (.4,1.5);
    \draw[thick] (0,1.5)--(0,-1.2);
    \node at (.3,-.3){$\varphi$};
}.$$
\end{proposition}
\end{tintedbox}

A morphism between two bimodules is an element
$t\in\Hom(\beta,\beta')$ satisfying 
\bea\label{bimorph}
(1_{\theta_1}\times t \times 1_{\theta_2})\scirc m = m' \scirc t:
\qquad 
\tikzmatht{
  \fill[black!10] (-1.8,-1.5) rectangle (1.8,1.2); 
    \draw[thick] (0,1.2)--(0,-1.5); 
    \draw[thick] (-1,1.2)--(-1,-.2) arc(180:270:.5)--(-.17,-.7);
    \draw[thick] (1,1.2)--(1,-.2) arc(0:-90:.5)--(.17,-.7);
    \draw[thick] (0,-.7) circle(.17);
  \fill[white] (-.5,0) rectangle (.5,.7);
    \draw[thick] (-.5,0) rectangle node{$t$} (.5,.7);
    \node at (.5,-1.1){$m$};
    \node at (-1.4,.7){$\theta_1$}; 
    \node at (1.4,.7){$\theta_2$};
} =
\tikzmatht{
  \fill[black!10] (-1.8,-1.5) rectangle (1.8,1.2); 
    \draw[thick] (0,1.2)--(0,-1.5); 
    \draw[thick] (-1,1.2)--(-1,.8) arc(180:270:.5)--(-.17,.3);
    \draw[thick] (1,1.2)--(1,.8) arc(0:-90:.5)--(.17,.3);
    \draw[thick] (0,.3) circle(.17);
  \fill[white] (-.5,-.5) rectangle (.5,-1.2);
    \draw[thick] (-.5,-.5) rectangle node{$t$} (.5,-1.2);
    \node at (.5,-.1){$m'$};
    \node at (-1.4,.7){$\theta_1$}; 
    \node at (1.4,.7){$\theta_2$};
} .
\eea

\begin{tintedbox}
\begin{proposition} \vskip-5mm \label{p:bimfunct}
Every morphism $t$ between two standard $\BA_1$-$\BA_2$-bimodules 
$(\ol\iota_1\varphi\iota_2,1_{\ol\iota_1}\times v_1\times 1_\varphi 
\times v_2\times 1_{\iota_2})$ and 
$(\ol\iota_1\varphi'\iota_2,1_{\ol\iota_1}\times v_1\times 1_\varphi'
\times v_2\times 1_{\iota_2})$ is of the form $t=1_{\ol\iota_1}\times
s\times 1_{\iota_2}$ where $s\in\Hom(\varphi,\varphi')$. This
establishes a bijective functor between the category of
$\BA_1$-$\BA_2$-bimodules and
the full subcategory of $\Hom(M_2,M_1)$ whose objects are the
homomorphisms $\prec\iota_1\rho\ol\iota_2$, $\rho\in\C$. 
\end{proposition}
\end{tintedbox}

Again, the homomorphism associated with a standard bimodule
$\mm=(\beta,m)$ can be 
computed. Namely, the formula for $m$ implies that 
$\ol\iota_1\varphi(v_2)=w_1^*m$ (corresponding to $\mm$ as a 
right $\BA_2$-module). Hence $\varphi(\iota_2(n)v_2) = \iota_1(k)v_1$ 
implies $\beta(n)w_1^*m = \theta_1(k)x_1$, hence $k=w_1^*\beta(n)w_1^*m$:
\be \label{phi-bim} 
\varphi(\iota_2(n)v_2) = \iota_1(w_1^*\beta(n)w_1^*m)v_1.
\ee 
In particular, $\varphi(v_2)=\iota_1(w_1^*w_1^*m)v_1$. 

The homomorphism associated with $\BA$ as an $\BA$-$\BA$-bimodule is
$\varphi=\id_M:M\to M$. 

\medskip

If $\varphi=\iota_1\rho\ol\iota_2$ (which is in general reducible), 
hence $\beta=\theta_1\rho\theta_2$ and $m=x_1\theta_1\rho(x_2)$, this 
simplifies to $\varphi(\iota_2(n))=\iota_1(\rho\theta_2(n))$ and
$\varphi(v_2)=\iota_1(\rho(x_2))$. Thus, $\varphi:M_2\to M_1$ happens 
to take values in $\iota_1(N)\subset M_1$. This property is, however, 
not intrinsic, as it is not stable under unitary equivalence in the target
algebra $M_1$. Also, the decomposition of $\varphi$ into irreducibles 
(which are unique only up to unitary equivalence within $M_1$) depends 
on the choice of isometries $s$, so that $\varphi_s=s^*\varphi(\cdot)s$. 
These may or may not be chosen in $\iota_1(N)$. As the 
\xref{x:Ibimodules} shows, there may be good reasons to choose the
homomorphisms {\em not} to take values in $\iota_1(N)$. 

\medskip

Also the analog of \pref{p:morita} holds for bimodules,
again with the same proof as for modules: 

\begin{tintedbox}
\begin{proposition} \vskip-5mm \label{p:moritabim} 
  There is a bijective functor
  between the category of $\BA_1$-$\BA_2$-bimodules and the category
  of $\BA_1'$-$\BA_2'$-bimodules, if and only if $\BA_1'$ is induced
  from $\BA_1$ by a standard module of $\BA_1$ (i.e.,
  $\iota_1\prec\iota_2\rho$, $\rho\in\C$), and $\BA_2'$ is induced
  from $\BA_2$ by a standard module of $\BA_2$.
\end{proposition}
\end{tintedbox}

In particular, the category of bimodules between a 
pair of Q-systems depends only on the Morita equivalence classes of
the latter. 

\begin{graybox}
\begin{example} \vskip-5mm \label{x:Ibimodules} (Bimodules in the Ising category)

Let $\BA=(\theta=\sigma^2,w=2^{\frac14} 
r,x=2^{-\frac14}(r+t))$ be the nontrivial Q-system as in \xref{x:IQ} and 
$M=N\vee\psi$ be the corresponding extension of $N$. The irreducible 
$\id$-$\id$-bimodules are just $\rho=\id,\sig,\tau$. 
The $\BA$-$\id$-bimodules are the same as the
modules of $\BA$, \xref{x:Imodules}. 

The $\id$-$\BA$-bimodules arising from $\varphi=\rho\ol\iota$ are: 
$\mm_\rho=(\beta=\rho\theta=\rho\sigma^2,m=\rho(x))$, where
$x=2^{-\frac14}(r+t)$. Thus, $\varphi$ maps $n\in N$ to $\rho\sigma^2(n)$ and
$v$ to $\rho(x)$. 

{\rm (i)} $\rho=\id$ and $\rho=\tau$: These are the same bimodules, because
$\tau\sig=\sig$ and $\tau(x)=x$. One finds $\varphi:n\mapsto
\sig^2(n),\psi\mapsto rt^*+tr^*$.

{\rm (ii)} $\rho=\sig$: $\beta=\sig^3$, $m=\sig(x)=2^{-\frac14}\sig(r+t) =
2^{-\frac34}(r+t+(r-t)u)$. $\varphi:n\mapsto\sig^3(n),\psi\mapsto rur^*-tut^*$.

The latter homomorphism $\varphi =\sig\ol\iota$ is reducible, with projections
$rr^*$ and $tt^*$ in the commutant. Then $\varphi_1=r^*\varphi(\cdot)r$
and $\varphi_2=t^*\varphi(\cdot)t$ give rise to 

{\rm (ii.1)} $\beta_1=\sig$, $m_1=2^{\frac14}r$, 
$\varphi_1:n\mapsto\sig(n),\psi\mapsto u$.

{\rm (ii.2)} $\beta_2=\tau\sig=\sig$, $m_2=2^{\frac14}t$,
$\varphi_2:n\mapsto\sig(n),\psi\mapsto -u$. One has
$\varphi_2=\varphi_1\circ\alpha=\tau\circ\varphi_1$.

The $\BA$-$\BA$-bimodules arising from $\varphi=\iota\rho\ol\iota$ are: 
$\mm_\rho=(\beta=\theta\rho\theta=\sigma^2\rho\sigma^2,m=x\theta\rho(x))$. 
Thus $\varphi$ maps $n$ to $\rho\sigma^2(n)$ and $v$ to $\rho(x)$.

{\rm (i)} $\rho=\id$ and $\rho=\tau$ are again the same bimodule.
$\varphi: n\mapsto \sig^2(n),\psi\mapsto rt^*+tr^*$.
$\varphi$ is reducible with projections $\frac12(1\pm\psi)=s_\pm s_\pm^*$,
$s_\pm=2^{-\frac12}(r\pm t\psi)$.
This gives the irreducible components $\varphi_\pm:n\mapsto
n,\psi\mapsto \pm\psi$, i.e., $\varphi_+=\id$ and $\varphi_-=\alpha$.

{\rm (ii)} $\rho=\sig$. $\varphi: n\mapsto \sig^3(n),\psi\mapsto
\sig(rt^*+tr^*)=rur^*-tut^*$. The commutant of $\varphi(M)$ contains 
$u$ and $\psi$. Thus $\varphi$ is the direct sum of two equivalent 
components, $[\varphi]=[\varphi']\oplus[\varphi']$. Choosing projections 
$rr^*$ and $tt^*$ to compute $\varphi_1=r^*\varphi(\cdot)r$ and 
$\varphi_2=t^*\varphi(\cdot)t$, one has $\varphi_1: n\mapsto
\sig(n),\psi\mapsto u$ and $\varphi_2: n\mapsto \sig(n),\psi\mapsto
-u$. These are equivalent to each other by $\psi\in\Hom(\varphi_1,\varphi_2)$. 
They are also equivalent to the $\alpha$-inductions $\alpha^\pm_\sig$
(cf.\ \sref{s:alpha}) by $U_\pm = 2^{-\frac12}(1\pm i\psi)$; with the 
former choice, $\varphi_i$ take values in $\iota_1(N)$, while 
$\Ad_{U_\pm}\varphi_1=\Ad_{U_\mp}\varphi_2=\alpha^\pm_\sig$ don't.)
\end{example}
\end{graybox}

\subsection{Tensor product of bimodules}
\setcounter{equation}{0}

The tensor product of bimodules is defined as follows. 
If $\mm_1=(\beta_1,m_1)$ is an $\BA$-$\BB$-bimodule 
and $\mm_2=(\beta_2,m_2)$ is an $\BB$-$\BC$-bimodule, then 
$$\wh m=\tikzmatht{
  \fill[black!10] (-2.2,-1.7) rectangle (2.2,1.2); 
    \draw[thick] (-.3,1.2)--(-.3,-1.7); 
    \draw[thick] (.3,1.2)--(.3,-1.7); 
    \draw[thick] (-1.3,1.2)--(-1.3,.2) arc(180:270:.5)--
          (-.47,-.3) (-.13,-.3)--(.13,-.3);
    \draw[thick] (1.3,1.2)--(1.3,.2) arc(0:-90:.5)--(.47,-.3);
    \draw[thick] (-.3,-.3) circle(.17);
    \draw[thick] (.3,-.3) circle(.17);
    \node at (-.8,-.7){$m_1$}; 
    \node at (.9,-.7){$m_2$};
    \node at (-1.8,.7){$\theta^A$}; 
    \node at (1.8,.7){$\theta^C$};
}\in\Hom(\beta_1\beta_2,\theta^\BA\beta_1\beta_2\theta^\BC)$$
satisfies the representation property of an $\BA$-$\BC$-bimodule, but
the unit property fails. Instead, we have

\begin{lemma} \label{l:bimprod}
The intertwiner 
$$p:=d_\BB\inv\cdot 
\tikzmatht{
  \fill[black!10] (-2,-1.7) rectangle (2,1.2); 
    \draw[thick] (-.3,1.2)--(-.3,-1.7); 
    \draw[thick] (.3,1.2)--(.3,-1.7); 
    \draw[thick] (-1.3,.6)--(-1.3,.2) arc(180:270:.5)--
          (-.47,-.3) (-.13,-.3)--(.13,-.3);
    \draw[thick] (1.3,.6)--(1.3,.2) arc(0:-90:.5)--(.47,-.3);
  \fill[white] (1.3,.6) circle(.2);
    \draw[thick] (1.3,.6) circle(.2);
  \fill[white] (-1.3,.6) circle(.2);
    \draw[thick] (-1.3,.6) circle(.2);
    \draw[thick] (-.3,-.3) circle(.17);
    \draw[thick] (.3,-.3) circle(.17);
    \node at (-.8,-.7){$m_1$}; 
    \node at (.9,-.7){$m_2$};
} \equiv d_\BB\inv\cdot 
\tikzmatht{
  \fill[black!10] (-1.5,-1.7) rectangle (1.5,1.2); 
    \draw[thick] (-.4,1.2)--(-.4,-1.7); 
    \draw[thick] (.4,1.2)--(.4,-1.7); 
    \draw[thick] (-.4,-.13) arc(90:-90:.17);
    \draw[thick] (.4,-.13) arc(90:270:.17);
    \draw[thick] (-.23,-.3)--(.23,-.3);
    \node at (-.8,.7){$\beta_1$}; 
    \node at (.9,.7){$\beta_2$};
}
\in\Hom(\beta_1\beta_2,\beta_1\beta_2)
$$
is a projection, and satisfies
\be \label{bimprod}
(1_{\theta^\BA}\times p \times 1_{\theta^\BC})\scirc \wh m = \wh m
= \wh m\scirc p: \;\;  d_\BB\inv\cdot
\tikzmatht{
  \fill[black!10] (-1.6,-1.2) rectangle (1.6,1.5); 
    \draw[thick] (-.5,1.5)--(-.5,-1.2); 
    \draw[thick] (.5,1.5)--(.5,-1.2); 
    \draw[thick] (-1.3,1.5)--(-1.3,.2) arc(180:270:.5)--
          (-.67,-.3) (-.33,-.3)--(.33,-.3);
    \draw[thick] (1.3,1.5)--(1.3,.2) arc(0:-90:.5)--(.67,-.3);
    \draw[thick] (-.5,-.3) circle(.17);
    \draw[thick] (.5,-.3) circle(.17);
    \draw[thick] (-.5,.77) arc(90:-90:.17);
    \draw[thick] (.5,.77) arc(90:270:.17);
    \draw[thick] (-.33,.6)--(.33,.6);
}  = 
\tikzmatht{
  \fill[black!10] (-2.3,-1.2) rectangle (2.3,1.5); 
    \draw[thick] (-.5,1.5)--(-.5,-1.2); 
    \draw[thick] (.5,1.5)--(.5,-1.2); 
    \draw[thick] (-1.3,1.5)--(-1.3,.5) arc(180:270:.5)--
          (-.67,0) (-.33,0)--(.33,0);
    \draw[thick] (1.3,1.5)--(1.3,.5) arc(0:-90:.5)--(.67,0);
    \draw[thick] (-.5,0) circle(.17);
    \draw[thick] (.5,0) circle(.17);
    \node at (-1.8,1.1) {$\theta^A$}; 
    \node at (1.8,1.1) {$\theta^C$};
    \node at (-1,-.5){$m_1$}; 
    \node at (1.1,-.5){$m_2$};
} =
d_\BB\inv\cdot
\tikzmatht{
  \fill[black!10] (-1.6,-1.2) rectangle (1.6,1.5); 
    \draw[thick] (-.5,1.5)--(-.5,-1.2); 
    \draw[thick] (.5,1.5)--(.5,-1.2); 
    \draw[thick] (-1.3,1.5)--(-1.3,1.1) arc(180:270:.5)--
          (-.67,.6) (-.33,.6)--(.33,.6);
    \draw[thick] (1.3,1.5)--(1.3,1.1) arc(0:-90:.5)--(.67,.6);
    \draw[thick] (-.5,.6) circle(.17);
    \draw[thick] (.5,.6) circle(.17);
    \draw[thick] (-.5,-.13) arc(90:-90:.17);
    \draw[thick] (.5,-.13) arc(90:270:.17);
    \draw[thick] (-.33,-.3)--(.33,-.3);
}.
\ee
\end{lemma}

{\em Proof:}
Idempotency of $p$ follows from the relation \eref{bimprod}. 
Self-adjointness of $p$ follows from \lref{l:conjmod}. To prove
\eref{bimprod}, we use the representation property, e.g.,  
$$\tikzmatht{
  \fill[black!10] (-2.3,-1.2) rectangle (2.3,1.5); 
    \draw[thick] (-.5,1.5)--(-.5,-1.2); 
    \draw[thick] (.5,1.5)--(.5,-1.2); 
    \draw[thick] (-1.3,1.5)--(-1.3,.2) arc(180:270:.5)--
          (-.67,-.3) (-.33,-.3)--(.33,-.3);
    \draw[thick] (1.3,1.5)--(1.3,.2) arc(0:-90:.5)--(.67,-.3);
    \draw[thick] (-.5,-.3) circle(.17);
    \draw[thick] (.5,-.3) circle(.17);
    \draw[thick] (-.5,.77) arc(90:-90:.17);
    \draw[thick] (.5,.77) arc(90:270:.17);
    \draw[thick] (-.33,.6)--(.33,.6);
    \node at (-1.8,1.1) {$\theta^A$}; 
    \node at (1.8,1.1) {$\theta^C$};
}  = 
\tikzmatht{
  \fill[black!10] (-2.2,-1.2) rectangle (2.2,1.5); 
    \draw[thick] (-1.1,1.5)--(-1.1,-1.2); 
    \draw[thick] (1.1,1.5)--(1.1,-1.2); 
    \draw[thick] (-2,1.5)--(-2,.2) arc(180:270:.5)--
          (-1.27,-.3); 
    \draw[thick] (2,1.5)--(2,.2) arc(0:-90:.5)--(1.27,-.3);
    \draw[thick] (-1.1,.47) arc(90:-90:.17);
    \draw[thick] (1.1,.47) arc(90:270:.17);
    \draw[thick] (1.1,-.13) arc(90:-90:.17);
    \draw[thick] (-1.1,-.13) arc(90:270:.17);
    \draw[thick] (-.93,.3)--(-.5,.3)(.5,.3)--(.93,.3);
    \draw[thick] (0,.3) circle(.5);
  \fill[black] (-.5,.3) circle(.12);
  \fill[black] (.5,.3) circle(.12);
    \node at (.2,-.6) {$\theta^B$};
} = d_\BB\cdot 
\tikzmatht{
  \fill[black!10] (-2.2,-1.2) rectangle (2.2,1.5); 
    \draw[thick] (-.5,1.5)--(-.5,-1.2); 
    \draw[thick] (.5,1.5)--(.5,-1.2); 
    \draw[thick] (-1.5,1.5)--(-1.5,.2) arc(180:270:.5)--
          (-.67,-.3) (-.33,-.3)--(.33,-.3);
    \draw[thick] (1.5,1.5)--(1.5,.2) arc(0:-90:.5)--(.67,-.3);
    \draw[thick] (-.5,-.3) circle(.17);
    \draw[thick] (.5,-.3) circle(.17);
}.$$ 
\qed

Then the bimodule tensor product is defined as the range of the
projection $p$: 

\begin{definition} \label{d:bimprod} 
Let $\mm_1=(\beta_1,m_1)$ be an $\BA$-$\BB$-bimodule
and $\mm_2=(\beta_2,m_2)$ a $\BB$-$\BC$-bimodule. Choose an isometry
$s\in N$ such that $ss^*=p$ and put $\beta(\cdot):=s^*\beta_1\beta_2(\cdot)s$ 
the range of $p$ in $\beta_1\beta_2$. Then the {\bf bimodule tensor product}
$$\mm_1\otimes_{\BB}\mm_2 = (\beta,m),$$
$$m:=d_\BB\inv\cdot
(1_{\theta^\BA}\times s^*\times 1_{\theta^\BC})\scirc \wh m \scirc s
=d_\BB\inv\cdot 
\tikzmatht{
  \fill[black!10] (-2.2,-1.7) rectangle (2.2,1.2); 
    \draw[thick] (-.3,.5)--(-.3,-1); 
    \draw[thick] (.3,.5)--(.3,-1); 
    \draw[thick] (0,-1.7)--(0,-1.4) (0,.9)--(0,1.2); 
  \fill[white] (0,-1.4)--(.3,-1)--(-.3,-1)--(0,-1.4);
\draw[thick] (0,-1.4)--(.3,-1)--(-.3,-1)--(0,-1.4);
  \fill[white] (0,.9)--(.3,.5)--(-.3,.5)--(0,.9);
    \draw[thick] (0,.9)--(.3,.5)--(-.3,.5)--(0,.9);
    \draw[thick] (-1.3,1.2)--(-1.3,.2) arc(180:270:.5)--
          (-.47,-.3) (-.13,-.3)--(.13,-.3);
    \draw[thick] (1.3,1.2)--(1.3,.2) arc(0:-90:.5)--(.47,-.3);
    \draw[thick] (-.3,-.3) circle(.17);
    \draw[thick] (.3,-.3) circle(.17);
    \node at (-.8,-.7){$m_1$}; 
    \node at (.9,-.7){$m_2$};
    \node at (-1.8,.7){$\theta^A$}; 
    \node at (1.8,.7){$\theta^C$};
    \node at (-.5,.9){$s^*$}; 
    \node at (-.5,-1.5){$s$}; 
}\in\Hom(\beta,\theta^\BA\beta\theta^\BC)
$$
is an $\BA$-$\BC$-bimodule. 
\end{definition}

\begin{tintedbox}
\begin{proposition} \vskip-5mm \label{p:bimprod}
Under the correspondence \pref{p:standbim}(ii), the bimodule tensor
product $\mm_1\otimes_\BB\mm_2$ corresponds to the composition of
homomorphisms $\varphi_1\circ\varphi_2:M^\BC\to M^\BA$. 
\end{proposition}
\end{tintedbox}

{\em Proof:} Using \pref{p:standbim}(ii), one computes 
$$p=d_\BB\inv\cdot 1_{\ol\iota^\BA\circ\varphi_1}\times w^\BB w^{\BB*} \times
1_{\varphi_2\circ\iota^\BC} = d_\BB\inv\cdot
\tikzmatht{
  \fill[black!10] (-2.4,-1) rectangle (2.4,1.5);
  \fill[black!18] (-1.6,-1)--(-1.6,1.5)--(1.6,1.5)--(1.6,-1);
    \draw[thick] (-1.6,-1)--(-1.6,1.5) (1.6,1.5)--(1.6,-1);
  \fill[black!25] (-.8,-1)--(-.8,1.5)--(.8,1.5)--(.8,-1);
    \draw[thick] (-.8,-1)--(-.8,1.5) (.8,1.5)--(.8,-1);
  \fill[black!10] (-.4,1.5)--(-.4,1) arc(180:360:.4)--(.4,1.5);
    \draw[thick] (-.4,1.5)--(-.4,1) arc(180:360:.4)--(.4,1.5);
  \fill[black!10] (-.4,-1)--(-.4,-.5) arc(180:0:.4)--(.4,-1);
    \draw[thick] (-.4,-1)--(-.4,-.5) arc(180:0:.4)--(.4,-1);
    \node at (-1.2,0) {$\varphi_1$};
    \node at (1.2,0) {$\varphi_2$};
    \node at (-2,.2) {$\ol\iota^A$};
    \node at (2,.2) {$\iota^C$};
}
,$$
hence (up to unitary equivalence) one may choose 
$$s=d_\BB^{-\frac12}\cdot 1_{\ol\iota^\BA\circ\varphi_1}\times w^\BB \times
1_{\varphi_2\circ\iota^\BC} \equiv d_\BB^{-\frac12}\cdot 
\tikzmatht{
  \fill[black!10] (-2,-1) rectangle (2,1.5);
  \fill[black!18] (-1.4,-1)--(-1.4,1.5)--(1.4,1.5)--(1.4,-1);
    \draw[thick] (-1.4,-1)--(-1.4,1.5) (1.4,1.5)--(1.4,-1);
  \fill[black!25] (-.8,-1)--(-.8,1.5)--(.8,1.5)--(.8,-1);
    \draw[thick] (-.8,-1)--(-.8,1.5) (.8,1.5)--(.8,-1);
  \fill[black!10] (-.4,1.5)--(-.4,1) arc(180:360:.4)--(.4,1.5);
    \draw[thick] (-.4,1.5)--(-.4,1) arc(180:360:.4)--(.4,1.5);
    \node at (0,.3) {$w^B$};
}. 
$$
With this, the claim is easily verified. The proper normalization is
fixed by \pref{p:standbim}(i).
\qed

In particular, we have equipped the category of $\BA$-$\BA$-bimodules
with the structure of a tensor category, such that the tensor product
corresponds to the composition of the corresponding endomorphisms in
$\End_0(M)$. By admitting bimodules between different Q-systems
$\BA_i$, one arrives naturally at a (non-strict) bicategory (with
1-objects $\BA_i$, 1-morphisms the 
bimodules and 2-morphisms the bimodule morphisms), corresponding to
homomorphisms among the associated extensions $M_i$. Fixing the von
Neumann algebra $N$
and some full subcategory $\C$ of $\End_0(N)$ in which the Q-systems,
bimodules and morphisms take their values, one obtains a full
sub-2-category of the latter 2-category. 

\medskip

In the tensor category of $\BA$-$\BA$-bimodules, the bimodule $\BA$ is
the tensor unit. Correspondingly, this category is simple iff $\BA$ is
irreducible as a $\BA$-$\BA$-bimodule. The following Lemma
characterizes the self-intertwiners of $\BA$:

\begin{lemma} \label{l:selfint} $t\in\Hom(\theta,\theta)$ is a
  self-morphism of $\BA$ as left (right) $\BA$-module if and only
  if $t$ satisfies the first (second) of \eref{1to3}. 
$t\in\Hom(\theta,\theta)$ is a self-morphism of $\BA$ as 
an $\BA$-$\BA$-bimodule if and only if 
$t\in\Hom_0(\theta,\theta)$.
\end{lemma}

{\em Proof:} The first statement is just the definition of morphisms. We prove only the last statement. ``If'': obvious. ``Only if'': by applying the unit
relation in several ways to the defining property of a bimodule morphism
$$\tikzmatht{
  \fill[black!10] (-1.4,-1) rectangle (1.4,1.4);
    \draw[thick] (-1,1.4)--(-1,.7) arc(180:360:1)--(1,1.4);
    \draw[thick] (0,-1)--(0,1.4);
  \fill[black] (0,-.3) circle(.12);
  \fill[white] (-.3,.8) rectangle (.3,.4);
    \draw[thick] (-.3,.8) rectangle (.3,.4);
}  =
\tikzmatht{
  \fill[black!10] (-1.4,-1) rectangle (1.4,1.4);
    \draw[thick] (-1,1.4) arc(180:360:1)--(1,1.4);
    \draw[thick] (0,-1)--(0,1.4);
  \fill[black] (0,.4) circle(.12);
  \fill[white] (-.3,-.1) rectangle (.3,-.5);
    \draw[thick] (-.3,-.1) rectangle (.3,-.5);
} .$$
\qed

\begin{tintedbox}
\begin{corollary} \vskip-5mm \label{c:simple=factorial} The following are
  equivalent. \\
{\rm (i)} A Q-system $\BA$ is simple. \\
{\rm (ii)} The corresponding extension $N\subset M$ is a factor. \\
{\rm (iii)} $\BA$ is irreducible as an $\BA$-$\BA$-bimodule. \\
{\rm (iv)} The tensor category of $\BA$-$\BA$-bimodules is simple.
\end{corollary}
\end{tintedbox}

Notice that (i) $\LRA$ (ii) is our \dref{d:simple} of
a simple Q-system. (iii) $\LRA$ (iv) is the definition of a simple
tensor category. Thus, \cref{c:simple=factorial} states the
equivalence of our \dref{d:simple} of simplicity with the
standard definition, which is given by the condition (iv). 

\medskip

{\em Proof:} It suffices to prove (ii) $\LRA$ (iii). The endomorphism
$\varphi:M\to M$ corresponding to the bimodule $\BA$ according to 
\pref{p:modulphi}, is $\varphi=\id_M$. Then, by \pref{p:bimfunct}, 
every self-intertwiner of $\BA$ as an $\BA$-$\BA$-bimodule is of the 
form $t=1_{\ol\iota}\times s\times 1_\iota\in\Hom(\theta,\theta)$, where $s\in
\Hom(\id_M,\id_M)$. But $\Hom(\id_M,\id_M)$ is the same as the centre
$M'\cap M$.  
\qed

\medskip

For later use, we mention

\begin{lemma} \label{l:moduleint}
If $\mm_i=(\beta_i,m_i)$ ($i=1,2$) are 
$\BA$-$\BB$-bimodules, and $t\in \Hom(\beta_1,\beta_2)$, then 
$$S:=m_2^*\scirc (1_{\theta^\BA}\times t\times 1_{\theta^\BB})\scirc m_1
=
\tikzmatht{
  \fill[black!10] (-1.5,-1.5) rectangle (1.5,1.5);
    \draw[thick] (-.17,.8) arc(90:270:.8) (.17,.8) arc(90:-90:.8);
    \draw[thick] (0,.8) circle(.17);
    \draw[thick] (0,-.8) circle(.17);
    \draw[thick] (0,-1.5)--(0,1.5); 
  \fill[white] (-.3,-.35) rectangle (.3,.35); 
    \draw[thick] (-.3,-.35) rectangle node {$t$} (.3,.35); 
    \node at (-.5,-1.2) {$\beta_1$};
    \node at (-.5,+1.2) {$\beta_2$};
}\in \Hom(\beta_1,\beta_2)$$
is a bimodule morphism $:\mm_1\to \mm_2$. 
\end{lemma}

The proof is rather easy in terms of the defining properties of
modules and module intertwiners, and actually becomes trivial if one
uses \pref{p:standbim}: namely $S\in 1_{\ol\iota^\BA}\times
\Hom(\varphi_1,\varphi_2)\times 1_{\iota^\BB}$. 

\smallskip 

This Lemma implies that if the two bimodules are irreducible and
inequivalent, then every intertwiner $S$ obtained in this way must be 
trivial. E.g., if $\mm_2$ is trivial $\BA$-$\BA$-bimodule 
$(\theta,x^{(2)})$, and $\mm=(\beta,m)$ is any nontrivial 
irreducible $\BA$-$\BA$-bimodule, then 
$x^*\scirc (1_{\theta}\times s^*\times 1_{\theta})\scirc m=
\tikzmatht{
  \fill[black!10] (-1.5,-1.3) rectangle (1.5,1.2);
    \draw[thick] (0,.7)--(-.17,.7) arc(90:270:.7) (0,.7)--(.17,.7)
    arc(90:-90:.7); 
  \fill[black] (0,.7) circle(.12);
    \draw[thick] (0,-.7) circle(.17);
    \draw[thick] (0,-1.3)--(0,-.1) (0,.7)--(0,1.2); 
  \fill[white] (0,.2)--(.2,-.1)--(-.2,-.1)--(0,.2);
    \draw[thick] (0,.2)--(.2,-.1)--(-.2,-.1)--(0,.2);
} = 0$ for
every $s\in\Hom(\id,\beta)$. This is a special case of the Lemma (with
$t=w\scirc s^*\in\Hom(\beta,\theta)$), that we shall make use of in
\sref{s:class}.

\section{Q-system calculus}
\label{s:Qsystem}
\setcounter{equation}{0}
Throughout this section, $N$ is an infinite factor, and
$\C\subset\End_0(N)$ with properties as specified in \sref{s:frob}.

Q-systems in $\C$ can be decomposed in several distinct ways. In the
first four subsections, we  
discuss various decompositions in turn, and characterize them 
in terms of suitable projections in the underlying category $\C$. 

In the remainder of this section, we discuss Q-systems in braided C*
tensor categories, introduce various operations with Q-systems (the 
centres, the braided products and the full centre), and compute the 
central decomposition of the extension corresponding to the braided
product of two full centres. The latter is motivated because this
decomposition gives the irreducible boundary conditions for phase
boundaries in local QFT \cite{BKLR}.

\subsection{Reduced Q-systems} 
\label{s:Qred}
\setcounter{equation}{0}
Let $(\theta,w,x)$ be a Q-system describing the extension $N\subset M$.
When the multiplicity $\dim\Hom(\id_N,\theta)=\dim\Hom(\iota,\iota)$ of 
$\id_N$ in $\theta$ is one, then the extension is irreducible 
($N'\cap M=\CC\eins$), and in particular $M$ is automatically a
factor. When the multiplicity is larger than one, then $M$ may or may 
not be a factor. 

Let $e$ be a nontrivial projection in $\Hom(\iota,\iota)$. If $M$ is 
a factor, then one can write $e=tt^*$ with an isometry $t\in M$,
define a sub-homomorphism $\iota_e\prec\iota$ by $\iota_e(\cdot)=
t^*\iota(\cdot)t$, and arrive at a decomposition $[\iota]=[\iota_e] \oplus
[\iota_{1-e}]$, cf.\ \cref{c:subQ}, where we shall characterize this
decomposition in terms of certain projections in $\Hom(\theta,\theta)$. 
In contrast, if $M$ is not a factor and $e\neq1$ belongs to the centre 
of $M$, such isometries $t$ do not exist in $M$. Namely, $t\in M$ and
$tt^*\in M'$ would imply $e=et^*t=t^*et=\eins_M$. 

One should therefore first perform a central decomposition of $M$ 
into factors $M_e=eM$ by the minimal central projections, and 
compute the reduced Q-systems (cf.\ \sref{s:Qcentral}) for the subfactors 
$N\simeq Ne\subset Me$. Each of these may still be reducible, and can 
be further reduced by decomposing $[\iota]=\bigoplus_e[\iota_e]$, as before. 

Finally, we also discuss in \sref{s:interm} the multiplicative 
``splitting'' decomposition of $\iota$, when there is an intermediate 
subfactor $N\subset L\subset M$, so that $\iota=\iota_2\circ\iota_1$. 
Whether an intermediate subfactor exists is independent of
reducibility of the subfactor, e.g., $[\iota] = [\id_N]\oplus [\id_N]$ 
is reducible but does not admit an intermediate subfactor, whereas 
$\iota_1\otimes\iota_2 =
(\iota_1\otimes\id_{N_2})\circ(\id_{N_1}\otimes\iota_2):N_1\otimes
N_2\to M_1\otimes M_2$ is an irreducible homomorphism (if $\iota_i$
are) but admits intermediate factors $N_1\otimes M_2$ and $M_1\otimes
N_2$. This example also shows that the splitting cannot be expected to
be unique. Also, even if both $N$ and $M$ are factors, the
intermediate algebra $L$ need not be a factor, as the example
$N\subset N\oplus N\subset \Mat_2(N)$ shows.

Although the three decompositions of a Q-system $(\theta,w,x)$ are of 
quite different nature, they all come with a projection 
$P\in\Hom(\theta,\theta)$ satisfying 
$$
(P\times P)\scirc x = (P\times 1_\theta)\scirc x\scirc P=(1_\theta\times
P)\scirc x\scirc P
$$
\be\label{2to3}
\tikzmatht{
  \fill[black!10] (-1.5,-1.4) rectangle (1.5,1.4);
    \draw[thick] (-1,1.4)--(-1,.8) arc(180:360:1)--(1,1.4);
    \draw[thick] (0,-1.4)--(0,-.2); 
    \draw[thick] (-1.2,.6)--(-0.8,1); 
    \draw[thick] (-1.2,1)--(-0.8,.6); 
    \draw[thick] (1.2,.6)--(0.8,1); 
    \draw[thick] (1.2,1)--(0.8,.6); 
  \fill[black] (0,-.2) circle (.12);
}=
\tikzmatht{
  \fill[black!10] (-1.5,-1.4) rectangle (1.5,1.4);
    \draw[thick] (-1,1.4)--(-1,.8) arc(180:360:1)--(1,1.4);
    \draw[thick] (0,-1.4)--(0,-.2); 
    \draw[thick] (-1.2,.6)--(-0.8,1); 
    \draw[thick] (-1.2,1)--(-0.8,.6); 
    \draw[thick] (-0.2,-0.6)--(0.2,-1); 
    \draw[thick] (-0.2,-1)--(0.2,-0.6); 
  \fill[black] (0,-.2) circle (.12);
}=
\tikzmatht{
  \fill[black!10] (-1.5,-1.4) rectangle (1.5,1.4);
    \draw[thick] (-1,1.4)--(-1,.8) arc(180:360:1)--(1,1.4);
    \draw[thick] (0,-1.4)--(0,-.2); 
    \draw[thick] (1.2,.6)--(0.8,1); 
    \draw[thick] (1.2,1)--(0.8,.6); 
    \draw[thick] (-0.2,-0.6)--(0.2,-1); 
    \draw[thick] (-0.2,-1)--(0.2,-0.6); 
  \fill[black] (0,-.2) circle (.12);
}
\ee
(``of three projection, any one is redundant'', or ``any two
projections imply the third''), and in each case
different further properties. One easily proves: 

\begin{lemma} \label{l:red}
\eref{2to3} alone implies  \\ 
{\rm (i)} The triple  
$$\theta_P:=S^*\theta(\cdot)S,\quad \wt w_P:=S^* \scirc
w,\quad \wt x_P:=(S^*\times S^*)\scirc x\scirc S
$$
is a C* Frobenius algebra, where $S\in N$ is any isometry such
that $SS^*=P$, i.e., $\theta_P\prec\theta$. \\
{\rm (ii)} 
$n_P:=\wt x^*_P{}\ncirc \wt x_P$ is a multiple of $1_{\theta_P}$ if and only if $x^*\scirc (P\times P) \scirc x$ is a multiple of $P$. 
\end{lemma}

(The notation emphasizes that the unitary equivalence class of
$(\theta_P,\wt w_P,\wt x_P)$ depends on $P$, but not on the choice of the
isometry $S$.) 

\medskip

{\em Proof:} (i) The unit property, associativity and Frobenius
property follow from the corresponding properties of $(\theta,w,x)$ by
``eliminating'' projections using \eref{2to3} and $PS=S$. \\
(ii) ``If'' is obvious. Conversely, $\wt x^*_P\ncirc \wt x_P=\mu\cdot
1_{\theta_P}$ implies that $P\scirc x^*\scirc (P\times P) \scirc x
\scirc P
= \mu\cdot P$. By \eref{2to3}, this equals $x^*\scirc (P\times P) \scirc x$. 
\qed

However, the property in (ii) may fail, in which case the C* Frobenius 
algebra fails to be special (and hence to be standard).  
Then, by \cref{c:special}, one can define the equivalent special 
C* Frobenius algebra $(\theta_P,\wh w_P,\wh x_P)$ with 
$$\wh w_P := n_P^{\frac12}\ncirc S^*\scirc w,\qquad \wh x_P :=
(n_P^{-\frac12}\ncirc S^*\times
n_P^{-\frac12}\ncirc S^*)\scirc x \scirc S\ncirc n_P^{\frac12}$$
with the ``normalization intertwiner'' $n_P=\wt x^*_P\ncirc \wt
x_P=S^*\ncirc X^*\ncirc (P\times P)\ncirc X\ncirc
S\in\Hom_0(\theta_P,\theta_P)$, such that  
$\wh x^*_P\ncirc \wh x_P = 1_{\theta_P}$. 

$\wh w_P\in\Hom(\id,\theta_P)$ is automatically a multiple of an
isometry, and 
$$\wh w^*_P\ncirc \wh w_P = w^*\scirc x^*\scirc (P\times P)\scirc x \scirc w= r^*\scirc (P\times P)\scirc r.$$


\begin{tintedbox}
\begin{corollary} \vskip-5mm \label{c:norm}
The appropriately rescaled triple  
\be\label{redQ}
\ba{c}
\theta_P=S^*\theta(\cdot)S, \\ 
w_P=\dim(\theta_P)^{-\frac14}\cdot 
n_P^{\frac12}\ncirc S^*\scirc w, \\ 
x_P=\dim(\theta_P)^{\frac14}\cdot (n_P^{-\frac12}S^*\times
n_P^{-\frac12}S^*) \scirc x \scirc Sn_P^{\frac12}
\ea
\ee
is a standard C* Frobenius algebra, i.e., a Q-system, called the 
{\bf reduced Q-system}, if and only if $\wh w_P$
has the correct normalization $\wh w^*_P\ncirc \wh w_P =r^*\scirc
(P\times P)\scirc r \stackrel !=\dim(\theta_P)$. 
\end{corollary}
\end{tintedbox}

In each of the three decompositions discussed in the subsequent
subsections, further properties of the characterizing projections
beyond \eref{2to3} will indeed ensure the correct normalization as 
required in \cref{c:norm}. 

\subsection{Central decomposition of Q-systems}
\label{s:Qcentral}
\setcounter{equation}{0}
In this subsection, we shall characterize decompositions of $\iota:
N\to M$ as a direct sum 
$$\iota = \iota_1 \oplus \iota_2$$
when $M=M_1\oplus M_2$ is not a factor, $\iota_i:N\to M_i$.

Let $N\subset M$ an inclusion of von Neumann algebras, where $N$ is an
infinite factor, and $M$ is properly infinite with a finite
centre. Every central projection  $e\in M'\cap M$ gives rise to an
inclusion $eN\subset eM$, where $eN$ is canonically isomorphic to
$N$. (Recall also the characterization of such projections as
$e=\iota(q)v$ with $q\in\Hom(\theta,\id_N)$, given in \lref{l:relc+cent}.) 
If $e$ is minimal, $eN\subset eM$ is a subfactor. 

We want to characterize the corresponding embeddings $eN\to eM$ in
terms of reduced Q-systems \eref{redQ}. Our starting observation is
that $P:=\ol\iota(e)\in N$ is a projection in $\Hom_0(\theta,\theta)$:
$$ (1_\theta\times P)\scirc x =
x\scirc P = (P\times 1_\theta)\scirc x: \qquad
\tikzmatht{
  \fill[black!10] (-1.5,-1.3) rectangle (1.5,1.3);
    \draw[thick] (-1,1.3)--(-1,.5) arc(180:360:1)--(1,1.3);
    \draw[thick] (0,-1.3)--(0,-.5); 
    \draw[thick] (1.2,.4)--(0.8,.8); 
    \draw[thick] (1.2,.8)--(0.8,.4); 
  \fill[black] (0,-.5) circle (.12);
}=
\tikzmatht{
  \fill[black!10] (-1.5,-1.3) rectangle (1.5,1.3);
    \draw[thick] (-1,1.3)--(-1,1) arc(180:360:1)--(1,1.3);
    \draw[thick] (0,-1.3)--(0,0); 
    \draw[thick] (-0.2,-.4)--(0.2,-.8); 
    \draw[thick] (-0.2,-.8)--(0.2,-.4); 
  \fill[black] (0,0) circle (.12);
}=
\tikzmatht{
  \fill[black!10] (-1.5,-1.3) rectangle (1.5,1.3);
    \draw[thick] (-1,1.3)--(-1,.5) arc(180:360:1)--(1,1.3);
    \draw[thick] (0,-1.3)--(0,-.5); 
    \draw[thick] (-1.2,.4)--(-0.8,.8); 
    \draw[thick] (-1.2,.8)--(-0.8,.4); 
  \fill[black] (0,-.5) circle (.12);
},$$
which is precisely \eref{1to3}.
This follows immediately by applying $\ol\iota$ to the equations
$\iota\ol\iota(e)v=ve$ (because $e\in M$) and $ev=ve$ (because $e\in
M'$). We now show the converse:

\begin{tintedbox}
\begin{proposition} \vskip-5mm \label{p:centraldec}
Let $\BA=(\theta,w,x)$ be a Q-system defining the extension $N\subset
M$. Let $P\in\Hom(\theta,\theta)$ be a projection satisfying
\eref{1to3}, hence also \eref{2to3}. Then \eref{redQ} with
normalization intertwiner $n_P=\sqrt{\dim(\theta)}\cdot 1_{\theta_p}$ defines a
reduced Q-system $\BA_P$. The reduced Q-system corresponds 
to the extension $eN\subset eM$ where $e\in M'\cap M$ and 
$\ol\iota(e)=P$.  

Along with $P$, also $1-P$ satisfies \eref{1to3}.
\end{proposition}
\end{tintedbox}

If $P$ is a {\em minimal} projection in $\Hom(\theta,\theta)$ with the
stated properties, we will also refer to the reduced Q-system as a
{\bf factor Q-system} of $\BA$. 

\medskip

{\em Proof:} 
Let us first compute the normalizations. Let $S\in N$ be any
isometry such that $P=SS^*$. Because $P\in\Hom_0(\theta,\theta)$, we 
have $n_P= S^*\ncirc x^*\ncirc (P\times P)\ncirc x\ncirc S = S^*\ncirc P \ncirc x^*\ncirc x\ncirc
P\ncirc S = \sqrt{\dim(\theta)}\cdot 1_\theta$, and $r^*\ncirc (P\times
P)\ncirc r = r^*\ncirc (1_\theta \times P)\ncirc r = \Tr_\theta(P) =
\dim(\theta_P)$ by \pref{p:traces}. Hence, by \cref{c:norm}, 
$\BA_P=(\theta_P,w_P,x_P)$ is a reduced Q-system.  

By \lref{l:selfint}, $P$ is a self-morphism
of $\BA$ considered as an $\BA$-$\BA$-bimodule. 
Hence, by \pref{p:bimfunct}, $P=1_{\ol\iota}\times e\times 1_\iota$
with $e\in\Hom(\id_M,\id_M)=M'\cap M$. If $P$ is a projection, so is
$e$. We claim that $n\mapsto en\equiv e\iota(n)$ is a *-isomorphism 
between $N$ and $eN$. Because $e$ is a central projection, the
*-homomorphism property is obvious, and so is
surjectivity. Injectivity follows because $e\iota(n)=0$ implies
$P\theta(n)=0$, hence $\theta_P(n):=S^*\theta(n)S=0$. Since $\theta_P$ is injective, $n=0$.

We now define a conjugate $\ol\iota_P$ for the embedding
$\iota_P:N\to eM$, $n\mapsto e\iota(n)$:
$$\ol\iota_P:em\mapsto S^*\ol\iota(m)S.$$
Then $\wt w_P:=S^*w \in\Hom(\id_{eN},\ol\iota_P\iota_P)$
and $\wt v_P:=e\iota(S^*)v\in\Hom(\id_{eM},\iota_P\ol\iota_P)$ are
intertwiners:  
$$\ol\iota_P\iota_P(n)\wt w_P = S^*\theta(n)SS^*w = S^*\theta(n)Pw  =
S^*P\theta(n)w = S^*\theta(n)w =S^*wn = \wt w_Pn,$$
$$\iota_P\ol\iota_P(em)\wt v_P =
e\iota(S^*)\iota\ol\iota(m)\iota(SS^*)v = e\iota(S^*P)\iota\ol\iota(m)v = 
e\iota(S^*)v m = \wt v_Pem,$$
because $\iota(SS^*)=\iota(P)=\iota\ol\iota(e)$ commutes with
$=\iota\ol\iota(m)$. $\wt w_P$ and $\wt v_P$ solve the conjugacy
relations \eref{conj}:  
$$\wt w_P^*\ol\iota_P(\wt v_P) = w^*SS^*\ol\iota(e\iota(S^*)v)S  
= w^*P\theta(S^*)xS = w^*\theta(S^*)xS = S^*w^*xS = S^*S=\eins_N,$$
because $P$ commutes with $\theta(S^*)$ and with $x$, and  
$$\ol\iota[\iota_P(\wt w_P^*)\wt v_P] =
\ol\iota(e\iota(w^*S)\iota(S^*)v) = P\theta(w^*P)x =  P\theta(w^*)xP =
P^2=P = \ol\iota(e), 
$$
which implies $\iota_P(\wt w_P^*)\wt v_P=e=1_{eM}$. Finally,
$\wt x_P:=\ol\iota_P(\wt v_P)$ equals $S^*\ol\iota(e)\theta(S^*)xS = 
(S^*\times S^*)\ncirc x\ncirc S$ because $\ol\iota(e)=P$. Thus, after 
the appropriate rescaling by $\dim(\theta_P)^{\mp\frac14}\cdot 
n_P^{\pm\frac12} = (\dim(\theta_P)/\dim(\theta))^{\mp\frac14}$, 
the Q-system for $eN\subset eM$ coincides with the reduced Q-system 
$\BA_P=(\theta_P,w_P,x_P)$. 

The last statement is obvious by linearity. 
\qed

\begin{corollary} \label{c:centraldec}
If $1_\theta=\sum P_i$ is the partition of unity into minimal projections
in $\Hom(\theta,\theta)$ satisfying \eref{1to3}, then $(\theta,w,x)$
is the direct sum of simple Q-systems as in \eref{Qsum}. The
corresponding partition $\eins_M=\sum e_i$ gives the decomposition of
$M'\cap M$ into minimal central projections, i.e., each simple
extension $e_iN\subset e_iM$ is a representation of the extension
$N\subset M$.   
\end{corollary}

\subsection{Irreducible decomposition of Q-systems}
\label{s:Qirred}
\setcounter{equation}{0}
In this subsection, we shall characterize decompositions of
$\iota:N\to M$ as a direct sum of sectors 
$$[\iota] = [\iota_1] \oplus [\iota_2], \qquad \hbox{i.e.}, \quad
\iota(\cdot) = s_1\iota_1(\cdot)s_1^* + s_2\iota_2(\cdot)s_2^*$$
for infinite subfactors $N\subset M$.

Thus, let both $N$ and $M$ be factors, i.e., there are no nontrivial
projections in $\Hom(\theta,\theta)$ satisfying \eref{1to3}.

If $\iota(N)' \cap M=\Hom(\iota,\iota)$ is nontrivial, then $\iota$ is a 
reducible homomorphism. If $e\in \iota(N)'\cap M$ is a projection,
then there is an isometry $s\in M$ such that $ss^*=e$, and
$\iota_s(n) = s^*\iota(n)s$ is a sub-homomorphism of $\iota$. Clearly 
$s\in\Hom(\iota_s,\iota)$. 

\begin{lemma}\label{l:subhom} The homomorphism $\iota_s:N\to M$ is
  isomorphic to the embedding $eN\equiv Ne\subset eMe$, i.e., the
  identical map $\iota_e:eN\to eMe$.  
\end{lemma}

{\em Proof:} We may write $eN\subset eMe$ as $ss^*\iota(N)ss^*\equiv
s\iota_s(N)s^*\subset sMs^*$. The claim follows because the map
$\Ad_s:M\to sMs^*\equiv eMe$ is an isomorphism. 
\qed

For $\iota_s\prec\iota$ one has a conjugate $\ol\iota_{\ol s}
\prec\ol\iota$, and an isometry $\ol s\in\Hom(\ol\iota_s,\ol\iota)$. 
Then $e=s\ncirc s^*\in\Hom(\iota,\iota)\subset M$ and $\ol e = \ol
s\ncirc\ol s^* \in\Hom(\ol\iota,\ol\iota)\subset N$ are projections
such that 
\be \label{conjug} 
(\ol e\times 1_{\iota})\scirc w = (1_{\ol\iota}\times e)\scirc
w,\qquad (e\times 1_{\ol\iota})\scirc v = (1_\iota\times \ol e)\scirc
v,\ee
and $w^*\scirc(\ol e\times e)\scirc w = \dim(\iota_s)\cdot\eins_N$,
$v^*\scirc(e\times \ol e)\scirc v = \dim(\iota_s)\cdot\eins_M$, 
Then $p=1_{\ol\iota}\times e$ and $\ol p=\ol e\times
1_{\iota}$ are a pair of commuting projections 
$\in\Hom(\theta,\theta)$ such that 
\be \label{1to1a}
\tikzmatht{
  \fill[black!10] (-1,-1.4) rectangle (1,1.2);
    \draw[thick] (0,1.2)--(0,-.5);
    \draw[thick] (-.2,.2)--(.2,.6); 
    \draw[thick] (.2,.2)--(-.2,.6); 
  \fill[white] (0,-.5) circle(.2);
    \draw[thick] (0,-.5) circle(.2);
    \node at (.5,.5) {$p$}; 
}=
\tikzmatht{
  \fill[black!10] (-1,-1.4) rectangle (1,1.2);
    \draw[thick] (0,1.2)--(0,-.5);
    \draw[thick] (-.2,.2)--(.2,.6); 
    \draw[thick] (.2,.2)--(-.2,.6); 
  \fill[white] (0,-.5) circle(.2);
    \draw[thick] (0,-.5) circle(.2);
    \node at (-.5,.5) {$\ol p$};
},
\ee
\be \label{1to1b}
\tikzmatht{
  \fill[black!10] (-1.5,-.8) rectangle (1.7,1.8);
    \draw[thick] (-.8,1.8)--(-.8,1) arc(180:360:.8)--(.8,1)--(.8,1.8);
    \draw[thick] (0,-.8)--(0,.2); 
    \draw[thick] (1,1.1)--(0.6,1.5); 
    \draw[thick] (1,1.5)--(0.6,1.1); 
  \fill[black] (0,.2) circle (.12);
    \node at (1.3,1.3) {$p$}; 
}=
\tikzmatht{
  \fill[black!10] (-1.5,-1.1) rectangle (1.5,1.5);
    \draw[thick] (-.8,1.5)--(-.8,1) arc(180:360:.8)--(.8,1.5);
    \draw[thick] (0,-1.1)--(0,.2); 
    \draw[thick] (-0.2,-.2)--(0.2,-.6); 
    \draw[thick] (-0.2,-.6)--(0.2,-.2); 
  \fill[black] (0,.2) circle (.12);
    \node at (.5,-.4) {$p$}; 
},\qquad
\tikzmatht{
  \fill[black!10] (-1.7,-.8) rectangle (1.5,1.8);
    \draw[thick] (-.8,1.8)--(-.8,1) arc(180:360:.8)--(.8,1.8);
    \draw[thick] (0,-.8)--(0,.2); 
    \draw[thick] (-1,1.1)--(-.6,1.5); 
    \draw[thick] (-1,1.5)--(-.6,1.1); 
  \fill[black] (0,.2) circle (.12);
    \node at (-1.3,1.3) {$\ol p$};
}=
\tikzmatht{
  \fill[black!10] (-1.5,-1.1) rectangle (1.5,1.5);
    \draw[thick] (-.8,1.5)--(-.8,1) arc(180:360:.8)--(.8,1.5);
    \draw[thick] (0,-1.1)--(0,.2); 
    \draw[thick] (-0.2,-.2)--(0.2,-.6); 
    \draw[thick] (-0.2,-.6)--(0.2,-.2); 
  \fill[black] (0,.2) circle (.12);
    \node at (-.5,-.4) {$\ol p$};
},\qquad
\tikzmatht{
  \fill[black!10] (-1.5,-.8) rectangle (1.5,1.8);
    \draw[thick] (-.8,1.8)--(-.8,1) arc(180:360:.8)--(.8,1.8);
    \draw[thick] (0,-.8)--(0,.2); 
    \draw[thick] (-1,1.1)--(-.6,1.5); 
    \draw[thick] (-1,1.5)--(-.6,1.1); 
  \fill[black] (0,.2) circle (.12);
    \node at (-.3,1.3) {$p$}; 
}=
\tikzmatht{
  \fill[black!10] (-1.5,-.8) rectangle (1.5,1.8);
    \draw[thick] (-.8,1.8)--(-.8,1) arc(180:360:.8)--(.8,1.8);
    \draw[thick] (0,-.8)--(0,.2); 
    \draw[thick] (1,1.1)--(.6,1.5); 
    \draw[thick] (1,1.5)--(.6,1.1); 
  \fill[black] (0,.2) circle (.12);
    \node at (.3,1.3) {$\ol p$};
}.
\ee

Conversely, if $\BA=(\theta,w,x)$ is a simple Q-system, and
$\dim\Hom(\id,\theta)>1$, then by Frobenius reciprocity also
$\dim\Hom(\iota,\iota)>1$, hence $\iota$ is reducible. In order to
decompose $\iota$, we want to characterize the projections in
$\Hom(\iota,\iota)$ in terms of projections in $\Hom(\theta,\theta)$.

Instead of characterizing $\jmath\prec\iota$ by the pair of
projections $p,\ol p$ satisfying \eref{1to1a} and \nref{1to1b}, we
observe that either $p$ or $\ol p$ suffices: namely, from the third
relation in \nref{1to1b}, one can express $p$ in
terms of $\ol p$, and vice versa:  
\be\label{ppbar}
\tikzmatht{
  \fill[black!10] (-1.3,-1) rectangle (1.3,1.8);
    \draw[thick] (0,-1)--(0,1.8); 
    \draw[thick] (-0.2,.6)--(0.2,.2); 
    \draw[thick] (-0.2,.2)--(0.2,.6); 
    \node at (.5,.4) {$p$};
} =\tikzmatht{
  \fill[black!10] (-1.7,-1) rectangle (1.7,1.8);
    \draw[thick] (-1,1.8)--(-1,1) arc(180:360:1)--(1,1)--(1,1.5);
    \draw[thick] (0,-1)--(0,0); 
    \draw[thick] (.8,0.6)--(1.2,1); 
    \draw[thick] (.8,1)--(1.2,0.6); 
    \node at (.5,.8) {$\ol p$};
  \fill[black] (0,0) circle (.12);
  \fill[white] (1,1.5) circle(.2);
    \draw[thick] (1,1.5) circle(.2);
} , \qquad
\tikzmatht{
  \fill[black!10] (-1.3,-1) rectangle (1.3,1.8);
    \draw[thick] (0,-1)--(0,1.8); 
    \draw[thick] (-0.2,.6)--(0.2,.2); 
    \draw[thick] (-0.2,.2)--(0.2,.6); 
    \node at (-.5,.4) {$\ol p$};
} =\tikzmatht{
  \fill[black!10] (-1.7,-1) rectangle (1.7,1.8);
    \draw[thick] (1,1.8)--(1,1) arc(360:180:1)--(-1,1)--(-1,1.5);
    \draw[thick] (0,-1)--(0,0); 
    \draw[thick] (-.8,0.6)--(-1.2,1); 
    \draw[thick] (-.8,1)--(-1.2,0.6); 
    \node at (-.5,.8) {$p$};
  \fill[black] (0,0) circle (.12);
  \fill[white] (-1,1.5) circle(.2);
    \draw[thick] (-1,1.5) circle(.2);
}. 
\ee
Expressing $p$ in terms of $\ol p$ as in \eref{ppbar}, turns
\eref{1to1a} into another relation for $\ol p$  
\be\label{pbar}
\tikzmatht{
  \fill[black!10] (-1.3,-1) rectangle (1.3,1.5);
    \draw[thick] (0,-.5)--(0,1.5); 
    \draw[thick] (-0.2,.6)--(0.2,.2); 
    \draw[thick] (-0.2,.2)--(0.2,.6); 
    \node at (-.5,.4) {$\ol p$};
  \fill[white] (0,-.5) circle(.2);
    \draw[thick] (0,-.5) circle(.2);
} =\tikzmatht{
  \fill[black!10] (-1.3,-1) rectangle (1.3,1.5);
    \draw[thick] (-.6,1.5)--(-.6,.2) arc(180:360:.6)--(.6,1)--(.6,1);
    \draw[thick] (.8,0.2)--(.4,.6); 
    \draw[thick] (.8,.6)--(.4,0.2); 
    \node at (.1,.4) {$\ol p$};
  \fill[white] (.6,1) circle(.2);
    \draw[thick] (.6,1) circle(.2);
}\ee
besides the second in \nref{1to1b}, while the first is automatically
satisfied. In the same way, expressing $\ol p$ 
in terms of $p$, turns \eref{1to1a} into another relation for $p$
\be\label{p}
\tikzmatht{
  \fill[black!10] (-1.3,-1) rectangle (1.3,1.5);
    \draw[thick] (0,-.5)--(0,1.5); 
    \draw[thick] (-0.2,.6)--(0.2,.2); 
    \draw[thick] (-0.2,.2)--(0.2,.6); 
    \node at (.5,.4) {$p$};
  \fill[white] (0,-.5) circle(.2);
    \draw[thick] (0,-.5) circle(.2);
} =\tikzmatht{
  \fill[black!10] (-1.3,-1) rectangle (1.3,1.5);
    \draw[thick] (.6,1.5)--(.6,.2) arc(360:180:.6)--(-.6,1)--(-.6,1);
    \draw[thick] (-.8,0.2)--(-.4,.6); 
    \draw[thick] (-.8,.6)--(-.4,0.2); 
    \node at (-.1,.4) {$p$};
  \fill[white] (-.6,1) circle(.2);
    \draw[thick] (-.6,1) circle(.2);
}
\ee
besides the first in \nref{1to1b}.  

\begin{lemma} \label{l:ppbar}
Let $\BA=(\theta,w,x)$ be a simple Q-system. 
Let either $p \in\Hom(\theta,\theta)$ be a projection
satisfying \eref{p} and the first of \eref{1to1b}, or $\ol
p\in\Hom(\theta,\theta)$ a projection satisfying \eref{pbar} and the
second of \eref{1to1b}. Defining $\ol p :=
\tikzmatht{
  \fill[black!10] (-1.3,-1) rectangle (1.3,1.5);
    \draw[thick] (.6,1.5)--(.6,.2) arc(360:180:.6)--(-.6,1);
    \draw[thick] (-.8,0.2)--(-.4,.6); 
    \draw[thick] (-.8,.6)--(-.4,0.2); 
    \node at (-.1,.4) {$p$};
  \fill[white] (-.6,1) circle(.2);
    \draw[thick] (-.6,1) circle(.2);
    \draw[thick] (0,-1)--(0,-.4); 
  \fill[black] (0,-.4) circle (.12);
}$ in the first case, and $p := 
\tikzmatht{
  \fill[black!10] (-1.3,-1) rectangle (1.3,1.5);
    \draw[thick] (-.6,1.5)--(-.6,.2) arc(180:360:.6)--(.6,1);
    \draw[thick] (.8,0.2)--(.4,.6); 
    \draw[thick] (.8,.6)--(.4,0.2); 
    \node at (.1,.4) {$\ol p$};
  \fill[white] (.6,1) circle(.2);
    \draw[thick] (.6,1) circle(.2);
    \draw[thick] (0,-1)--(0,-.4);
  \fill[black] (0,-.4) circle (.12);
}$
in the second case, gives another projection such that $p$ and $\ol p$
satisfy the system \eref{1to1a}, \nref{1to1b}.
\end{lemma}

{\em Proof:}
We first establish that $\ol p$ defined from $p$ is a projection: 
$$\ol p^* =\tikzmatht{ 
  \fill[black!10] (-1.5,-2) rectangle (1.5,.8);
    \draw[thick] (1,-2)--(1,-1) arc(0:180:1)--(-1,-1)--(-1,-1.5);
    \draw[thick] (0,.8)--(0,0); 
    \draw[thick] (-.8,-.6)--(-1.2,-1); 
    \draw[thick] (-.8,-1)--(-1.2,-.6); 
    \node at (-.5,-.8) {$p$};
  \fill[black] (0,0) circle (.12);
  \fill[white] (-1,-1.5) circle(.2);
    \draw[thick] (-1,-1.5) circle(.2);
} = 
\tikzmatht{ 
  \fill[black!10] (-1.7,-1.8) rectangle (1.7,1);
    \draw[thick] (1.2,-1.8)--(1.2,-.4) arc(0:180:.6)
    arc(360:180:.6)--(-1.2,.5); 
    \draw[thick] (0.6,1)--(0.6,.2); 
    \draw[thick] (-1,-.6)--(-1.4,-.2); 
    \draw[thick] (-1,-.2)--(-1.4,-.6); 
    \node at (-.7,-.4) {$p$};
  \fill[black] (.6,.2) circle (.12);
  \fill[white] (-1.2,.5) circle(.2);
    \draw[thick] (-1.2,.5) circle(.2);
}
=\tikzmatht{
  \fill[black!10] (-1.5,-1) rectangle (1.5,1.8);
    \draw[thick] (1,1.8)--(1,1) arc(360:180:1)--(-1,1)--(-1,1.5);
    \draw[thick] (0,-1)--(0,0); 
    \draw[thick] (-.8,0.6)--(-1.2,1); 
    \draw[thick] (-.8,1)--(-1.2,0.6); 
    \node at (-.5,.8) {$p$};
  \fill[black] (0,0) circle (.12);
  \fill[white] (-1,1.5) circle(.2);
    \draw[thick] (-1,1.5) circle(.2);
} = \ol p,$$
$$\ol p^2 = 
\tikzmatht{
  \fill[black!10] (-1.7,-1) rectangle (1.7,1.8);
    \draw[thick] (1.2,1.8)--(1.2,1) arc(360:180:.6)--(0,1.5);
    \draw[thick] (.6,.4) arc(360:180:.9)--(-1.2,1.5);
    \draw[thick] (-.3,-.5)--(-.3,-1); 
    \draw[thick] (-1,0.6)--(-1.4,1); 
    \draw[thick] (-1,1)--(-1.4,0.6); 
    \node at (-.7,.8) {$p$};
    \draw[thick] (-.2,0.8)--(.2,1.2); 
    \draw[thick] (-.2,1.2)--(.2,0.8); 
    \node at (.5,1) {$p$};
  \fill[black] (.6,.4) circle (.12);
  \fill[black] (-.3,-.5) circle (.12);
  \fill[white] (-1.2,1.5) circle(.2);
    \draw[thick] (-1.2,1.5) circle(.2);
  \fill[white] (0,1.5) circle(.2);
    \draw[thick] (0,1.5) circle(.2);
} = 
\tikzmatht{
  \fill[black!10] (-1.7,-1) rectangle (1.7,1.8);
    \draw[thick] (0,1.5)--(0,1) arc(360:180:.6)--(-1.2,1.5);
    \draw[thick] (1.2,1.8)--(1.2,.4) arc(360:180:.9);
    \draw[thick] (.3,-.5)--(.3,-1); 
    \draw[thick] (-1,0.8)--(-1.4,1.2); 
    \draw[thick] (-1,1.2)--(-1.4,0.8); 
    \node at (-.7,1) {$p$};
    \draw[thick] (-.2,0.8)--(.2,1.2); 
    \draw[thick] (-.2,1.2)--(.2,0.8); 
    \node at (.5,1) {$p$};
  \fill[black] (-.6,.4) circle (.12);
  \fill[black] (.3,-.5) circle (.12);
  \fill[white] (-1.2,1.5) circle(.2);
    \draw[thick] (-1.2,1.5) circle(.2);
  \fill[white] (0,1.5) circle(.2);
    \draw[thick] (0,1.5) circle(.2);
} = 
\tikzmatht{
  \fill[black!10] (-1.7,-1) rectangle (1.7,1.8);
    \draw[thick] (0,1.5)--(0,1) arc(360:180:.6)--(-1.2,1.5);
    \draw[thick] (1.2,1.8)--(1.2,.4) arc(360:180:.9);
    \draw[thick] (.3,-.5)--(.3,-1); 
    \draw[thick] (-1,0.8)--(-1.4,1.2); 
    \draw[thick] (-1,1.2)--(-1.4,0.8); 
    \node at (-.7,1) {$p$};
    \draw[thick] (-.3,.2)--(-.7,-.2); 
    \draw[thick] (-.3,-.2)--(-.7,.2); 
    \node at (0,0) {$p$};
  \fill[black] (-.6,.4) circle (.12);
  \fill[black] (.3,-.5) circle (.12);
  \fill[white] (-1.2,1.5) circle(.2);
    \draw[thick] (-1.2,1.5) circle(.2);
  \fill[white] (0,1.5) circle(.2);
    \draw[thick] (0,1.5) circle(.2);
} = 
\tikzmatht{
  \fill[black!10] (-1.5,-1) rectangle (1.5,1.8);
    \draw[thick] (1,1.8)--(1,1) arc(360:180:1)--(-1,1)--(-1,1.5);
    \draw[thick] (0,-1)--(0,0); 
    \draw[thick] (-.8,0.6)--(-1.2,1); 
    \draw[thick] (-.8,1)--(-1.2,0.6); 
    \node at (-.5,.8) {$p$};
  \fill[black] (0,0) circle (.12);
  \fill[white] (-1,1.5) circle(.2);
    \draw[thick] (-1,1.5) circle(.2);
} = \ol p, 
$$
where we have used $p^*=p$ and the first defining property \eref{p} of $p$, and
$p^2=p$ and the second defining property \eref{1to1b} of $p$, respectively.

That $\ol p$ satisfies \eref{1to1a}, is an immediate
consequence of the defining property \eref{p} of $p$. The second of 
\eref{1to1b} is trivial by associativity. It remains to verify the last
of \eref{1to1b}:
$$\tikzmatht{
  \fill[black!10] (-1.7,-1) rectangle (1.7,1.8);
    \draw[thick] (-1,1.8)--(-1,1) arc(180:360:1)--(1,1)--(1,1.8);
    \draw[thick] (0,-1)--(0,0); 
    \draw[thick] (1.2,1.1)--(0.8,1.5); 
    \draw[thick] (1.2,1.5)--(0.8,1.1); 
  \fill[black] (0,0) circle (.12);
    \node at (.5,1.3) {$\ol p$};
}=\tikzmatht{
  \fill[black!10] (-1.7,-1) rectangle (1.7,1.8);
    \draw[thick] (1.2,1.8)--(1.2,1) arc(360:180:.6)--(0,1.5);
    \draw[thick] (.6,.4) arc(360:180:.9)--(-1.2,1.8);
    \draw[thick] (-.3,-.5)--(-.3,-1); 
    \draw[thick] (-.2,0.8)--(.2,1.2); 
    \draw[thick] (-.2,1.2)--(.2,0.8); 
    \node at (.5,1) {$p$};
  \fill[black] (.6,.4) circle (.12);
  \fill[black] (-.3,-.5) circle (.12);
  \fill[white] (0,1.5) circle(.2);
    \draw[thick] (0,1.5) circle(.2);
} = 
\tikzmatht{
  \fill[black!10] (-1.7,-1) rectangle (1.7,1.8);
    \draw[thick] (0,1.5)--(0,1) arc(360:180:.6)--(-1.2,1.8);
    \draw[thick] (1.2,1.8)--(1.2,.4) arc(360:180:.9);
    \draw[thick] (.3,-.5)--(.3,-1); 
    \draw[thick] (-.2,0.8)--(.2,1.2); 
    \draw[thick] (-.2,1.2)--(.2,0.8); 
    \node at (.5,1) {$p$};
  \fill[black] (-.6,.4) circle (.12);
  \fill[black] (.3,-.5) circle (.12);
  \fill[white] (0,1.5) circle(.2);
    \draw[thick] (0,1.5) circle(.2);
} = 
\tikzmatht{
  \fill[black!10] (-1.7,-1) rectangle (1.7,1.8);
    \draw[thick] (0,1.5)--(0,1.2) arc(360:180:.6)--(-1.2,1.8);
    \draw[thick] (1.2,1.8)--(1.2,.4) arc(360:180:.9)--(-.6,.6);
    \draw[thick] (.3,-.5)--(.3,-1); 
    \draw[thick] (-.4,.4)--(-.8,0); 
    \draw[thick] (-.4,0)--(-.8,.4); 
    \node at (0,.1) {$p$};
  \fill[black] (-.6,.6) circle (.12);
  \fill[black] (.3,-.5) circle (.12);
  \fill[white] (0,1.5) circle(.2);
    \draw[thick] (0,1.5) circle(.2);
} = 
\tikzmatht{
  \fill[black!10] (-1.7,-1) rectangle (1.7,1.8);
    \draw[thick] (-1,1.8)--(-1,1) arc(180:360:1)--(1,1)--(1,1.8);
    \draw[thick] (0,-1)--(0,0); 
    \draw[thick] (-1.2,1.1)--(-0.8,1.5); 
    \draw[thick] (-1.2,1.5)--(-0.8,1.1); 
  \fill[black] (0,0) circle (.12);
    \node at (-.5,1.3) {$p$}; 
}.
$$
The properties of $p$ defined from $\ol p$ follow similarly. 
\qed

\begin{lemma}\label{l:ppbarcomm}
If projections $p\in\Hom(\theta,\theta)$ and $\ol p
\in\Hom(\theta,\theta)$ satisfy \eref{1to1a} and \nref{1to1b}, 
then $p$ and $\ol p$ commute, and $P=\ol pp$ is a projection
satisfying \eref{2to3}.  
\end{lemma}

{\em Proof:}
Using in turn \eref{ppbar}, the second and 
the last relation of \eref{1to1b}, one finds
\be\label{commute}
\tikzmatht{
  \fill[black!10] (-1.3,-1.3) rectangle (1.3,1.8);
    \draw[thick] (0,-1.3)--(0,1.8); 
    \draw[thick] (-0.2,0.6)--(0.2,1); 
    \draw[thick] (-0.2,1)--(0.2,0.6); 
    \node at (-.5,.8) {$\ol p$};
    \draw[thick] (-0.2,-0.6)--(0.2,-.2); 
    \draw[thick] (-0.2,-.2)--(0.2,-0.6); 
    \node at (.5,-.4) {$p$};
} =\tikzmatht{
  \fill[black!10] (-1.7,-1.3) rectangle (1.7,1.8);
    \draw[thick] (-1,1.8)--(-1,1) arc(180:360:1)--(1,1)--(1,1.5);
    \draw[thick] (0,-1.3)--(0,0); 
    \draw[thick] (.8,0.6)--(1.2,1); 
    \draw[thick] (.8,1)--(1.2,0.6); 
    \node at (.5,.8) {$\ol p$};
    \draw[thick] (-.8,0.6)--(-1.2,1); 
    \draw[thick] (-.8,1)--(-1.2,0.6); 
    \node at (-1.5,.8) {$\ol p$};
  \fill[black] (0,0) circle (.12);
  \fill[white] (1,1.5) circle(.2);
    \draw[thick] (1,1.5) circle(.2);
} =\tikzmatht{
  \fill[black!10] (-1.7,-1.3) rectangle (1.7,1.8);
    \draw[thick] (-1,1.8)--(-1,1) arc(180:360:1)--(1,1)--(1,1.5);
    \draw[thick] (0,-1.3)--(0,0); 
    \draw[thick] (.8,0.6)--(1.2,1); 
    \draw[thick] (.8,1)--(1.2,0.6); 
    \node at (.5,.8) {$\ol p$};
    \draw[thick] (-0.2,-0.6)--(0.2,-1); 
    \draw[thick] (-0.2,-1)--(0.2,-0.6); 
    \node at (-.5,-.8) {$\ol p$};
  \fill[black] (0,0) circle (.12);
  \fill[white] (1,1.5) circle(.2);
    \draw[thick] (1,1.5) circle(.2);
} =\tikzmatht{
  \fill[black!10] (-1.7,-1.3) rectangle (1.7,1.8);
    \draw[thick] (-1,1.8)--(-1,1) arc(180:360:1)--(1,1)--(1,1.5);
    \draw[thick] (0,-1.3)--(0,0); 
    \draw[thick] (-0.2,-0.6)--(0.2,-1); 
    \draw[thick] (-0.2,-1)--(0.2,-0.6); 
    \node at (-.5,-.8) {$\ol p$};
    \draw[thick] (-.8,0.6)--(-1.2,1); 
    \draw[thick] (-.8,1)--(-1.2,0.6); 
    \node at (-.5,.8) {$p$};
  \fill[black] (0,0) circle (.12);
  \fill[white] (1,1.5) circle(.2);
    \draw[thick] (1,1.5) circle(.2);
} = 
\tikzmatht{
  \fill[black!10] (-1.3,-1.3) rectangle (1.3,1.8);
    \draw[thick] (0,-1.3)--(0,1.8); 
    \draw[thick] (-0.2,0.6)--(0.2,1); 
    \draw[thick] (-0.2,1)--(0.2,0.6); 
    \node at (.5,.8) {$p$};
    \draw[thick] (-0.2,-0.6)--(0.2,-.2); 
    \draw[thick] (-0.2,-.2)--(0.2,-0.6); 
    \node at (-.5,-.4) {$\ol p$};
}.
\ee
It follows that $P=\ol pp$ is a projection, and the relations
\eref{1to1b} immediately imply \eref{2to3} for $P$.
\qed
 
Instead of characterizing either of the projections $p$ or $\ol
p\in\Hom(\theta,\theta)$ as in \lref{l:ppbar}, it is also possible to
characterize directly the projection $P=\ol p\ncirc
p\in\Hom(\theta,\theta)$. 

\begin{proposition} \label{p:subQP}
A projection $P\in\Hom(\theta,\theta)$ is of the form $P=\ol p\ncirc p$
with $p$ and $\ol p$ as in \lref{l:ppbar}, if and only if $P$
satisfies 
\be\label{subQP}
\tikzmatht{
  \fill[black!10] (-1.3,-.8) rectangle (1.3,1.7);
    \draw[thick] (0,1.7)--(0,-.8);
    \draw[thick] (-.8,1.2)--(-.8,.5) arc(180:360:.8) -- (.8,1.2); 
  \fill[white] (-1.1,.9) rectangle (-.5,.4); 
    \draw[thick] (-1.1,.9) rectangle (-.5,.4); 
  \fill[white] (1.1,.9) rectangle (.5,.4); 
    \draw[thick] (1.1,.9) rectangle (.5,.4); 
  \fill[white] (-.8,1.3) circle(.2);
    \draw[thick] (-.8,1.3) circle(.2);
  \fill[white] (.8,1.3) circle(.2);
    \draw[thick] (.8,1.3) circle(.2);
  \fill[black] (0,-.3) circle (.12);
} = 
\tikzmatht{
  \fill[black!10] (-1,-.8) rectangle (1,1.7);
    \draw[thick] (0,1.7)--(0,-.8);
  \fill[white] (-.3,.7) rectangle (.3,.2); 
    \draw[thick] (-.3,.7) rectangle (.3,.2); 
    \node at (-.6,.8) {$P$};
} 
\ee
\end{proposition}

{\em Proof:} ``Only if'': we show that $P=\ol p\ncirc p$ satisfies
\eref{subQP}: 
$$\tikzmatht{
  \fill[black!10] (-1.7,-1) rectangle (1.7,1.7);
    \draw[thick] (0,1.7)--(0,-1);
    \draw[thick] (-1,1.2)--(-1,.4) arc(180:360:1) -- (1,1.2); 
  \fill[white] (-1,1.3) circle(.2);
    \draw[thick] (-1,1.3) circle(.2);
  \fill[white] (1,1.3) circle(.2);
    \draw[thick] (1,1.3) circle(.2);
    \draw[thick] (1.2,1)--(.8,.6);
    \draw[thick] (1.2,.6)--(.8,1);
    \draw[thick] (1.2,0)--(.8,.4);
    \draw[thick] (1.2,.4)--(.8,0);
    \draw[thick] (-1.2,1)--(-.8,.6);
    \draw[thick] (-1.2,.6)--(-.8,1);
    \draw[thick] (-1.2,0)--(-.8,.4);
    \draw[thick] (-1.2,.4)--(-.8,0);
  \fill[black] (0,-.6) circle (.12);
    \node at (-1.4,.2) {$\ol p$};
    \node at (.6,.8) {$\ol p$};
    \node at (1.4,0) {$p$};
    \node at (-.6,.8) {$p$};
} \stapel{\nref{1to1b}}= 
\tikzmatht{
  \fill[black!10] (-1.5,-1) rectangle (1.5,1.7);
    \draw[thick] (0,1.7)--(0,-1);
    \draw[thick] (-1,1.3) arc(180:360:1) -- (1,1.3); 
  \fill[white] (-1,1.4) circle(.2);
    \draw[thick] (-1,1.4) circle(.2);
  \fill[white] (1,1.4) circle(.2);
    \draw[thick] (1,1.4) circle(.2);
    \draw[thick] (-.2,1)--(.2,.6);
    \draw[thick] (-.2,.6)--(.2,1);
    \draw[thick] (-.2,1.2)--(.2,1.6);
    \draw[thick] (-.2,1.6)--(.2,1.2);
    \draw[thick] (-.2,-.8)--(.2,-.4);
    \draw[thick] (-.2,-.4)--(.2,-.8);
    \draw[thick] (-.2,.1)--(.2,-.3);
    \draw[thick] (-.2,-.3)--(.2,.1);
  \fill[black] (0,.3) circle (.12);
    \node at (-.4,-.6) {$\ol p$};
    \node at (-.4,.8) {$\ol p$};
    \node at (.4,1.4) {$p$};
    \node at (.4,-.2) {$p$};
} =
\tikzmatht{
  \fill[black!10] (-1,-1) rectangle (1,1.7);
    \draw[thick] (0,1.7)--(0,-1);
    \draw[thick] (-.2,1)--(.2,.6);
    \draw[thick] (-.2,.6)--(.2,1);
    \draw[thick] (-.2,1.2)--(.2,1.6);
    \draw[thick] (-.2,1.6)--(.2,1.2);
    \draw[thick] (-.2,-.8)--(.2,-.4);
    \draw[thick] (-.2,-.4)--(.2,-.8);
    \draw[thick] (-.2,.1)--(.2,-.3);
    \draw[thick] (-.2,-.3)--(.2,.1);
    \node at (-.4,-.6) {$\ol p$};
    \node at (-.4,.8) {$\ol p$};
    \node at (.4,1.4) {$p$};
    \node at (.4,-.2) {$p$};
} \stapel{\nref{commute}}=  
\tikzmatht{
  \fill[black!10] (-1,-1) rectangle (1,1.7);
    \draw[thick] (0,1.7)--(0,-1);
    \draw[thick] (-.2,1)--(.2,.6);
    \draw[thick] (-.2,.6)--(.2,1);
    \draw[thick] (-.2,.2)--(.2,-.2);
    \draw[thick] (-.2,-.2)--(.2,.2);
    \node at (-.4,0) {$\ol p$};
    \node at (.4,.8) {$p$};
} 
.$$

``If'': We first show that \eref{subQP} implies further
identities. Namely, we obviously get by the unit property: 
\be \label{idem}
\tikzmatht{
  \fill[black!10] (-1.3,-.8) rectangle (1.3,1.7);
    \draw[thick] (0,1.7)--(0,-.4);
    \draw[thick] (-.8,1.2)--(-.8,.4) arc(180:360:.8) -- (.8,1.2); 
  \fill[white] (-1.1,.8) rectangle (-.5,.3); 
    \draw[thick] (-1.1,.8) rectangle (-.5,.3); 
  \fill[white] (1.1,.8) rectangle (.5,.3); 
    \draw[thick] (1.1,.8) rectangle (.5,.3); 
  \fill[white] (-.8,1.3) circle(.2);
    \draw[thick] (-.8,1.3) circle(.2);
  \fill[white] (.8,1.3) circle(.2);
    \draw[thick] (.8,1.3) circle(.2);
  \fill[black] (0,-.4) circle (.12);
} = 
\tikzmatht{
  \fill[black!10] (-1,-.8) rectangle (1,1.7);
    \draw[thick] (0,1.7)--(0,0);
  \fill[white] (-.3,1) rectangle (.3,.5); 
    \draw[thick] (-.3,1) rectangle (.3,.5); 
  \fill[white] (0,0) circle(.2);
    \draw[thick] (0,0) circle(.2);
},\quad\hbox{and}\quad
\tikzmatht{
  \fill[black!10] (-1.3,-.8) rectangle (1.3,1.7);
    \draw[thick] (0,-.3)--(0,-.8);
    \draw[thick] (-.8,1.2)--(-.8,.5) arc(180:360:.8) -- (.8,1.2); 
  \fill[white] (-1.1,.9) rectangle (-.5,.4); 
    \draw[thick] (-1.1,.9) rectangle (-.5,.4); 
  \fill[white] (1.1,.9) rectangle (.5,.4); 
    \draw[thick] (1.1,.9) rectangle (.5,.4); 
  \fill[white] (-.8,1.3) circle(.2);
    \draw[thick] (-.8,1.3) circle(.2);
  \fill[white] (.8,1.3) circle(.2);
    \draw[thick] (.8,1.3) circle(.2);
  \fill[white] (-.8,1.3) circle(.2);
    \draw[thick] (-.8,1.3) circle(.2);
  \fill[white] (.8,1.3) circle(.2);
    \draw[thick] (.8,1.3) circle(.2);
  \fill[black] (0,-.3) circle (.12);
} = 
\tikzmatht{
  \fill[black!10] (-1,-.8) rectangle (1,1.7);
    \draw[thick] (0,1)--(0,-.8);
  \fill[white] (-.3,-.1) rectangle (.3,.4); 
    \draw[thick] (-.3,-.1) rectangle (.3,.4); 
  \fill[white] (0,1) circle(.2);
    \draw[thick] (0,1) circle(.2);
} . 
\ee
Moreover, by \pref{p:dotstar}, we have
\be\label{aux}
\tikzmatht{
  \fill[black!10] (-1.3,-1.8) rectangle (1.3,.7);
    \draw[thick] (0,-.3) arc(180:0:.4); 
  \fill[white] (-.3,-.3) rectangle (.3,-.8); 
    \draw[thick] (-.3,-.3) rectangle (.3,-.8); 
    \draw[thick] (0,-.8) arc(360:180:.4)--(-.8,.7)
    (.8,-.3)--(.8,-1.8); 
}=
\tikzmatht{
  \fill[black!10] (-1.3,-.7) rectangle (1.3,1.8);
    \draw[thick] (0,.3) arc(180:360:.4); 
  \fill[white] (-.3,.3) rectangle (.3,.8); 
    \draw[thick] (-.3,.3) rectangle (.3,.8); 
    \draw[thick] (0,.8) arc(0:180:.4)--(-.8,-.7)
    (.8,.3)--(.8,1.8); 
}\quad \RA\quad 
\tikzmatht{
  \fill[black!10] (-1.9,-1.8) rectangle (1.2,.7);
    \draw[thick] (0,0) arc(180:0:.4); 
  \fill[white] (-.3,0) rectangle (.3,-.5); 
    \draw[thick] (-.3,0) rectangle (.3,-.5); 
    \draw[thick] (0,-.5) arc(360:180:.4)--(-.8,.7)
    (.8,0)--(.8,-.4) arc(360:180:1.1) -- (-1.4,.7); 
}=
\tikzmatht{
  \fill[black!10] (-.7,-.7) rectangle (1.3,1.8);
    \draw[thick] (0,.3) arc(180:360:.4); 
  \fill[white] (-.3,.3) rectangle (.3,.8); 
    \draw[thick] (-.3,.3) rectangle (.3,.8); 
    \draw[thick] (0,1.8)--(0,.8) (.8,1.8)--(.8,.3);
}  \quad \RA\quad
\tikzmatht{
  \fill[black!10] (-2,-.8) rectangle (1.4,1.7);
    \draw[thick] (.2,1.2) arc(180:0:.4); 
    \draw[thick] (.2,.7)--(.2,.2) arc(360:180:.5);
    \draw[thick] (1,1.2)--(1,.7) arc(360:180:1.4) -- (-1.8,1.7); 
    \draw[thick] (-1.2,1.7)--(-1.2,.6) arc(180:360:.4)--(-.4,1.7);
  \fill[white] (-1.5,1.2) rectangle (-.9,.7); 
    \draw[thick] (-1.5,1.2) rectangle (-.9,.7); 
  \fill[white] (-.1,1.2) rectangle (.5,.7); 
    \draw[thick] (-.1,1.2) rectangle (.5,.7); 
  \fill[black] (-.8,.2) circle (.12);}=
\tikzmatht{
  \fill[black!10] (-2.5,-.8) rectangle (.2,1.7);
    \draw[thick] (-.8,.3) arc(360:180:.6) -- (-2,1.7); 
    \draw[thick] (-1.2,1.7)--(-1.2,.7) arc(180:360:.4)--(-.4,1.7);
  \fill[white] (-2.3,1.2) rectangle (-1.7,.7); 
    \draw[thick] (-2.3,1.2) rectangle (-1.7,.7); 
  \fill[white] (-1.5,1.2) rectangle (-.9,.7); 
    \draw[thick] (-1.5,1.2) rectangle (-.9,.7); 
  \fill[black] (-.8,.3) circle (.12);},
\ee
implying
$$
\tikzmatht{
  \fill[black!10] (-1,-.8) rectangle (1,1.7);
    \draw[thick] (0,1.7)--(0,0);
  \fill[white] (-.3,1) rectangle (.3,.5); 
    \draw[thick] (-.3,1) rectangle (.3,.5); 
  \fill[white] (0,0) circle(.2);
    \draw[thick] (0,0) circle(.2);
}=
\tikzmatht{
  \fill[black!10] (-1.3,-.8) rectangle (1.3,1.7);
    \draw[thick] (0,1.7)--(0,-.4);
    \draw[thick] (-.8,1.2)--(-.8,.4) arc(180:360:.8) -- (.8,1.2); 
  \fill[white] (-1.1,.8) rectangle (-.5,.3); 
    \draw[thick] (-1.1,.8) rectangle (-.5,.3); 
  \fill[white] (1.1,.8) rectangle (.5,.3); 
    \draw[thick] (1.1,.8) rectangle (.5,.3); 
  \fill[white] (-.8,1.3) circle(.2);
    \draw[thick] (-.8,1.3) circle(.2);
  \fill[white] (.8,1.3) circle(.2);
    \draw[thick] (.8,1.3) circle(.2);
  \fill[black] (0,-.4) circle (.12);
}=
\tikzmatht{
  \fill[black!10] (-2.1,-.8) rectangle (1.4,1.7);
    \draw[thick] (.2,1.1) arc(180:0:.4); 
    \draw[thick] (.2,.7)--(.2,.2) arc(360:180:.5);
    \draw[thick] (1,1.2)--(1,.7) arc(360:180:1.4) -- (-1.8,1.5); 
    \draw[thick] (-1.2,1.5)--(-1.2,.6) arc(180:360:.4)--(-.4,1.7);
  \fill[white] (-1.5,1.1) rectangle (-.9,.6); 
    \draw[thick] (-1.5,1.1) rectangle (-.9,.6); 
  \fill[white] (-.1,1.1) rectangle (.5,.6); 
    \draw[thick] (-.1,1.1) rectangle (.5,.6); 
  \fill[white] (-1.8,1.4) circle(.2);
    \draw[thick] (-1.8,1.4) circle(.2);
  \fill[white] (-1.2,1.4) circle(.2);
    \draw[thick] (-1.2,1.4) circle(.2);
  \fill[black] (-.8,.2) circle (.12);
} \stapel{\nref{aux}}= 
\tikzmatht{
  \fill[black!10] (-1.5,-.8) rectangle (1.5,1.7);
    \draw[thick] (0,1.2)--(0,-.6);
    \draw[thick] (-1,1.2)--(-1,.4) arc(180:360:1) -- (1,1.7); 
  \fill[white] (-1.3,.8) rectangle (-.7,.3); 
    \draw[thick] (-1.3,.8) rectangle (-.7,.3); 
  \fill[white] (-.3,.8) rectangle (.3,.3); 
    \draw[thick] (-.3,.8) rectangle (.3,.3); 
  \fill[white] (-1,1.3) circle(.2);
    \draw[thick] (-1,1.3) circle(.2);
  \fill[white] (0,1.3) circle(.2);
    \draw[thick] (0,1.3) circle(.2);
  \fill[black] (0,-.6) circle (.12);
}\stapel{\nref{idem}}= 
\tikzmatht{
  \fill[black!10] (-1,-.8) rectangle (1,1.7);
    \draw[thick] (-.4,1.2)--(-.4,.2) arc(180:360:.4) -- (.4,1.7); 
  \fill[white] (-.1,.8) rectangle (-.7,.3); 
    \draw[thick] (-.1,.8) rectangle (-.7,.3); 
  \fill[white] (-.4,1.3) circle(.2);
    \draw[thick] (-.4,1.3) circle(.2);
},
$$
and similarly
$$
\tikzmatht{
  \fill[black!10] (-1,-.8) rectangle (1,1.7);
    \draw[thick] (0,1.7)--(0,0);
  \fill[white] (-.3,1) rectangle (.3,.5); 
    \draw[thick] (-.3,1) rectangle (.3,.5); 
  \fill[white] (0,0) circle(.2);
    \draw[thick] (0,0) circle(.2);
}=
\tikzmatht{
  \fill[black!10] (-1,-.8) rectangle (1,1.7);
    \draw[thick] (-.4,1.7)--(-.4,.2) arc(180:360:.4) -- (.4,1.2); 
  \fill[white] (.1,.8) rectangle (.7,.3); 
    \draw[thick] (.1,.8) rectangle (.7,.3); 
  \fill[white] (.4,1.3) circle(.2);
    \draw[thick] (.4,1.3) circle(.2);
}.
$$
Now, given $P$, we define $\ol p:= (w^*\ncirc P\times
1_\theta)\scirc x = 
\tikzmatht{
  \fill[black!10] (-1.3,-.8) rectangle (.5,1.7);
    \draw[thick] (0,1.7)--(0,-.8);
    \draw[thick] (-.8,1.2)--(-.8,.4) arc(180:270:.8); 
  \fill[white] (-1.1,.8) rectangle (-.5,.3); 
    \draw[thick] (-1.1,.8) rectangle (-.5,.3); 
  \fill[white] (-.8,1.3) circle(.2);
    \draw[thick] (-.8,1.3) circle(.2);
  \fill[black] (0,-.4) circle (.12);
}$ 
and $p:= (1_\theta\times w^*\ncirc
P)\scirc x = 
\tikzmatht{
  \fill[black!10] (-.5,-.8) rectangle (1.3,1.7);
    \draw[thick] (0,1.7)--(0,-.8);
    \draw[thick] (0,-.4) arc(270:360:.8) -- (.8,1.2); 
  \fill[white] (1.1,.8) rectangle (.5,.3); 
    \draw[thick] (1.1,.8) rectangle (.5,.3); 
  \fill[white] (.8,1.3) circle(.2);
    \draw[thick] (.8,1.3) circle(.2);
  \fill[black] (0,-.4) circle (.12);
}$.   
Then, obviously $P=\ol p\ncirc p$, and $p$ and $\ol p$ satisfy the
last of \eref{1to1b}:
\be\label{self}
\tikzmatht{
  \fill[black!10] (-2,-1.5) rectangle (.5,1.5);
    \draw[thick] (0,1.5)--(0,-.3) arc(360:180:.8)--(-1.6,1.5);
    \draw[thick] (-.8,1.2)--(-.8,.4) arc(360:270:.8); 
  \fill[white] (-1.1,.8) rectangle (-.5,.3); 
    \draw[thick] (-1.1,.8) rectangle (-.5,.3); 
  \fill[white] (-.8,1.2) circle(.2);
    \draw[thick] (-.8,1.2) circle(.2);
  \fill[black] (-1.6,-.4) circle (.12);
  \fill[black] (-.8,-1.1) circle (.12);
    \draw[thick] (-.8,-1.5)--(-.8,-1.1); 
} = 
\tikzmatht{
  \fill[black!10] (-2,-1.5) rectangle (.5,1.5);
    \draw[thick] (0,1.5)--(0,-.3) arc(360:180:.8)--(-1.6,1.5);
    \draw[thick] (-.8,1.2)--(-.8,.4) arc(180:270:.8); 
  \fill[white] (-1.1,.8) rectangle (-.5,.3); 
    \draw[thick] (-1.1,.8) rectangle (-.5,.3); 
  \fill[white] (-.8,1.2) circle(.2);
    \draw[thick] (-.8,1.2) circle(.2);
  \fill[black] (0,-.4) circle (.12);
  \fill[black] (-.8,-1.1) circle (.12);
    \draw[thick] (-.8,-1.5)--(-.8,-1.1); 
},\ee
hence $p$ and $\ol p$  are related to each other by \eref{ppbar}. 
Moreover, $\ol p$ obviously satisfies \eref{pbar} and the second of
\eref{1to1b} by associativity and the unit property. 
In view of \lref{l:ppbar}, it remains to verify that $\ol p$ is a
projection. Idempotency of $\ol p$ is just \eref{idem}, and
selfadjointness of follows from \eref{self}:   
$$
\tikzmatht{
  \fill[black!10] (-1.3,.8) rectangle (.5,-1.7);
    \draw[thick] (0,-1.7)--(0,.8);
    \draw[thick] (-.8,-1.2)--(-.8,-.4) arc(180:90:.8); 
  \fill[white] (-1.1,-.8) rectangle (-.5,-.3); 
    \draw[thick] (-1.1,-.8) rectangle (-.5,-.3); 
  \fill[white] (-.8,-1.3) circle(.2);
    \draw[thick] (-.8,-1.3) circle(.2);
  \fill[black] (0,.4) circle (.12);
} = 
\tikzmatht{
  \fill[black!10] (-1.3,-.8) rectangle (1.3,1.7);
    \draw[thick] (.7,1.7)--(.7,-.8);
    \draw[thick] (-.8,1.2)--(-.8,.4) arc(180:270:.5) arc(270:360:.5) arc(180:90:.5); 
  \fill[white] (-1.1,.8) rectangle (-.5,.3); 
    \draw[thick] (-1.1,.8) rectangle (-.5,.3); 
  \fill[white] (-.8,1.3) circle(.2);
    \draw[thick] (-.8,1.3) circle(.2);
  \fill[black] (.7,.9) circle (.12);
} = 
\tikzmatht{
  \fill[black!10] (-1.3,-.8) rectangle (.5,1.7);
    \draw[thick] (0,1.7)--(0,-.8);
    \draw[thick] (-.8,1.2)--(-.8,.4) arc(180:270:.8); 
  \fill[white] (-1.1,.8) rectangle (-.5,.3); 
    \draw[thick] (-1.1,.8) rectangle (-.5,.3); 
  \fill[white] (-.8,1.3) circle(.2);
    \draw[thick] (-.8,1.3) circle(.2);
  \fill[black] (0,-.4) circle (.12);
}, \qquad
\tikzmatht{
  \fill[black!10] (-1.3,-.8) rectangle (1.3,1.7);
    \draw[thick] (.8,1.7)--(.8,-.8);
    \draw[thick] (-.8,1.2)--(-.8,.4) arc(180:270:.8)--(.8,-.4); 
    \draw[thick] (0,1.2)--(0,.4) arc(180:270:.4)--(.8,0); 
  \fill[white] (-1.1,1) rectangle (-.5,.5); 
    \draw[thick] (-1.1,1) rectangle (-.5,.5); 
  \fill[white] (-.8,1.4) circle(.2);
    \draw[thick] (-.8,1.4) circle(.2);
  \fill[white] (-.3,1) rectangle (.3,.5); 
    \draw[thick] (-.3,1) rectangle (.3,.5); 
  \fill[white] (0,1.4) circle(.2);
    \draw[thick] (0,1.4) circle(.2);
  \fill[black] (.8,0) circle (.12);
  \fill[black] (.8,-.4) circle (.12);
} =
\tikzmatht{
  \fill[black!10] (-1.3,-.8) rectangle (1.3,1.7);
    \draw[thick] (.8,1.7)--(.8,-.8);
    \draw[thick] (-.8,1.2)--(-.8,.4) arc(180:360:.4)--(0,1.2); 
    \draw[thick] (-.4,0) arc(180:270:.4)--(.8,-.4); 
  \fill[white] (-1.1,1) rectangle (-.5,.5); 
    \draw[thick] (-1.1,1) rectangle (-.5,.5); 
  \fill[white] (-.8,1.4) circle(.2);
    \draw[thick] (-.8,1.4) circle(.2);
  \fill[white] (-.3,1) rectangle (.3,.5); 
    \draw[thick] (-.3,1) rectangle (.3,.5); 
  \fill[white] (0,1.4) circle(.2);
    \draw[thick] (0,1.4) circle(.2);
  \fill[black] (-.4,0) circle (.12);
  \fill[black] (.8,-.4) circle (.12);
} =
\tikzmatht{
  \fill[black!10] (-1.3,-.8) rectangle (.5,1.7);
    \draw[thick] (0,1.7)--(0,-.8);
    \draw[thick] (-.8,1.2)--(-.8,.4) arc(180:270:.8); 
  \fill[white] (-1.1,.8) rectangle (-.5,.3); 
    \draw[thick] (-1.1,.8) rectangle (-.5,.3); 
  \fill[white] (-.8,1.3) circle(.2);
    \draw[thick] (-.8,1.3) circle(.2);
  \fill[black] (0,-.4) circle (.12);
}.$$
\qed

The following is the main result of this subsection.

\begin{tintedbox}
\begin{proposition} \vskip-5mm \label{p:subQ}
Let $\BA=(\theta,w,x)$ be a simple Q-system, describing a subfactor
$N\subset M$. Let either $p$ or $\ol p$
or $P$ be a projection in $\Hom(\theta,\theta)$ with properties as
specified in \lref{l:ppbar} resp.\ \pref{p:subQP}, thus defining the
respective other two projections. Then \eref{redQ} with
normalization factor 1 (i.e., $n_P=\sqrt{\dim(\theta_P)}\cdot
1_{\theta_p}$) defines a reduced Q-system $\BA_P$. The sub-Q-system $\BA_P$ 
is associated with a homomorphism $\iota_p\prec \iota$ which is the
range of $p\in\Hom(\iota,\iota)$ such that $\iota$ is a direct sum
$\iota_p\oplus\iota_{1-p}$.  
\end{proposition}
\end{tintedbox}

We will also refer to this reduced Q-system as a sub-Q-system of $\BA$.

\medskip

{\em Proof:} From \lref{l:ppbar} we know that $p,\ol p$ satisfy
the system \eref{1to1a}, \nref{1to1b}.

The first and second relations of \eref{1to1b} state that
$p$ (resp.\ $\ol p$) are self-intertwiners of $(\theta,x)$ as a left (resp.\
right) $\BA$-module. By \pref{p:modulint} for left and right 
modules, we conclude that $p=1_{\ol\iota}\times e$ 
and $\ol p= \ol e \times 1_{\ol\iota}$ with projections 
$e\in\Hom(\iota,\iota)$, $\ol e\in\Hom(\ol\iota,\ol\iota)$. In terms
of the projections $e$ and $\ol e$, \eref{1to1a} and the last relation
of \eref{1to1b} read
\be\label{ewev}
\tikzmatht{
  \fill[black!10] (-1.5,-1.2) rectangle (1.5,1.2);
  \fill[black!18] (-1,1.2)--(-1,.5) arc(180:360:1)--(1,1.2);
    \draw[thick] (-1,1.2)--(-1,.5) arc(180:360:1)--(1,1.2);
    \draw[thick] (.8,.4)--(1.2,.8) (.8,.8)--(1.2,.4) ; 
    \node at (.4,.6) {$e$};
    \node at (0,-.9) {$w$};  
\node at (-.8,.8) {$\ol\iota$};
} = 
\tikzmatht{
  \fill[black!10] (-1.5,-1.2) rectangle (1.5,1.2);
  \fill[black!18] (-1,1.2)--(-1,.5) arc(180:360:1)--(1,1.2);
    \draw[thick] (-1,1.2)--(-1,.5) arc(180:360:1)--(1,1.2);
    \draw[thick] (-.8,.4)--(-1.2,.8) (-.8,.8)--(-1.2,.4) ; 
    \node at (-.4,.6) {$\ol e$};
    \node at (0,-.9) {$w$};  
    \node at (.8,.8) {$\iota$};
}, \qquad 
\tikzmatht{
  \fill[black!18] (-1.5,-1.2) rectangle (1.5,1.2);
  \fill[black!10] (-1,1.2)--(-1,.5) arc(180:360:1)--(1,1.2);
    \draw[thick] (-1,1.2)--(-1,.5) arc(180:360:1)--(1,1.2);
    \draw[thick] (-.8,.4)--(-1.2,.8) (-.8,.8)--(-1.2,.4) ; 
    \node at (-.4,.6) {$e$};
    \node at (0,-.9) {$v$}; 
    \node at (.8,.8) {$\ol\iota$};
} = 
\tikzmatht{
  \fill[black!18] (-1.5,-1.2) rectangle (1.5,1.2);
  \fill[black!10] (-1,1.2)--(-1,.5) arc(180:360:1)--(1,1.2);
    \draw[thick] (-1,1.2)--(-1,.5) arc(180:360:1)--(1,1.2);
    \draw[thick] (.8,.4)--(1.2,.8) (.8,.8)--(1.2,.4) ; 
    \node at (.4,.6) {$\ol e$};
    \node at (0,-.9) {$v$}; 
    \node at (-.8,.8) {$\iota$};
}.
\ee
Because $M$ is a factor, one can write $e=ss^*$ with an isometry $s\in
M$ and define $\iota_s=s^*\iota(\cdot)s\prec\iota$ as the range of $e$. Similarly, 
$\ol\iota_{\ol s}=\ol s^*\ol\iota(\cdot)\ol s\prec\ol\iota$ is the range of $\ol e$. Then $\iota_s$ and 
$\ol\iota_{\ol s}\prec\ol\iota$ are conjugate homomorphisms, because
$w_P=(\ol s^*\times s^*)\scirc w$, $v_P=(s^*\times\ol s^*)\scirc v$
solve the conjugacy relations \eref{conj}: 
$$(1_{\ol\iota_s}\times
v_P^*)\scirc (w_P\times 1_{\iota_s}) =\tikzmatht{          
  \fill[black!10] (-1.5,-1.5) rectangle (1.5,1.5);
  \fill[black!18] (-1.2,1.5)--(-1.2,-.3) arc(180:360:.6)--(0,.3)
arc(180:0:.6)--(1.2,-1.5)--(1.5,-1.5)--(1.5,1.5);
    \draw[thick] (-1.2,1.5)--(-1.2,-.3) arc(180:360:.6)--(0,.3)
arc(180:0:.6)--(1.2,-1.5);
  \fill[white] (-1.2,0)--(-1.4,-.3)--(-1,-.3)--(-1.2,0);
    \draw[thick] (-1.2,0)--(-1.4,-.3)--(-1,-.3)--(-1.2,0);
  \fill[white] (0,0)--(-.2,-.3)--(.2,-.3)--(0,0);
    \draw[thick] (0,0)--(-.2,-.3)--(.2,-.3)--(0,0);
  \fill[white] (0,0)--(-.2,.3)--(.2,.3)--(0,0);
    \draw[thick] (0,0)--(-.2,.3)--(.2,.3)--(0,0);
  \fill[white] (1.2,0)--(1.4,.3)--(1,.3)--(1.2,0);
    \draw[thick] (1.2,0)--(1.4,.3)--(1,.3)--(1.2,0);
    \node at (-.8,1) {$\ol\iota_s$};
    \node at (.8,-1) {$\ol\iota_s$};
    \node at (-.6,-1.2) {$w$};
    \node at (.6,1.2) {$v^*$};
    \node at (.4,0) {$e$};
} \stapel{\nref{ewev}}= 
\tikzmatht{          
  \fill[black!10] (-1.5,-1.5) rectangle (1.5,1.5);
  \fill[black!18] (-1.2,1.5)--(-1.2,-.3) arc(180:360:.6)--(0,.3)
  arc(180:0:.6)--(1.2,-1.5)--(1.5,-1.5)--(1.5,1.5);
    \draw[thick] (-1.2,1.5)--(-1.2,-.3) arc(180:360:.6)--(0,.3)
    arc(180:0:.6)--(1.2,-1.5); 
  \fill[white] (-1.2,1)--(-1.4,.7)--(-1,.7)--(-1.2,1);
    \draw[thick] (-1.2,1)--(-1.4,.7)--(-1,.7)--(-1.2,1);
  \fill[white] (-1.2,0)--(-1.4,-.3)--(-1,-.3)--(-1.2,0);
    \draw[thick] (-1.2,0)--(-1.4,-.3)--(-1,-.3)--(-1.2,0);
  \fill[white] (-1.2,0)--(-1.4,.3)--(-1,.3)--(-1.2,0);
    \draw[thick] (-1.2,0)--(-1.4,.3)--(-1,.3)--(-1.2,0);
  \fill[white] (1.2,0)--(1.4,.3)--(1,.3)--(1.2,0);
    \draw[thick] (1.2,0)--(1.4,.3)--(1,.3)--(1.2,0);
    \node at (.3,0) {$\iota$};
} =
\tikzmatht{          
  \fill[black!10] (-1,-1.5) rectangle (1,1.5);
  \fill[black!18] (0,-1.5) rectangle (1,1.5);
    \draw[thick] (0,-1.5)--(0,1.5);
  \fill[white] (0,-.8)--(-.2,-.5)--(.2,-.5)--(0,-.8);
    \draw[thick] (0,-.8)--(-.2,-.5)--(.2,-.5)--(0,-.8);
  \fill[white] (0,.8)--(-.2,.5)--(.2,.5)--(0,.8);
    \draw[thick] (0,.8)--(-.2,.5)--(.2,.5)--(0,.8);
    \node at (.3,0) {$\ol\iota$};
    \node at (.4,1) {$\ol\iota_s$};
}
=\tikzmatht{          
  \fill[black!10] (-1,-1.5) rectangle (1,1.5);
  \fill[black!18] (0,-1.5) rectangle (1,1.5);
    \draw[thick] (0,-1.5)--(0,1.5);
} = 1_{\ol\iota_s}, 
$$
where we have used the first of \eref{ewev}, and similarly $(1_{\iota_s}\times
w_P^*)\scirc (v_P\times 1_{\ol\iota_s}) =1_{\iota_s}$. 

Now let $S=\ol s^*\ol\iota(s)$, hence $w_P=S^*\circ w$. We compute 
$$\ol\iota_s(v_P)= \ol s^*\ol\iota[s^*\iota(\ol s^*)v]\ol s= \ol
s^*\ol\iota[s^*\iota(\ol s^*\ol\iota(s^*))vs]\ol s = S^*\theta(\ol
s^*\ol\iota(s^*)) \ol\iota(vs)\ol s = (S^*\times
S^*)\scirc x \scirc S=:x_P.$$
It remains to show that $(\theta_P=\ol\iota_s\iota_s,w_P,x_P)$ is the
reduced Q-system. Clearly, $\theta_P=S^*\theta(\cdot)S$. 
With $SS^*=P=\ol p p\in\Hom(\theta,\theta)$, we compute
$$ w_P^*\ncirc w_P = w^*\scirc P\scirc w = \Tr_{\iota}(e) = \Tr_{\ol\iota}(\ol e) 
 = \dim(\iota_s) = \dim(\ol\iota_{\ol s}) = \sqrt{\dim(\theta_P)}
$$
by \pref{p:traces}, and 
$$P\scirc X^*\scirc (P\times P)\scirc X\scirc P = 
\tikzmatht{          
  \fill[black!10] (-1.5,-1.8) rectangle (1.5,1.8);
  \fill[black!18] (-.2,1.8)--(-.2,1.55) arc(360:300:.9)
  arc(120:240:.9) arc(60:0:.9)--(-.2,-1.8)--(.2,-1.8)--(.2,-1.55)
  arc(180:120:.9) arc(-60:60:.9) arc(240:180:.9)--(.2,1.8); 
    \draw[thick] (-.2,1.8)--(-.2,1.55) arc(360:300:.9)
    arc(120:240:.9) arc(60:0:.9)--(-.2,-1.8);
    \draw[thick] (.2,-1.8)--(.2,-1.55) arc(180:120:.9)
    arc(-60:60:.9) arc(240:180:.9)--(.2,1.8);
  \fill[black!10] (0,0) circle (.6);
    \draw[thick] (0,0) circle (.6);
    \draw[thick] (-1.3,-.2)--(-.9,.2) (-1.3,.2)--(-.9,-.2);
    \draw[thick] (1.3,-.2)--(.9,.2) (1.3,.2)--(.9,-.2);
    \draw[thick] (-.8,-.2)--(-.4,.2) (-.8,.2)--(-.4,-.2);
    \draw[thick] (.8,-.2)--(.4,.2) (.8,.2)--(.4,-.2);
    \draw[thick] (-.4,1.2)--(0,1.6) (-.4,1.6)--(0,1.2);
    \draw[thick] (.4,1.2)--(0,1.6) (.4,1.6)--(0,1.2);
    \draw[thick] (-.4,-1.2)--(0,-1.6) (-.4,-1.6)--(0,-1.2);
    \draw[thick] (.4,-1.2)--(0,-1.6) (.4,-1.6)--(0,-1.2);
} = 
\tikzmatht{          
  \fill[black!10] (-1.5,-1.8) rectangle (1.5,1.8);
  \fill[black!18] (-.2,1.8)--(-.2,1.55) arc(360:300:.9) arc(120:240:.9)
  arc(60:0:.9)--(-.2,-1.8)--(.2,-1.8)--(.2,-1.55) 
  arc(180:120:.9) arc(-60:60:.9) arc(240:180:.9)--(.2,1.8); 
    \draw[thick] (-.2,1.8)--(-.2,1.55) arc(360:300:.9)
          arc(120:240:.9) arc(60:0:.9)--(-.2,-1.8);
    \draw[thick] (.2,-1.8)--(.2,-1.55) arc(180:120:.9)
          arc(-60:60:.9) arc(240:180:.9)--(.2,1.8);
  \fill[black!10] (0,0) circle (.6);
    \draw[thick] (0,0) circle (.6);
    \draw[thick] (-.8,-.2)--(-.4,.2) (-.8,.2)--(-.4,-.2);
    \draw[thick] (.8,-.2)--(.4,.2) (.8,.2)--(.4,-.2);
    \draw[thick] (-.4,1.2)--(0,1.6) (-.4,1.6)--(0,1.2);
    \draw[thick] (.4,1.2)--(0,1.6) (.4,1.6)--(0,1.2);
} = 
\tikzmatht{          
  \fill[black!10] (-1.5,-1.8) rectangle (1.5,1.8);
  \fill[black!18] (-.2,1.8)--(-.2,1.55) arc(360:300:.9)
          arc(120:240:.9) arc(60:0:.9)--(-.2,-1.8)--(.2,-1.8)--
          (.2,-1.55) arc(180:120:.9) arc(-60:60:.9) arc(240:180:.9)--
          (.2,1.8); 
    \draw[thick] (-.2,1.8)--(-.2,1.55) arc(360:300:.9)
          arc(120:240:.9) arc(60:0:.9)--(-.2,-1.8);
    \draw[thick] (.2,-1.8)--(.2,-1.55) arc(180:120:.9)
          arc(-60:60:.9) arc(240:180:.9)--(.2,1.8);
  \fill[black!10] (0,0) circle (.6);
    \draw[thick] (0,0) circle (.6);
    \draw[thick] (.8,-.2)--(.4,.2) (.8,.2)--(.4,-.2);
    \draw[thick] (-.4,1.2)--(0,1.6) (-.4,1.6)--(0,1.2);
    \draw[thick] (.4,1.2)--(0,1.6) (.4,1.6)--(0,1.2);
} = \dim(\iota_s)\cdot P,
$$
hence $n_P= x_P^*\ncirc x_P= \dim(\iota_s)\cdot 1_{\theta_P} =
\sqrt{\dim(\theta_P)}\cdot 1_{\theta_P}$. (Contact with \cref{c:norm}
is made by noting that  
$\wt w_P^*\ncirc \wt w_P = w_P^* \ncirc n_P \ncirc w_P = 
\sqrt{\dim(\theta_P)}\cdot w_P^*\ncirc w_P= \dim(\theta_P)$.)
\qed

\begin{corollary} \label{c:subQ}
If $p$ or $\ol p$ satisfy the conditions in \pref{p:subQ}, the same is true
for $1-p$ resp.\ $1-\ol p$. Thus, 
every simple Q-system with $\dim\Hom(\id,\theta)>1$ has a
decomposition into irreducible Q-systems $\BA_P$ with
$\dim\Hom(\id,\theta_P)=1$. 
\end{corollary}

Namely, if $\BA_P$ is reducible, one can just continue the
decomposition. 

(Notice, however, that unlike in \sref{s:Qcentral} the decomposition
corresponds to a partition of unity by $p_i$, not by $P_i=\ol p_i
p_i$! This reflects the obvious fact that 
$[\theta]=(\bigoplus_n[\ol\iota_n])(\bigoplus_m[\iota_m])$ is different from 
$\bigoplus[\theta_P] \simeq \bigoplus[\ol\iota_n\iota_n]$. 
)

\medskip

Finally, instead of characterizing the projection
$p=\ol\iota(e)\in\Hom(\theta,\theta)$ satisfying the pair of
relations as in \pref{p:subQ}, 
one may also write $e=\iota(q)v$ which is in $\Hom(\iota,\iota)$ iff 
$q\in\Hom(\theta,\id)$, and characterize the operator $q$. Indeed, by
\lref{l:relc+cent}, $e$ is idempotent iff    
$q=(q\times q)\scirc x$, and $e$ is selfadjoint iff
$q^*=(1_\theta\times q)\scirc x\scirc w$. In view of these properties,
the first of the two conditions on $p=\theta(q)x$ is equivalent to  
$q^*= (q\times 1_\theta)\scirc x\scirc w$, whereas the second one is
automatic. Therefore, $q\in \Hom(\theta,\id)$ satisfying 
\be \label{q} q\equiv 
\tikzmatht{
  \fill[black!10] (-1,-1) rectangle (1,1.2);
    \draw[thick] (0,-1)--(0,.5);
  \fill[white] (0,.8)--(-.2,.5)--(.2,.5)--(0,.8);
    \draw[thick] (0,.8)--(-.2,.5)--(.2,.5)--(0,.8);
    \node at (.3,-.5) {$\theta$};
} =
\tikzmatht{
  \fill[black!10] (-1.2,-1) rectangle (1.2,1.2);
    \draw[thick] (0,-1)--(0,-.1);
    \draw[thick] (-.6,.5) arc(180:360:.6);
  \fill[black] (0,-.1) circle(.12);
  \fill[white] (-.6,.8)--(-.4,.5)--(-.8,.5)--(-.6,.8);
    \draw[thick] (-.6,.8)--(-.4,.5)--(-.8,.5)--(-.6,.8);
  \fill[white] (.6,.8)--(.8,.5)--(.4,.5)--(.6,.8);
    \draw[thick] (.6,.8)--(.8,.5)--(.4,.5)--(.6,.8);
}=
\tikzmatht{
  \fill[black!10] (-1.2,-1.2) rectangle (1.2,1);
    \draw[thick] (-.6,-.5)--(-.6,-.1) arc(180:0:.6)--(.6,-1.2);
  \fill[white] (-.6,-.8)--(-.8,-.5)--(-.4,-.5)--(-.6,-.8);
    \draw[thick] (-.6,-.8)--(-.8,-.5)--(-.4,-.5)--(-.6,-.8);
} =
\tikzmatht{
  \fill[black!10] (-1.2,-1.2) rectangle (1.2,1);
    \draw[thick] (.6,-.5)--(.6,-.1) arc(0:180:.6)--(-.6,-1.2);
  \fill[white] (.6,-.8)--(.8,-.5)--(.4,-.5)--(.6,-.8);
    \draw[thick] (.6,-.8)--(.8,-.5)--(.4,-.5)--(.6,-.8);
}
\ee
give rise to projections $e=\iota(q)v\in\Hom(\iota,\iota)$, hence
$p=\ol\iota(e)=\theta(q)x\in\Hom(\theta,\theta)$, hence also $\ol
p\in\Hom(\theta,\theta)$ as in the proposition, hence the sub-Q-system. 

Notice that the last equality in \eref{q} is an instance of
\pref{p:dotstar}, which applies since $M$ is a factor ($\BA$ is
simple). 

\subsection{Intermediate Q-systems}
\label{s:interm}
\setcounter{equation}{0}
In this subsection, we shall characterize decompositions of
$\iota:N\to M$ as 
$$\iota = \iota_2\circ\iota_1$$
when $M$ is a factor, i.e., intermediate von Neumann algebras
$\iota_1(N)$ between $N$ and $M$. 

Let $N\subset L\subset M$ be an intermediate extension with
$\iota=\iota_2\circ\iota_1$, hence $\theta=\ol\iota_1\theta_2\iota_1$.
Let $\BA=(\theta,w,x)$ and $\BA_2=(\theta_2,w_2,x_2)$ be the Q-systems for 
$N\subset M$ and $N\subset L$, respectively. The projection 
$e_2=d_2\inv \cdot w_2\ncirc w_2^*\in\Hom(\theta_2,\theta_2)$ onto 
$\id_L\prec\theta_2$ defines a projection
$P=\ol\iota_1(e_2) = 
\tikzmatht{
  \fill[black!10] (-1,-.7) rectangle (1,.7);
  \fill[black!18] (-.7,-.7) rectangle (.7,.7);
    \draw[thick] (-.7,.7)--(-.7,-.7) (.7,.7)--(.7,-.7);
  \fill[black!25] (-.3,.7)--(-.3,.5) arc(180:360:.3)--(.3,.7);
    \draw[thick] (-.3,.7)--(-.3,.5) arc(180:360:.3)--(.3,.7);
  \fill[black!25] (-.3,-.7)--(-.3,-.5) arc(180:0:.3)--(.3,-.7) ;
    \draw[thick] (-.3,-.7)--(-.3,-.5) arc(180:0:.3)--(.3,-.7) ;
}
\in\Hom(\theta,\theta)$. The projection $P$
satisfies the relations \eref{2to3} and 
\be\label{pw}
P\scirc w=w: \qquad 
\tikzmatht{
  \fill[black!10] (-1,-1.3) rectangle (1,1);
    \draw[thick] (0,1)--(0,-.5);
    \draw[thick] (-.2,.1)--(.2,.5); 
    \draw[thick] (.2,.1)--(-.2,.5); 
  \fill[white] (0,-.5) circle(.2);
    \draw[thick] (0,-.5) circle(.2);
    \node at (.5,.4) {$P$}; 
} = d_2\inv\cdot
\tikzmatht{
  \fill[black!10] (-1.2,-1.3) rectangle (1.2,1);
  \fill[black!18] (-.7,1)--(-.7,-.5) arc(180:360:.7)--(.7,1);
    \draw[thick] (-.7,1)--(-.7,-.5) arc(180:360:.7)--(.7,1);
  \fill[black!25] (-.3,1)--(-.3,.7) arc(180:360:.3)--(.3,1);
    \draw[thick] (-.3,1)--(-.3,.7) arc(180:360:.3)--(.3,1);
  \fill[black!25] (.3,-.2) arc(0:180:.3)--(-.3,-.5) arc(180:360:.3);
    \draw[thick] (.3,-.2) arc(0:180:.3)--(-.3,-.5)
          arc(180:360:.3)--(.3,-.2);
} =
\tikzmatht{
  \fill[black!10] (-1,-1.3) rectangle (1,1);
    \draw[thick] (0,1)--(0,-.5);
  \fill[white] (0,-.5) circle(.2);
    \draw[thick] (0,-.5) circle(.2);
},
\ee
hence $w^*\scirc P \scirc w=w^*\scirc w=d\cdot \eins_N$. It also
satisfies 
$$x^*\scirc(P\times P)\scirc x = d_2^{-2} \cdot 
\tikzmatht{
  \fill[black!10] (-2,-2) rectangle (2,2);
  \fill[black!18] (-.6,2) arc(0:-54:.8) arc(126:180:1)--(-1.35,-.55)
          arc(180:234:1) arc(54:0:.8)--(.6,-2) arc(180:126:.8)
          arc(-54:0:1)--(1.35,.55) arc(0:54:1) arc(234:180:.8)--
          (-.6,2); 
    \draw[thick] (-.6,2) arc(0:-54:.8) arc(126:180:1)--
          (-1.35,-.55) arc(180:234:1) arc(54:0:.8);
    \draw[thick] (.6,-2) arc(180:126:.8) arc(-54:0:1)--
          (1.35,.55) arc(0:54:1) arc(234:180:.8);
  \fill[black!25] (-.2,2) arc(0:-53:1) arc(126:180:1) arc(180:360:.2)
          arc(180:0:.6) arc(180:360:.2)  arc(0:54:1) arc(234:180:1);
    \draw[thick] (-.2,2) arc(0:-53:1) arc(126:180:1)
          arc(180:360:.2) arc(180:0:.6) arc(180:360:.2)  arc(0:54:1)
          arc(234:180:1); 
  \fill[black!25] (-.2,-2) arc(0:53:1) arc(234:180:1) arc(180:0:.2)
          arc(180:360:.6) arc(180:0:.2) arc(0:-54:1) arc(126:180:1);
    \draw[thick] (-.2,-2) arc(0:53:1) arc(234:180:1)
          arc(180:0:.2) arc(180:360:.6) arc(180:0:.2) arc(0:-54:1)
          arc(126:180:1); 
  \fill[black!10] (-.3,.2) arc(180:0:.3)--(.3,-.2)
          arc(360:180:.3)--(-.3,.2);
    \draw[thick] (-.3,.2) arc(180:0:.3)--(.3,-.2)
          arc(360:180:.3)--(-.3,.2);
} = d_\BA d_2^{-2}\cdot P.$$

Conversely, the intermediate extension is characterized by the 
projection $P$:

\begin{tintedbox}
\begin{proposition} \vskip-5mm \label{p:intermQ}
Let $\BA=(\theta,w,x)$ be a Q-system in $\C\subset\End_0(N)$, defining
$\iota:N\to M$ of dimension $\dim(\iota)=d_\BA$. Let $P\in\Hom(\theta,\theta)$ 
be a projection satisfying \eref{2to3} and \eref{pw}. Then \eref{redQ} 
defines a reduced Q-system $\BA_P$. The intermediate Q-system corresponds to the intermediate von Neumann algebra $N\subset
L_P\subset M$ given by
\be L_P:= \iota(NP)v . 
\ee
\end{proposition}
\end{tintedbox}

We will also refer to this reduced Q-system as 
{\bf intermediate Q-system} of $\BA$.

\begin{remark}\label{r:bisch}
A similar characterization of intermediate subfactors by projections has 
been given for the type $I\!I$ case in \cite{B}.  
\end{remark}

\begin{remark} \label{r:automatic}
The ``normalization intertwiner'' $n_P\in\Hom_0(\theta_P,\theta_P)$
as in \lref{l:red} will in general not be a multiple of $1_{\theta_P}$, or equivalently, 
$x^*\scirc(P\times P)\scirc x$ will not be a multiple of $P$. Because
of \cref{c:special} and \lref{l:relc+cent}, this can only occur when 
$L_P$ is not a factor. We shall present an example below 
(\xref{x:counter}). On the other hand, when $\dim\Hom(\id,\theta_P)=1$, 
then $n_P\in\Hom_0(\theta_P,\theta_P)$ is trivially a multiple of 
$1_{\theta_P}$. In particular, when $N\subset M$ is irreducible, hence 
$\dim\Hom(\id,\theta)=1$, then $N\subset L_P$ is irreducible, and 
$L_P$ is a factor. We also have: if $n_P\in\Hom_0(\theta_P,\theta_P) 
= \mu \cdot 1_{\theta_P}$, then $\mu = \dim(\theta_P)/d_\BA$, because
by \eref{pr}, $r^*(P\times P)r = r^*(1_{\theta_P}\times P)r = \Tr_{\theta} (P) =
\dim(\theta_P)$, while on the other hand, by \eref{pw}, 
$r^*(P\times P)r = w^*x^*(P\times P)xw = \mu\cdot w^*Pw = 
\mu\cdot w^*w=\mu\cdot d_\BA$.  
\end{remark}

{\em Proof of \pref{p:intermQ}:} We first observe that by the assumed
relations, 
\be\label{pr} 
\tikzmatht{
  \fill[black!10] (-1,-1.2) rectangle (1,1.2);
     \draw[thick] (-.5,1.2)--(-.5,.3) arc(180:360:.5)--(.5,1.2);
     \draw[thick] (.3,.5)--(.7,.9) (.3,.9)--(.7,.5);
} \stapel{\nref{pw}} =
\tikzmatht{
  \fill[black!10] (-1,-1.2) rectangle (1,1.2);
     \draw[thick] (-.5,1.2)--(-.5,.6) arc(180:360:.5)--(.5,1.2);
     \draw[thick] (.3,.5)--(.7,.9) (.3,.9)--(.7,.5);
     \draw[thick] (0,.1)--(0,-.8);
     \draw[thick] (-.2,-.5)--(.2,-.1) (.2,-.5)--(-.2,-.1);
   \fill[black] (0,.1) circle(.12);
   \fill[white] (0,-.8) circle(.2);
     \draw[thick] (0,-.8) circle(.2);
} \stapel{\nref{2to3}} =
\tikzmatht{
  \fill[black!10] (-1,-1.2) rectangle (1,1.2);
     \draw[thick] (-.5,1.2)--(-.5,.6) arc(180:360:.5)--(.5,1.2);
     \draw[thick] (-.3,.5)--(-.7,.9) (-.3,.9)--(-.7,.5);
     \draw[thick] (0,.1)--(0,-.8);
     \draw[thick] (-.2,-.5)--(.2,-.1) (.2,-.5)--(-.2,-.1);
   \fill[black] (0,.1) circle(.12);
   \fill[white] (0,-.8) circle(.2);
     \draw[thick] (0,-.8) circle(.2);
} \stapel{\nref{pw}}= 
\tikzmatht{
  \fill[black!10] (-1,-1.2) rectangle (1,1.2);
     \draw[thick] (-.5,1.2)--(-.5,.3) arc(180:360:.5)--(.5,1.2);
     \draw[thick] (-.3,.5)--(-.7,.9) (-.3,.9)--(-.7,.5);
}.
\ee
Thus, by \pref{p:traces},
$$r^*\scirc (P\times P)\scirc r = r^*\scirc (1_\theta\times P)
\scirc r = \Tr_\theta(P) =\dim(\theta_P).$$ 
Hence, by \cref{c:norm}, $\BA_P=(\theta_P,w_P,x_P)$
is a reduced Q-system.  

We write $\iota(n)\equiv n$ in the following. 

To show that $L_P=NPv$ is a subalgebra of $M$, we compute
$(n_1Pv)(n_2Pv)=n_1P\theta(n_2P)xv = n_1\theta(n_2)P\theta(P)xv 
 = n_1P\theta(n_2)xPv$, using \eref{2to3} in the last step. To show
that $L_P$ is a *-algebra, we compute $(nPv)^*=r^*vPn^*=r^*\theta(Pn^*)v = 
r^*P\theta(n^*)v=r^*\theta(n^*)Pv$, using \eref{pr} in the third step. 
$L_P=N\cdot Pv$ is clearly weakly closed, and it is contained in
$N\cdot v=M$. 

We now compute the Q-system for $N\subset L_P$. Let $P=SS^*$ with $S\in N$,
$S^*S=\eins_N$, and put $\wt w_P:=S^*w$ and $\wt v_P:= S^*v\in
L_P$. Then the embedding $\iota_P:N\to L_P$ is given by  
$$\iota_P(n)\equiv n=nw^*v=nw^*Pv = nw^*SS^*v = n\wt w_P^*\wt v_P.$$
The conjugate map 
$$\ol\iota_P(\cdot):= S^*\ol\iota(\cdot)S$$
is a homomorphism by \eref{2to3}, because every element of $L_P$ is of
the form $nPv=nS\wt v_P$ with $n\in N$. 

We claim that the pair $(\wt w_P,\wt v_P)$ solves the conjugacy
relations \eref{conj} for $(\iota_P,\ol\iota_P)$. Certainly,
$\wt w_P\in\Hom(\id_{N},\ol\iota_P\iota_P)$, because
$$\ol\iota_P\iota_P(n) = S^*\ol\iota(n)S = S^*\theta(n)S=\theta_P(n).
$$
Furthermore, $\wt v_P\in\Hom(\id_{L_P},\iota_P\ol\iota_P)$ because $\wt v_Pn=
S^*vn=S^*\theta(n)v=S^*\theta(n)SS^*v=\theta_P(n)S^*v=\theta_P(n)\wt v_P$,
and $\wt v_P\wt v_P=S^*vS^*v = S^*\theta(S^*)xv = S^*\theta(S^*)xSS^*v =
\wt x_P\wt v_P=\ol\iota_P(\wt v_P)\wt v_P$, using \eref{2to3} in the
third step. The conjugacy relations then follow from \eref{2to3}.  

Finally,  
$\ol\iota(\wt v_P)=S^*\ol\iota(S^*v)S = S^*\theta(S^*)xS=\wt x_P$. 
Thus, after the appropriate rescaling as in \cref{c:norm}, the
Q-system for $N\subset L_P$ coincides with the reduced Q-system $\BA_P=(\theta_P,w_P,x_P)$. 
\qed


\begin{graybox}
\begin{example} \vskip-5mm \label{x:counter} 
We give here a counterexample, showing that $n_P$ is {\em not
necessarily} a multiple of $1_{\theta_P}$. 

Let $N\subset L\subset M$, where $N$ and $M$ are factors, 
and $L=\bigoplus_i L_i$ a finite direct sum of factors. 
Let $\iota:N\to M$ given by $[\iota]=\bigoplus_{i}[\iota_{2i}\iota_{1i}]$ 
where $\iota_{1i}:N\to L_i$ and $\iota_{2i}:L_i\to M$. Similarly, 
$[\ol\iota]=\bigoplus_{i}[\ol\iota_{1i}\ol\iota_{2i}]$ where 
$\ol\iota_{1i}:L_i\to N$ and $\ol\iota_{2i}:M\to L_1$. We choose 
orthonormal isometries $s_i\in\Hom(\iota_{2i}\iota_{1i},\iota)$ and 
$t_i\in\Hom(\ol\iota_{1i}\ol\iota_{2i},\ol\iota)$. The canonical
endomorphism is $[\theta]=[\ol\iota\iota] = \bigoplus_{ij}
[\ol\iota_{1i}\ol\iota_{2i}\iota_{2j}\iota_{1j}]$. 

The intermediate embedding is described by $\iota_1 =
\bigoplus_i\iota_{1i}:N\to L$, as in \sref{s:nonfact}, with canonical
endomorphism $[\theta_1]=\bigoplus_i [\ol\iota_{1i}\iota_{1i}]\prec
[\theta]$. 

For $N\subset M$ we construct a standard Q-system as usual (cf.\
\lref{l:addmult}): with 
standard pairs $(w_{1i},\ol w_{1i})$ for $\iota_{1i}(N)\subset L_i$ and 
$(w_{2i},\ol w_{2i})$ for $\iota_{2i}(L_i)\subset M$, we have the
``composite'' standard pairs as in \lref{l:addmult}(i) 
$$(w_i=\ol\iota_{1i}(w_{2i})w_{1i},\quad \ol w_i=\iota_{2i}(\ol w_{1i})\ol
w_{2i})$$
for $\iota_i(N)=\iota_{2i}\iota_{1i}(N)\subset M$. Then $w\in \Hom(\id_N,\ol\iota\iota)$ and
$\ol w\in \Hom(\id_M,\iota\ol\iota)$ given by
$$w=\sum_i(t_i\times s_i)\scirc w_i =\sum_i 
\tikzmatht{
  \fill[black!3] (-1.5,-1) rectangle (1.5,1.3);
  \fill[black!14] (-.6,1.3)--(-.6,.4) arc(180:360:.6)--(.6,1.3);
    \draw[thick] (-.6,1.3)--(-.6,.4) arc(180:360:.6)--(.6,1.3);
  \fill[black!25] (.2,1.3)--(.2,.4) arc(360:180:.2)--(-.2,1.3);
    \draw[thick] (.2,1.3)--(.2,.4) arc(360:180:.2)--(-.2,1.3);
  \fill[white] (-.4,.4)--(-.1,.7)--(-.7,.7);
    \draw[thick] (-.4,.4)--(-.1,.7)--(-.7,.7)--(-.4,.4);
  \fill[white] (.4,.4)--(.1,.7)--(.7,.7);
    \draw[thick] (.4,.4)--(.1,.7)--(.7,.7)--(.4,.4);
    \node at (-1.1,.7) {$t_i$};     \node at (1.1,.7) {$s_i$};
    \node at (0,-.6) {$w_i$};
},\quad
\ol w = \sum_i(s_i\times t_i)\scirc \ol w_i = \sum_i 
\tikzmatht{
  \fill[black!25] (-1.5,-1) rectangle (1.5,1.3);
  \fill[black!14] (-.6,1.3)--(-.6,.4) arc(180:360:.6)--(.6,1.3);
    \draw[thick] (-.6,1.3)--(-.6,.4) arc(180:360:.6)--(.6,1.3);
  \fill[black!3] (.2,1.3)--(.2,.4) arc(360:180:.2)--(-.2,1.3);
    \draw[thick] (.2,1.3)--(.2,.4) arc(360:180:.2)--(-.2,1.3);
  \fill[white] (-.4,.4)--(-.1,.7)--(-.7,.7);
    \draw[thick] (-.4,.4)--(-.1,.7)--(-.7,.7)--(-.4,.4);
  \fill[white] (.4,.4)--(.1,.7)--(.7,.7);
    \draw[thick] (.4,.4)--(.1,.7)--(.7,.7)--(.4,.4);
    \node at (-1.1,.7) {$s_i$}; 
    \node at (1.1,.7) {$t_i$};
    \node at (0,-.6) {$\ol w_i$};
}$$
form a standard pair for $\iota:N\to M$, hence $(\theta,
w,x)$ is the Q-system for $\iota(N)\subset M$, where
$$x=\ol\iota(\ol w)=\sum_{i,j,k}\tikzmatht{
  \fill[black!3] (-2.2,-1.3) rectangle (2.2,2.3);
  \fill[black!14] (-1.4,2.3)--(-1.4,1.4) arc(180:225:1.4)
  arc(45:0:1.4)--(-.6,-1.3)--(.6,-1.3)--(.6,-.6) arc(180:135:1.4)
  arc(315:360:1.4)--(1.4,2.3); 
    \draw[thick] (-1.4,2.3)--(-1.4,1.4) arc(180:225:1.4)
    arc(45:0:1.4)--(-.6,-1.3); 
    \draw[thick] (.6,-1.3)--(.6,-.6) arc(180:135:1.4)
    arc(315:360:1.4)--(1.4,2.3); 
  \fill[black!25] (-.2,-1.3)--(-.2,-.2) arc(0:53:1)
  arc(233:180:1)--(-1,2.3)--(1,2.3)--(1,1.4) arc(360:307:1)
  arc(127:180:1)--(.2,-1.3); 
    \draw[thick] (-.2,-1.3)--(-.2,-.2) arc(0:53:1)
    arc(233:180:1)--(-1,2.3); 
    \draw[thick] (1,2.3)--(1,1.4) arc(360:307:1)
    arc(127:180:1)--(.2,-1.3); 
  \fill[black!14] (-.6,2.3)--(-.6,1.4) arc(180:360:.6)--(.6,2.3);
    \draw[thick] (-.6,2.3)--(-.6,1.4) arc(180:360:.6)--(.6,2.3);
  \fill[black!3] (.2,2.3)--(.2,1.4) arc(360:180:.2)--(-.2,2.3);
    \draw[thick] (.2,2.3)--(.2,1.4) arc(360:180:.2)--(-.2,2.3);
  \fill[white] (-1.2,1.4)--(-1.5,1.7)--(-.9,1.7);
    \draw[thick] (-1.2,1.4)--(-1.5,1.7)--(-.9,1.7)--(-1.2,1.4);
  \fill[white] (-.4,1.4)--(-.1,1.7)--(-.7,1.7);
    \draw[thick] (-.4,1.4)--(-.1,1.7)--(-.7,1.7)--(-.4,1.4);
  \fill[white] (.4,1.4)--(.1,1.7)--(.7,1.7);
    \draw[thick] (.4,1.4)--(.1,1.7)--(.7,1.7)--(.4,1.4);
  \fill[white] (1.2,1.4)--(1.5,1.7)--(.9,1.7);
    \draw[thick] (1.2,1.4)--(1.5,1.7)--(.9,1.7)--(1.2,1.4);
  \fill[white] (-.4,-.6)--(-.1,-.9)--(-.7,-.9);
    \draw[thick] (-.4,-.6)--(-.1,-.9)--(-.7,-.9)--(-.4,-.6);
  \fill[white] (.4,-.6)--(.1,-.9)--(.7,-.9);
    \draw[thick] (.4,-.6)--(.1,-.9)--(.7,-.9)--(.4,-.6);
    \node at (-.4,2) {\scriptsize $i$}; 
    \node at (.4,2) {\scriptsize $i$};
    \node at (-1.9,1.7) {$t_j$};  
    \node at (1.9,1.6) {$s_k$};
    \node at (-1.1,-.8) {$t_j^*$}; 
    \node at (1.1,-.8) {$s_k^*$};
}.$$
The projection $P\in\Hom(\theta,\theta)$ onto $\theta_1\prec \theta$
is given by 
$$P=\sum_i (t_i\times s_i)\scirc (1_{\ol\iota_{1i}}\times E_{2i}\times
1_{\iota_{1i}})\scirc (t_i\times s_i) = 
\sum_{i}\dim(\iota_{2i})\inv\cdot
\tikzmatht{
  \fill[black!10] (-1.5,-1.3) rectangle (1.5,1.3);
  \fill[black!18] (-.6,1.3)--(-.6,-1.3)--(.6,-1.3)--(.6,1.3);
    \draw[thick] (-.6,1.3)--(-.6,-1.3)(.6,-1.3)--(.6,1.3);
  \fill[black!25] (.2,1.3)--(.2,.4) arc(360:180:.2)--(-.2,1.3);
    \draw[thick] (.2,1.3)--(.2,.4) arc(360:180:.2)--(-.2,1.3);
  \fill[black!25] (.2,-1.3)--(.2,-.4) arc(0:180:.2)--(-.2,-1.3);
    \draw[thick] (.2,-1.3)--(.2,-.4) arc(0:180:.2)--(-.2,-1.3);
  \fill[white] (-.4,.4)--(-.1,.7)--(-.7,.7);
    \draw[thick] (-.4,.4)--(-.1,.7)--(-.7,.7)--(-.4,.4);
  \fill[white] (.4,.4)--(.1,.7)--(.7,.7);
    \draw[thick] (.4,.4)--(.1,.7)--(.7,.7)--(.4,.4);
  \fill[white] (-.4,-.4)--(-.1,-.7)--(-.7,-.7);
    \draw[thick] (-.4,-.4)--(-.1,-.7)--(-.7,-.7)--(-.4,-.4);
  \fill[white] (.4,-.4)--(.1,-.7)--(.7,-.7);
    \draw[thick] (.4,-.4)--(.1,-.7)--(.7,-.7)--(.4,-.4);
    \node at (-1,.8) {$t_i$}; 
    \node at (1,.8) {$s_i$};
    \node at (-1,-.7) {$t_i^*$}; 
    \node at (1,-.7) {$s_i^*$};
}$$
where $E_{2i}=\dim(\iota_{2i})\inv\cdot w_{2i}w_{2i}^* \in
\Hom(\ol\iota_{2i}\iota_{2i},\ol\iota_{2i}\iota_{2i})$ is the
projection onto $\id_{L_i}\prec \ol\iota_{2i}\iota_{2i}$. Then, one
computes 
$$x^*\scirc(P\times P)\scirc x=\sum_i
\frac{\dim(\iota_{1i})}{\dim(\iota_{2i})}\cdot (t_i\times s_i)\scirc
(1_{\ol\iota_{1i}}\times E_{2i}\times 1_{\iota_{1i}})\scirc (t_i\times
s_i).$$
Since $\frac{\dim(\iota_{1i})}{\dim(\iota_{2i})}$ in general depends
on $i$, this is not a multiple of $P$ in general. In contrast, 
the normalization condition in \cite{BDH} (cf.\ \rref{r:BDH}) would be
satisfied. 
\end{example}
\end{graybox}

The following Lemma states how modules restrict to modules of
intermediate Q-systems: 
\begin{lemma} \label{l:modrest}
If $\BA$ is a Q-system, and $\BA_P$ is an intermediate Q-system, then
a (left) $\BA$-module $\mm=(\beta,m)$ restricts to a (left) $\BA_P$-module 
$$\mm_P = (\beta,m_P) \quad\hbox{with}\quad m_P:=\dim(\theta_P)^{\frac14}\cdot
(n_P^{-\frac12}S^*\times 1_\beta)\scirc m,$$
where $S^*S=1$, $SS^*=P$. If $n_P\in\Hom_0(\theta_P,\theta_P)$ is a
multiple of $1_{\theta_P}$, then the normalization factor equals 
$\dim(\theta_P)^{\frac14}\cdot n_P^{-\frac12} = (d_{\BA}/d_{\BA_P})^{\frac12}$. 
If $\mm$ is standard, then $\mm_P$ is standard. The analogous
statements hold for right modules and bimodules. 
\end{lemma}

{\em Proof:}
One easily verifies, using \eref{2to3}, that the defining unit and 
representation properties of a module are satisfied. As for
standardness of $\mm_P$, one has 
$$
\tikzmatht{
  \fill[black!10] (-2,-1.8) rectangle (.8,1.8);
    \draw[thick] (.5,1.8)--(.5,-1.8);
    \draw[thick] (.33,1.1)--(.3,1.1) arc(90:180:.6)--(-.3,-.6)
    arc(180:270:.6)--(.33,-1.1); 
    \draw[thick] (.5,1.27) arc(90:270:.17); 
    \draw[thick] (.5,-.93) arc(90:270:.17); 
  \fill[black] (-.1,-.2) rectangle (-.5,.2);
  \fill[white] (-.3,-.4)--(-.5,-.7)--(-.1,-.7)--(-.3,-.4);
    \draw[thick] (-.3,-.4)--(-.5,-.7)--(-.1,-.7)--(-.3,-.4);
  \fill[white] (-.3,.4)--(-.5,.7)--(-.1,.7)--(-.3,.4);
    \draw[thick] (-.3,.4)--(-.5,.7)--(-.1,.7)--(-.3,.4);
    \node at (-1.2,0) {$n\inv$};
} \stapel{\nref{repn}}= 
\tikzmatht{
  \fill[black!10] (-.9,-2.1) rectangle (1.8,1.5);
    \draw[thick] (1.5,1.5)--(1.5,-2.1);
    \draw[thick] (.9,.6) arc(0:180:.6)--(-.3,-.6)
    arc(180:360:.6)--(.9,.6); 
    \draw[thick] (.3,-1.2) arc(180:270:.5)--(1.33,-1.7); 
    \draw[thick] (1.5,-1.53) arc(90:270:.17); 
  \fill[black] (-.1,-.2) rectangle (-.5,.2);
  \fill[white] (-.3,-.4)--(-.5,-.7)--(-.1,-.7)--(-.3,-.4);
    \draw[thick] (-.3,-.4)--(-.5,-.7)--(-.1,-.7)--(-.3,-.4);
  \fill[white] (-.3,.4)--(-.5,.7)--(-.1,.7)--(-.3,.4);
    \draw[thick] (-.3,.4)--(-.5,.7)--(-.1,.7)--(-.3,.4);
}, $$ with
$$
\tikzmatht{
  \fill[black!10] (-.8,-2.1) rectangle (1.3,1.5);
    \draw[thick] (.9,.6) arc(0:180:.6)--(-.3,-.6)
    arc(180:360:.6)--(.9,.6); 
  \fill[black] (-.1,-.2) rectangle (-.5,.2);
  \fill[white] (-.3,-.4)--(-.5,-.7)--(-.1,-.7)--(-.3,-.4);
    \draw[thick] (-.3,-.4)--(-.5,-.7)--(-.1,-.7)--(-.3,-.4);
  \fill[white] (-.3,.4)--(-.5,.7)--(-.1,.7)--(-.3,.4);
    \draw[thick] (-.3,.4)--(-.5,.7)--(-.1,.7)--(-.3,.4);
    \draw[thick] (.3,-2.1)--(.3,-1.2);
} \stapel{\nref{pr},\nref{2to3}}=
\tikzmatht{
  \fill[black!10] (-.8,-2.1) rectangle (1.3,1.5);
    \draw[thick] (.9,.6) arc(0:180:.6)--(-.3,-.3)
    arc(180:360:.6)--(.9,.6); 
  \fill[black] (-.1,-.1) rectangle (-.5,.3);
  \fill[white] (-.3,-.2)--(-.5,-.5)--(-.1,-.5)--(-.3,-.2);
    \draw[thick] (-.3,-.2)--(-.5,-.5)--(-.1,-.5)--(-.3,-.2);
  \fill[white] (-.3,.4)--(-.5,.7)--(-.1,.7)--(-.3,.4);
    \draw[thick] (-.3,.4)--(-.5,.7)--(-.1,.7)--(-.3,.4);
  \fill[white] (.9,-.2)--(1.1,-.5)--(.7,-.5)--(.9,-.2);
    \draw[thick] (.9,-.2)--(1.1,-.5)--(.7,-.5)--(.9,-.2);
  \fill[white] (.9,.4)--(1.1,.7)--(.7,.7)--(.9,.4);
    \draw[thick] (.9,.4)--(1.1,.7)--(.7,.7)--(.9,.4);
    \draw[thick] (.3,-2.1)--(.3,-.9);
  \fill[black] (.3,-.9) circle(.12);
  \fill[white] (.3,-1.5)--(.5,-1.8)--(.1,-1.8)--(.3,-1.5);
    \draw[thick] (.3,-1.5)--(.5,-1.8)--(.1,-1.8)--(.3,-1.5);
  \fill[white] (.3,-1.4)--(.5,-1.1)--(.1,-1.1)--(.3,-1.4);
    \draw[thick] (.3,-1.4)--(.5,-1.1)--(.1,-1.1)--(.3,-1.4);
} =
\tikzmatht{
  \fill[black!10] (-1.5,-2.1) rectangle (1.3,1.5);
    \draw[thick] (.9,.6) arc(0:180:.6)--(-.3,.1)
    arc(180:360:.6)--(.9,.6); 
   \fill[white] (-.3,.3)--(-.5,0)--(-.1,0)--(-.3,.3);
    \draw[thick] (-.3,.3)--(-.5,0)--(-.1,0)--(-.3,.3);
  \fill[white] (-.3,.4)--(-.5,.7)--(-.1,.7)--(-.3,.4);
    \draw[thick] (-.3,.4)--(-.5,.7)--(-.1,.7)--(-.3,.4);
  \fill[white] (.9,.3)--(1.1,0)--(.7,0)--(.9,.3);
    \draw[thick] (.9,.3)--(1.1,0)--(.7,0)--(.9,.3);
  \fill[white] (.9,.4)--(1.1,.7)--(.7,.7)--(.9,.4);
    \draw[thick] (.9,.4)--(1.1,.7)--(.7,.7)--(.9,.4);
    \draw[thick] (.3,-2.1)--(.3,-.5);
  \fill[black] (.3,-.5) circle(.12);
  \fill[white] (.3,-1)--(.5,-.7)--(.1,-.7)--(.3,-1);
    \draw[thick] (.3,-1)--(.5,-.7)--(.1,-.7)--(.3,-1);
  \fill[black] (.1,-1.5) rectangle (.5,-1.1);
  \fill[white] (.3,-1.6)--(.5,-1.9)--(.1,-1.9)--(.3,-1.6);
    \draw[thick] (.3,-1.6)--(.5,-1.9)--(.1,-1.9)--(.3,-1.6);
    \node at (-.6,-1.2) {$n\inv$};} = 
\tikzmatht{
  \fill[black!10] (-1.5,-2.1) rectangle (1.3,1.5);
  \fill[white] (.3,.7) circle(.2);
    \draw[thick] (.3,.7) circle(.2);
    \draw[thick] (.3,-2.1)--(.3,.5);
  \fill[white] (.3,0)--(.5,.3)--(.1,.3)--(.3,0);
    \draw[thick] (.3,0)--(.5,.3)--(.1,.3)--(.3,0);
  \fill[black] (.1,-.2) rectangle (.5,-.6);
  \fill[black] (.1,-1.4) rectangle (.5,-1);
  \fill[white] (.3,-1.6)--(.5,-1.9)--(.1,-1.9)--(.3,-1.6);
    \draw[thick] (.3,-1.6)--(.5,-1.9)--(.1,-1.9)--(.3,-1.6);
    \node at (-.4,-.4) {$n$};
    \node at (-.6,-1.2) {$n\inv$};
} = 
\tikzmatht{
  \fill[black!10] (-.8,-2.1) rectangle (1.3,1.5);
  \fill[white] (.3,.7) circle(.2);
    \draw[thick] (.3,.7) circle(.2);
    \draw[thick] (.3,-2.1)--(.3,.5);
  \fill[white] (.3,-.2)--(.5,.1)--(.1,.1)--(.3,-.2);
    \draw[thick] (.3,-.2)--(.5,.1)--(.1,.1)--(.3,-.2);
  \fill[white] (.3,-.6)--(.5,-.9)--(.1,-.9)--(.3,-.6);
    \draw[thick] (.3,-.6)--(.5,-.9)--(.1,-.9)--(.3,-.6);
} \stapel{\nref{pw}}= 
\tikzmatht{
  \fill[black!10] (-.6,-2.1) rectangle (1.2,1.5);
  \fill[white] (.3,.7) circle(.2);
    \draw[thick] (.3,.7) circle(.2);
    \draw[thick] (.3,-2.1)--(.3,.5);
},
$$
where in the second step, we have used that
$n\inv\in\Hom_0(\theta_P,\theta_P)$, and the
definition of $n_P$ and \eref{pw} in the third step.
Thus, 
$m_P^*\ncirc m_P = \dim(\theta_P)^{\frac12}\cdot 1_\beta$ by
\eref{m-unit}. Because $\dim(\theta_P)^{\frac12}=d_{\BA_P}$, this is
the proper normalization of a standard module in accord with \lref{l:normmod}.

If $n_P\in\Hom_0(\theta_P,\theta_P)$ is a
multiple of $1_{\theta_P}$, then
$n_P=\dim(\theta_P)/\dim(\theta)^{\frac12}\cdot 1_{\theta_P}$ by \rref{r:automatic}, giving
the stated normalization factor.

The right module and bimodule cases are proven similarly.
\qed

\subsection{Q-systems in braided tensor categories}
\label{s:braid}
\setcounter{equation}{0}
Let now $\C\subset \End_0(N)$ be in addition be {\em braided}. 
The braiding is denoted by
$$\eps_{\rho,\sig}\equiv 
\tikzmatht{
  \fill[black!10] (-1.2,-1) rectangle (1.2,1);
    \draw[thick] (1,1)--(-1,-1) (1,-1)--(.2,-.2) (-.2,.2)--(-1,1);
    \node at (-1,-.6) {$\rho$};  
    \node at (1,-.6) {$\sig$}; 
}
\in\Hom(\rho\sig,\sig\rho).$$ 
We also write  
$\eps^+_{\rho,\sig}\equiv \eps_{\rho,\sig}$, and
$\eps^-_{\rho,\sig}\equiv \eps_{\sig,\rho}^*$ for the opposite braiding. 

\begin{definition} \label{d:Copp}
If $\C$ is a braided C* tensor category with braiding
$\eps\equiv\eps^+$, then $\C\opp$ is the braided C* tensor category, 
which coincides with $\C$ as a C* tensor category, equipped with 
the opposite braiding $\eps^-$. 
\end{definition}

\begin{remark} \label{r:opp}
This definition is tantamount to the
  more fundamental definition (as in \sref{s:canonical}), according to which the monoidal product
  is regarded as a functor $\times:\C\times\C\to\C$, and $\C\opp$ is 
  the category equipped with the opposite monoidal product
  $\sig\times\opp\rho=\rho\times\sigma$. The
  braiding is a natural transformation between the functors $\times$
  and $\times\opp$, and its inverse $:\times\opp\to\times$ is the
  opposite braiding.
The equivalence can be seen ``by left-right reflection of every
diagram''. 
\end{remark}

The cases of interest in QFT are $\C=\DHR(\A)$, the categories of DHR 
endomorphisms of local quantum field nets. These categories are
braided categories, where the DHR braiding is defined in terms of
unitary ``charge transporters'' changing the localization of DHR
endomorphisms, as exposed \sref{s:DHR-braiding}. In low dimensions, the
braiding and the opposite braiding arise, depending on the choice of a
connected component of the spacelike complement. In particular, for
a two-dimensional conformal net $\A_2=\A_+\otimes \A_-$ arising as
a product of its two chiral subnets, we have  
$\DHR(\A_2)=\DHR(\A_+)\boxtimes\DHR(\A_-)\opp$, cf.\ \sref{s:two-d}. 

\begin{definition} \label{d:twist}
If $\rho\in\C$, then the operator 
$$\LTr_\rho(\eps_{\rho,\rho})\equiv 
\tikzmatht{
  \fill[black!10] (-1.5,-1) rectangle (1.5,1);
    \draw[thick] (1,1)--(-.4,-.4) arc(315:45:.58)--(-.2,.2) (.2,-.2)--(1,-1);
    \node at (1,.6) {$\rho$}; 
} = 
\RTr_\rho(\eps_{\rho,\rho})\in\Hom(\rho,\rho)$$ 
is called the {\bf twist}. The twist is a unitary self-intertwiner
\cite[after Lemma 4.3]{LRo},\cite[Prop.\ 2.4]{M00}; in particular, 
it is a complex phase denoted $\kappa_\rho$ if $\rho$ is irreducible. 
\end{definition}

\begin{graybox}
\begin{example} \vskip-5mm \label{x:Ibraid} (Braiding of the Ising category)

The tensor category \xref{x:Icat} can be equipped with four
inequivalent braidings.

The braiding of the DHR category of the Ising model is specified by 
$$\eps_{\tau,\tau}=-1, \quad \eps_{\sig,\sig}
=\kappa_\sig^{-1}\cdot rr^*+\kappa_\sig^3\cdot tt^*, \quad 
\eps_{\sig,\tau}=\eps_{\tau,\sig}=-iu,$$
where $\kappa_\sig = \exp\frac{2\pi i}{16}$. 

\newpage

(This braiding and its opposite, and a second pair of braidings
obtained by replacing $\kappa_\sig$ by $-\kappa_\sig$, exhaust all
possibilities. The second tensor category mentioned in \xref{x:Icat}
also admits four inequivalent braidings.) 
\end{example}
\end{graybox}

\begin{definition} \label{d:commutative}
A Q-system $(\theta,w,x)$ in a braided tensor category is called
{\bf commutative} if 
\be \label{comm}
\eps_{\theta,\theta}\scirc x= x: \qquad 
\tikzmatht{
  \fill[black!10] (-1.4,-1.3) rectangle (1.4,1.5);
    \draw[thick] (.8,1.5)--(-.43,.27) arc(-225:45:.6)--(.2,.5)
    (-.2,.9)--(-.8,1.5);  
    \draw[thick] (0,-1.3)--(0,-.75); 
  \fill[black] (0,-.75) circle (.12);
} = 
\tikzmatht{
  \fill[black!10] (-1.4,-1.3) rectangle (1.4,1.5);
    \draw[thick] (-.8,1.5)--(-.8,.4) arc(180:360:.8)--(.8,1.5);
    \draw[thick] (0,-1.3)--(0,-.4); 
  \fill[black] (0,-.4) circle (.12);
}
\ee
\end{definition}

\begin{tintedbox}
\begin{proposition} \vskip-5mm \label{p:cancomm} \cite{LR95} 
The canonical Q-system (cf.\ \pref{p:canonical}) of a braided C*
tensor category is a commutative Q-system in the Deligne product
$\C\boxtimes \C\opp$.
\end{proposition}
\end{tintedbox}

In local quantum field theory, commutative Q-systems describe local
extensions of a given local quantum field theory \cite{LR95}, cf.\
\sref{s:Qnets}. 

Recall that the DHR category of a two-dimensional QFT which is the
tensor product $\A_2=\A_+\otimes\A_-$ of two chiral QFTs, is 
$\DHR(\A_2)= \DHR(\A_+)\boxtimes \DHR(\A_-)\opp$ as a braided
category. Therefore, if $\A_+$ and $\A_-$ are isomorphic, the 2D
extension associated with the canonical Q-system in $\DHR(\A)\boxtimes
\DHR(\A)\opp$ is always a local QFT.

\subsection{$\alpha$-induction}
\label{s:alpha}
\setcounter{equation}{0}
If $\BA=(\theta,w,x)$ is a Q-system in a braided category,
then $\mm=(\beta=\theta\rho, m=\theta^2(\eps^\pm_{\theta,\rho})x^{(2)})$ is a
standard $\BA$-$\BA$-bimodule. The formula \eref{phi-bim} for the associated endomorphism
$\varphi:M\to M$ becomes 
$$\varphi(\iota(n)v)=\iota(\rho(n)\eps^\pm_{\theta,\rho})v,$$ 
which is known as the $\alpha$-induction of $\rho\in\End_0(N)$ to
$\alpha^\mp_\rho\in\End_0(M)$, originally defined by 
$\ol\iota\circ\alpha^\pm_\rho = \Ad_{\eps_{\rho,\theta}}\circ\rho\circ\ol\iota$ 
\cite{LR95,BE}. 

The endomorphisms $\alpha_\rho^\pm$ {\em extend} the
endomorphism $\rho\in\End(N)$: 
\be \label{extend}
\alpha^\pm_\rho\circ\iota =
\iota\circ\rho,
\ee
and the mappings $\rho\mapsto\alpha^\pm_\rho$, $t\to\iota(t)$ are
functorial, namely if $t\in\Hom(\rho_1,\rho_2)$, then
\be\label{t-ind}
\iota(t)\in\Hom(\alpha^\pm_{\rho_1},\alpha^\pm_{\rho_2}).\ee 
However, $\iota: \Hom(\rho_1,\rho_2)\to
\Hom(\alpha^\pm_{\rho_1},\alpha^\pm_{\rho_2})$ is in general not
surjective. E.g., $\alpha^\pm_\rho$ may possess self-intertwiners
(i.e., $\alpha^\pm_\rho$ is reducible), while $\rho$ is irreducible. 

\begin{corollary}\label{c:alphaind} 
{\rm (i)} One has $\alpha^\pm_{\ol\rho}=\ol\alpha^\pm_\rho$ and
$\dim(\alpha^\pm_\rho) = \dim(\rho)$. \\
{\rm (ii)} If $(\theta,w,x)$ is a Q-system in $\End_0(N)$, then  
$(\alpha^\pm_\theta,\iota(w),\iota(x))$ is a Q-system in $\End_0(M)$.
\end{corollary}

{\em Proof:} Since conjugacy and dimension are defined in terms of
intertwiners and their algebraic relations, (i) follows from
\eref{t-ind}. Similarly, (ii) follows because also Q-systems are
defined in terms of intertwiners and their algebraic relations. 
\qed

If the category $\C$ is modular (cf.\ \sref{s:modular}), then the
matrices 
\be\label{Z}
Z_{\rho,\sig}=\dim\Hom(\alpha^-_\rho,\alpha^+_\sig)
\ee
are ``modular invariants'', i.e., they commute with the unitary
representation of the modular group $SL(2,\ZZ)$ defined by the
braiding \cite{R00,BEK99,BEK00}, and have many other remarkable
properties \cite{BEK99,BEK00,EP03} that can, not least, be exploited for
classifications and actual computations.

\subsection{Mirror Q-systems}
\label{s:mirror}
\setcounter{equation}{0}
Let $N\otimes \wt N\subset M$ be an irreducible finite-index
subfactor, and $\BA=(\Theta,W,X)$ its Q-system. The subfactor is called a
canonical tensor product subfactor (CTPS), if $\Theta$ has the form 
$$[\Theta] = \bigoplus Z_{\rho,\wt\sig}[\rho]\otimes[\wt\sig],$$
where $\rho\in\End_0(N)$ and $\wt\sig\in\End_0(\wt N)$ are
irreducible, and $Z_{\rho,\wt\sig}$ are multiplicities. 

The following proposition was derived in \cite[Thm.\ 3.6]{R00}:

\begin{proposition} \label{p:normal} 
The following are equivalent: \\
{\rm (i)} $[\id]\otimes[\wt\sig]\prec\Theta$ implies $[\wt\sig]=[\id_{\wt N}]$, and  
$[\sig]\otimes[\id]\prec\Theta$ implies $[\sig]=[\id_N]$.  \\
{\rm (ii)} It holds  
$$(N\otimes \eins)'\cap M = (\eins\otimes \wt N),\qquad (\eins\otimes
\wt N)'\cap M = (N\otimes \eins).$$
\noindent
{\rm (iii)} There is a bijection $F$ between the set $\Delta$ of sectors $[\rho]$ and
the set $\wt\Delta$ of sectors $[\wt\sig]$ contributing to $\Theta$ such that 
$$[\Theta] = \bigoplus [\rho]\otimes F[\rho];$$
$\Delta$ and $\wt\Delta$ are closed under fusion (i.e., the product
$[\rho_1][\rho_2]$ decomposes into irreducibles in $\Delta$ resp.\
$\wt\Delta$), and $F$ is an isomorphism of fusion rings. 
\end{proposition}

Under stronger conditions 
(the tensor categories generated by the endomorphisms
$\rho\in[\rho]\in\Delta$ and $\wt\rho\in[\wt\rho]\in\wt\Delta$ are
braided and modular, and $\BA$ is commutative), the
isomorphism of fusion rings is even an isomorphism of braided
tensor categories \cite{DMNO}. 

The canonical Q-systems in $\C\boxtimes \C\opp$ with $[\Theta]=\bigoplus
[\rho]\otimes[\ol\rho]$ and $F[\rho]=[\ol\rho]$, cf.\ \cref{c:opp}, are
examples fulfilling the properties in \pref{p:normal}.  

Xu \cite{X} has strengthened the statement:

\begin{proposition} \label{p:mirror} 
Assume that the equivalent conditions of \pref{p:normal}
are fulfilled. Let $\C\subset \End_0(N)$ be the full tensor
subcategory generated by endomorphisms $\rho\in[\rho]\in\Delta$,
and similarly $\wt\C\subset \End_0(\wt N)$. If $\C$ and $\wt\C$ are
braided categories, hence $\BA$ is a Q-system in the braided category
$\C\boxtimes\wt\C$, the $\alpha$-induction of
$\rho\otimes\wt\sig\in\C\otimes\wt\C$ is well-defined (choosing
$\alpha^+$ for definiteness). Then, one has 
\be\label{homalpha}
\iota(\Hom(\rho_1,\rho_2)\otimes \eins) =
\Hom(\alpha_{\rho_1\otimes\id},\alpha_{\rho_2\otimes\id})
\ee
(and similar for $\wt\rho$), rather than just the inclusion $\subset$
acccording to \eref{t-ind}. Moreover, if
$[\wt\rho]=F[\rho]$, then $\alpha_{\id\otimes \wt\rho}$ and
$\alpha_{\ol\rho\otimes\id}$ are unitarily equivalent. 
If $\BA$ is commutative, the unitary
$u\in\Hom(\alpha_{\rho\otimes\id},\alpha_{\id\otimes\wt\rho})$
can be chosen such that 
\be\label{mirrcomm}
(u\times u)\scirc \iota(\eps_{\rho,\rho}\otimes\eins) =
\iota(\eins\otimes\eps_{\wt\rho,\wt\rho}^*)\scirc(u\times u).
\ee
\end{proposition}

From this, he concludes the existence of the ``mirror extension'' 
defined by a ``mirror Q-system'' in $\wt\C$ associated with a Q-system
in $\C$, as follows. 

Assume that the equivalent conditions of \pref{p:normal}
are fulfilled. If $(\theta,w,x)$ is a Q-system in $\C$, there is
$\wt\theta$ such that $\alpha_{\id\otimes \wt\theta}$ and
$\alpha_{\theta\otimes\id}$ are unitarily equivalent, i.e., 
$[\wt\theta]=F[\theta]$. Let
$u\in\Hom(\alpha_{\theta\otimes\id},\alpha_{\id\otimes\wt\theta})$
unitary. Then, by \eref{homalpha}, 
$$u\scirc\iota(w\otimes \eins)\in u\scirc\Hom(\id,\alpha_{\theta\otimes\id}) =
\Hom(\id,\alpha_{\id\otimes\wt\theta}) =
\iota(\eins\otimes\Hom(\id,\wt\theta)),$$
and similarly 
$$(u\times u)\scirc\iota(x\otimes \eins)\scirc u^*\in
\iota(\eins\otimes\Hom(\wt\theta,\wt\theta^2)).$$ 
This defines $\wt w$ and $\wt x$ such that $u\scirc\iota(w\otimes
\eins)=\iota(\eins\otimes \wt w)$ and
$(u\times u)\scirc\iota(x\otimes \eins)\scirc u^* =
\iota(\eins\otimes \wt x)$. 

\begin{corollary} \cite[Thm.\ 3.8]{X} \label{c:mirror} 
$(\wt\theta,\wt w,\wt x)$ is a Q-system in $\wt\C$. If $\BA=(\Theta,W,X)$ is
commutative, then $(\wt\theta,\wt w,\wt x)$ is commutative iff
$(\theta,w,x)$ is commutative.
\end{corollary}

{\em Proof:} The defining relations for $(\wt\theta,\wt w,\wt x)$ to
be a Q-system are satisfied because by \cref{c:alphaind}
$(\alpha_{\theta\otimes\id},\iota(w\otimes\eins),\iota(x\otimes\eins))$ 
is a Q-system in $\End_0(M)$, and hence $(\alpha_{\id\otimes\wh\theta}, u\scirc\iota(w\otimes
\eins),(u\times u)\scirc\iota(x\otimes\eins)\scirc u^*)$ is an
equivalent Q-system in $\End_0(M)$. If $\BA$ is commutative, then
\eref{mirrcomm} proves the second statement.
\qed

\subsection{Centre of Q-systems}
\label{s:centre}
\setcounter{equation}{0}
Let $\BA=(\theta,w,x)$ be a Q-system of dimension $d_\BA$ in a braided C*
category $\C$, $r=x\scirc w$, and $\mm=(\beta,m)$ an
$\BA$-$\BA$-bimodule. Define $Q_\mm^\pm\in\Hom(\beta,\beta)$ by
$$Q_\mm^\pm:=(r^*\times 1_\beta)\scirc
(1_\theta\times\eps_{\beta,\theta}^\pm)\scirc m = (1_\beta\times 
r^*)\scirc (\eps_{\theta,\beta}^\mp\times 1_\theta)\scirc m:
\quad Q^+_\mm=
\tikzmatht{
  \fill[black!10] (-1.7,-2.2) rectangle (1.7,.3);
    \draw[thick] (-.3,-.2)--(-.5,-.2) arc(90:270:.4)--
           (-.17,-1) (.17,-1)--
           (.5,-1) arc(-90:90:.4)--(.3,-.2);
    \draw[thick] (0,-2.2)--(0,.3);
    \draw[thick] (0,-1) circle (.17);
    \node at (-.1,-1.9) [right] {$\beta$}; 
    \node at (-1.1,-1.4) [right] {$m$}; 
    \node at (-1.6,-.3) [right] {$\theta$};
}.$$

\begin{lemma} \label{l:bimproj} (cf.\ \cite{FFRS06})
$P_\mm^\pm:=d_\BA\inv\cdot Q_\mm^\pm$ are projections. For $\mm=\BA$ the trivial
$\BA$-$\BA$-bimodule, the projections $P^\pm\equiv P_\BA^\pm$ satisfy the relations
\be\label{centp}
\tikzmatht{
  \fill[black!10] (-1.7,-1.5) rectangle (1.7,2);
    \draw[thick] (-.8,2)--(-.8,.5) arc(180:360:.8)--(.8,2);
    \draw[thick] (0,-1.5)--(0,-.3); 
    \draw[thick] (-.6,1)--(-1,.6); 
    \draw[thick] (-.6,.6)--(-1,1); 
    \node at (0,1) {$P^+$}; 
  \fill[black] (0,-.3) circle (.12);
} = 
\tikzmatht{
  \fill[black!10] (-1.7,-2) rectangle (1.7,1.5);
    \draw[thick] (0,-1.4) arc(-90:45:.8)--(-.8,1.5); 
    \draw[thick] (0,-1.4) arc(270:135:.8)--(-.2,.4) (.2,.8)--(.8,1.5); 
    \draw[thick] (0,-2)--(0,-1.4); 
    \draw[thick] (1,-.4)--(.6,-.8); 
    \draw[thick] (1,-.8)--(.6,-.4); 
  \fill[black] (0,-1.4) circle (.12);
    \node at (0,-.5) {$P^+$}; 
},\qquad 
\tikzmatht{
  \fill[black!10] (-1.7,-1.5) rectangle (1.7,2);
    \draw[thick] (-.8,2)--(-.8,.5) arc(180:360:.8)--(.8,2); 
    \node at (.2,1) {$P^-$}; 
    \draw[thick] (0,-1.5)--(0,-.3); 
    \draw[thick] (-.6,1)--(-1,.6); 
    \draw[thick] (-.6,.6)--(-1,1); 
  \fill[black] (0,-.3) circle (.12);
} = 
\tikzmatht{
  \fill[black!10] (-1.7,-2) rectangle (1.7,1.5);
    \draw[thick] (0,-1.4) arc(-90:45:.8)--(.2,.4) (-.2,.8)--(-.8,1.5); 
    \draw[thick] (0,-1.4) arc(270:135:.8)--(.8,1.5); 
    \draw[thick] (0,-2)--(0,-1.4); 
    \draw[thick] (1,-.4)--(.6,-.8); 
    \draw[thick] (1,-.8)--(.6,-.4); 
  \fill[black] (0,-1.4) circle (.12);
    \node at (0,-.4) {$P^-$}; 
}.
\ee
\end{lemma}

{\em Proof:} We prove idempotency and selfadjointness of $P_\mm^+$,
using the representation property of the bimodule, the associativity of 
the Q-system, and the unitarity of the twist (cf.\ \dref{d:twist})
in the last step: 
$$
\tikzmatht{
  \fill[black!10] (-1.5,-1.7) rectangle (1.5,1.8);
    \draw[thick] (-.3,0)--(-.5,0) arc(90:270:.4)--(-.17,-.8)
    (.17,-.8)--(.5,-.8) arc(-90:90:.4)--(.3,0);
    \draw[thick] (-.3,1.3)--(-.5,1.3) arc(90:270:.4)--(-.17,.5)
    (.17,.5)--(.5,.5) arc(-90:90:.4)--(.3,1.3);
    \draw[thick] (0,-1.7)--(0,1.8);
    \draw[thick] (0,-.8) circle (.17);
    \draw[thick] (0,.5) circle (.17);
    \node at (-.1,-1.4) [right] {$\beta$};
} =
\tikzmatht{
  \fill[black!10] (-1.5,-1.7) rectangle (1.5,1.8);
    \draw[thick] (-.2,0.5)--(-.5,0.5) arc(90:270:.3)--(-.17,-.1)
    (.17,-.1)--(.5,-.1) arc(-90:90:.3)--(.2,0.5);
    \draw[thick] (-.2,.9)--(-.5,.9) arc(90:270:.7)--(-.17,-.5)
    (.17,-.5)--(.5,-.5) arc(-90:90:.7)--(.2,.9);
    \draw[thick] (0,-1.7)--(0,1.8);
    \draw[thick] (0,-.1) circle (.17);
    \draw[thick] (0,-.5) circle (.17);
} =
\tikzmatht{
  \fill[black!10] (-1.5,-1.7) rectangle (1.5,1.8);
    \draw[thick] (-.2,.3)--(-.7,.3) arc(90:270:.25) arc(90:-90:.25)
    arc(270:90:.75)--(-.2,.8); 
    \draw[thick] (.2,.3)--(.7,.3) arc(90:-90:.25) arc(90:270:.25)
    arc(-90:90:.75)--(.2,.8); 
    \draw[thick] (0,-1.7)--(0,1.8);
    \draw[thick] (-.4,-.45)--(-.17,-.45) (.17,-.45)--(.4,-.45);
    \draw[thick] (0,-.45) circle (.17);
  \fill[black] (-.45,-.45) circle (.12); 
  \fill[black] (.45,-.45) circle (.12);
} =
\tikzmatht{
  \fill[black!10] (-1.5,-1.7) rectangle (1.5,1.8);
    \draw[thick] (-.2,.3)--(-.5,.3) arc(270:90:.25)--(-.2,.8);
    \draw[thick] (.2,.3)--(.5,.3) arc(-90:90:.25)--(.2,.8);
    \draw[thick] (-.75,.55) arc(90:270:.5)--(-.17,-.45); 
    \draw[thick] (.75,.55) arc(90:-90:.5)--(.17,-.45);
    \draw[thick] (0,-1.7)--(0,1.8);
    \draw[thick] (0,-.45) circle (.17);
  \fill[black] (-.75,.55) circle (.12);   
  \fill[black] (.75,.55) circle (.12);
}= 
d_\BA\cdot 
\tikzmatht{
  \fill[black!10] (-1.5,-1.7) rectangle (1.5,1.8);
    \draw[thick] (-.3,.6)--(-.5,.6) arc(90:270:.4)--(-.17,-.2)
    (.17,-.2)--(.5,-.2) arc(-90:90:.4)--(.3,.6);
    \draw[thick] (0,-1.7)--(0,1.8);
    \draw[thick] (0,-.2) circle (.17);
},$$
and
$$
\tikzmatht{
  \fill[black!10] (-1.5,-1.2) rectangle (1.5,1.8);
    \draw[thick] (-.3,0)--(-.5,0) arc(270:90:.4)--(-.17,.8) (.17,.8)--(.5,.8)
arc(90:-90:.4)--(.3,0);
    \draw[thick] (0,-1.2)--(0,1.8);
    \draw[thick] (0,.8) circle (.17);
}=
\tikzmatht{
  \fill[black!10] (-1.5,-1.2) rectangle (1.5,1.8);
    \draw[thick] (.2,1.2)--(.4,1.2) arc(90:0:.4) (.8,.3)--(.8,.2)
arc(-180:90:.3)--(.17,.5);
    \draw[thick] (-.2,1.2)--(-.4,1.2) arc(90:180:.4) (-.8,.3)--(-.8,.2)
arc(360:90:.3)--(-.17,.5);
    \draw[thick] (0,-1.2)--(0,1.8);
    \draw[thick] (0,.5) circle (.17);
}=
\tikzmatht{
  \fill[black!10] (-1.5,-1.2) rectangle (1.5,1.8);
    \draw[thick] (-.3,.8)--(-.5,.8) arc(90:270:.4)--(-.17,0) (.17,0)--(.5,0)
arc(-90:90:.4)--(.3,.8);
    \draw[thick] (0,-1.2)--(0,1.8);
    \draw[thick] (0,0) circle (.17);
}.
$$
We then prove the relation for $P^+\equiv P_\BA^+$:
$$\tikzmatht{
  \fill[black!10] (-1.5,-2) rectangle (2,1.6);
    \draw[thick] (-.8,1.6)--(.7,0) arc(40:-230:.8)--(-.1,.5)  
          (.2,.9)--(.8,1.6); 
    \draw[thick] (.7,.6) arc(100:-200:.7);
    \draw[thick] (0,-2)--(0,-1.35); 
  \fill[black] (0,-1.35) circle (.12);
  \fill[black] (.8,-.8) circle (.12);
}=
\tikzmatht{
  \fill[black!10] (-1.7,-1.6) rectangle (1.8,2);
    \draw[thick] (.2,-1.6)--(.2,.3) arc(0:40:1)--(-1,2);
    \draw[thick] (0,2)--(-.3,1.65) (-.7,1.25)--(-.8,1.1) 
          arc(135:315:.4) (.3,1.1) arc(135:0:.5)--(1.2,0) arc
          (360:180:.95)--(-.7,.5);
  \fill[black] (.2,-.95) circle (.12);
  \fill[black] (-.7,.45) circle (.12);
}=
\tikzmatht{
  \fill[black!10] (-1.7,-1.6) rectangle (1.8,2);
    \draw[thick] (.2,-1.6)--(.2,-.6) arc(0:25:1.2)--(-1,2);
    \draw[thick] (-.1,1.1) arc(110:-210:1) (.4,1.2)--(.4,2);
  \fill[black] (.2,-.85) circle (.12);
  \fill[black] (.4,1.2) circle (.12); 
}=
\tikzmatht{
  \fill[black!10] (-1.7,-1.8) rectangle (1.8,1.8);
    \draw[thick] (0,-1.8)--(0,0) arc(0:30:1.2)--(-1,1.8); 
    \draw[thick] (-.2,1.1) arc(110:-210:.8) (.2,.8);
    \draw[thick] (.4,-.45) arc(180:360:.5)--(1.4,1.8);
  \fill[black] (.4,-.4) circle (.12);
  \fill[black] (0,-.45) circle (.12); 
}=
\tikzmatht{
  \fill[black!10] (-1.7,-1.6) rectangle (1.8,2);
    \draw[thick] (-.2,0) arc(180:360:.7)--(1.2,2); 
    \draw[thick] (-.2,0)--(-.2,.3) arc(0:30:1.2)--(-.8,2);          
    \draw[thick] (-.3,1.4) arc(90:-220:.6) (.45,-1.6)--(.45,-.9);
  \fill[black] (.45,-.75) circle (.12);
  \fill[black] (-.2,.2) circle (.12); 
},$$
where we have several times used associativity of the Q-system. The
proofs for $P^-$ are similar.
\qed

\begin{lemma} \label{l:centreproj}
(cf.\ \cite{FFRS06}) The projections $P_\BA^\pm$ satisfy \eref{2to3}
and \eref{pw}. Hence, they define intermediate extensions 
by \pref{p:intermQ}; the corresponding reduced Q-systems 
$(\theta^\pm_P,w^\pm_P,x^\pm_P)$ are called {\bf left resp.\ right centre} 
$C^\pm[\BA]$.  
\end{lemma}

{\em Proof:} We prove \eref{pw} by
$$
\tikzmatht{
  \fill[black!10] (-1.2,-1.2) rectangle (1.2,1.2);
    \draw[thick] (0,1.2)--(0,-1);
    \draw[thick] (0,-.2) arc(270:110:.5) (0,-.2) arc(270:430:.5);
  \fill[white] (0,-.8) circle(.2);
    \draw[thick] (0,-.8) circle(.2);
  \fill[black] (0,-.2) circle(.12);
} 
= 
\tikzmatht{
  \fill[black!10] (-1.2,-1.2) rectangle (1.2,1.2);
    \draw[thick] (0,1.2)--(0,-.8);
    \draw[thick] (0,.8) arc(90:-70:.5) (0,.8) arc(90:250:.5);
  \fill[white] (0,-.8) circle(.2);
    \draw[thick] (0,-.8) circle(.2);
  \fill[black] (0,.8) circle(.12);} 
= 
\tikzmatht{
  \fill[black!10] (-1.2,-1.2) rectangle (1.2,1.2);
    \draw[thick] (0,1.2)--(0,.8);
    \draw[thick] (0,.8) arc(90:-90:.5) (0,.8) arc(90:270:.5);
  \fill[black] (0,.8) circle(.12);} 
\equiv 
\tikzmatht{  \fill[black!10] (-1.2,-1.2) rectangle (1.2,1.2);
    \draw[thick] (0,1.2)--(0,.8) (0,-.2)--(0,-.8);
    \draw[thick] (0,.8) arc(90:-90:.5) (0,.8) arc(90:270:.5);
  \fill[white] (0,-.8) circle(.2);
    \draw[thick] (0,-.8) circle(.2);
  \fill[black] (0,.8) circle(.12);
  \fill[black] (0,-.2) circle(.12);} 
= d_\BA\cdot
\tikzmatht{
  \fill[black!10] (-1,-1.2) rectangle (1,1.2);
    \draw[thick] (0,1.2)--(0,-.5);
  \fill[white] (0,-.5) circle(.2);
    \draw[thick] (0,-.5) circle(.2);},
$$
using selfadjointness of $P^\pm$, and the unit property and
standardness of $\BA$. In order to establish 
\eref{2to3} (for $P_\BA^+$), we compute  
$$ 
\tikzmatht{
  \fill[black!10] (-2,-2) rectangle (2,2);
    \draw[thick] (-1,2)--(-1,.7) arc(180:360:1)--(1,2);
    \draw[thick] (-1,.7) arc(270:120:.5) (-1,.7) arc(270:420:.5);
    \draw[thick] (1,.7) arc(270:120:.5) (1,.7) arc(270:420:.5);
    \draw[thick] (0,-1.7) arc(270:120:.5) (0,-1.7) arc(270:420:.5);
    \draw[thick] (0,-2)--(0,-.3);
  \fill[black] (-1,.7) circle(.12); 
  \fill[black] (1,.7) circle(.12);
  \fill[black] (0,-.3) circle(.12); 
  \fill[black] (0,-1.7) circle(.12);
}=
\tikzmatht{
  \fill[black!10] (-2,-2) rectangle (2,2);
    \draw[thick] (-1,2)--(-1,.7) arc(180:360:1)--(1,2);
    \draw[thick] (-1,.7) arc(270:120:.5) (-1,.7) arc(270:420:.5);
    \draw[thick] (1,.7) arc(270:120:.5) (1,.7) arc(270:420:.5);
    \draw[thick] (0,-1.7) arc(270:150:1.1) (0,-1.7) arc(270:390:1.1);
    \draw[thick] (0,.5) arc(90:120:1.1) (0,.5) arc(90:60:1.1);
    \draw[thick] (0,-2)--(0,-.3);
  \fill[black] (-1,.7) circle(.12); 
  \fill[black] (1,.7) circle(.12);
  \fill[black] (0,-.3) circle(.12); 
  \fill[black] (0,-1.7) circle(.12);
}=
\tikzmatht{
  \fill[black!10] (-2,-2) rectangle (2,2);
    \draw[thick] (-1,.7) arc(270:120:.5) (-1,.7) arc(270:420:.5);
    \draw[thick] (1,.7) arc(270:120:.5) (1,.7) arc(270:420:.5);
    \draw[thick] (-1,2)--(-1,-.2) arc(360:300:1); 
    \draw[thick] (1,2)--(1,-.2) arc(180:240:1); 
    \draw[thick] (-1.2,.2) arc(90:270:.7)--(1.2,-1.2) arc(-90:90:.7);
    \draw[thick] (-.8,.2)--(.8,.2); 
    \draw[thick] (0,-2)--(0,-1.2); 
  \fill[black] (-1,.7) circle(.12); 
  \fill[black] (1,.7) circle(.12);
  \fill[black] (-1.5,-1.1) circle(.12); 
  \fill[black] (1.5,-1.1) circle(.12);
} \stapel{\nref{centp}}=
\tikzmatht{
  \fill[black!10] (-2,-2) rectangle (2,2);
    \draw[thick] (-1.3,.7) arc(270:120:.5) (-1.3,.7) arc(270:420:.5);
    \draw[thick] (1.3,.7) arc(270:120:.5) (1.3,.7) arc(270:420:.5);
    \draw[thick] (-1.3,2)--(-1.3,0) arc(180:360:.4) arc(180:0:.5)
arc(180:360:.4)--(1.3,2);
    \draw[thick] (-.9,-.4) arc(180:360:.9);
    \draw[thick] (0,-2)--(0,-1.3);
  \fill[black] (-1.3,.7) circle(.12); 
  \fill[black] (1.3,.7) circle(.12);
  \fill[black] (-.9,-.4) circle(.12); 
  \fill[black] (.9,-.4) circle(.12);
}= d_\BA\cdot
\tikzmatht{
  \fill[black!10] (-2,-2) rectangle (2,2);
    \draw[thick] (-1,2)--(-1,.7) arc(180:360:1)--(1,2);
    \draw[thick] (-1,.7) arc(270:120:.5) (-1,.7) arc(270:420:.5);
    \draw[thick] (1,.7) arc(270:120:.5) (1,.7) arc(270:420:.5);
    \draw[thick] (0,-2)--(0,-.3);
  \fill[black] (-1,.7) circle(.12); 
  \fill[black] (1,.7) circle(.12);
  \fill[black] (0,-.3) circle(.12); 
},$$
using associativity in the second step, \eref{centp} in the third step, and
the Frobenius property and standardness in the last step. Thus, one of
the three projections is redundant. Redundancy of the other two is
obtained similarly. The other statements follow from \pref{p:intermQ}.
\qed

The left and right centre projections can be characterized as the
maximal ones satisfying \eref{centp}: 

\begin{tintedbox}
\begin{proposition} \vskip-5mm \label{p:maximal} \cite{FFRS06} 
Among all projections $p\in\Hom(\theta,\theta)$ satisfying \eref{centp},
$P_\BA^\pm$ are the maximal ones. 
\end{proposition}
\end{tintedbox}

{\em Proof:} For $P_\BA^+$: 
$$\tikzmatht{
  \fill[black!10] (-1.6,-1.6) rectangle (1.6,1.6);
    \draw[thick] (0,-1) arc(270:120:.8) (0,-1) arc(270:420:.8);
    \draw[thick] (-.2,1)--(.2,1.4);
    \draw[thick] (-.2,1.4)--(.2,1);
    \draw[thick] (0,-1.6)--(0,1.6);
  \fill[black] (0,-1) circle(.12);
} =
\tikzmatht{
  \fill[black!10] (-1.6,-1.6) rectangle (1.6,1.6);
    \draw[thick] (0,1.6)--(0,.2) arc(360:90:.6)--(-.2,.8); 
    \draw[thick] (-.6,-.4) arc(180:360:.8)--(1,0) arc(0:90:.8)--(.2,.8);
    \draw[thick] (.2,-1.6)--(.2,-1.2);
    \draw[thick] (-.2,1)--(.2,1.4);
    \draw[thick] (-.2,1.4)--(.2,1);
  \fill[black] (-.6,-.4) circle(.12);
  \fill[black] (.2,-1.2) circle(.12);
} \stapel{\nref{centp}}=
\tikzmatht{
  \fill[black!10] (-1.6,-1.6) rectangle (1.6,1.6);
    \draw[thick] (-1.2,1.6)--(-1.2,.6) arc(180:360:.6) arc(180:0:.6)--
(1.2,0) arc(360:180:.9); 
    \draw[thick] (.3,-1.6)--(.3,-.9);
    \draw[thick] (-1,1)--(-1.4,1.4);
    \draw[thick] (-1,1.4)--(-1.4,1);
  \fill[black] (-.6,0) circle(.12);
  \fill[black] (.3,-.9) circle(.12);
} = d_\BA\cdot 
\tikzmatht{
  \fill[black!10] (-1,-1.6) rectangle (1,1.6);
    \draw[thick] (0,-1.6)--(0,1.6);
    \draw[thick] (-.2,.2)--(.2,.6);
    \draw[thick] (-.2,.6)--(.2,.2);
}
$$
Thus, $p<P_\BA^+$, concluding the proof. 
\qed

\begin{tintedbox}
\begin{corollary} \vskip-5mm \label{c:centres}
The left and right centres of a Q-system are maximal commutative
intermediate Q-systems. A Q-system $\BA$ is commutative iff  
$P_\BA^+=1_\theta$ iff $P_\BA^+=1_\theta$ (i.e., $C^\pm[\BA]=\BA$).
\end{corollary}
\end{tintedbox}

{\em Proof:} Follows from \pref{p:maximal} and \pref{p:intermQ}
because by definition, a Q-system is commutative iff $1_\theta$ satisfies
\eref{centp}. 
\qed 

This result is of interest in the applications to local QFT, where the
intermediate extension associated with the centre projections can be
identified as certain relative commutants of local algebras
\cite{BKLR}, cf.\ \sref{s:wedges}. 

\subsection{Braided product of Q-systems}
\label{s:Qprod}
\setcounter{equation}{0}
\begin{definition} \label{d:Qprod}
Let $\BA=(\theta^\BA,w^\BA,x^\BA)$ and 
$\BB=(\theta^\BB,w^\BB,x^\BB)$ be two Q-systems in a braided C* tensor category
$\C$. Then there are two natural product Q-systems, called {\bf braided
products} and denoted as $\BA\times^\pm\BB$, given by the object
$\theta=\theta^\BA\theta^\BB$ and the interwiners
$$w=w^\BA\times w^\BB \equiv
\tikzmatht{
  \fill[black!10] (-1.4,-1.5) rectangle (1.5,1.2);
    \draw[thick] (-.8,1.2)--(-.8,-.2) (.6,1.2)--(.6,-.2);
  \fill[white] (-.8,-.2) circle(.2) (.6,-.2) circle(.2);
    \draw[thick] (-.8,-.2) circle(.2) (.6,-.2) circle(.2);
    \node at (-.7,-.7) {$w^A$}; 
    \node at (.7,-.7) {$w^B$};
    \node at (-.3,.8) {$\theta^A$}; 
    \node at (1.1,.8) {$\theta^B$};
},\quad 
x^\pm= (1_{\theta^\BA}\times\eps^\pm_{\theta^\BA,\theta^\BB}\times
1_{\theta^\BB})\scirc (x^\BA\times x^\BB):\;\; x^+=
\tikzmatht{
  \fill[black!10] (-1.5,-1) rectangle (1.5,1.7);
    \draw[thick] (-1.2,1.7)--(-1.2,1) arc(180:360:.8)--(.4,1.7);
    \draw[thick] (-.4,1.7)--(-.4,1) arc(180:220:.8) (1.2,1.7)--(1.2,1)
    arc(360:260:.8); 
    \draw[thick] (-.4,-1)--(-.4,.2) (.4,-1)--(.4,.2); 
  \fill[black] (-.4,.2) circle (.12) (.4,.2) circle (.12);
    \node at (-1,-.3) {$x^A$}; 
    \node at (1,-.3) {$x^B$}; 
}
.$$
The extension $N\subset M^\pm$ corresponding to the braided product of two
Q-systems is called the {\bf braided product of extensions}.
\end{definition}

Notice that
$\dim\Hom(\id_N,\theta^\BA\theta^\BB)=\dim\Hom(\theta^\BA,\theta^\BB)$ 
can in general be larger than 1, even if
$\dim\Hom(\id_N,\theta^\BA)=\dim\Hom(\id_N,\theta^\BB)=1$. Thus, the
braided product of extensions is in general not irreducible,
and not even a factor, even if both extensions are irreducible. We
shall return to this issue below. 

\medskip

One can easily see that the braided product $\BA\times^\pm\BB$ contains
both $\BA$ and $\BB$ as intermediate Q-systems, via the natural
projections $d_\BA\inv\cdot (w^\BA w^{\BA*}\times 1_{\theta^\BB})$ onto
$\theta^\BA\prec\theta^\BA\theta^\BB$ and $d_\BB\inv\cdot
(1_{\theta^\BA}\times w^\BB w^{\BB*})$ onto $\theta^\BB\prec\theta^\BA\theta^\BB$,
respectively. 

Expressed in terms of the corresponding extensions, the braided
products $N\subset M^\pm$ of extensions $N\subset M^\BA$, $N\subset
M^\BB$ contain both $M^\BA$ and $M^\BB$ as intermediate extensions: 
\be\label{interm}
\begin{array}{ccccc}
&&M^\BA&&\\[-3mm]
&\rotatebox[origin=c]{20}{$\subset$}&&\rotatebox[origin=c]{-20}{$\subset$}&\\[-2.5mm]
N&&&&M^\pm\\[-2.5mm]
&\rotatebox[origin=c]{-20}{$\subset$}&&\rotatebox[origin=c]{20}{$\subset$}&\\[-3mm]
&&M^\BB&&\end{array}.
\ee
More precisely, we have

\begin{lemma} \label{l:braidext}
The braided products $N\subset M^\pm$ of two extensions $N\subset
M^\BA=\iota^\BA(N)v^\BA$, 
$N\subset M^\BB=\iota^\BB(N)v^\BB$ are generated by the subalgebra $N$
and the generator $v^\pm=v^\BA v^\BB$, where $v^\BA$ and $v^\BB$ are
embedded into $M^\pm$ as
$$v^\BA=\iota^\pm(\theta^\BA(w^{\BB*}))v^\pm = 
\tikzmatht{
  \fill[black!10] (-1.8,-2) rectangle (1.2,1);
  \fill[black!18] (-1.8,-2)--(-.5,-2)--(-.5,-1.5)--(-1.3,-.7)--
(-1.3,1)--(-1.8,1)--(-1.8,-2) ;
  \fill[white] (-.5,-1.5)--(-1.3,-.7)--(.3,-.7)--(-.5,-1.5);
    \draw[thick] (-.5,-1.5)--(-1.3,-.7)--(.3,-.7)--(-.5,-1.5);
    \draw[thick] (.3,-.7)--(.3,.3); 
    \draw[thick] (-1.3,-.7)--(-1.3,1); 
    \draw[thick] (-.5,-.7)--(-.5,1); 
    \draw[thick] (-.5,-1.5)--(-.5,-2); 
  \fill[white] (.3,.3) circle(.2);
    \draw[thick] (.3,.3) circle(.2);
    \node at (-.5,-1.1) {$v$};     
    \node at (-1,.5) {$\iota$}; 
    \node at (.6,-.3) {$B$};
    \node at (-.8,-.3) {$A$};
} ,\qquad
v^\BB=\iota^\pm(w^{\BA*})v^\pm = 
\tikzmatht{
  \fill[black!10] (-1.8,-2) rectangle (1.2,1);
  \fill[black!18] (-1.8,-2)--(-.5,-2)--(-.5,-1.5)--(-1.3,-.7)--
(-1.3,1)--(-1.8,1)--(-1.8,-2) ;
  \fill[white] (-.5,-1.5)--(-1.3,-.7)--(.3,-.7)--(-.5,-1.5);
    \draw[thick] (-.5,-1.5)--(-1.3,-.7)--(.3,-.7)--(-.5,-1.5);
    \draw[thick] (.3,-.7)--(.3,1); 
    \draw[thick] (-1.3,-.7)--(-1.3,1); 
    \draw[thick] (-.5,-.7)--(-.5,.3); 
    \draw[thick] (-.5,-1.5)--(-.5,-2); 
  \fill[white] (-.5,.3) circle(.2);
    \draw[thick] (-.5,.3) circle(.2);
    \node at (-.5,-1.1) {$v$};     
    \node at (-1,.5) {$\iota$}; 
    \node at (.6,-.3) {$B$};
    \node at (-.8,-.3) {$A$};
}.$$ 
Thus $M^\pm$ contain both $M^\BA=\iota^\pm(N)v^\BA$ and
$M^\BB=\iota^\pm(N)v^\BB$ as intermediate algebras. 
In $M^\pm$, the generators $v^\BA$ and $v^\BB$ satisfy the relations 
$$v^\BB v^\BA=\iota(\eps^\pm_{\theta^\BA,\theta^\BB})\cdot v^\BA v^\BB.$$
\end{lemma}

We can relate the braided product of Q-systems with the
$\alpha$-induction of Q-systems, \cref{c:alphaind}, as follows. 

\begin{tintedbox}
\begin{proposition} \vskip-5mm \label{p:alpha-Q}
Let $\iota^\BA:N\to M^\BA$ and $\iota^\BB:N\to M^\BB$, and 
$\BA=(\theta^\BA,w^\BA,x^\BA)$ and $\BB=(\theta^\BB,w^\BB,x^\BB)$
the associated Q-systems in a braided C* tensor category
$\C\subset \End_0(N)$.   
Denote by
$$\alpha^\pm(\BB)=(\alpha_{\theta^\BB}^\pm,\iota^\BA(w^\BB),\iota^\BA(x^\BB))$$
the Q-system in $\End_0(M^\BA)$ obtained from $\BB$ by $\alpha$-induction 
along $\BA$ (\cref{c:alphaind}(ii)). Then $\alpha^\pm(\BB)$ is the Q-system for the extension
$M^\BA\subset M^\mp$ in the diagram \eref{interm}. 

More precisely, if we write the extensions corresponding to
the braided products $\BA\times^\pm\BB$ as
$\iota^\pm:N\to M^\pm$, and the extension corresponding to
$\alpha^\pm(\BB)$ as $\jmath^\BB{}^\pm:M^\BA\to M^{\alpha\pm}$, 
such that $\alpha_{\theta^\BB}^\pm=\ol\jmath^\BB{}^\pm\jmath^\BB{}^\pm$,
then we have $M^{\alpha\pm}=M^\mp$ and
$$\iota^\mp=\jmath^\BB{}^\pm\circ\iota^\BA.$$
\end{proposition}
\end{tintedbox}

{\em Proof:} It suffices to verify that the composite Q-system 
according to \lref{l:addmult}(i) arising by the composition of embeddings 
$\iota^\BA:N\to M^\BA$ and $\jmath^\BB{}^\pm:M^\BA\to
M^{\alpha\pm}$, coincides with $\BA\times^\mp\BB=(\Theta,W,X^\mp)$. Indeed,
by the definitions and \eref{extend} we have 
$$\ol\iota^\BA\ol\jmath^\BB{}^\pm\circ\jmath^\BB{}^\pm\iota^\BA = 
\ol\iota^\BA\alpha^\pm_{\theta^\BB}\iota^\BA =
\ol\iota^\BA\iota^\BA\theta^\BB = \theta^\BA\theta^\BB = \Theta
,$$
$$\ol\iota^\BA(\iota^\BA(w^\BB))w^\BA= \theta_1(w^\BB)w^\BA=W,$$ 
and, denoting the generator of $\alpha^\pm(\BB)$ by $v^\pm$,
such that $\ol\jmath^\BB{}^\pm(v^\pm)=\iota^\BA(x^\BB)$: 
\bea \notag
\ol\iota^\BA\ol\jmath^\BB{}^\pm
\big[\jmath^\BB{}^\pm(v^\BA)v^\pm\big] =  
\ol\iota^\BA\big[\alpha^\pm_{\theta^\BB}(v^\BA)\iota^\BA(x^\BB)\big] = \qquad
\\ \notag 
\ol\iota^\BA\big[\iota^\BA(\eps_{\theta^\BA,\theta^\BB}^\mp)v^\BA
\iota^\BA(x^\BB)\big] = 
\theta^\BA(\eps_{\theta^\BA,\theta^\BB}^\mp)x^\BA \theta^\BA(x^\BB) = X^\mp.
\eea
\qed

Of course, a similar result is true for the $\alpha$-induction of
the Q-system $\BA$ to a Q-system in $M^\BB$, namely $\alpha^\pm(\BA)$
is the Q-system for $M^\BB$ in the braided product of extensions
corresponding to $\BB\times^\mp \BA$, which is in turn unitarily
equivalent to the braided product of extensions corresponding to 
$\BA\times^\pm \BB$. 

\medskip

As mentioned before, the braided product of two extensions
may fail to be irreducible, or to be a factor, even if both extensions
are irreducible. For the braided product of two {\em commutative}
extensions, the centre equals the relative commutant. This result is
of particular interest in the applications to local QFT, where phase
boundaries are described by the braided product of two local 
extensions \cite{BKLR}. 

\begin{tintedbox}
\begin{proposition} \vskip-5mm \label{p:commprod}
Let $\BA=(\theta^\BA,w^\BA,x^\BA)$ and $\BB=(\theta^\BB,w^\BB,x^\BB)$ be two 
commutative Q-systems in a braided category, and 
$\BA\times^\pm \BB = (\theta,w,x)$ the product Q-system 
(with either braiding). Let $N\subset M$ be the corresponding 
braided product of extensions. Then the centre $M'\cap M$ of $M$
equals the relative commutant $\iota(N)'\cap M$. 
\end{proposition}
\end{tintedbox}

{\em Proof:} In view of \lref{l:relc+cent}, we have to show that 
every $q\in\Hom(\theta^\BA\theta^\BB,\id_N)$ satisfies \eref{hom-c}. 
Let $q\in\Hom(\theta^\BA\theta^\BB,\id_N)$. Then 
$$\tikzmatht{
  \fill[black!10] (-1.6,-1.8) rectangle (1.2,1.5);
    \draw[thick] (-.6,1.5)--(-.6,0);
    \draw[thick] (-1.4,1.5)--(-1.4,0);
  \fill[white] (1,1) rectangle (.2,.2);
    \draw[thick] (1,1) rectangle node{$q$} (.2,.2);
    \draw[thick] (1,.2)--(1,0) arc(360:260:.8);
    \draw[thick] (-.6,.2)--(-.6,0) arc(180:220:.8);
    \draw[thick] (-1.4,.2)--(-1.4,0) arc(180:360:.8)--(.2,.2);
    \node at (-.9,-1.5) {$A$};
    \node at (.5,-1.5) {$B$};
    \draw[thick] (-.6,-1.8)--(-.6,-.8);           
  \fill[black] (-.6,-.8) circle (.12); 
    \draw[thick] (.2,-1.8)--(.2,-.8);
  \fill[black] (.2,-.8) circle (.12);
} =
\tikzmatht{
  \fill[black!10] (-1.4,-1.8) rectangle (1.4,1.5);
  \fill[white] (-.4,1) rectangle (.4,.2);
    \draw[thick] (-.4,1) rectangle node{$q$} (.4,.2);
    \draw[thick] (-1.2,1.5)--(-1.2,-.5) arc(180:360:.4)--(-.4,.2);
    \draw[thick] (.4,.2)--(.9,-.3) arc(405:135:.35)--(.5,-.2)
    (.8,.1)--(1,.3) arc(-45:0:.6)--(1.2,1.5);
    \node at (-1.1,-1.5) {$A$};
    \node at (.9,-1.5) {$B$};
    \draw[thick] (-.8,-1.8)--(-.8,-.9);           
  \fill[black] (-.8,-.9) circle (.12); 
    \draw[thick] (.6,-1.8)--(.6,-.9);
  \fill[black] (.6,-.9) circle (.12);
},\qquad 
\tikzmatht{
  \fill[black!10] (-1.6,-1.8) rectangle (1.2,1.5);
    \draw[thick] (.2,1.5)--(.2,0);
    \draw[thick] (1,1.5)--(1,0);
  \fill[white] (-1.4,1) rectangle (-.6,.2);
    \draw[thick] (-1.4,1) rectangle node{$q$} (-.6,.2);
    \draw[thick] (1,.2)--(1,0) arc(360:260:.8);
    \draw[thick] (-.6,.2)--(-.6,0) arc(180:220:.8);
    \draw[thick] (-1.4,.2)--(-1.4,0) arc(180:360:.8)--(.2,.2);
    \node at (-.9,-1.5) {$A$};
    \node at (.5,-1.5) {$B$};
    \draw[thick] (-.6,-1.8)--(-.6,-.8);           
  \fill[black] (-.6,-.8) circle (.12); 
    \draw[thick] (.2,-1.8)--(.2,-.8);
  \fill[black] (.2,-.8) circle (.12);
} = 
\tikzmatht{
  \fill[black!10] (1.4,-1.8) rectangle (-1.4,1.5);
  \fill[white] (.4,1) rectangle (-.4,.2);
    \draw[thick] (.4,1) rectangle node{$q$} (-.4,.2);
    \draw[thick] (1.2,1.5)--(1.2,-.5) arc(360:180:.4)--(.4,.2);
    \draw[thick] (-.4,.2)--(-.5,.1) (-.8,-.2)--(-.9,-.3)
    arc(135:415:.35)--(-1,.4) arc(225:180:.6)--(-1.2,1.5);
    \node at (1.1,-1.5) {$B$};
    \node at (-1,-1.5) {$A$};
    \draw[thick] (.8,-1.8)--(.8,-.9);           
  \fill[black] (.8,-.9) circle (.12); 
    \draw[thick] (-.7,-1.8)--(-.7,-.9);
  \fill[black] (-.7,-.9) circle (.12);
}.
$$
If both Q-systems are commutative, the two expressions are the same.
\qed

\subsection{The full centre}
\label{s:Qfullc}
\setcounter{equation}{0}
\begin{definition} \label{d:fullcent} (\cite{FFRS06})
Let $\BA=(\theta,w,x)$ be a Q-system in $\C$. It trivially gives
rise to a Q-system $\BA\otimes \eins = (\theta_i\otimes \id_N, w\otimes
\eins_N,x\otimes \eins_N)$ in $\C\boxtimes \C\opp$. Let $\BR$ be the
canonical Q-system in $\C\boxtimes \C\opp$. Then the {\bf full centre} 
of $\BA$ is defined as the commutative Q-system in $\C\boxtimes\C\opp$ 
given by the left centre of the $\times^+$-product
\be 
Z[\BA]= C^+[(\BA\otimes \eins)\times^+ \BR].
\ee
\end{definition}

Because $\dim\Hom(\id,(\theta\otimes\id)\Theta\can) =
\dim\Hom(\id,\theta)$, and the centre projection can only decrease
multiplicities, the full centre is irreducible if $\BA$ is
irreducible. For a stronger statement, see \pref{p:fullsimple}.

\begin{tintedbox}
\begin{proposition} \vskip-5mm \label{p:fullc=alpha} \cite[Prop.\ 4.18]{BKL} 
The full centre equals the $\alpha$-induction construction in \cite{R00}.
\end{proposition}
\end{tintedbox}

This result was conjectured in \cite{KR10}, and proven in \cite{BKL}. 
In fact, it is rather easy to show that both the $\alpha$-induction
construction and the full centre give 
$$[\Theta]=\bigoplus Z_{\rho,\sig}\;[\rho]\otimes[\ol\sig]$$
with the multiplicities $Z_{\rho,\sig}$ given by \eref{Z}; 
whereas the equality of the respective intertwiners $X$ is more
difficult to establish. 


\begin{remark} 
\label{r:alpha}
The $\alpha$-induction construction \cite{R00} was originally found as a
  construction of two-dimensional local conformal QFT models out of
  chiral data, cf.\ \sref{s:b-ext}. It is in fact a 
  construction of commutative Q-systems in $\C\boxtimes\C\opp$ out of
  a Q-system in $\C$, using the $\alpha$-induction (\sref{s:alpha}) to 
  extend $\rho\in\End(N)$ to $\alpha_\rho^\pm\in \End(M)$. In the
  simplest case, when the Q-system in $\C$ is trivial or Morita
  equivalent to the trivial Q-system, then one obtains the canonical
  Q-system \pref{p:canonical} in $\C\boxtimes\C\opp$. A more general
  analysis is given in \cite{I,BEK99,BEK00}. 
\end{remark}

\begin{tintedbox}
\begin{proposition} \vskip-5mm \label{p:fullsimple} 
Let $\BA$ be a Q-system in a braided C* tensor category. 
The full centre $Z[\BA]$ is irreducible iff $\BA$ is simple, i.e., iff
the extension described by $\BA$ is a factor 
(cf.\ \cref{c:simple=factorial}). More generally, the following linear 
spaces have equal dimension: \\
{\rm (i)} $\Hom(\id\otimes\id,Z[\BA])$ \\
{\rm (ii)} $\Hom(\id,C^+[\BA])$ and $\Hom(\id,C^-[\BA])$ \\
{\rm (iii)} The centre $M'\cap M$ of the extension described by $\BA$. 
\end{proposition}
\end{tintedbox}

{\em Proof:} The projection defining the full centre is a multiple of 
$$\bigoplus_\sig 
\tikzmatht{
  \fill[black!10] (-1.1,-1) rectangle (1.7,1.4);
    \draw[thick] (0,-.5) arc(270:110:.7);
    \draw[thick] (0,-.5)--(.6,-.5) arc(270:430:.7) (.2,.9)--(.4,.9);
  \fill[black] (0,.-.5) circle (.12);
    \draw[thick] (0,-1)--(0,1.4);
    \draw[thick] (.6,-1)--(.6,-.7) (.6,-.3)--(.6,1.4);
    \node at (.9,-.8) {$\sig$};
} \otimes 
\tikzmatht{
  \fill[black!10] (-1,-1) rectangle (1,1.4);
    \draw[thick] (0,-1)--(0,1.4);
    \node at (.6,-.7) {$\ol\sig$};
}.$$ 
Therefore, for the multiplicity of the identity in $Z[\BA]$, we can
replace the canonical Q-system $\BR$ by the trivial Q-system
$\id\otimes\id$ in $\C \boxtimes\C\opp$. Then trivially, 
$\dim\Hom(\id\otimes\id,Z[\BA])=\dim\Hom(\id\otimes\id,C^+[\BA\otimes\id])
=\dim\Hom(\id,C^+[\BA])$. Writing the centre projection as
$p_\BA^+=SS^*$, we have $t\in \Hom(\id,C^+[\BA])\subset \Hom(\id,\theta_{P})$
iff $tn = \theta_{P}(n)t = S^*\theta(n)St$ iff $q=St\in \Hom(\id,\theta)$
satisfies $qn = P\theta(n)q=Pqn$ for all $n\in N$, i.e.,
$q=Pq$. Then \lref{l:relc+cent}(iii) together with the following Lemma 
prove the claim. 
\qed

\begin{lemma} Let $\BA$ be a Q-system in a braided C* tensor category,
  and $P^\pm\equiv p_\BA^+$ its centre projections. Then $q\in\Hom(\id,\theta)$
satisfies $qP^+=q$ iff $qP^-=q$ iff $q$ satisfies \eref{hom-c}. 
\end{lemma}

{\em Proof:} We have 
$$qP^+ = d_\BA\inv\cdot
\tikzmatht{
  \fill[black!10] (-1.3,-1) rectangle (1.3,1.5);
    \draw[thick] (0,-.5) arc(270:110:.6);
    \draw[thick] (0,-.5) arc(270:430:.6);
  \fill[white] (0,1.4)--(-.2,1.1)--(.2,1.1)--(0,1.4);
    \draw[thick] (0,1.4)--(-.2,1.1)--(.2,1.1)--(0,1.4);
    \draw[thick] (0,-1)--(0,1.1);
  \fill[black] (0,.-.5) circle (.12);
} = d_\BA\inv\cdot
\tikzmatht{
  \fill[black!10] (-1.3,-1) rectangle (1.3,1.5);
    \draw[thick] (0,.2) circle(.7);
  \fill[white] (0,.4)--(-.2,.1)--(.2,.1)--(0,.4);
    \draw[thick] (0,.4)--(-.2,.1)--(.2,.1)--(0,.4);
    \draw[thick] (0,-1)--(0,.1);
  \fill[black] (0,.-.5) circle (.12);
} = d_\BA\inv\cdot
\tikzmatht{
  \fill[black!10] (-1.3,-1) rectangle (1.3,1.5);
    \draw[thick] (0,.1) circle(.6);
  \fill[white] (0,1.4)--(-.2,1.1)--(.2,1.1)--(0,1.4);
    \draw[thick] (0,1.4)--(-.2,1.1)--(.2,1.1)--(0,1.4);
    \draw[thick] (0,-1)--(0,.5) (0,.9)--(0,1.1);
  \fill[black] (0,.-.5) circle (.12);
} = 
 qP^-.$$
If $qP^\pm = q$, then $q$ satisfies \eref{hom-c} (using associativity):
$$
\tikzmatht{
  \fill[black!10] (-1.3,-1) rectangle (1.3,1.5);
    \draw[thick] (.6,1.5)--(.6,.2) arc(360:180:.6)--(-.6,.8);
  \fill[white] (-.6,1)--(-.8,.7)--(-.4,.7)--(-.6,1);
    \draw[thick] (-.6,1)--(-.8,.7)--(-.4,.7)--(-.6,1);
    \draw[thick] (0,-1)--(0,-.4); 
  \fill[black] (0,-.4) circle (.12);
    \node at (-1,1.1) {$q$};
} = d_\BA\inv\cdot
\tikzmatht{
  \fill[black!10] (-1.3,-1) rectangle (1.5,1.5);
    \draw[thick] (0,.5) circle(.7);
  \fill[white] (0,.7)--(-.2,.4)--(.2,.4)--(0,.7);
    \draw[thick] (0,.7)--(-.2,.4)--(.2,.4)--(0,.7);
    \draw[thick] (0,-.6) arc(270:360:1.1) -- (1.1,1.5); 
    \draw[thick] (0,-1)--(0,.4);
  \fill[black] (0,-.6) circle (.12);
  \fill[black] (0,-.2) circle (.12);
} = d_\BA\inv\cdot
\tikzmatht{
  \fill[black!10] (-1.3,-1) rectangle (1.3,1.5);
    \draw[thick] (0,.2) circle(.7);
  \fill[white] (0,.4)--(-.2,.1)--(.2,.1)--(0,.4);
    \draw[thick] (0,.4)--(-.2,.1)--(.2,.1)--(0,.4);
    \draw[thick] (0,.9)--(0,1.5);
    \draw[thick] (0,-1)--(0,.1);
  \fill[black] (0,.9) circle (.12);
  \fill[black] (0,.-.5) circle (.12);
} = d_\BA\inv\cdot
\tikzmatht{
  \fill[black!10] (-1.5,-1) rectangle (1.3,1.5);
    \draw[thick] (0,.5) circle(.7);
  \fill[white] (0,.7)--(-.2,.4)--(.2,.4)--(0,.7);
    \draw[thick] (0,.7)--(-.2,.4)--(.2,.4)--(0,.7);
    \draw[thick] (0,-.6) arc(270:180:1.1) -- (-1.1,1.5); 
    \draw[thick] (0,-1)--(0,.4);
  \fill[black] (0,-.6) circle (.12);
  \fill[black] (0,-.2) circle (.12);
} = 
\tikzmatht{
  \fill[black!10] (-1.3,-1) rectangle (1.3,1.5);
    \draw[thick] (-.6,1.5)--(-.6,.2) arc(180:360:.6)--(.6,.8);
  \fill[white] (.6,1)--(.8,.7)--(.4,.7)--(.6,1);
    \draw[thick] (.6,1)--(.8,.7)--(.4,.7)--(.6,1);
    \draw[thick] (0,-1)--(0,-.4);
  \fill[black] (0,-.4) circle (.12);
    \node at (1,1.1) {$q$};
}.
$$
Conversely, if $q$ satisfies \eref{hom-c}, then 
$$qP^+ = d_\BA\inv\cdot
\tikzmatht{
  \fill[black!10] (-1.3,-1) rectangle (1.3,1.5);
    \draw[thick] (0,-.5) arc(270:110:.6);
    \draw[thick] (0,-.5) arc(270:430:.6);
  \fill[white] (0,1.4)--(-.2,1.1)--(.2,1.1)--(0,1.4);
    \draw[thick] (0,1.4)--(-.2,1.1)--(.2,1.1)--(0,1.4);
    \draw[thick] (0,-1)--(0,1.1);
  \fill[black] (0,.-.5) circle (.12);
} = d_\BA\inv\cdot
\tikzmatht{
  \fill[black!10] (-1.3,-1) rectangle (1.3,1.5);
    \draw[thick] (0,-.5) arc(270:180:.8) arc(180:0:.4)--(0,-.5);
    \draw[thick] (0,-.5) arc(270:360:.8)--(.8,1.1);
  \fill[white] (.8,1.2)--(.6,.9)--(1,.9)--(.8,1.2);
    \draw[thick] (.8,1.2)--(.6,.9)--(1,.9)--(.8,1.2);
    \draw[thick] (0,-1)--(0,.1);
  \fill[black] (0,.-.5) circle (.12);
} = 
\tikzmatht{
  \fill[black!10] (-1.3,-1) rectangle (1.3,1.5);
  \fill[white] (0,.6)--(-.2,.3)--(.2,.3)--(0,.6);
    \draw[thick] (0,.6)--(-.2,.3)--(.2,.3)--(0,.6);
    \draw[thick] (0,-1)--(0,.3);
} = q
.$$
\qed

\subsection{Modular tensor categories}
\label{s:modular}
\setcounter{equation}{0}
A C* tensor category with finitely many inequivalent irreducible 
objects (denoted $\rho,\sig,\tau$, etc.), all of finite dimension, 
is called {\em rational}. In a {\em braided} rational C* tensor
category, one can introduce the matrices 
$$S_{\sig,\tau}:= \dim(\C)^{-\frac12}\cdot
\tikzmatht{
  \fill[black!10] (-1.5,-1) rectangle (1.5,1);
    \draw[thick] (-1.3,0) arc(180:70:.8); 
    \draw[thick] (-1.3,0) arc(180:390:.8); 
    \draw[thick] (1.3,0) arc(0:210:.8); 
    \draw[thick] (1.3,0) arc(360:250:.8);
    \node at (-.6,0) {$\sig$}; 
    \node at (.6,0) {$\tau$}; 
} = S_{\tau,\sig}, \quad 
T^0_{\sig,\tau}:= \frac{\delta_{\sig,\tau}}{\dim(\tau)}\cdot 
\tikzmatht{
  \fill[black!10] (-1.5,-1) rectangle (1.5,1);
    \draw[thick] (0,0)--(.4,.4) arc(135:-135:.58)--(.2,-.2);
    \draw[thick] (0,0)--(-.4,-.4) arc(315:45:.58)--(-.2,.2);
    \node at (.3,.6) {$\tau$}; 
} \equiv \delta_{\sig,\tau}\cdot \kappa_\tau,
$$
where $\dim(\C)=\sum_\rho \dim(\rho)^2$ is the global dimension
\eref{globaldim}, and $\kappa_\tau$ is the twist (\dref{d:twist}).   

\begin{definition}\label{d:modular}
A braiding of a tensor category $\C$ is called {\bf non-degenerate} if
there is no nontrivial sector $[\rho]$ such that
$\eps_{\rho,\sig}^+=\eps_{\rho,\sig}^-$ for all $\sig\in\C$. A braided 
rational C* tensor category is called {\bf modular}, if the symmetric
matrix $S$ is invertible.  
\end{definition}

\begin{proposition}\label{p:modular} \cite{R89} 
A braided rational C* tensor category is modular if and only if it is
non-degenerate. In this case, the matrix $S$ is unitary, and there is
a complex phase $\omega$ (unique up to a third root of unity) such that the
matrices $S$ and $T:=\omega\cdot T^0$ form a unitary representation 
of the modular group $SL(2,\ZZ)$:
$$(ST\inv)^3=S^2, \quad S^4=E.$$
Moreover, 
$S_{\sig,\ol\tau}=\ol{S_{\sig,\tau}}=S_{\ol\sig,\tau}$, i.e., the
central element $S^2$ of $SL(2,\ZZ)$ is represented by the
conjugation matrix $C$.   
\end{proposition}

Recall that $\dim(\C)^{\frac12}=d_\BR$ is also the dimension of
the canonical Q-system in $\C\boxtimes \C\opp$
(\pref{p:canonical}). By considering the $\id$-$\id$-component of the
equation $T\inv ST\inv ST\inv =S$, one finds that
$\omega^3=\sum_\tau\kappa_\tau\inv\dim(\tau)^2/d_\BR$.

All the braidings mentioned in \xref{x:Ibraid} are non-degenerate,
giving rise to eight inequivalent modular categories associated with
the same ``fusion rules'' of three irreducible sectors. 

\begin{lemma}\label{l:ring} For $\tau$ and $\sig$ irreducible, one has
  in a modular category
$$\RTr_\sig(\eps_{\sig,\tau}\eps_{\tau,\sig}) \equiv 
\tikzmatht{
  \fill[black!10] (-1.2,-1) rectangle (1.5,1);
    \draw[thick] (-.3,1)--(-.1,.8) (.2,.5) arc(35:-35:.8)--(-.3,-1); 
    \draw[thick] (1.3,0) arc(0:210:.8); 
    \draw[thick] (1.3,0) arc(360:250:.8);
    \node at (-.6,.7) {$\tau$}; 
    \node at (.6,0) {$\sig$}; 
} = 
\frac{d_\BR\cdot S_{\tau,\sig}}{\dim(\tau)}\cdot 1_\tau =
\LTr_\sig(\eps^*_{\sig,\ol\tau}\eps^*_{\tau,\ol\sig}).$$ 
\end{lemma}

{\em Proof:} Clearly, $\RTr_\sig(\eps_{\sig,\tau}\eps_{\tau,\sig})$ is
a multiple of $1_\tau$. Thus, one can compute the coefficient by
applying $\Tr_\tau$, where $\Tr_\tau(1_\tau)=\dim(\tau)$. Similar for
the second equation.
\qed

\begin{proposition} \label{p:killing} (The ``killing ring'')
For $\rho$ an object of $\C$, consider $\rho\otimes\id$ as
an object of $\C\boxtimes \C\opp$. If $\C$ is modular, then 
$$\tikzmatht{
  \fill[black!10] (-2,-1) rectangle (1.5,1);
    \draw[thick] (-.3,1)--(-.1,.8) (.2,.5) arc(35:-35:.8)--(-.3,-1); 
    \draw[thick] (1.3,0) arc(0:210:.8); 
    \draw[thick] (1.3,0) arc(360:250:.8);
    \node at (-1.2,.7) {$\rho\!\otimes\!\id$}; 
    \node at (-.7,-.2) {$\Theta$}; 
} 
= d_\BR^2\cdot E_{\id} = 
\tikzmatht{
  \fill[black!10] (-1.5,-1) rectangle (2,1);
    \draw[thick] (.3,-1)--(.1,-.8) (-.2,-.5) arc(215:145:.8)--(.3,1); 
    \draw[thick] (-1.3,0) arc(180:390:.8); 
    \draw[thick] (-1.3,0) arc(180:70:.8);
    \node at (1.2,.7) {$\rho\!\otimes\!\id$}; 
    \node at (.7,-.2) {$\Theta$}; 
},
$$ 
where $\Theta$ is the endomorphism of the canonical Q-system, \cref{c:opp},
and $E_\id= \tikzmatht{
  \fill[black!10] (-1,-1) rectangle (1,1);
    \draw[thick] (0,1)--(0,.6) (0,-1)--(0,-.6); 
  \fill[white] (0,.2)--(.3,.6)--(-.3,.6)--(0,.2); 
    \draw[thick] (0,.2)--(.3,.6)--(-.3,.6)--(0,.2); 
  \fill[white] (0,-.2)--(.3,-.6)--(-.3,-.6)--(0,-.2); 
    \draw[thick] (0,-.2)--(.3,-.6)--(-.3,-.6)--(0,-.2); 
}
\in \Hom(\rho,\rho)$ is the projection on the identity component
$\id\prec\rho$ (which is zero if $\id$ is not contained in $\rho$). 
\end{proposition}

{\em Proof:} If $\tau$ is irreducible, then 
$$\tikzmatht{
  \fill[black!10] (-2,-1) rectangle (1.5,1);
    \draw[thick] (-.3,1)--(-.1,.8) (.2,.5) arc(35:-35:.8)--(-.3,-1); 
    \draw[thick] (1.3,0) arc(0:210:.8); 
    \draw[thick] (1.3,0) arc(360:250:.8);
    \node at (-1.2,.7) {$\tau\!\otimes\!\id$}; 
    \node at (-.7,-.2) {$\Theta$}; 
} = \sum_\sig
\tikzmatht{
  \fill[black!10] (-1.2,-1) rectangle (1.5,1);
    \draw[thick] (-.3,1)--(-.1,.8) (.2,.5) arc(35:-35:.8)--(-.3,-1); 
    \draw[thick] (1.3,0) arc(0:210:.8); 
    \draw[thick] (1.3,0) arc(360:250:.8);
    \node at (-.6,.7) {$\tau$}; 
    \node at (-.6,0) {$\sig$}; 
}\otimes 
\tikzmatht{
  \fill[black!10] (-1,-1) rectangle (1,1);
    \draw[thick] (0,0) circle(.8);
    \node at (-.5,0) {$\ol\sig$};
} = \sum_\sig \frac{d_\BR\cdot
  S_{\tau,\sig}}{\dim(\tau)}\cdot 1_\tau \cdot \dim(\sig).
$$
Then, $\dim(\sig)=d_\BR\cdot S_{\sig,\id}=d_\BR\cdot
\ol{S_{\sig,\id}}$ and unitarity of $S$ yield $d_\BR^2$ if $\tau=\id$
and zero otherwise. If $\rho$ is reducible, then write $1_\rho=\sum_\tau
E_\tau$ where $E_\tau\in\Hom(\rho,\rho)$ are the projections on the
irreducible $\tau\prec\rho$. Under the ``killing ring'', only
$\tau=\id$ survives. 
\qed

We can now see that it was essential to choose matching signs in the
definition \dref{d:fullcent} of the full centre:

\begin{corollary} \label{c:cR=R} For $\BA$ an irreducible Q-system in
  $\C$ and $\BR$ the canonical Q-system, one has
$$C^-[(\BA\otimes \eins)\times^+ \BR] = C^+[(\BA\otimes \eins)\times^-
\BR] = \BR.$$ 
\end{corollary}

{\em Proof:} By using \pref{p:killing} and the fact that $\BR$ is
commutative, one can compute the trace $\Tr(p^\pm)$ of the respective
  centre projections of $(\BA\otimes \eins)\times^{\mp} \BR$. The result is 
  $\Tr(p^\pm)=d_\BR^2$. On the other hand, the projection $p_\BR$ onto the 
  intermediate Q-system $\BR=(1\otimes\eins)\times^{\mp} \BR \prec
  (\BA\otimes\eins) 
  \times^{\mp} \BR$ satisfies \eref{centp}, hence $p_\BR< p^\pm$ by
  \pref{p:maximal}. Since by \pref{p:traces},
  $\Tr(p_\BR)=\dim(\Theta\can) = d_\BR^2$, the claim follows. 
\qed

\subsection{The braided product of two full centres}
\label{s:class}
\setcounter{equation}{0}
We assume $\C$ to be modular. 

The following \tref{t:centchar} provides the minimal central
projections for the braided product of two commutative Q-systems which
arise as full centres. By way of preparation of this result, let us
compile several equivalent ways of describing the centre. 

Recall that the centre $M'\cap M$ of the extension corresponding to 
the braided product of two commutative Q-systems equals the relative
commutant $\iota(N)'\cap M = \iota(\Hom(\Theta^\BA\Theta^\BB,\id))V$ by
\lref{l:relc+cent} and \pref{p:commprod}. The space
$\Hom(\Theta^\BA\Theta^\BB,\id)$ is isomorphic to
$\Hom(\Theta^\BB,\Theta^\BA)$ by Frobenius reciprocity. Thus, there is a
linear bijection 
\be \label{chi}
\chi:\Hom(\Theta^\BB,\Theta^\BA)\to M'\cap M,\quad 
\chi(T):=\iota\big(R^{\BA*}\scirc
(1_{\Theta^\BA}\times T)\big)V =
\tikzmatht{
  \fill[black!10] (-1.8,-2) rectangle (1.2,1.8);
  \fill[black!18] (-1.8,-2)--(-.5,-2)--(-.5,-1.7)--(-1.5,-.7)--
(-1.5,1.8)--(-1.8,1.8)--(-1.8,-2) ;
    \draw[thick] (.5,1.2) arc(0:180:.5)--(-.5,-.7); 
    \draw[thick] (.5,.4)--(.5,-.7);
  \fill[white] (1,1.2) rectangle (0,.4);
    \draw[thick] (1,1.2) rectangle node{$T$} (0,.4);
  \fill[white] (-.5,-1.7)--(-1.5,-.7)--(.5,-.7)--(-.5,-1.7);
    \draw[thick] (-.5,-1.7)--(-1.5,-.7)--(.5,-.7)--(-.5,-1.7);
    \draw[thick] (-1.5,-.7)--(-1.5,1.8); 
    \draw[thick] (-.5,-1.7)--(-.5,-2); 
    \node at (-.5,-1.1) {$V$};     \node at (-1.2,1.1) {$\iota$}; 
    \node at (.9,-.2) {$B$};
    \node at (-.8,-.2) {$A$};
}
\ee
with inverse
$$\chi\inv(\cdot) = [1_{\Theta^\BA}\times
(W^*\scirc\ol\iota(\cdot))]\scirc R^\BA.$$
Notice also that $\ol\iota$ maps the centre into
$\Hom(\Theta^\BA\Theta^\BB,\Theta^\BA\Theta^\BB)$: 
$$\ol\iota\chi(T) = \Big(1_{\Theta^\BA\Theta^\BB}\times \big(R^{\BA*}\scirc
(1_{\Theta^\BA}\times T)\big)\Big)\scirc X = 
\tikzmatht{
  \fill[black!10] (-1.6,-1.8) rectangle (1.6,1.8);
    \draw[thick] (1,.9)--(1,1.1) arc(0:180:.4)--(0.2,0); 
    \draw[thick] (-.6,1.8)--(-.6,0);
    \draw[thick] (-1.4,1.8)--(-1.4,0);
  \fill[white] (1.5,.9) rectangle (.5,.1);
    \draw[thick] (1.5,.9) rectangle node{$T$} (.5,.1);
    \draw[thick] (1,.1)--(1,0) arc(360:260:.8);
    \draw[thick] (-.6,.4)--(-.6,0) arc(180:220:.8);
    \draw[thick] (-1.4,0) arc(180:360:.8);
    \node at (-1.1,1.5) {$A$};
    \node at (-.3,1.5) {$B$};
    \draw[thick] (-.6,-1.8)--(-.6,-.8);           
  \fill[black] (-.6,-.8) circle (.12); 
    \draw[thick] (.2,-1.8)--(.2,-.8);
  \fill[black] (.2,-.8) circle (.12);
    \node at (-.3,-1.5) {$A$};
    \node at (.6,-1.5) {$B$};
}  \stapel{\nref{comm}}= 
\tikzmatht{
  \fill[black!10] (-1.6,-1.8) rectangle (1.6,1.8);
    \draw[thick] (.4,.9)--(.4,1.1) arc(0:180:.4)--(-.4,0); 
    \draw[thick] (-1.4,1.8)--(-1.4,0);
  \fill[white] (.9,.9) rectangle (-.1,.1);
    \draw[thick] (.9,.9) rectangle node{$T$} (-.1,.1);
    \draw[thick] (-1.4,0) arc(180:360:.5);
    \draw[thick] (1.4,1.8)--(1.4,0) arc(360:180:.5)--(.4,.1);
    \node at (-1.1,1.5) {$A$};
    \node at (1,1.5) {$B$};
    \draw[thick] (-.9,-1.8)--(-.9,-.5); 
  \fill[black] (-.9,-.5) circle (.12);
    \draw[thick] (.9,-1.8)--(.9,-.5);
  \fill[black] (.9,-.5) circle (.12);
    \node at (-.6,-1.5) {$A$};
    \node at (1.2,-1.5) {$B$};
} \equiv 
\tikzmatht{
  \fill[black!10] (-1.6,-1.8) rectangle (1.6,1.8);
    \draw[thick] (-1.2,1.8)--(-1.2,-1.8);
    \draw[thick] (1.2,1.8)--(1.2,-1.8);
    \draw[thick] (-1.2,.8)--(-.3,.8) arc(90:0:.4)--(0,-.4)
    arc(180:270:.4)--(1.2,-.8); 
  \fill[black] (-1.2,.8) circle (.12);
  \fill[black] (1.2,-.8) circle (.12);
  \fill[white] (-.4,-.4) rectangle (.4,.4);
    \draw[thick] (-.4,-.4) rectangle node{$T$} (.4,.4);
    \node at (.5,-1.2) {$B$};
    \node at (-.5,1.2) {$A$};
}$$
where we have used commutativity of $\BB$, and are freely appealing to 
Frobenius reciprocity in the last way of drawing the diagram. 

Then, one easily verifies that
$$\chi (T_1)\scirc\chi(T_2) = \chi(T_1\ast T_2), \qquad 
\ol\iota\chi (T_1)\scirc\ol\iota\chi(T_2) = \ol\iota\chi(T_1\ast T_2),$$
where $T_1\ast T_1$ is the commutative ``convolution'' product on
$\Hom(\Theta^\BB,\Theta^\BA)$ with
unit $W^\BA W^{\BB*}$: 
\be\label{convol} 
T_1\ast T_2:= 
\tikzmatht{
  \fill[black!10] (-1.4,-1.8) rectangle (1.4,1.8);
    \draw[thick] (-.6,.5) arc(180:0:.6); 
    \draw[thick] (-.6,-.5) arc(180:360:.6); 
    \draw[thick] (0,1.8)--(0,1.1); 
  \fill[black] (0,1.1) circle (.12);
    \draw[thick] (0,-1.8)--(0,-1.1);
  \fill[black] (0,-1.1) circle (.12);
  \fill[white] (-1.1,-.5) rectangle (-.1,.5);
    \draw[thick] (-1.1,-.5) rectangle node{$T_1$} (-.1,.5);
  \fill[white] (1.1,-.5) rectangle (.1,.5);
    \draw[thick] (1.1,.5) rectangle node{$T_2$} (.1,-.5);
}\stapel{\nref{comm}}=
\tikzmatht{
  \fill[black!10] (-1.4,-1.8) rectangle (1.4,1.8);
    \draw[thick] (0,-1.5) arc(270:130:.3)--(.5,-.5);
    \draw[thick] (0,-1.5) arc(-90:50:.3) (-.2,-.7)--(-.5,-.5);
    \draw[thick] (0,1.5) arc(90:230:.3)--(.5,.5);
    \draw[thick] (0,1.5) arc(90:-50:.3) (-.2,.7)--(-.5,.5);
    \draw[thick] (0,1.8)--(0,1.5); 
  \fill[black] (0,1.5) circle (.12);
    \draw[thick] (0,-1.8)--(0,-1.5);
  \fill[black] (0,-1.5) circle (.12);
  \fill[white] (-1.1,-.5) rectangle (-.1,.5);
    \draw[thick] (-1.1,-.5) rectangle node{$T_1$} (-.1,.5);
  \fill[white] (1.1,-.5) rectangle (.1,.5);
    \draw[thick] (1.1,.5) rectangle node{$T_2$} (.1,-.5);
}=
\tikzmatht{
  \fill[black!10] (-1.4,-1.8) rectangle (1.4,1.8);
    \draw[thick] (-.6,.5) arc(180:0:.6); 
    \draw[thick] (-.6,-.5) arc(180:360:.6); 
    \draw[thick] (0,1.8)--(0,1.1); 
  \fill[black] (0,1.1) circle (.12);
    \draw[thick] (0,-1.8)--(0,-1.1);
  \fill[black] (0,-1.1) circle (.12);
  \fill[white] (-1.1,-.5) rectangle (-.1,.5);
    \draw[thick] (-1.1,-.5) rectangle node{$T_2$} (-.1,.5);
  \fill[white] (1.1,-.5) rectangle (.1,.5);
    \draw[thick] (1.1,.5) rectangle node{$T_1$} (.1,-.5);
},\qquad
\tikzmatht{
  \fill[black!10] (-1.4,-1.8) rectangle (1.4,1.8);
    \draw[thick] (-.6,.5) arc(180:0:.6); 
    \draw[thick] (-.6,-.5) arc(180:360:.6); 
    \draw[thick] (0,1.8)--(0,1.1); 
  \fill[black] (0,1.1) circle (.12);
    \draw[thick] (0,-1.8)--(0,-1.1);
  \fill[black] (0,-1.1) circle (.12);
  \fill[white] (-1.1,-.5) rectangle (-.1,.5);
    \draw[thick] (-1.1,-.5) rectangle node{$T$} (-.1,.5);
  \fill[white] (.6,.5) circle(.2) (.6,-.5) circle(.2);
    \draw[thick] (.6,.5) circle(.2) (.6,-.5) circle(.2);
} = 
\tikzmatht{
  \fill[black!10] (-1,-1.8) rectangle (1,1.8);
    \draw[thick] (0,1.8)--(0,.5); 
    \draw[thick] (0,-1.8)--(0,-.5);
  \fill[white] (-.5,-.5) rectangle (.5,.5);
    \draw[thick] (-.5,-.5) rectangle node{$T$} (.5,.5);
} = 
\tikzmatht{
  \fill[black!10] (-1.4,-1.8) rectangle (1.4,1.8);
    \draw[thick] (-.6,.5) arc(180:0:.6); 
    \draw[thick] (-.6,-.5) arc(180:360:.6); 
    \draw[thick] (0,1.8)--(0,1.1); 
  \fill[black] (0,1.1) circle (.12);
    \draw[thick] (0,-1.8)--(0,-1.1);
  \fill[black] (0,-1.1) circle (.12);
  \fill[white] (1.1,-.5) rectangle (.1,.5);
    \draw[thick] (1.1,.5) rectangle node{$T$} (.1,-.5);
  \fill[white] (-.6,.5) circle(.2) (-.6,-.5) circle(.2);
    \draw[thick] (-.6,.5) circle(.2) (-.6,-.5) circle(.2);
}.
\ee
Likewise, the adjoint is given by 
$$\chi(T)^* = \chi(F(T)),\qquad\ol\iota\chi(T)^* =
\ol\iota\chi(F(T)),$$
where $F$ is the antilinear Frobenius conjugation on
$\Hom(\Theta^{\BB},\Theta^{\BA})$ 
$$F(T) = 
\tikzmatht{
  \fill[black!10] (-1.4,-1.6) rectangle (1.4,1.6);
    \draw[thick] (-1,-1.6)--(-1,.7) arc(180:0:.5)--(0,.5); 
    \draw[thick] (1,1.6)--(1,-.7) arc(360:180:.5)--(0,-.5); 
  \fill[white] (-.5,-.5) rectangle (.5,.5);
    \draw[thick] (-.5,-.5) rectangle node{$T^*$} (.5,.5);
} \stapel{\pref{p:dotstar}}= 
\tikzmatht{
  \fill[black!10] (-1.4,-1.6) rectangle (1.4,1.6);
    \draw[thick] (-1,1.6)--(-1,-.7) arc(180:360:.5)--(0,-.5); 
    \draw[thick] (1,-1.6)--(1,.7) arc(0:180:.5)--(0,.5); 
  \fill[white] (-.5,-.5) rectangle (.5,.5);
    \draw[thick] (-.5,-.5) rectangle node{$T^*$} (.5,.5);
}
\in \Hom(\Theta^{\BB},\Theta^{\BA}).$$
Therefore, finding the minimal projections $E_m\in M'\cap M$ is
equivalent to finding the minimal projections 
$I_m\in \Hom(\Theta^{\BB},\Theta^{\BA})$ w.r.t.\ the convolution
product, i.e., to solving the system 
\bea\label{projD}
\hbox{self-adjointness} && I_m^*= F(I_m), \nonumber \\ 
\hbox{idempotency} && I_m\ast I_m'=\delta_{mm'}\cdot I_m,
\qquad\qquad \\
\hbox{completeness} &&\sum_m I_m= W^{\BA}W^{\BB*}. \nonumber \eea
Minimality is ensured if the number of $I_m$ exhausts the dimension 
of $\Hom(\Theta^\BB,\Theta^\BA)$. We therefore have to solve these 
equations by a basis $I_m$ of $\Hom(\Theta^\BB,\Theta^\BA)$, and put
$E_m=\chi(I_m)$. 
Obviously, then also 
$P_m=\ol\iota(E_m) = 
\tikzmatht{
  \fill[black!10] (-1.6,-1.5) rectangle (1.6,1.5);
    \draw[thick] (-1.2,1.5)--(-1.2,-1.5);
    \draw[thick] (1.2,1.5)--(1.2,-1.5);
    \draw[thick] (-1.2,.8)--(-.3,.8) arc(90:0:.4)--(0,-.4)
    arc(180:270:.4)--(1.2,-.8); 
  \fill[black] (-1.2,.8) circle (.12);
  \fill[black] (1.2,-.8) circle (.12);
  \fill[white] (-.5,-.4) rectangle (.5,.4);
    \draw[thick] (-.5,-.4) rectangle node{$I_m$} (.5,.4);
    \node at (.5,-1.2) {$B$};
    \node at (-.5,1.2) {$A$};
} 
\in\Hom(\Theta^{\BA}\Theta^{\BB},\Theta^{\BA}\Theta^{\BB})$ will be projections.

The following theorem gives the solution to \eref{projD}, where
$I_\mm$ are labelled by the irreducibe $\BA$-$\BB$-bimodules $\mm$ in
$\C$. This result is of great interest for boundary conformal
QFT: it provides a bijection between chiral bimodules
and phase boundaries \cite{BKLR}. It therefore establishes the
link between our algebraic QFT approach to phase boundaries, and the TFT
approach by \cite{TFT1,TFT}. The fact that the central projections for
the braided product extension of two full centre Q-systems in 
$\C \boxtimes\C\opp$ are labelled by bimodules in $\C$, means in 
physical terms that the boundary conditions between two maximal 
local two-dimensional extensions is fixed by chiral data. 

\begin{tintedbox}
\begin{theorem} \vskip-5mm \label{t:centchar}
Let $\BA$ and $\BB$ be two simple Q-systems in a modular tensor
category $\C$, and let $Z[\BA]=(\Theta^\BA,W^\BA,X^\BA)$ and
$Z[\BB]=(\Theta^\BB,W^\BB,X^\BB)$ be their full 
centre Q-systems in $\C\boxtimes\C\opp$. Let $N\subset M$ be the extension
defined by either of the product Q-systems $Z[\BA]\times^\pm Z[\BB]$. 
Then $M$ has a centre given by $M'\cap M = \iota(N)'\cap M=
\Hom(\iota,\iota)$. The minimal central projections $E_\mm$ can be
characterized as follows. 

The irreducible $\BA$-$\BB$-bimodules $\mm=(\beta,m)$ naturally give rise to 
intertwiners $D_{R[\mm]\vert_Z}\in\Hom(\Theta^\BB,\Theta^\BA)$ 
(to be defined in the proof). Then 
$$I_\mm = \frac{\dim(\beta)}{d_\BA^2 d_\BB^2 d_\BR^2} \cdot D_{R[\mm]\vert_Z}$$
solve the system \eref{projD}. Thus $E_\mm =\chi(I_\mm)$ are the
minimal central projections. 
\end{theorem}
\end{tintedbox}

As a byproduct, we shall also prove:

\begin{proposition} \label{p:dimcent}
Let $\BA$ be a simple Q-system in a modular tensor category, so that
its centre $Z[\BA]$ is irreducible (\pref{p:fullsimple}). Then 
$d_{Z[\BA]}=d_\BR=\dim(\C)^{\frac12}$ equals the dimension of the
canonical Q-system. In particular, all irreducible full centres have
the same dimension.  
\end{proposition}

(This is not a new result, cf.\ \cite{KR08}, but the proof seems to be new.)

\bigskip

The proof of \tref{t:centchar} is rather lengthy, but it is worthwhile
because it paves the way to an efficient computation of the centre,
with ensuing classification results. The operators $I_\mm$ first
appeared in \cite{FFRS06}, but their 
idempotent property is not manifest there. It was proven in a more 
special case in \cite{KR08} (with the hindsight that the general case 
can be reduced to the special case by highly nontrivial properties of 
modular tensor categories). We attempt to give here a streamlined 
version of the proof that does not require the general theory of modular 
tensor categories. The use of the C*-structure of the DHR category 
allows for some substantial simplification as compared to \cite{KR08}. 

\medskip

{\em Proof of \tref{t:centchar}}: The statement about the centre is just
an instance of \lref{l:relc+cent} and \pref{p:commprod}, because the 
full centres are commutative. 

To prepare the solution of \eref{projD}, we associate with every
$\BA$-$\BB$-bimodule $\mm =(\beta,m\in\Hom(\beta,\theta^\BA\beta\theta^\BB))$ an intertwiner
$D_\mm\in\Hom(\theta^\BB,\theta^\BA)$ as follows
\cite{KR08}:
$$D_\mm:= \Tr_\beta \big(\eps_{\theta^\BA,\beta}\scirc
(1_{\theta^\BA\beta}\times r^{\BB*}) 
\scirc (m\times 1_{\theta^\BB})\big) : 
\qquad
\tikzmatht{
  \fill[black!10] (-.8,-1.5) rectangle (2.3,1.5);
    \draw[thick] (1.5,1.5)--(.5,.5) arc(135:225:.7) (.5,-.5)--(.87,-.87)
(1.13,-1.13)--(1.5,-1.5);
    \draw[thick] (.8,1.2) arc(45:180:.7)--(-.4,-.7) arc(180:315:.7)--
(1.5,-.5) arc(-45:45:.7)--(1.2,.8);
    \draw[thick] (1,-1) circle (.17);
    \node at (1.8,-1) {$m$};
    \node at (1.8,1.1) {$\theta^A$}; 
    \node at (1.3,0) {$\beta$};
} \equiv
\tikzmatht{
  \fill[black!10] (-1.5,-1.5) rectangle (1.5,1.5);
    \draw[thick] (0,1.5)--(0,-.43) (0,-.77)--(0,-1.5); 
    \draw[thick] (-.3,.6)--(-.5,.6) arc(90:270:.6)--(.5,-.6)
arc(-90:90:.6)--(.3,.6);
    \draw[thick] (0,-.6) circle (.17);
    \node at (.5,-1.1) {$\theta^B$};
    \node at (.5,1.1) {$\theta^A$};
    \node at (-1,.8) {$\beta$};}.
$$
(We freely use Frobenius reciprocity in the graphical
representations.) One easily sees 

\begin{lemma} \label{l:Dm} (cf.\ \cite{KR08})
The following statements hold. \\ 
{\rm (i)} $D_\mm$ depends only on the unitary equivalence
class of $\mm=(\beta,m)$. \\
{\rm (ii)} $D_\mm^*=D_{\ol \mm}$. \\
{\rm (iii)} $D_{\mm_1\oplus \mm_2}=D_{\mm_1}+D_{\mm_2}$. \\
{\rm (iv)} If $\mm=(\beta,m)$ is an $\BA$-$\BB$-bimodule and 
$\mm'=(\beta',m')$ a $\BB$-$\BC$-bimodule, hence $\mm\otimes_\BB \mm'$ 
an $\BA$-$\BC$-bimodule, then $D_{\mm}D_{\mm'}=d_\BB\cdot D_{\mm\otimes_\BB \mm'}$. \\
{\rm (v)} $w^{\BA*}D_\mm w^\BB = \dim(\mm)\equiv \dim(\beta)$ for
$\mm = (\beta,m)$. 
\end{lemma}

{\em Proof:} (i) follows because the ``closed $\beta$-line''
represents a trace, absorbing a unitary bimodule morphism $:\mm\to \mm'$. 
(ii) is proven in the same way as \lref{l:bimproj}, using the unitarity of
the twist. (iii) follows by  
$$\sum_i
\tikzmatht{
  \fill[black!10] (-2,-1.5) rectangle (2,1.5);
    \draw[thick] (0,1.5)--(0,-.63) (0,-.97)--(0,-1.5); 
    \draw[thick] (-.2,.8)--(-.8,.8) arc(90:270:.8)--(0,-.8);
    \draw[thick] (.2,.8)--(.8,.8) arc(90:-90:.8)--(0,-.8);
    \draw[thick] (0,-.8) circle (.17);
  \fill[white] (-.6,-.8)--(-.9,-.6)--(-.9,-1)--(-.6,-.8);
    \draw[thick] (-.6,-.8)--(-.9,-.6)--(-.9,-1)--(-.6,-.8);
  \fill[white] (.6,-.8)--(.9,-.6)--(.9,-1)--(.6,-.8);
    \draw[thick] (.6,-.8)--(.9,-.6)--(.9,-1)--(.6,-.8);
    \node at (-.7,-.2) {$s_i^*$};    
    \node at (.7,-.3) {$s_i$};
    \node at (-1.4,1) {$\beta$};
} = \sum_i
\tikzmatht{
  \fill[black!10] (-2,-1.5) rectangle (2,1.5);
    \draw[thick] (0,1.5)--(0,-.63) (0,-.97)--(0,-1.5); 
    \draw[thick] (-.2,.8)--(-.8,.8) arc(90:270:.8)--(0,-.8);
    \draw[thick] (.2,.8)--(.8,.8) arc(90:-90:.8)--(0,-.8);
    \draw[thick] (0,-.8) circle (.17);
    \node at (-.6,-.5) {$m_i$};
    \node at (-1.5,1) {$\beta_i$};
}.
$$
(iv) follows from 
$$
\tikzmatht{
  \fill[black!10] (-1.5,-1.7) rectangle (1.5,1.8);
    \draw[thick] (-.3,0)--(-.5,0) arc(90:270:.4)--(.5,-.8)
arc(-90:90:.4)--(.3,0);
    \draw[thick] (-.3,1.3)--(-.5,1.3) arc(90:270:.4)--(.5,.5)
arc(-90:90:.4)--(.3,1.3);
    \draw[thick] (0,1.8)--(0,.67) (0,.33)--(0,-.63) (0,-.97)--(0,-1.7);
    \draw[thick] (0,-.8) circle (.17);
    \draw[thick] (0,.5) circle (.17);
}=
\tikzmatht{
  \fill[black!10] (-1.5,-1.7) rectangle (1.5,1.8);
    \draw[thick] (-.2,0.5)--(-.5,.5) arc(90:270:.3)
--(.5,-.1)
arc(-90:90:.3)--(.2,0.5);
    \draw[thick] (-.2,.9)--(-.5,.9) arc(90:270:.7)--(.5,-.5)
arc(-90:90:.7)--(.2,.9);
    \draw[thick] (0,1.8)--(0,.07) (0,-.27)--(0,-.33) (0,-.67)--(0,-1.7);
    \draw[thick] (0,-.1) circle (.17);
    \draw[thick] (0,-.5) circle (.17);
    \node at (.5,-1.3){$\theta^C$};
    \node at (.5,1.4){$\theta^A$};
}
$$ 
in combination with the property \eref{bimprod} in \lref{l:bimprod}.
(v) follows from the unit property and \eref{Tr=dim}. 
\qed

In particular, taking $\BA\equiv(\theta^\BA,x^{\BA(2)})$ as the trivial
$\BA$-$\BA$-bimodule, one has 

\begin{corollary} \label{c:Dm}
{\rm (i)} $D_\BA=d_\BA\cdot P_\BA^l$, hence (by \lref{l:Dm}(iv)) \\
{\rm (ii)} $P_\BA^l\cdot D_\mm \cdot P_\BB^l= D_\mm$. 
\end{corollary}

We also define more generally, for any $\rho\in\C$, 
$$D_\mm(\rho):= 
\tikzmatht{
  \fill[black!10] (-2,-1.5) rectangle (2,2);
    \draw[thick] (-.5,-1.5)--(-.5,-.57) (-.5,-.23)--(-.5,.8) (-.5,1.2)--(-.5,2);
    \draw[thick] (.5,-1.5)--(.5,.8) (.5,1.2)--(.5,2);
    \draw[thick] (.7,1) arc(90:-90:.7)--(.67,-.4) (.33,-.4)--
          (-.67,-.4) arc(270:90:.7)--(.7,1); 
    \draw[thick] (-.5,-.4) circle (.17);
    \node at (-1,-.8) {$m$};
    \node at (.8,1.6) {$\rho$}; 
    \node at (-1.3,1.2) {$\beta$};
} \in \Hom(\theta^\BB\circ\rho,\theta^\BA\circ\rho).$$
Clearly, $D_\mm\equiv D_\mm(\id)$. The properties (i)--(iv) hold as well for 
$D_\mm(\rho)$. 

Next, consider $\BA\otimes \eins\equiv (\theta\otimes\id,x\otimes 1_\id)$ as
a Q-system in $\C\boxtimes \C\opp$, and $\mm\otimes \eins\equiv 
(\beta\otimes\id,m\otimes 1_\id)$ as an $\BA\otimes \eins$-$\BB\otimes \eins$-bimodule. 
Taking the product 
$$R[\BA]:= (\BA\otimes \eins)\times^+ \BR$$
where $\BR=(\Theta\can,W\can,X\can)$ is the canonical Q-system in
$\C\boxtimes\C\opp$, we get
$$R[\mm] = ((\beta\otimes\id)\circ\Theta\can,R[\mm]), \qquad R[\mm]=
\tikzmatht{
  \fill[black!10] (-3.7,-1.1) rectangle (3,2);
    \draw[thick] (-1.3,2)--(-1.3,.7) arc(180:260:1) (.7,2)--(.7,.7)
arc(360:280:1);
    \draw[thick] (-.3,-.3) circle(.17);
    \draw[thick] (-.3,-1.1)--(-.3,2);
    \draw[thick] (.3,-1.1)--(.3,-.3) (.3,.2)--(.3,2) (.3,.7) arc(270:285:1) 
    (.3,.7) arc(270:250:1) (-.7,2)--(-.7,1.7) arc(180:220:1) (1.3,2)--(1.3,1.7) arc(360:310:1);
  \fill[black] (.3,.7) circle (.12);
    \node at (.6,-.7) {$R$}; 
    \node at (1.6,1.6) {$R$};
    \node at (-1.5,-.7) {$\beta\otimes\id$}; 
    \node at (-2.5,.4) {$A\otimes\id$};
    \node at (1.8,.4) {$B\otimes\id$};
}
$$
as an $R[\BA]$-$R[\BB]$-bimodule. Because $\BR$ is commutative, one
has $D_\BR=d_\BR\cdot 1_{\Theta\can}$ by \cref{c:Dm}, and hence 
\be\label{DRm}
D_{R[\mm]} = 
\tikzmatht{
  \fill[black!10] (-3,-1.8) rectangle (2.5,1.8);
    \draw[thick] (-.2,1) arc(90:270:1)--(1.2,-1) arc(-90:90:1) (.8,1)--(.2,1);
    \draw[thick] (0,-1.8)--(0,-1.17) (0,-.83)--(0,1.8);
    \draw[thick] (1,-1.8)--(1,-1.2) (1,-.8)--(1,1.8);
    \draw[thick] (1,-.8) (.2,-.7) (-.2,-.7) arc(270:180:.7)--(-.9,.4);
    \draw[thick] (-.2,1.3) arc(90:150:.7) (.2,1.3)--(.8,1.3) (1.2,1.3)
    arc(90:30:.7); 
    \draw[thick] (.2,-.7)--(1.2,-.7) arc(-90:0:.7)--(1.9,.4); 
    \draw[thick] (0,-1) circle(.17);
  \fill[black] (1,-.7) circle(.12);
    \node at (-1.3,-1.5) {$A\otimes \id$};
    \node at (-2,.7) {$\beta\otimes \id$};
    \node at (1.4,-.1) {$R$};
} = 
\tikzmatht{
  \fill[black!10] (-1.5,-1.8) rectangle (2.5,1.8);
    \draw[thick] (-.2,1) arc(90:270:1)--(1.2,-1) arc(-90:90:1) (.8,1)--(.2,1);
    \draw[thick] (0,-1.8)--(0,-1.17) (0,-.83)--(0,1.8);
    \draw[thick] (1,-1.8)--(1,-1.2) (1,-.8)--(1,1.8);
    \draw[thick] (1,-.6) arc(270:110:.6) (1,-.6) arc(-90:70:.6);
    \draw[thick] (0,-1) circle(.17);
  \fill[black] (1,-.6) circle(.12);
} = d_\BR\cdot
\tikzmatht{
  \fill[black!10] (-1.5,-1.8) rectangle (2.5,1.8);
    \draw[thick] (-.2,1) arc(90:270:1)--(1.2,-1) arc(-90:90:1) (.8,1)--(.2,1);
    \draw[thick] (0,-1.8)--(0,-1.17) (0,-.83)--(0,1.8);
    \draw[thick] (1,-1.8)--(1,-1.2) (1,-.8)--(1,1.8);
    \draw[thick] (0,-1) circle(.17);
    \node at (1.4,-.1) {$R$};
}
= d_\BR\cdot D_{\mm\otimes \eins}(\Theta\can).
\ee

As the full centre $Z[\BA]=(\Theta^\BA,W^\BA,X^\BA)$ is an irreducible
(by \pref{p:fullsimple}) intermediate Q-system to $R[\BA]$, the
$R[\BA]$-$R[\BB]$-bimodule $R[\mm]$ restricts to a $Z[\BA]$-$Z[\BB]$-bimodule
according to \lref{l:modrest}
\bea \notag 
R[\mm]\vert_Z =
\sqrt{\frac{d_{R[\BA]}d_{R[\BB]}}{d_{Z[\BA]}d_{Z[\BB]}}}\cdot
(S^{\BA*}\times 1_{(\beta\otimes\id)\times\Theta\can}\times S^{\BB})\scirc R[\mm]
\equiv \\ 
\equiv \sqrt{\frac{d_{R[\BA]}d_{R[\BB]}}{d_{Z[\BA]}d_{Z[\BB]}}} \cdot
\tikzmatht{
  \fill[black!10] (-3.7,-1) rectangle (3,2.6);
    \draw[thick] (-1.3,1.7)--(-1.3,.7) arc(180:260:1) (.7,1.7)--
    (.7,.7) arc(360:280:1);
    \draw[thick] (-.3,-.3) circle(.17);
    \draw[thick] (-.3,-1)--(-.3,2.6);
    \draw[thick] (.3,-1)--(.3,-.3) (.3,.2)--(.3,2.6) (.3,.7)
    arc(270:285:1) (.3,.7) arc(270:250:1) (-.7,1.7) arc(180:220:1)
    (1.3,1.7) arc(360:310:1); 
  \fill[black] (.3,.7) circle (.12);
    \node at (.6,-.7) {$R$}; 
    \node at (-1.5,-.7) {$\beta\otimes\id$}; 
    \node at (-2.5,.4) {$A\otimes\id$}; 
    \node at (-2.1,2.2) {$Z[A]$};
    \node at (1.8,.4) {$B\otimes\id$}; 
    \node at (2.1,2.2) {$Z[B]$};
  \fill[white] (-1,2.1)--(-1.3,1.7)--(-.7,1.7)--(-1,2.1);
    \draw[thick] (-1,2.1)--(-1.3,1.7)--(-.7,1.7)--(-1,2.1);
  \fill[white] (1,2.1)--(1.3,1.7)--(.7,1.7)--(1,2.1);
    \draw[thick] (1,2.1)--(1.3,1.7)--(.7,1.7)--(1,2.1);
    \draw[thick] (-1,2.1)--(-1,2.6) (1,2.1)--(1,2.6);
},
\eea
where $S^\BA\in\Hom(Z[\BA],R[\BA])$, $S^\BB\in\Hom(Z[\BB],R[\BB])$ are isometric
intertwiners such that $S^\BA S^{\BA*}=P^l_{R[\BA]}$, $S^\BB
S^{\BB*}=P^l_{R[\BB]}$, cf.\ \lref{l:modrest}. 

Next, we consider the intertwiners 
$$D_{R[\mm]\vert_Z}\in \Hom(\Theta^\BB,\Theta^\BA).$$
In particular, for the trivial $\BA$-$\BA$-bimodule $\BA$,
one has, using \cref{c:Dm} 
\be\label{RAZ}
D_{R[\BA]\vert_Z} = \frac{d_{R[\BA]}}{d_{Z[\BA]}} \cdot S^*D_{R[\BA]}S =
\frac{d_{R[\BA]}^2}{d_{Z[\BA]}} \cdot S^*P^l_{R[\BA]}S = \frac{d_\BA^2
  d_\BR^2}{d_{Z[\BA]}}\cdot 1_{\Theta^\BA}.
\ee

We introduce the positive-definite inner product on
$\Hom(\Theta^\BB,\Theta^\BA)$ w.r.t.\ the trace: 
\be \label{inner}
(T_2,T_1):= \Tr_{\Theta^{\BB}}
(T_1^*T_2) = 
\tikzmatht{
  \fill[black!10] (-1.2,-1.8) rectangle (1.4,1.8);
    \draw[thick] (.5,1.1) arc(0:180:.6)--(-.7,-1.1) arc(180:360:.6); 
    \draw[thick] (.5,.2)--(.5,-.2);
  \fill[white] (1,1.1) rectangle (0,.2);
    \draw[thick] (1,1.1) rectangle node{$T_1^*$} (0,.2);
  \fill[white] (1,-1.1) rectangle (0,.-.2);
    \draw[thick] (1,-1.1) rectangle node{$T_2$} (0,-.2);
} \equiv \Tr_{\Theta^{\BA}} (T_2T_1^*) = 
\tikzmatht{
  \fill[black!10] (-1.2,-1.8) rectangle (1.4,1.8);
    \draw[thick] (.5,1.1) arc(0:180:.6)--(-.7,-1.1) arc(180:360:.6); 
    \draw[thick] (.5,.2)--(.5,-.2);
  \fill[white] (1,1.1) rectangle (0,.2);
    \draw[thick] (1,1.1) rectangle node{$T_2$} (0,.2);
  \fill[white] (1,-1.1) rectangle (0,.-.2);
    \draw[thick] (1,-1.1) rectangle node{$T_1^*$} (0,-.2);
}.
\ee
Then we compute
\bea (D_{R[\mm']\vert_Z},D_{R[\mm]\vert_Z})\hspace{-32mm} && =\quad 
\frac{d_{R[\BA]}d_{R[\BB]}}{d_{Z[\BA]}d_{Z[\BB]}} \cdot 
\Tr_{\Theta^\BA}(S^*D_{R[\mm']}SS^*D_{R[\ol \mm]}S)
\\ \nonumber \hspace{-10mm}&\stapel{\cref{c:Dm}, \nref{DRm}}=&
\frac{d_{R[\BA]}d_{R[\BB]}}{d_{Z[\BA]}d_{Z[\BB]}} \cdot d_\BR^2\cdot  \Tr_{(\theta^\BA\otimes \id)\Theta\can} 
(D_{\mm'\otimes \eins}(\Theta\can)D_{\ol \mm\otimes \eins}(\Theta\can) )
 \\\nonumber \hspace{-10mm}& \stapel{\lref{l:Dm}(iii),(iv)}=&
 \frac{d_{R[\BA]}d_{R[\BB]}}{d_{Z[\BA]}d_{Z[\BB]}} \cdot d_\BR^2 
\cdot d_\BB \sum_\mn
N_{\mm'\ol {\mm}}^\mn\cdot \LTr_{\theta^\BA\otimes
  \id}\RTr_{\Theta\can}(D_{\mn\otimes \eins}(\Theta\can)), 
\eea
where $\mm'\otimes_\BB \ol \mm \simeq \bigoplus_\mn N_{\mm'\ol
  \mm}^\mn\cdot \mn$ as
$\BA$-$\BA$-bimodules, $\mn = (\alpha,n)$ irreducible.  

At this point, modularity comes to bear through \pref{p:killing}:
namely $\RTr_{\Theta\can}$ projects on the contribution
$\id\prec\alpha$:
$$\RTr_{\Theta\can}(D_{\mn\otimes \eins}(\Theta\can)) = 
\tikzmatht{
  \fill[black!10] (-2.9,-1.5) rectangle (3.7,1.5);
    \draw[thick] (0,-.5) circle (.17);
    \draw[thick] (0,-1.5)--(0,1.5);
    \draw[thick] (-.17,-.5) arc(270:90:.5)  (.17,-.5) arc(270:395:.5)
    (.17,.5) arc(90:70:.5); 
    \draw[thick] (1.9,0) arc(0:190:.8) (1.9,0) arc(360:230:.8); 
    \node at (-1.8,0) {$\alpha\otimes \id$};
    \node at (2.8,0) {$\Theta\can$}; 
    \node at (-1.2,-1.1) {$\theta^A\!\otimes\id$};
    \node at (-1.2,1.1) {$\theta^A\!\otimes\id$};
} = d_{\BR}^2\cdot 
\tikzmatht{
  \fill[black!10] (-1.5,-1.5) rectangle (1.5,1.5);
    \draw[thick] (0,-.5) circle (.17);
    \draw[thick] (0,-1.5)--(0,1.5);
    \draw[thick] (-.17,-.5) arc(270:180:.5)  (.17,-.5) arc(270:360:.5);
    \node at (1,-.3) {$\alpha$};
    \node at (-.5,-1.1) {$\theta^A$};
    \node at (-.5,1.1) {$\theta^A$};
  \fill[white] (.67,.3)--(.87,0)--(.47,0)--(.67,.3);
    \draw[thick] (.67,.3)--(.87,0)--(.47,0)--(.67,.3);
  \fill[white] (-.67,.3)--(-.87,0)--(-.47,0)--(-.67,.3);
    \draw[thick] (-.67,.3)--(-.87,0)--(-.47,0)--(-.67,.3);
    \node at (1.1,.6) {$s^*$};
}\otimes\id,$$
where $s$ is an isometry such that
$ss^*=E_\id\in\Hom(\alpha,\alpha)$. 
After taking also $\LTr_{\theta^\BA\otimes
  \id}$, one obtains 
$$\LTr_{\theta^\BA\otimes
  \id}\RTr_{\Theta\can}(D_{\mn\otimes \eins}(\Theta\can))=d_\BR^2\cdot 
\tikzmatht{
  \fill[black!10] (-1.5,-1.5) rectangle (1.5,1.5);
    \draw[thick] (0,.4) circle(.8);
    \draw[thick] (0,-.4) circle(.17);
    \draw[thick] (0,-.57)--(0,-.9)(0,-.23)--(0,.1);
  \fill[white] (0,.4)--(-.2,.1)--(.2,.1)--(0,.4);
    \draw[thick] (0,.4)--(-.2,.1)--(.2,.1)--(0,.4);
  \fill[white] (0,-1.2)--(-.2,-.9)--(.2,-.9)--(0,-1.2);
    \draw[thick] (0,-1.2)--(-.2,-.9)--(.2,-.9)--(0,-1.2);
    \node at (1,-.3) {$\theta^A$};
}$$
Now, \lref{l:moduleint} applies, and accordingly the value vanishes
unless $\mn$ is the trivial $\BA$-$\BA$-bimodule $\mn=\BA$. In this
case $s=w/\sqrt{d_\BA}$, hence the value of the diagram is 
$\dim(\theta)/d_\BA=d_\BA$. Since $N_{\mm'\ol
  \mm}^\BA=\delta_{\mm\mm'}$, we arrive at the orthogonality relation 

\begin{corollary} \label{c:ortho}
In a modular category, 
$$(D_{R[\mm']\vert_Z},D_{R[\mm]\vert_Z})=
\frac{d_\BA^2d_\BB^2d_\BR^6}{d_{Z[\BA]}d_{Z[\BB]}}\cdot \delta_{\mm\mm'}.$$
\end{corollary}

For $\mm=\BA$, $D_{R[\BA]\vert_Z}=D_{R[\BA]\vert_Z}^*$, one can compute 
$(D_{R[\BA]\vert_Z},D_{R[\BA]\vert_Z})$ in two diff\-erent ways: By
\eref{RAZ}, the result is $d_\BA^4 d_\BR^4$. By \cref{c:ortho}, it is 
$\frac{d_\BA^2d_\BB^2d_\BR^6}{d_{Z[\BA]}d_{Z[\BB]}}$. By comparison, 
$d_{Z[\BA]}=d_\BA$ for all simple Q-systems $\BA$. This proves 
\pref{p:dimcent}. Moreover, the coefficient in \cref{c:ortho} equals 
$d_\BA^2d_\BB^2d_\BR^4$. 

By \cref{c:ortho}, the intertwiners $D_{R[\mm]\vert_Z}$ are linearly
independent. It is also known that the number of inequivalent
irreducible $\BA$-$\BB$-bimodules equals
$\dim\Hom(\Theta^\BA,\Theta^\BB)$ \cite{EP03}, hence the intertwiners
$I_\mm$ span $\Hom(\Theta^\BA,\Theta^\BB)$. Since both $(d_\BA d_\BB d_\BR^2)\inv
\cdot D_{R[\mm]\vert_Z}$ and $T:=(\dim(\sig)\dim(\tau))^{-\frac12}
\cdot T^\BA_{\sigma\otimes\tau}T^{\BB*}_{\sigma\otimes\tau}$ (where $T^\BA$ and
$T^\BB$ are isometric bases of $\Hom(\sigma\otimes\tau,\Theta^\BA)$ resp.\
$\Hom(\sigma\otimes\tau,\Theta^\BB)$ for all irreducible common
sub-endomorphism $\sigma\otimes\tau$ of $\Theta^\BA$ and $\Theta^\BB$) form orthonormal bases w.r.t.\
trace, the matrix 
$$S_{\mm T}:= \frac1{d_\BA d_\BB d_\BR^2\sqrt{\dim(\sig\otimes\tau)}} \cdot
\Tr_{\Theta^\BA}(D_{R[\mm]\vert_Z}T^\BB_{\sigma\otimes\tau}T^{\BA*}_{\sigma\otimes\tau})
$$
is unitary. In particular, for $\sig=\tau=\id$, $T=T_0\equiv d_\BR\inv\cdot
W^\BA W^{\BB*}$, one finds
$$S_{\mm0} = \frac{\dim(\beta)}{d_\BA d_\BB},$$
hence 
\be\label{W0} 
W^\BA W^{\BB*} = d_\BR\cdot T_0 = d_\BR \sum_\mm \frac{S_{\mm0}}{d_\BA d_\BB
d_\BR^2} \cdot
D_{R[\mm]\vert_Z} = \sum_\mm \frac{\dim(\beta)}{d_\BA^2d_\BB^2d_\BR^2}\cdot
D_{R[\mm]\vert_Z}.
\ee

Now, we define 
\be \label{Im} 
I_\mm:= \frac{\dim(\beta)}{d_\BA^2d_\BB^2d_\BR^2}\cdot
D_{R[\mm]\vert_Z}\in\Hom(\Theta^\BB,\Theta^\BA).
\ee
From the definition and properties (i) and (ii) in \lref{l:Dm}, one
can see that $D_{R[\mm]\vert_Z}$ and hence $I_\mm$ satisfy the
selfadjointness condition in \eref{projD}. Because $W^\BA W^{\BB*}$ is
the unit w.r.t.\ the convolution product \eref{convol}, \eref{W0} is the
completeness relation, i.e., $\sum_\mm E_\mm =\eins_M$ under the
isomorphism $\chi$. It remains to prove the
idempotency relation in \eref{projD}. Using \eref{W0}, it
suffices to show that  
$$I_\mm\ast I_{\mm'}=0$$
for $\mm\neq \mm'$, in order to conclude that $I_\mm\ast I_\mm =
I_\mm\ast \sum_{\mm'}I_{\mm'}=I_\mm\ast (W^\BA W^{\BB*})=I_\mm$. 

Let $\mm$ and $\mm'$ be two $\BA$-$\BB$-bimodules. Define 
$$Q_{\mm,\mm'}:=
\tikzmatht{
  \fill[black!10] (-2,-1.5) rectangle (2,2);
    \draw[thick] (-.5,-1.5)--(-.5,.8) (-.5,1.2)--(-.5,2);
    \draw[thick] (.5,-1.5)--(.5,.8) (.5,1.2)--(.5,2);
    \draw[thick] (.7,1) arc(90:-90:.7)--(.67,-.4) (.33,-.4)--
    (-.33,-.4) (-.67,-.4) arc(270:90:.7)--(.7,1); 
    \draw[thick] (-.5,-.4) circle (.17);
    \draw[thick] (.5,-.4) circle (.17);
    \node at (-1,-.8) {$m$}; 
    \node at (1.1,-.7) {$m'$};
    \node at (1.6,1.3) {$\theta^\BA$}; 
}
\in\Hom(\beta\ol\beta',\beta\ol\beta').$$
By a similar computation as for the projection property of the
left and right centre, \lref{l:bimproj}, one sees that 
$(d_\BA d_\BB)\inv\cdot Q_{\mm,\mm'}$ is a projection. Now,
$$\LTr_{\beta}\RTr_{\ol\beta'}(Q_{\mm,\mm'}) =
\Tr_{\theta^\BA}(D_\mm D_{\mm'}^*).$$
Thus, replacing $\BA$ and $\BB$ by $Z[\BA]$ and $Z[\BB]$, $\mm$ and $\mm'$ by 
$R[\mm]\vert_Z$ and $R[\mm']\vert_Z$, and $\beta$ and $\beta'$ by
$R[\beta]=(\id\otimes\beta)\Theta\can$ and $R[\beta']$, we conclude that
$$\LTr_{R[\beta]}\RTr_{R[\ol\beta']}(Q_{R[\mm]\vert_Z,R[\mm']\vert_Z})=0$$
for $\mm\neq \mm'$ by the orthogonality of $D_{R[\mm]\vert_Z}$,  
\cref{c:ortho}. Since $Q_{R[\mm]\vert_Z,R[\mm']\vert_Z}$ is a multiple of
a projection, hence a positive operator, and because the
traces are faithful positive maps, it follows that 
$Q_{R[\mm]\vert_Z,R[\mm']\vert_Z}=0$ for $\mm\neq \mm'$.  

Now in order to conclude that $X^{\BA*}(I_\mm\times I_{\mm'})X^\BB=0$ for
$\mm\neq \mm'$, it suffices to compute 
$$\tikzmatht{
  \fill[black!10] (-2,-2) rectangle (2,2);
    \draw[thick] (0,2)--(0,1.2) arc (90:197:1.2);
    \draw[thick] (0,1.2) arc (90:-17:1.2);
    \draw[thick] (0,-2)--(0,-1.2) arc (270:215:1.2);
    \draw[thick] (0,-1.2) arc (270:325:1.2);
    \draw[thick] (1.4,.4) arc(60:-240:.5); 
    \draw[thick] (-.9,.4) arc(60:-240:.5); 
  \fill[black] (0,1.2) circle(.12);
  \fill[black] (0,-1.2) circle(.12);
    \draw[thick] (1.1,-.52) circle (.17);
    \draw[thick] (-1.1,-.52) circle (.17);
    \node at (-1.5,-.9) {$m$}; 
    \node at (1.6,-.8) {$m'$};
    \node at (-.6,1.6) {$A$}; 
    \node at (-.6,-1.6) {$B$}; 
} \stapel{\nref{bimrep}}=
\tikzmatht{
  \fill[black!10] (-2,-2) rectangle (2,2);
    \draw[thick] (-1.2,2)--(-1.2,-.38) (-1.2,-.72)--(-1.2,-2); 
    \draw[thick] (-.2,-.36)--(-.2,.5) arc(180:0:.7)--(1.2,-.36)
    (1.2,-.70) arc(0:-180:.7); 
    \draw[thick] (1.4,.4) arc(60:-240:.5); 
    \draw[thick] (0,.4) arc(60:-90:.5)--(-1.2,-.55) arc(270:120:.5)
    (-.9,.45)--(-.5,.45); 
    \draw[thick] (-1.2,-.55) circle (.17);
    \draw[thick] (-.2,-.53) circle (.17);
    \draw[thick] (1.2,-.53) circle (.17);
    \node at (1,1.6) {$A$}; 
    \node at (1,-1.6) {$B$}; 
} =
\tikzmatht{
  \fill[black!10] (-2,-2) rectangle (2,2);
    \draw[thick] (-1.2,2)--(-1.2,-1.08) (-1.2,-1.42)--(-1.2,-2); 
    \draw[thick] (.8,.6) arc(90:-75:.6) (.63,-.6)--(.17,-.6);
    \draw[thick] (.8,.6)--(0,.6) arc(90:255:.6);
    \draw[thick] (-1.2,-1.25) arc(270:180:.5)--(-1.7,.6)
    arc(180:110:.5) (-.9,1.1)--(-.5,1.1) arc(90:25:.5)
    (0,.4)--(0,-.75) arc (0:-90:.5)--(-1.2,-1.25);  
    \draw[thick] (.8,-.75)--(.8,.4) (.8,.8) arc(170:0:.5)--(1.8,-.72)
    arc(0:-180:.5); 
    \draw[thick] (-1.2,-1.25) circle (.17);
    \draw[thick] (0,-.59) circle (.17);
    \draw[thick] (.8,-.59) circle (.17);          
    \draw[thin,dashed] (-.9,-.9)--(-.9,.9)--(1.6,.9)--(1.6,-.9)--(-.9,-.9);
}.
$$
Inside the dashed box, there appears the intertwiner $Q_{R[\mm]\vert_Z,R[\mm']\vert_Z}$,
which we have just shown to be zero if $\mm\neq \mm'$. In step (r), the
representation property of $\mm$ as a left $\BA$-module has been used. This
concludes the proof that $I_\mm$ solve \eref{projD}. 

\tref{t:centchar} now follows from the considerations before \eref{projD}.
\qed 

The minimal projections $E_\mm\in M'\cap M$ define representations
$m\mapsto E_\mm m$ as in \sref{s:Qcentral}. In these representations,
the generators $V^\BA$ and $V^\BB$ of the intermediate algebras
$M^\BA$ and $M^\BB$ (defined as in \lref{l:braidext}) are no longer 
independent. Let us describe the nature of these relations. 

\begin{lemma} \label{l:chi}
The bijection $\chi$, \eref{chi}, can be equivalently written as 
$$\chi(T)=V^{BA*}\iota(T)V^\BB.$$
Therefore, in particular, $E_\mm =V^{\BA*}\iota(I_\mm)V^\BB$.
\end{lemma}

{\em Proof:} By \lref{l:braidext} and $V^*=\iota(R^*)V$, we have   
$$V^{BA*}\iota(T)V^\BB =
V^*\iota\Theta^\BA(W^\BB)\iota(T)\iota(W^{\BA*})V = 
\iota\Big(R^*\Theta(\Theta^{\BA}(W^\BB)TW^{\BA*})X\Big)V,
$$
and the argument of $\iota$ equals 
$$\tikzmatht{
  \fill[black!10] (-2,-1.5) rectangle (2,2);
    \draw[thick] (0,0)--(0,-.2) arc(360:180:.8)--(-1.6,.6)
    arc(180:0:1.1)--(.6,-.2) arc(360:270:.8); 
    \draw[thick] (1.2,.4)--(1.2,.6) arc(0:90:1.1); 
    \draw[thick] (-1,-.2)--(-1,.6) arc(180:130:1.1); 
    \draw[thick] (-1,-.2) arc(180:230:.8); 
    \draw[thick] (-.8,-1.5)--(-.8,-1);
    \draw[thick] (-.2,-1.5)--(-.2,-1);
  \fill[black] (-.8,-1) circle(.12);
  \fill[black] (-.2,-1) circle(.12);
  \fill[white] (1.2,.4) circle(.2);
    \draw[thick] (1.2,.4) circle(.2);
  \fill[white] (0,0) circle(.2);
    \draw[thick] (0,0) circle(.2);
  \fill[white] (.4,0) rectangle (.8,.4);
    \draw[thick] (.4,0) rectangle (.8,.4);
    \node at (-1.3,-1.2) {$A$};
    \node at (.3,-1.2) {$B$};
}
=
\tikzmatht{
  \fill[black!10] (-1.2,-1.5) rectangle (1.2,2);
    \draw[thick] (-.8,-1.5)--(-.8,.6) arc(180:0:.8)--(.8,-1.5);
  \fill[white] (.6,0) rectangle (1,.4);
    \draw[thick] (.6,0) rectangle (1,.4);
}.
$$
This coincides with \eref{chi}. The last statement follows from
$E_\mm=\chi(I_\mm)$. 
\qed

Expanding a general element $T\in\Hom(\Theta^\BB,\Theta^\BA)$ in the
basis $I_\mm$, such that $T=\sum_\mm c_\mm(T)\cdot I_\mm$, we get
$$V^{\BA*}\iota(T)V^\BB = \sum_\mm c_\mm(T) \cdot E_\mm,$$
i.e., in the representation defined by each $E_\mm$, the central
elements $V^{\BA*}\iota(T)V^\BB$ take the values $c_\mm(T)$. In
particular, for $\sigma\otimes\tau$ an irreducible common
sub-endomorphism of $\Theta^\BA$ and $\Theta^\BB$, and
$T=(\dim(\sig)\dim(\tau))^{-\frac12} \cdot
T^\BA_{\sigma\otimes\tau}T^{\BB*}_{\sigma\otimes\tau}$ as above, 
these values are 
$$c_\mm(T) = \frac{d_\BA d_\BB}{\dim(\beta)} \cdot \ol{S_{\mm T}}.$$
Since on the other hand, the charged intertwiners
$\Psi^\BA_{\sigma\otimes\tau} = \iota(T^{\BA*}_{\sigma\otimes\tau})
V^\BA\in\Hom(\iota^\BA,\iota^\BA\circ (\sigma\otimes\tau))$ and  
$\Psi^\BB_{\sigma\otimes\tau} = \iota(T^{\BB*}_{\sigma\otimes\tau}) V^\BB$ are
multiples of isometries because $\iota^\BA$ and $\iota^\BB$ are
irreducible, the numerical values for
$\Psi^{\BA*}_{\sigma\otimes\tau}\Psi^\BB_{\sigma\otimes\tau}$ define 
``angles'' between them \cite{BKLR}. 

\begin{graybox}
\begin{example} \vskip-5mm \label{x:cardy}
Let $\BA$ and $\BB$ be the trivial Q-system (or Morita equivalent),
such that the full centres coincide with $\BR$. The irreducible
bimodules of the trivial Q-system are just the irreducible
endomorphisms $\sig\in\C$, $\mm = (\sig,1_\sig)$. The irreducible
sub-endomorphism of $\Theta^\BA=\Theta^\BB=\Theta\can$ are
$\tau\otimes\ol\tau$. The operators $I_\mm$, \eref{Im}, simplify to 
$\frac{\dim(\sig)}{d_\BR}\cdot 
\tikzmatht{
  \fill[black!3] (-2.2,-1) rectangle (1.5,1);
    \draw[thick] (0,-.4)--(-.2,-.4) arc(270:90:.4); 
    \draw[thick] (0,-.4)--(.2,-.4) arc(270:450:.4); 
    \draw[thick] (0,-1)--(0,-.6) (0,-.2)--(0,1);
    \node at (-1.2,-.6) {$\sig\otimes\id$};
    \node at (.7,.7) {$\Theta\can$}; 
} \sim \bigoplus_\tau
\tikzmatht{
  \fill[black!3] (-1.2,-1) rectangle (1.2,1);
    \draw[thick] (0,-.4)--(-.2,-.4) arc(270:90:.4); 
    \draw[thick] (0,-.4)--(.2,-.4) arc(270:450:.4); 
    \draw[thick] (0,-1)--(0,-.6) (0,-.2)--(0,1);
    \node at (-.7,-.6) {$\sig$};
    \node at (.3,.8) {$\tau$}; 
} \otimes
\tikzmatht{
  \fill[black!3] (-.6,-1) rectangle (.6,1);
    \draw[thick] (0,-1)--(0,1);
    \node at (.4,.7) {$\ol\tau$}; 
}
$. 
The matrix $(S_{\mm,T})_{\mm,T}$ determining the angles turns out to
coincide with the ``modular'' matrix
$(S_{\sigma,\tau})_{\sigma,\tau}$, cf.\ \dref{d:modular}.  
In particular, if $S_{\sig,\tau}$ happens to equal a
complex phase $\omega$ times $\dim(\sig)\dim(\tau)\cdot\big(\sum_\rho 
  \dim(\rho)^2\big)^{-\frac12}$ (this is always the case whenever
$\sig$ has dimension $\dim(\sig)=1$), it follows that the
generators $\Psi^\BA_{\tau\otimes\ol\tau}=\omega\cdot
\Psi^\BB_{\tau\otimes\ol\tau}$ are linearly dependent in the
representation given by $E_\mm$.   
\end{example}
\end{graybox}

\section{Applications in QFT}
\label{s:QFT}
\setcounter{equation}{0}
We review some applications of the above abstract theory in the
context of local quantum field theory. In a nutshell, Q-systems
provide a complete characterization of (finite index) extensions of
local quantum field theories, and the notions and operations discussed
in the main body of this work have counterparts in conformal QFT that are of
interest for the construction and classification of local extensions 
and of boundary conditions. More details can be found in
\cite{LR95,LR09,CKL13} and \cite{BKL,BKLR}. 

\subsection{Basics of algebraic quantum field theory}
\label{s:AQFT}
\setcounter{equation}{0}

\subsubsection{Local nets}
\label{s:nets}

The additional feature in quantum field theory is the local structure:
quantum fields are operator-valued distributions in spacetime, such
that the support of the test function specifies the localization of
field operators. In the algebraic approach \cite{H} one rather considers 
local algebras $\A(O)$ of bounded operators generated by quantum fields 
evaluated on test functions with a given spacetime support $O$.  

In fact, it is not necessary to assume that the local algebras are
generated by actual quantum fields. It is sufficient to assume that
the net of local algebras is isotonous, i.e., 
$O_1\subset O_2$ implies $\A(O_1)\subset \A(O_2)$. 

Thus, rather than with a single von Neumann algebra, one deals with a
directed net of von Neumann algebras $\A(O)$, where $O$ runs over a
suitable family of connected open regions in spacetime, and $\A(O)$
are the von Neumann algebras generated by local observables localized 
in $O$. (If $O$ has a sufficiently large causal complement, then the algebra
$\A(O)$ does not depend on the representation in which the weak closure
is taken.) The regions $O$ can be chosen to be doublecones (intersections
of a future and a past lightcone) 
$$\tikzmatht{
  \fill[black!10] (-1.5,-1.5) rectangle (1.5,1.5);
  \fill[white] (.2,.6)--(.8,0)--(0,-.8)--(-.6,-.2)--(.2,.6);
    \draw[thick] (.2,.6)--(.8,0)--(0,-.8)--(-.6,-.2)--(.2,.6);
    \draw[thin] (-1.5,-1.1)--(.2,.6)--(1.5,-.7) (-1.5,.7)--(0,-.8)--(1.5,.7);
    \node at (0,1.2) {future}; 
    \node at (0,-1.2) {past}; 
    \node at (.1,-.1) {$O$}; 
}
$$
in Minkowski spacetime $M^D=(\RR^D,\eta)$, or
intervals $I\subset \RR$ where $\RR$ is the ``spacetime'' of a chiral quantum
field theory. 

The use of $\RR$ as the ``spacetime'' of a chiral quantum field
theory is due to the feature of conformal quantum field theory in two
dimensions that it necessarily contains fields (notably the stress-energy
tensor which is the local generator of diffeomorphism covariance) that
depend only on either lightcone coordinate $t+x$ or $t-x$. The net of local
algebras generated by chiral fields is therefore indexed by the
intervals $I\subset \RR$. By virtue of the conformal symmetry, a
conformal net on $\RR$ actually extends to a net (more precisely: a
pre-cosheaf) on $S^1$ by identifying $\RR$ with $S^1$ minus a point;
but this feature will not be essential for the applications that we
are going to review.  

In two spacetime dimensions, the spacelike complement of a doublecone
has two connected components (``wedges''). In chiral theories, the 
spacelike complement of an interval is just its complement in $\RR$,
which is a disconnected union of two halfrays.  

The quasilocal algebra $\A\ql$ associated with a net $\A: O\mapsto \A(O)$ is
the C* algebra defined as the inductive C* limit as $O$ exhausts the
entire spacetime. A group $G$ of spacetime symmetries (Poincar\'e group,
conformal group) is assumed to act on $\A\ql$ as automorphisms
$\alpha_g$ such that $\alpha_g(\A(O))=\A(gO)$.  

The principle of causality (locality) expresses the absence of
superluminal causal influences. In quantum theory, it asserts that
observables localized at spacelike distance must commute with each
other. Thus, if two regions $O_1$, $O_2$ are spacelike separated (in
the chiral case: disjoint), then $[\A(O_1),\A(O_2)]=\{0\}$ (as
subalgebras of $\A\ql$), or equivalently  
$$\A(O)\subset \A(O')',$$ 
where $O'$ is the causal complement of $O$, and $\A(O')$ the C* algebra
generated by $\A(\wh O)$, $\wh O\subset O'$, and $\A(O')'$ its commutant 
in $\A\ql$. 

An overview of the consequences of these axioms (isotony, covariance,
locality, vacuum representation) in chiral conformal QFT can be found in
\cite{GL96} and in \cite{CKL08}.

The enormous benefit of the approach lies in the fact that, once the
validity of the formalism is established, one does not need any
dynamical details of the quantum field theory at hand, except the
knowledge of its representation category as a braided C* tensor
category. In turn, as \xref{x:Icat} and \xref{x:Ibraid} show, this
information typically requires very few data (like the fusion rules
and the twist parameters $\kappa_\rho$) which in many cases uniquely
fix the category.

There is a variety of methods to construct local conformal nets. Free
field nets can be constructed as CAR or CCR algebras, equipped with a
vacuum state. Local nets associated with affine Kac-Moody algebras 
can be constructed from unitary implementers of local gauge
transformations acting as automorphisms of the CAR algebra, giving
rise to projective representations of loop groups \cite{PS}. Local nets
associated with a chiral stress-energy tensor can similarly be
obtained from unitary implementers of local diffeomorphisms acting as
automorphisms of the CAR algebra; an alternative, more explicit
construction from a given heighest-weight representation of the
Virasoro algebra is given in \cite{CKL08}. By orbifold (fixed point)
and coset (relative commutants) constructions, one can construct new
nets from given ones. Finally, the extension of local nets by
commutative Q-systems will be described in \sref{s:qftext}.

\subsubsection{Representations and DHR endomorphisms}
\label{s:DHRendo}

Nets of local algebras possess inequivalent Hilbert space representations. 
A representation $\pi$ is covariant if the automorphisms $\alpha_g$ 
are implemented by a unitary representation $U_\pi(g)$ on the
representation Hilbert space, $\pi(\alpha_g(a))=U_\pi(g)\pi(a)U_\pi(g)^*$.  
A representation $\pi$ is said to have positive energy if the
generator of the unitary one-parameter group $U_\pi(t)$ corresponding to 
the subgroup of time translations has positive spectrum. We assume
that there is a unique vacuum representation $\pi_0$, i.e., a faithful
positive-energy representation with an invariant ground state $\Omega$,
$U(g)\Omega=\Omega$, and we assume that in the vacuum representation a
stronger version of locality holds, namely Haag duality: 
$$\pi_0(\A(O))=\pi_0(\A(O'))'.$$

Under these standard assumptions, one can show that the local algebras
$A(O)$ are infinite factors. Moreover \cite{DHR}, an important class of
positive-energy representations (in two-dimensional conformal QFT: all
positive-energy representations) can be described in terms of
DHR endomorphisms $\rho$ of the quasilocal algebra $\A\ql$ such that 
$$\pi = \pi_0\circ\rho.$$ 
DHR endomorphisms are {\em localized} in some region $O$ in the sense that 
the restriction of $\rho$ to the algebra $\A(O')$ of the causal
complement acts like the identity; and {\em transportable} in the sense that
for every other region $\wh O$, there is an endomorphism $\wh\rho$ 
localized in $\wh O$ which is unitarily equivalent, namely, there is a
unitary {\bf charge transporter} $u\in \A\ql$ (actually localized in any
doublecone that contains $O$ and $\wh O$) such that
$\wh\rho=\Ad_u\circ\rho$: 
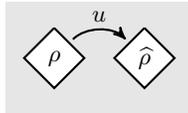
$$\tikzmatht{
  \fill[black!10] (-2.5,-1.5) rectangle (2.5,1.5);
  \fill[white] (-1.2,.8)--(-.4,0)--(-1.2,-.8)--(-2,0)--(-1.2,.8);
    \draw[thick] (-1.2,.8)--(-.4,0)--(-1.2,-.8)--(-2,0)--(-1.2,.8);
  \fill[white] (1.2,.8)--(.4,0)--(1.2,-.8)--(2,0)--(1.2,.8);
    \draw[thick] (1.2,.8)--(.4,0)--(1.2,-.8)--(2,0)--(1.2,.8);
    \node at (1.2,0) (r) {$\wh\rho$}; 
    \node at (-1.2,0) (s) {$\rho$} edge[arrow, bend left=45]
    node[auto] {$u$} (r);   
}.$$

By Haag duality it follows that
$\rho(\A(O))\subset \A(O)$ if $\rho$ is localized in $O$, i.e., $\rho$
restricts to an endomorphism of the von Neumann algebra $N=\A(O)$. 

The composition $\rho_1\circ\rho_2$ of DHR endomorphisms is again a
DHR endomorphism. 
Inertwiners between DHR endomorphisms are defined as operators $t\in
\A\ql$ satisfying $t\rho_1(a)=\rho_2(a)t$ for all $a\in \A\ql$. By Haag
duality, it follows that $t\in \A(O)$ if $\rho_i$ are localized in
$O_i$ and $O_1\cup O_2\subset O$. In particular, all intertwiners
among DHR endomorphisms localized in the same region $O$ are elements
of $\A(O)$.

In this way, picking any fixed region $O$ and putting $N=\A(O)$, 
the restrictions of DHR endomorphisms localized in $O$ form
a C* tensor subcategory of $\End(N)$. One can show \cite[Thm.\ 2.3]{GL96} 
that this subcategory is full, i.e., every intertwiner between
$\rho_1$ and $\rho_2$ regarded as endomorphisms of the von 
Neumann algebra $N$ is also an intertwiner between $\rho_1$ and
$\rho_2$ regarded as endomorphisms of the C* algebra $\A\ql$ (``local
intertwiners = global intertwiners''). In particular, notions like
sector, conjugates and dimension have the same meaning for DHR
endomorphisms as endomorphisms of $\A\ql$ and as endomorphisms of $N$. 

We denote by $\DHR(\A)\vert_O$ the full subcategory of $\End_0(N)$,
whose objects are the DHR endomorphisms of finite dimension, localized
in $O$, and by $\DHR(\A)$ the C* tensor category of all DHR
endomorphisms of finite dimension. \xref{x:Icat} specifies the DHR
category of the chiral Ising model.

\subsubsection{DHR braiding}
\label{s:DHR-braiding}

The C* tensor category $\DHR(\A)$ is equipped with a distinguished
unitary braiding $\eps_{\rho,\sig}\in\Hom(\rho\sig,\sig\rho)$. 
It is defined using unitary charge transporters
$u_\rho\in\Hom(\rho,\wh\rho)$ and $u_\sig\in\Hom(\sig,\wh\sig)$, such
that $\wh\rho$ is localized to the spacelike right (in the chiral case:
in the future) of $\wh\sig$:
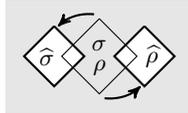
$$\tikzmatht{
  \fill[black!10] (-2.5,-1.5) rectangle (2.5,1.5);
  \fill[white] (-1.2,.8)--(-.4,0)--(-1.2,-.8)--(-2,0)--(-1.2,.8);
    \draw[thick] (-1.2,.8)--(-.4,0)--(-1.2,-.8)--(-2,0)--(-1.2,.8);
  \fill[white] (1.2,.8)--(.4,0)--(1.2,-.8)--(2,0)--(1.2,.8);
    \draw[thick] (1.2,.8)--(.4,0)--(1.2,-.8)--(2,0)--(1.2,.8);
    \draw[thin] (0,1)--(1,0)--(0,-1)--(-1,0)--(0,1);
    \node at (1.4,0) {$\wh\rho$}; \node at (-1.4,0) {$\wh\sig$};
    \node at (0,-.3) {$\rho$}; \node at (0,.3) {$\sig$};
    \draw[thick] (0,1.1) edge[arrow, bend right=25] (-1.2,.7);
    \draw[thick] (0,-1.1) edge[arrow, bend right=25] (1.2,-.7);
}.
$$
One shows with Haag duality that the auxiliary endomorphisms $\wh\rho$
and $\wh\sig$, being localized at spacelike distance, commute with
each other, and defines  
\be\label{DHR-b}
\eps_{\rho,\sig}:= (u_\sig\times u_\rho)^*\scirc (u_\rho\times u_\sig)
\equiv \sig(u_\rho^*)u_\sig^*u_\rho\rho(u_\sig)\in\Hom(\rho\sig,\sig\rho). 
\ee
This unitary does not depend on the choice of $\wh\rho$, $\wh\sig$
with the specified relative localization, nor on the choice of the
charge transporters $u_\rho$, $u_\sig$. It satisfies the defining
properties of a braiding. By construction, if $\rho$ is localized to
the spacelike right (in the chiral case: in the future) of $\sig$, then 
$$\eps_{\rho,\sig}=1,$$
because one may just choose $\wh\rho=\rho$, $\wh\sig=\sig$, and
$u_\rho=u_\sig=1$. In contrast, if $\rho$ is localized to
the spacelike left (past) of $\sig$, one will have
$\eps_{\sig,\rho}=1$ but $\eps_{\rho,\sig}\neq 1$ in general, because
the braiding $\eps^+_{\rho,\sig}\equiv\eps_{\rho,\sig}$ and its
opposite $\eps^-_{\rho,\sig}\equiv\eps_{\sig,\rho}^*$ differ in
low-dimensional QFT, due to the two connected components of the causal
complement. 

In four dimensional QFT, the braiding is degenerate:
$\eps_{\rho,\sig}\eps_{\sig,\rho}=\eins$, i.e., it is a permutation symmetry, and
the twist parameter $\kappa_\rho=\pm1$ distinguishes fermionic and
bosonic sectors \cite{DHR}. In chiral conformal QFT, the conformal 
spin-statistics theorem \cite{GL96} relates the twist parameter
$\kappa_\rho=e^{2\pi i h_\rho}$ of a sector to the lowest eigenvalue
$h_\rho$ of $L_0$. 

If both $\rho$ and $\sig$ are localized in $O$, then
$\eps_{\rho,\sig}\in\A(O)$, hence the DHR braiding restricts to each
$\DHR(\A)\vert_O$. The local structure of a QFT net therefore provides
us intrinsically with a braided C* tensor category, the arena of the
abstract theory of the previous sections. 

Of particular interest in the context of the present work is the case
when the quantum field theory $\A$ possesses only finitely many
irreducible DHR sectors of finite dimension. In chiral conformal QFT,
this property (referred to as ``completely rationality''), is known to
follow from the split property and Haag duality for intervals. Many
models of interest, including the chiral Virasoro models with central
charge $c<1$, are completely rational. The case $c=\frac12$ is the
chiral Ising model, \xref{x:Icat} and \xref{x:Ibraid}. 

(Complete rationality should be regarded, however, rather as a
technically useful regularity condition with far-reaching
consequences, than an axiom based on physical principles -- since
important models, like the $u(1)$ current algebra, do not share this
property.)  

In completely rational chiral models, the DHR braiding is non-degenerate
\cite{KLM}, making the braided category $\DHR(\A)$ a modular category, 
cf.\ \sref{s:modular}. Moreover, the global dimension of $\DHR(\A)$ 
(i.e., the quantity $\sum_{[\rho]\,\mathrm{irr}} \dim(\rho)^2$, 
\eref{globaldim}) coincides with the $\mu$-index of $\A$ which
measures the violation of Haag duality for pairs of disconnected
intervals \cite{KLM,LX}. Thus, the presence of DHR sectors can be
``detected'' by inspection of the two-interval subfactor
$$\pi_0\big(\A(I_1\cup I_2)\big)\subset \pi_0\big(\A((I_1\cup I_2)')\big)'$$ 
where $I_1, I_2$
are any pair of non-touching intervals.   
Recall that the global dimension also is the common dimension of
$\Theta$ in all irreducible full centre Q-systems by \pref{p:dimcent}. 
In particular, the two-interval subfactor is isomorphic with the
subfactor described by the canonical Q-system \cite{LR95}.  

\medskip 

We now turn to the interpretations of Q-systems and the various
operations on them, in the QFT context.

\subsection{Local and nonlocal extensions}
\label{s:qftext}
\setcounter{equation}{0}

\subsubsection{Q-systems for quantum field theories}
\label{s:Qnets}

A Q-system in a C* tensor category $\C\subset \End_0(N)$ describes an
extension $N\subset M$. A Q-system $\BA=(\theta,w,x)$ in the category
$\C=\DHR(\A)$ describes a family of extensions 
$$\A(O)\subset \B(O)$$
in very much the same way. Namely, let $\B$ be the * algebra generated
by $\A\ql$ and $v$ subject to the relations 
$$v\cdot a=\theta(a)\cdot v,\quad v^2=x\cdot v,\quad v^*=w^*x^*\cdot v,$$
such that $\B=\A\ql\cdot v$ as a vector space. Embed $\A\ql$ by
$\iota(a)=aw^*\cdot v$ as a * subalgebra. Define * subalgebras 
$$\B(O):=\A(O)u\cdot v$$
where $u\in\A\ql$ is a unitary such that $\wh\theta=\Ad_u\circ\theta$ is
localized in $O$. Because $(u\times u)\scirc x \scirc
u^*\in\Hom(\wh\theta,\wh\theta^2)\subset \A(O)$ and $(u\times u)\scirc
x\scirc w\in\Hom(\id,\wh\theta^2)\subset \A(O)$, $\B(O)$ are indeed *
algebras. In fact, $\A(O)\subset\B(O)$ is precisely the von Neumann
algebra extension of $\A(O)$ by the Q-system 
$(\wh\theta=\Ad_u\circ\theta,u\scirc w,(u\times u)\scirc x \scirc u^*)$. 

One obtains a net of von Neumann algebras $O\mapsto \B(O)$ extending
$\A(O)$, and $\B$ is its inductive limit as the regions $O$ exhaust
the entire spacetime. 
 
Charged intertwiners $\psi_\rho$, defined for $\rho\prec\theta$ as in
\rref{r:charged}, are elements of $\B(O)$ whenever $\rho$ is localized
in $O$, and these operators together with $\A(O)$ generate $\B(O)$. As
$O$ varies, these operators are the substitute of charged ``fields''
in the language of algebraic QFT. 

The charged intertwiners create charged states from the vacuum as
follows  \cite{LR95}. The positive map 
$$\mu:b\mapsto d_\BA\inv\cdot w^*\ol\iota(b)w$$ 
is a conditional expectation $\mu:\B\to\A$. It allows to extend the
vacuum state $\omega_0$ on $\A$ to a vacuum state
$\omega:=\omega_0\circ\mu$ on $\B$. Since
$\mu(\psi_\rho\psi_\rho^*)\in \Hom(\rho,\rho)$ is a multiple of 
$\eins$ if $\rho$ is irreducible, we may assume it to be $=\eins$ by
normalizing $\psi_\rho$. Then one has 
$$\omega_0\circ\mu(\psi_\rho \iota(a)\psi_\rho^*) =
\omega_0(\rho(a)\mu(\psi_\rho \psi_\rho^*)) = \omega_0\circ\rho(a).$$
Thus, in the GNS representation $\pi$ of the state
$\omega$, the vector $\pi(\psi_\rho^*)\Omega_{\omega}$ belongs to the
DHR representation $\pi_0\circ\rho$ of $\A$. Indeed, upon
restriction to $\A$, the GNS representation of $\omega$ is
equivalent to the DHR representation $\pi_0\circ\theta$ of $\A$.

\medskip

The net $\B$ is by construction relatively local w.r.t.\ the subnet $\A$: 
if $b=auv\in \B(O)$ with $a\in \A(O)$, and $a'\in \A(O')$, then 
$$b\cdot a' = auv\cdot a' =au\theta(a')v=a\wh\theta(a')uv=aa'uv =
a'\cdot auv = a'\cdot b,$$
where we have used the localization of $\wh\theta$ and the local
commutativity of $a$ with $a'$. In fact, every relatively local net of
extensions of finite index arises this way \cite{LR95}.  

An extension $\B$ of $\A$ is in general not local. It is local iff
$u_1v$ commutes with $u_2v$ whenever $\theta_1=\Ad_{u_1}\circ\theta$
and $\theta_2=\Ad_{u_2}\circ\theta$ are localized in spacelike
separated regions $O_1$, $O_2$. But 
$$u_1v \cdot u_2v = u_1\theta(u_2)xv,\qquad u_2v \cdot u_1v = u_2\theta(u_1)xv,$$
which are equal iff $\eps_{\theta,\theta}\scirc x=x$, by the
definition of the DHR braiding \eref{DHR-b}. Thus, $\B$ is
local iff the Q-system $(\theta,w,x)$ is commutative. 

\medskip

Given a local extension $\A\subset \B$, one may apply
$\alpha$-induction (cf.\ \sref{s:alpha}) to the DHR 
endomorphisms $\rho$ of $\A$, defining endomorphisms $\alpha^\pm_\rho$
of $\B$. These are, however, in general not DHR endomorphisms of $\B$,
since they act trivially only on one of the two components of the
causal complement of the localization region of $\rho$. DHR
endomorphisms of $\B$ are obtained as sub-endomorphisms which are
contained in both $\alpha^+_{\rho}$ and $\alpha^-_{\sig}$ for some
$\rho,\sig\in\DHR(\A)$. The common (``ambichiral'') sub-endomorphisms
are counted by the numbers
$Z_{\rho,\sig}=\dim\Hom(\alpha^-_\rho,\alpha^+_\sig)$, cf.\ \eref{Z}. 

\medskip 

By classifying (commutative) Q-systems within the DHR category of a
given completely rational quantum field theory, one obtains a
classification of its (local) extensions. This program has been
completed (profiting from existence and uniqueness results of
\cite{KO} and the previous classifications of modular invariant
matrices in \cite{CIZ}) 
for the local extensions of chiral nets associated with the 
stress-energy tensor with central charge $c<1$, which are known to be
completely rational \cite{KL04}. All models in this classification can
be realized by coset constructions, except one which arises as a mirror
extension (cf.\ \sref{s:mirror}) of a coset extension. The
classification of relatively 
local extensions with $c<1$ (which is of interest in the presence of
boundaries, \sref{s:hardb}) can be found in \cite{KLPR}; and the
classification of local two-dimensional extensions (\sref{s:two-d}) with
$c<1$ was achieved in \cite{KL04-2}. 

\subsubsection{Two-dimensional extensions}
\label{s:two-d}
The chiral observables of a two-dimensional conformal QFT are
given by a tensor product of two chiral nets $ \A_2:=\A_+\otimes \A_-$
such that 
$$\A_2(O)=\A_+(I)\times \A_-(J)$$
if 
$$O=I\times J = \{(t,x):t+x\in I, t-x\in J\}:\qquad 
\tikzmatht{
  \fill[black!10] (-2.5,-1.5) rectangle (1.5,1.5);
  \fill[white] (.7,-.8)--(-.1,0)--(.5,.6)--(1.3,-.2)--(.7,-.8);
    \draw[thick] (.7,-.8)--(-.1,0)--(.5,.6)--(1.3,-.2)--(.7,-.8);
    \draw[thick] (-.6,.4)--(0,1) (-.6,-.4)--(.2,-1.2);
    \draw[dotted,thin,arrow] (-2.5,-1.5)--(.5,1.5); 
    \draw[dotted,thin,arrow] (.5,-1.5)--(-2.5,1.5);
    \node at (-.5,1) {$I$}; \node at (-.5,-1.1) {$J$}; 
    \node at (.6,-.1) {$O$};
}.
$$
Its DHR endomorphisms are direct sums of
$\rho_+\otimes\rho_-\in\DHR(\A_+)\otimes\DHR(\A_-)$. 
From the definition of the DHR braidings, and because $O_1=I_1\times
J_1$ is in the right spacelike complement of $O_2=I_2\times
J_2$ if and only if $I_1$ is in the future of $I_2$ and $J_1$ is in
the past of $J_2$, it follows that the braiding of $\A_2$ is given by
$\eps^+\otimes \eps^-$. Therefore, as a braided category,
$\DHR(\A_2)=\DHR(\A_+)\boxtimes \DHR(\A_-)\opp$.  

In particular, if the chiral nets $\A_+$ and $\A_-$ are isomorphic, 
then the canonical Q-system gives rise to a local two-dimensional
extension $\B_2$ of $\A_2=\A\otimes\A$, which is also 
known as the ``Cardy type'' extension. Its charged fields carry
conjugate charges w.r.t.\ the $+$ and $-$ chiral observables.
For the construction of this extension, it is actually not essential
that $\A_+$ and $\A_-$ are isomorphic, but it is sufficient that they
have isomorphic DHR categories. Obviously, one may as well construct a
Cardy type extension based on any pair of isomorphic subcategories of 
$\DHR(\A_+)$ and of $\DHR(\A_-)$. 

A more general class of local two-dimensional extensions of $\A_2$ was
constructed in \cite{R00}, by exhibiting the numerical coefficients of
the Q-system 
$(\Theta,W,X)$ in $\DHR(\A)\boxtimes\DHR(\A)\opp$ by a method involving
chiral $\alpha$-induction along a possibly noncommutative chiral
Q-system (the ``$\alpha$-induction construction''). The multiplicities
of the irreducible subsectors $\rho\otimes\ol\sig\prec\Theta$ coincide
with the matrix elements 
$Z_{\rho,\sig}=\dim\Hom(\alpha^-_\rho,\alpha^+_\sig)$ of the modular
invariants, mentioned before, cf.\ \sref{s:Qnets}. 

\subsubsection{Left and right centre}
\label{s:wedges}
In general, the extension $\A\subset \B$ described by a Q-system $\BA$
in $\DHR(\A)$ will be nonlocal. Since the left and right centres
$C^\pm[\BA]$ of $\BA$ are commutative Q-systems, cf.\ \sref{s:centre},
they correspond to
local extensions $\B^\pm\loc$ intermediate between $\A$ and $\B$. 

In \cite{BKL}, we have identified these local intermediate extensions
with relative commutants 
$$\B^+\loc(O):=\B(W_L)'\cap \B(W_R'), \quad\hbox{resp.}\quad 
\B^-\loc(O):=\B(W_R)'\cap \B(W_L').$$
Here, the wedges $W_L$ and $W_R$ are the left and right components of
the spacelike complement of the doublecone $O$ (resp.\ the past and future
complements of an interval in the chiral case):
$$\tikzmatht{
  \fill[black!10] (-2.5,-1.5) rectangle (2.5,1.5);
  \fill[white] (0,1)--(1,0)--(0,-1)--(-1,0)--(0,1);
    \draw[thick] (-2.5,1.5)--(0,-1)--(2.5,1.5); 
    \draw[thick] (-2.5,-1.5)--(0,1)--(2.5,-1.5);
    \node at (-1.8,0) {$W_L$}; \node at (1.8,0) {$W_R$}; 
    \node at (0,0) {$O$};
}, \qquad 
\tikzmatht{
  \fill[black!10] (-2.5,-1.5) rectangle (2.5,1.5);
    \fill[white] (0,1)--(1,0)--(0,-1)--(-1,0)--(0,1);
    \draw[thick] (-2.5,1.5)--(-1,0)--(-2.5,-1.5);
    \draw[thick] (-.5,1.5)--(1,0)--(-.5,-1.5);  
    \draw[thin] (0,1)--(-1,0)--(0,-1);
    \node at (-1.8,0) {$W_L$}; \node at (-.9,.9) {$W_R'$}; 
    \node at (0,0) {$O$};
}.
$$

In order to establish this result, one has to verify that the relative
commutant $\B^+\loc(O)=\B(W_L)'\cap \B(W_R')$ is intermediate between
$\A(O) \subset \B(O)$, and that the projection $p^+\loc$ corresponding to
the intermediate extension coincides with the right centre projection
$p^+$ of the Q-system for $\A(O)\subset \B(W_L)'\cap \B(W_R')\subset
\B(O)$. Thanks to \pref{p:maximal}, it is sufficient to prove that
$p^+\loc$ satisfies the relation \eref{centp}, and that $\B^+\loc(O)$
is maximal with this property. 
 
The inclusion $\A(O)\subset \B^+\loc(O)$ is obvious by isotony of $\B$ and
relative locality of $\B$ w.r.t.\ $\A$. The inclusion
$\B^+\loc(O)\subset \B(O)$ can be established with the help of Haag
duality for wedges, which was assumed to be valid for the net
$\A$. The intersection $\B^+\loc(O)$ is therefore the maximal
subalgebra of $\B(O)$ commuting with $\B(W_L)$. One has 
$\B^+\loc(O)=\A(O)v^+\loc=\A(O)p^+\loc v$. This algebra commutes with
$\B(W_L)$ iff the generator $p^+\loc v$ of $\B^+\loc(O)$  
commutes with the generator $\wh v=uv$ of $\B(W_L)$, where 
$u\in\Hom(\theta,\wh\theta)$ is a unitary charge transporter taking 
$\theta$ to $\wh\theta$ localized in $W_L$. Now, $p^+\loc v\cdot
uv=p^+\loc\theta(u)xv = \theta(u)p^+\loc xv$, whereas $uv\cdot
p^+\loc v=u\theta(p^+\loc)xv$. Commutativity is therefore equivalent to 
$p^+\loc x=\theta(u)^*u\theta(p^+\loc)x$, which is \eref{centp} by the
definition of the DHR braiding. The claim, $p^+\loc=p^+$, then follows by
the maximality of $\B^+\loc(O)$ and the characterization of $p^+$ in
\pref{p:maximal}. 

In establishing this result, it is again essential that the braiding
is the DHR braiding \eref{DHR-b}, defined in terms of unitary charge
transporters. This explains why a similar interpretation of the left and
right centre as the Q-system of some relative commutant cannot be
given in a general braided subcategory of $\End_0(N)$ for a single von
Neumann algebra $N$. It would be interesting to have such a theory,
which would require -- in addition to a braided tensor category
$\C\subset\End_0(N)$ -- as additional data a splitting of $N'$ into two
commuting subalgebras $N'=N_L\vee N_R$, and unitary intertwiners between
endomorphisms $\rho\in \End_0(N)$ and endomorphisms of $N_L$ and of
$N_R$, connected to the given braiding by (versions of) \eref{DHR-b}. 

\subsubsection{Braided product of extensions}
\label{s:b-ext}

According to \lref{l:braidext}, the braided products 
$\BA\times^\pm\BB$ of two Q-systems $\BA=(\theta^\BA,w^\BA,x^\BA)$ and
$\BB=(\theta^\BB,w^\BB,x^\BB)$ describe extensions $M^\pm$ which are
generated by the algebra $N$ and the generators $v^\BA$ and $v^\BB$
such that $Nv^\BA=M^\BA$ and $Nv^\BB=M^\BB$ are intermediate algebras,
and the generators $v^\BA$, $v^\BB$ satisfy the commutation relations 
\be\label{CRv}
v^\BB v^\BA=\iota(\eps^\pm_{\theta^\BA,\theta^\BB})\cdot v^\BB v^\BA.
\ee
These properties uniquely specify $M^\pm$. 

The same holds true in the QFT setting, that is, the braided product $\B^\pm$
of two extensions $\B^\BA$ and $\B^\BB$ of a net $\A$ is generated by
$\A$ and the generators $v^\BA$, $v^\BB$ subject to the same
commutation relation \eref{CRv}.

It follows that for any unitary charge transporters
$u_1\in\Hom(\theta^\BA,\wh\theta^\BA)$,
$u_2\in\Hom(\theta^\BB,\wh\theta^\BB)$, 
$$(u_2v^\BB)(u_1v^\BA)=\iota(\eps^\pm_{\wh\theta^\BA,\wh\theta^\BB})\cdot
(u_1v^\BA) (u_2v^\BB).$$
Now, if $\wh\theta^\BA$ is localized to the right (left) of  
$\wh\theta^\BB$, then $\eps^+_{\wh\theta^\BA,\wh\theta^\BB}=1$
($\eps^-_{\wh\theta^\BA,\wh\theta^\BB}=1$), hence in these cases the
generators $u_1v^\BA$ and $u_2v^\BB$ commute. Since the local algebras
of an extension are generated by $\A(O)$ and $uv$ such that
$\wh\theta=\Ad_u\circ\theta$ is localized in $O$, it follows (using
relative locality w.r.t.\ $\A$) that, as subalgebras of $\B^+$ resp.\ 
of $\B^-$, the algebras $\B^\BA(O_1)$ and $\B^\BB(O_2)$ commute with
each other if $O_1$ is located to the spacelike right resp.\ left of
$O_2$ -- but in general not in the converse order.

We paraphrase these situations by saying that the net $\B^\BA$ is ``right
local'' resp.\ ``left local'' w.r.t.\ the net $\B^\BB$. 
Thus, the braided products of of two extensions $\B^\BA$ and $\B^\BB$
can be regarded as a quotient of the free product by the
relations that identify the common subnets $\A\subset\B^\BA$ and 
$\A\subset\B^\BB$, and by the relations expressing that $\B^\BA$ is ``right
local'' resp.\ ``left local'' w.r.t.\ $\B^\BB$. 

(The same is true in the chiral case, replacing ``right'' by ``future'' and
``left'' by ``past''.)

Let $\A$ be a chiral QFT and $\BA$ a Q-system in $\DHR(\A)$,
describing a (nonlocal) chiral extension $\B$. The full centre 
$Z[\BA]$ is a Q-system in $\DHR(\A\otimes\A)$, i.e., it describes 
a two-dimensional extension $\B_2$ of $\A\otimes\A$. Because the full
centre is the right centre of the braided product $(\BA\otimes
1)\times^+ \BR$, we recognize the corresponding extension $\B_2$ as
the relative commutant of right wedges of the nonlocal extension,
obtained by the right-local braided product of the possibly nonlocal
chiral extension $\B\otimes 1$ with the local canonical
extension $\B_2^\BR$.

By \pref{p:fullc=alpha}, the full centre coincides with the 
$\alpha$-induction construction which was originally found as a
construction of two-dimensional local conformal QFT models out of
chiral data. This result therefore not only gives a more satisfactory,
purely algebraic interpretation of the $\alpha$-induction construction
in terms of braided products of nets and relative commutants of wedge
algebras, cf.\ \sref{s:wedges}; it also explains the fact (known
before) that the latter depends only on the Morita equivalence class
of the chiral Q-system in $\C$ \cite{LR04}; namely two Q-systems in a modular
tensor category $\C$ have the same full centre if and only if they are
Morita equivalent \cite{KR08}.  

Since moreover, every irreducible extension $\B_2$ of $\A_2$ is
intermediate between $\A_2$ and an $\alpha$-induction extension
\cite{LR04,LR09,CKL13}, it follows that full centre extensions are
precisely the maximal irreducible extensions (if the underlying chiral
theory $\A$ is completely rational).

\subsection{Hard boundaries}
\label{s:hardb} 
\setcounter{equation}{0}
A conformal quantum field theory with a ``hard boundary'' arises, when
Minkowski spacetime $M^2$ is restricted to a halfspace, say the right
halfspace $M^2_R=\{(t,x): x>0\}$. The stress-energy tensor defined on
$M^2_R$ still splits into two chiral components, but if one imposes
conservation of energy at the boundary, the two components are no
longer independent fields, but instead they {\em coincide} as
operator-valued distributions on $\RR$ \cite{LR04}. Thus, for $I$ and $J$
intervals such that $O=I\times J=\{(t,x):t+x\in I,\, t-x\in J\}$ lies
inside $M^2_R$ ($\LRA$ $I>J$ elementwise as subsets of $\RR$), the
local algebra of chiral observables is   
$$\A_R(O)=\A(I)\vee \A(J)$$
rather than the tensor product $\A_+(I)\otimes\A_-(J)$. Here, the
chiral algebras are generated by the stress-energy tensor and possibly
further chiral fields whose boundary conditions might also impose an
identification of the fields. 

An analysis of local extensions $\A_R(O)\subset \B_R(O)$ on the
halfspace was given in \cite{LR04}. One finds (assuming $\A$ to be
completely rational) that the local algebras of every maximal such
extension are of the form  
\be\label{holo}
\B_R(O)=\B(K)'\cap \B(L)\qquad (O=I\times J\subset M^2_R)
\ee
where $\A\subset\B$ is a possibly nonlocal chiral extension
(given by a Q-system in $\DHR(\A)$), and $K\subset L$ is the unique
pair of open intervals such that $I\cup J = L\setminus \ol K$:  
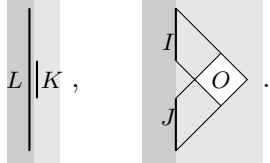
$$
\tikzmatht{
  \fill[black!10] (-.7,-2.2) rectangle (.7,2.2);
  \fill[black!20] (-.7,-2.2) rectangle (0,2.2);
    \draw[thick] (-.1,-1.9)--(-.1,1.9) (.1,-.5)--(.1,.5);
    \node at (-.5,0) {$L$}; \node at (.5,0) {$K$};
},\qquad 
\tikzmatht{
  \fill[black!10] (-.9,-2.2) rectangle (2.3,2.2);
  \fill[black!20] (-.9,-2.2) rectangle (0,2.2);
  \fill[white] (.5,0)--(1.2,.7)--(1.9,0)--(1.2,-.7)--(.5,0);
    \draw[thick] (0,-1.9)--(0,-.5) (0,1.9)--(0,.5);
    \draw[thin] (0,-1.9)--(1.9,0)--(0,1.9);
    \draw[thin] (0,-.5)--(1.2,.7) (0,.5)--(1.2,-.7);
    \node at (-.2,1) {$I$}; \node at (-.2,-1) {$J$}; 
    \node at (1.2,0) {$O$};
}.
$$
This formula
is ``holographic'' in the sense that the local observables in a region
$O\subset M^2_R$ are given in terms of operators in a chiral net that
can be thought of as a net on the boundary. 

The simplest case is the trivial chiral extension $\B=\A$. In this
case, $\B_R(O)$ is generated by $\A_R(O)= \A(I)\vee \A(J)$ and charge
transporters in $\A(L)$, namely unitary intertwiners transporting a DHR
endomorphisms localized in $J$ to an equivalent DHR endomorphism
localized in $J$. Accordingly, the charged generators for the
subfactor $\A_R(O)\subset\B_R(O)$ ``carry a charge $\rho$ in $I$ and a
charge $\ol\rho$ in $J$''. Indeed, under the split isomorphism between
the von Neumann algebras $\A(I)\vee \A(J)$ and $\A(I)\otimes \A(J)$,
the subfactor turns out to be isomorphic with the subfactor associated
with the canonical Q-system \pref{p:canonical} with $[\Theta_\BR] = 
\bigoplus_{[\rho]\,\mathrm{irr}} \rho\otimes\ol\rho$.    

For general chiral extensions $\A\subset\B$ with irreducible Q-system
$\BA$, the local subfactor $\A_R(O)\subset \B_R(O)$ for any bounded doublecone
$O\subset M^2_R$ not touching the boundary is isomorphic to the subfactor
obtained from the full centre Q-system $Z[\BA]$, and hence 
depends (up to isomorphism) only on the Morita equivalence class 
(cf.\ \sref{s:indQ}) of the chiral Q-system $\BA$.  

As mentioned in \sref{s:b-ext}, the full centre gives also a local net
on the full two-dimensional Minkowski spacetime as an extension
$\B_2\supset\A\otimes\A$ of the tensor product of a 
pair of isomorphic chiral nets. Indeed, this net can be recovered from
the maximal boundary net $\B_R$ by a procedure called 
``removing the boundary''. It proceeds by taking the limit of a
sequence of states on right wedge algebras $\B_R(W_R+a)$ as $a\in W_R$
tends to infinity (``far away from the boundary''). The net $\B_2$ can then
be defined in the GNS Hilbert space of this state, which carries two
commuting unitary representations of the Möbius group. First defining
the local algebra $\B_2(W_R)$ of a single right wedge, and
$\B_2(W_R'):=\B_s(W_R)'$ as its commutant, the two unitary representations
of the translations
are used to define the local algebras for general wedge regions, and
the local algebras for doublecones by intersections of algebra for wedges. 

The converse procedure of ``adding a boundary'' can also be performed
algebraically \cite{CKL13}. Starting from an extension
$\A_2\subset\B_2$ defined on Minkowski spacetime, one can redefine the
representation of its restriction to $M^2_R$, obtaining a reducible
representation. Its decomposition yields a direct sum of boundary extensions
$\A_R\subset\B_R$ related to chiral extensions $\A\subset\B$ by the
``holographic formula'' \eref{holo}, which all give back 
$\A_2\subset\B_2$ when the ``boundary is removed''. In particular, 
for every right wedge $W_R\subset M^2_R$ not touching the boundary, 
the subnets $\A_R(O)\subset\B_R(O)$ indexed by $O\subset W_R$ are all 
isomorphic to the subnets $\A_2(O)\subset\B_2(O)$, so that the
boundary nets in the decomposition can be interpreted as different
boundary conditions imposed on the original net $\A_2\subset\B_2$.  

The procedure of ``adding a boundary'' amounts, in the language of
C*-tensor categories, to the tensor functor $T:\C\boxtimes\C\opp$,
$\rho\otimes\sig\mapsto\rho\sig$, taking Q-systems in
$\C\boxtimes\C\opp$ to Q-systems in $\C$. This functor 
is adjoint \cite{KR09} to the full centre, taking Q-systems in $\C$ to
Q-systems in $\C\boxtimes\C\opp$. In is proven in \cite{KR08} that the
image of the full centre Q-system $Z[\BA]$ under $T$ is the direct sum 
(in the sense of \sref{s:Qcentral}) of Q-systems given by the
irreducible $\BA$-modules. Thus, the hard boundary conditions are
classified in 1:1 correspondence with the irreducible modules of the
underlying chiral Q-system $\BA$.

\subsection{Transparent boundaries}
\label{s:transb} 
\setcounter{equation}{0}
Whereas a hard boundary describing a QFT on a halfspace identifies
the left- and right-moving chiral observables in the halfspace, a
transparent boundary separates two possibly different quantum field
theories $\B^L$ and $\B^R$ in the halfspaces $M^2_L$, $M^2_R$ on
either side of the boundary: 
$$\tikzmatht{
  \fill[black!10] (-2.5,-1.7) rectangle (2.5,1);
  \fill[black!15] (-2.5,-1.7) rectangle (0,1);
  \fill[white] (-1.2,.8)--(-.4,0)--(-1.2,-.8)--(-2,0)--(-1.2,.8);
    \draw[thick] (-1.2,.8)--(-.4,0)--(-1.2,-.8)--(-2,0)--(-1.2,.8);
  \fill[white] (1.2,.8)--(.4,0)--(1.2,-.8)--(2,0)--(1.2,.8);
    \draw[thick] (1.2,.8)--(.4,0)--(1.2,-.8)--(2,0)--(1.2,.8);
    \draw[thin,dotted] (0,-1.7)--(0,1);
    \node at (-1.4,-1.2) {$\scriptstyle{\B^L(O_1)}$}; 
    \node at (1.4,-1.2) {$\scriptstyle{\B^R(O_2)}$}; 
} 
\qquad (O_1\subset M^2_L, \, O_2\subset M^2_R).
$$
Physically speaking, the boundary is thought to separate regions with
different dynamics, e.g., two different phases of a relativistic  
system with a phase transition. For the example of the Ising model,
cf.\ \cite{ST78} and \xref{x:ising-ph}.    

The two theories are defined on the same Hilbert space, and share a
tensor product $\A_+\otimes\A_-$ of common chiral subtheories. The latter
property arises from the physical assumption that energy {\em and}
momentum are conserved at the boundary, which identifies the chiral
stress-energy tensors on either side of the boundary \cite{BKLR}. 

Because the presence of the boundary cannot violate the principle of
causality, quantum observables of $\B^L$ localized in the left
halfspace $M^2_L$ must commute with observables of $\B^R$ localized in
the right halfspace $M^2_R$ at spacelike separation. 

Because the stress-energy tensor is the local generator of
diffeomorphisms, the common chiral subtheory $\A_+\otimes\A_-$ can be
used to extend both theories to the full Minkowski spacetime. 

Motivated by these two (heuristic) observations, one should define a
transparent boundary as a pair of quantum field theories on
two-dimensional Minkowski spacetime, sharing a common chiral
subtheory, such that $\B^L$ is left-local w.r.t.\ $\B^R$. As we have
seen before, such a pair is described by the braided product of two
extensions of the common chiral subtheory, and every irreducible such
pair is a quotient of the braided product. 

The mathematical issue is therefore the central decomposition of the
braided product of a pair of commutative Q-systems in
$\DHR(\A_+)\boxtimes\DHR(\A_-)$. The centre of the braided product of
extensions is given by \pref{p:commprod} as a linear space isomorphic
to $\Hom(\Theta^L,\Theta^R)$. In order to know its central
projections, it must be computed as an algebra. This 
is precisely what we have achieved in \tref{t:centchar}, provided
$\DHR(\A_+)$ and $\DHR(\A_-)$ are isomorphic as modular braided
categories, and the pair of commutative Q-systems are full centres of
chiral Q-systems $\BA^L$ and $\BA^R$. Namely, \tref{t:centchar}
classifies the transparent boundary conditions in 1:1 correspondence
with the irreducible chiral $\BA^L$-$\BA^R$-bimodules.  

In \cite{BKLR}, this classification is further elaborated. 
As discussed in \sref{s:class}, the space $\Hom(\Theta^L,\Theta^R)$ has two
distinguished bases, orthogonal w.r.t.\ the inner product
\eref{inner}: one arising by ``diagonalizing'' the left and right
compositions with $\Hom(\Theta^L,\Theta^L)$ and
$\Hom(\Theta^R,\Theta^R)$, the other corresponding to the minimal
central projections of the braided product, i.e., the minimal
projections in $\Hom(\Theta^L,\Theta^R)$ w.r.t.\ the convolution
product \eref{convol}. The unitary transition matrix is a generalized
Verlinde matrix, and can be computed by its distinguishing property
that it ``diagonalizes'' the bimodule fusion rules. Its matrix
elements finally turn out to determine the specific identifications
between charged fields of $\B^L$ and charged fields of $\B^R$, that
make up the specific boundary conditions. 

\begin{graybox}
\begin{example}
\label{x:ising-ph}
The special case where both Q-systems are the canonical one, i.e., the
boundary between two conformal quantum field theories both isomorphic
to the Cardy extension, has been given in \xref{x:cardy}.
For the Ising model (i.e., the chiral net is given by the
Virasoro net with central charge $c=1$), one obtains three boundary
conditions given by the three sets of linear dependencies between the
charged generators $\Psis,\Psit$:
\bea (\mathrm{i}) & \Psit^L=\Psit^R,\quad \Psis^L=\Psis^R; \nonumber \\
(\mathrm{ii}) & \Psit^L=\Psit^R,\quad \Psis^L=-\Psis^R; \nonumber \\
(\mathrm{iii}) & \Psit^L=-\Psit^R. \nonumber
\eea
The first case is the trivial boundary; the second the ``fermionic''
boundary where the field $\Psis$ changes sign, and the third the
``dual'' boundary, in which there are two independent fields 
$\Psis^R$ and $\Psis^L$ (corresponding to the order and 
disorder parameter $\sig$ and $\mu$ in \cite{ST78}). 
\end{example}
\end{graybox}

\subsection{Further directions} 
\setcounter{equation}{0}
We have outlined the remarkably tight links between the abstract
theory of Q-systems in braided C* tensor categories and the
representation theory of conformal quantum field theories in two
dimensions. Notably the classifications of ``hard'' and
``transparent'' boundary conditions have very natural counterparts in
the abstract setting. 

Thinking of systems with several transparent boundaries, some
immediate questions arise: the juxtaposition of two boundaries is
described by the (associative) braided product of three Q-systems. The
individual boundary conditions are classified as $\BA$-$\BB$-bimodules
and as $\BB$-$\BC$-bimodules. Thus, it is expected that the
juxtaposition of boundary conditions is described in terms of the
bimodule tensor product. 

It is much less clear which mathematical structure should be expected to
describe situations where two transparent boundaries {\em intersect}
each other.  

Finally, hard and transparent boundaries are only two ``opposite
extremes'' in a wide spectrum of possible behaviour of chiral fields
at a boundary \cite{BKLR}. It would be rewarding to describe also more
general boundaries in terms of the present unifying framework.

\section{Conclusions}

Q-systems are a tool to describe extensions $N\subset M$ of an 
infinite von Neumann factor $N$ in terms of ``data'' referring only to
$N$. We have extended this notion, well-known for subfactors, to the
case when $M$ is admitted to be a finite direct sum of factors. 
Modules and bimodules of Q-systems are equivalent to homomorphisms 
between extensions. Decompositions of Q-systems and other operations 
defined in {\em braided} C* tensor categories: the centres, braided 
products and the full centre -- which are known in the setting of
abstract tensor categories -- are interpreted in terms of the associated
extensions of von Neumann algebras. 

The meaning of these operations in the context of local quantum field
theory is elaborated in \cite{BKLR}. Especially the determination of
the centre of the von Neumann algebra which arises as the braided
product of two commutative extensions, is a problem motivated by these
applications. We have completely solved this task for the braided
product of two full centres in {\em modular} C* tensor categories.  

In the last section, we have given a brief outline of this and other
applications of the theory of braided and modular C* tensor categories
in the context of quantum field theory. It is here, where the interpretation 
in terms of endomorphisms of von Neumann algebras is most substantial,
since local quantum observables are (selfadjoint) elements of von Neumann
algebras. This application was not only our original motivation for
the analysis presented in the main body of this work; it is also not
an exaggeration to say that the (rather natural) appearance of {\em modular}
C* categories in chiral conformal QFT, as an offspring of the original
DHR theory designed for massive QFT in four spacetime dimensions, has
triggered many of the developments described in this work.

\bigskip

{\large\bf Acknowledgments} 

\medskip

We are very much indebted to J. Fuchs, I. Runkel, and C. Schweigert 
for their hospitality and enlightening explanations of their work,
which were most beneficial for the results presented in \sref{s:class}. 
Y.K. thanks M. Izumi for an interesting question.


\begin{thebibliography}{99} \itemsep-.8mm 
\bibitem{BK} B. Bakalov, A. Kirillov Jr., {\it Lectures on
  tensor categories and modular functors}, University Lecture
  Series 21, A.M.S., Providence, RI (2001). 

\bibitem{BDH} A. Bartels, C.L. Douglas, A. Henriques, {\it
  Dualizability and index of subfactors}, arXiv:1110.5671.

\bibitem{B} D. Bisch, {\it A note on intermediate subfactors},
  Pac.\ J. Math.\ 163, 201--216 (1994).

\bibitem{BKL} M. Bischoff, Y. Kawahigashi, R. Longo, {\it
  Characterization of 2D rational local conformal nets and its
  boundary conditions: the maximal case}, arXiv:1410.8848.  

\bibitem{BKLR} M. Bischoff, Y. Kawahigashi, R. Longo, K.-H. Rehren, 
  {\it Phase boundaries in algebraic conformal QFT}, arXiv:1405.7863.  

\bibitem{BE} J. B\"ockenhauer, D. Evans, {\it Modular invariants,
  graphs and $\alpha$-induction for nets of subfactors, I}, Commun.\
  Math.\ Phys.\ 197, 361--186 (1998), and {II}, Commun.\
  Math.\ Phys.\ 200, 57--103 (1999), and {III} Commun.\
  Math.\ Phys.\ 205, 183--228 (1999). 
\bibitem{BEK99} J. B\"ockenhauer, D. Evans, Y. Kawahigashi, {\it On
  $\alpha$-induction, chiral projectors and modular invariants for
  subfactors}, Commun.\ Math.\ Phys.\ 208, 429--487 (1999). 

\bibitem{BEK00} J. B\"ockenhauer, D. Evans, Y. Kawahigashi, 
  {\it Chiral structure of modular invariants for subfactors}, Commun.\
  Math.\ Phys.\ 210, 733--784 (2000). 

\bibitem{CIZ} A. Cappelli, C. Itzykson, J.-B. Zuber, {\it The $A$-$D$-$E$
  classification of minimal and $A^{(1)}_1$ conformal invariant
  theories}, Commun.\ Math.\ Phys.\ 113, 1--26 (1987).
\bibitem{CKL08} S. Carpi, Y. Kawahigashi, R. Longo, {\it Structure and
  classification of superconformal nets}, Ann.\ H. Poincar\'e 9,
  1069--1121 (2008). 
\bibitem{CKL13} S. Carpi, Y. Kawahigashi, R. Longo, {\it How to
  add a boundary condition}, Commun.\ Math.\ Phys.\ 322, 149--166 (2013).
\bibitem{DMNO} A. Davydov, M. M\"uger, D. Nikshych, V. Ostrik, {\it The
  Witt group of non-degenerate braided fusion categories}, arXiv:1009.2117.

\bibitem{D91} P. Deligne, {\it Cat\'egories tannakiennes}, in:
  P. Cartier et al.\ (eds.), Grothendieck Festschrift, vol.\ II,
  Birkh\"auser, Basel, 1991, pp. 111--195. 

\bibitem{D02} P. Deligne, {\it Cat\'egories tensorielles},
  Mosc.\ Math.\ J. 2, 227--248 (2002).

\bibitem{DHR} S. Doplicher, R. Haag, J.E. Roberts, {\it Local
  observables and particle statistics. I}, Commun.\ Math.\
  Phys.\ 23, 199--230 (1971).
\bibitem{ENO} P. Etingof, D. Nikshych, V. Ostrik, {\it On fusion
  categories}, Ann.\ Math.\ 162, 581--642 (2005). 

\bibitem{EP03} D. Evans, P. Pinto, {\it Subfactor realizations
  of modular invariants}, Commun.\ Math.\ Phys.\ 237, 309--363 (2003).
\bibitem{FI} F. Fidaleo, T. Isola, {\it Minimal expectations for
  inclusions with atomic centre}, Int.\ J. Math.\ Phys.\ 7, 307--327 (1996). 

\bibitem{FRS} K. Fredenhagen, K.-H. Rehren, B. Schroer, {\it Superselection 
  sectors with braid group statistics and exchange
  algebras I}, Commun.\ Math.\ Phys.\ 125, 201--226 (1989). 
\bibitem{FFRS04} J. Fr\"ohlich, J. Fuchs, I. Runkel,
  C. Schweigert, {\it Kramers-Wannier duality from conformal defects},
  Phys.\ Rev.\ Lett.\ 93, 070601 (2004).

\bibitem{FFRS06} J. Fr\"ohlich, J. Fuchs, I. Runkel, C. Schweigert,
  {\it Correspondences of ribbon categories}, Ann.\ Math.\ 199
  192--329, (2006).  
\bibitem{FFRS07} J. Fr\"ohlich, J. Fuchs, I. Runkel, C. Schweigert,
  {\it Duality and defects in rational conformal field theory}, Nucl.\
  Phys.\ B 763, 354--430 (2007). 
\bibitem{TFT1} J. Fuchs, I. Runkel, C. Schweigert, {\it TFT construction 
  of RCFT correlators I: partition functions}, Nucl.\ Phys.\ B 646
  [PM], 353--497 (2002).  
\bibitem{TFT} J. Fuchs, I. Runkel, C. Schweigert, {\it TFT construction of 
  RCFT correlators, II--IV}, Nucl.\ Phys.\ B 678 [PM] 511--637 (2004), 
  Nucl.\ Phys.\ B 694 [PM] 277--353 (2004), Nucl.\ Phys.\ B 715 [PM],
  539--638 (2005). 
\bibitem{GL96} D. Guido, R. Longo, {\it The conformal spin and
  statistics theorem},  Commun.\ Math.\ Phys.\ 181, 11--35 (1996).
\bibitem{H} R. Haag: Local Quantum Physics, Springer Verlag,
  Berlin -- Heidelberg -- New York, 1996.

\bibitem{I} M. Izumi, {\it The structure of sectors associated with
  Longo-Rehren inclusions. I. General theory}, Commun.\ Math.\
  Phys.\ 213, 127--179 (2000). 
\bibitem{IK} M. Izumi, H. Kosaki, {\it On a subfactor analogue of the second
  cohomology}, Rev.\ Math.\ Phys.\ 14, 733--737 (2002).
\bibitem{J} V.F.R. Jones, {\it Index for subfactors}, Invent.\
  Math.\ 72, 1--25 (1983).
\bibitem{J98} V.F.R. Jones, {\it The planar algebra of a bipartite
    graph}, in: Knots in Hellas '98, Ser.\ Knots Everything 24, World
  Scientific, River Edge, NJ, 2000; pp.\ 94--117. 
\bibitem{JS} A. Joyal, R. Street, {\it Braided tensor categories},
  Adv.\ Math.\ 102, 20--78 (1993). 
\bibitem{KL04} Y. Kawahigashi, R. Longo, {\it Classification of local
  conformal nets. Case $c<1$}, Ann.\ Math.\ 160, 493--522 (2004). 
\bibitem{KL04-2} Y. Kawahigashi, R. Longo, {\it Classification of
  two-dimensional local conformal nets with $c<1$ and 2-cohomology
  vanishing for tensor categories}, Commun.\ Math.\ Phys.\ 244,
  63--97 (2004). 
\bibitem{KLM} Y. Kawahigashi, R. Longo, M. M\"uger, {\it
  Multi-interval subfactors and modularity of representations in
  conformal field theory}, Commun.\ Math.\ Phys.\ 219, 631--669
(2001).
\bibitem{KLPR} Y. Kawahigashi, R. Longo, U. Pennig, K.-H. Rehren, {\it
  The classification of non-local chiral CFT with $c<1$}, Commun.\
  Math.\ Phys.\ 271, 375--385 (2007).  
\bibitem{KO} A. Kirillov Jr., V. Ostrik, {\it On $q$-analog of McKay
    correspondence and ADE classification of $sl(2)$ conformal field
    theories}, Adv.\ Math.\ 171, 183--227 (2002). 
\bibitem{Ko} H. Kosaki, {\it Extension of Jones' theory on index
  to arbitrary factors}, J. Funct.\ An.\ 66, 123--140 (1986).
\bibitem{KL} H. Kosaki, R. Longo, {\it A remark on the minimal
  index of subfactors}, J. Funct.\ An.\ 107, 458--470 (1992).
\bibitem{KR08} L. Kong, I. Runkel, {\it Morita classes of algebras
  in modular tensor categories}, Adv.\ Math.\ 219, 1548--1576 (2008).
\bibitem{KR09} L. Kong, I. Runkel, {\it Cardy algebras and sewing
  constraints, I}, Commun.\ Math.\ Phys.\ 292, 871--912 (2009). 
\bibitem{KR10} L. Kong, I. Runkel, {\it Algebraic structures in
  Euclidean and Minkowskian two-dimensional conformal field theory}, 
  arXiv:0902.3829. 
\bibitem{L89} R. Longo, {\it Index of subfactors and statistics 
  of quantum fields I}, Commun.\ Math.\ Phys.\ 126, 217--247 (1989).
\bibitem{L94} R. Longo, {\it A duality for Hopf algebras and for
  subfactors}, Commun.\ Math.\ Phys.\ 159, 133--150 (1994).
\bibitem{LR95} R. Longo, K.-H. Rehren, {\it Nets of subfactors}, 
  Rev.\ Math.\ Phys.\ 7, 567--597 (1995).
\bibitem{LR04} R. Longo, K.-H. Rehren, {\it Local fields in
  boundary CFT}, Rev.\ Math.\ Phys.\ 16, 909--960 (2004). 
\bibitem{LR09} R. Longo, K.-H. Rehren, {\it How to remove the boundary 
  in CFT -- an operator algebraic procedure}, Commun.\ Math.\ Phys.\
  285, 1165--1182 (2009). 
\bibitem{LRo} R. Longo, J.E. Roberts, {\it A theory of dimension},
  K-Theory 11, 103--159 (1997), notably Chap.\ 3 and 4. 
\bibitem{LX} R. Longo, F. Xu, {\it Topological sectors and a
  dichotomy in conformal field theory}, Commun.\ Math.\ Phys.\ 251,
  321--364 (2004). 
\bibitem{ML} S. Mac Lane, Categories for the Working
  Mathematician, second ed., Springer-Verlag, Berlin, 1998. 
\bibitem{M00} M. M\"uger, {\it Galois theory for braided tensor categories
  and the modular closure}, Adv.\ Math.\ 150, 151--201 (2000). 
\bibitem{M03} M. M\"uger, {\it From subfactors to categories and
  topology I. Frobenius algebras in and Morita equivalence of tensor
  categories}, J. Pure Appl.\ Algebra 180, 81--157 (2003), and
  {\it II. The quantum double of tensor categories and subfactors},
  J. Pure Appl.\ Algebra 180, 159--219 (2003).
\bibitem{M04} M. M\"uger, {\it Galois extensions of braided tensor categories
  and braided crossed G-categories}, J. Algebra 277, 256--281 (2004).
\bibitem{O} V. Ostrik, {\it Module categories, weak Hopf algebras and
    modular invariants}, Transform.\ Groups 8, 177--206 (2003).
\bibitem{P} S. Popa, {\it Classification of subfactors and their
  endomorphisms}, CBMS Regional Conference Series in Mathematics,
  86, American Mathematical Society, Providence, RI, 1995.
\bibitem{PS} A. Pressley, G. Segal, Loop Groups, Oxford
  University Press, Oxford, 1986.  
\bibitem{R89} K.-H. Rehren, {\it Braid group statistics and their
    superselection rules}, In: The Algebraic Theory of Superselection
  Sectors, D. Kastler (ed.), World Scientific 1990, pp.\ 333--355.
\bibitem{R00} K.-H. Rehren, {\it Canonical tensor product subfactors}, 
  Commun.\ Math.\ Phys.\ 211, 395--406 (2000).
\bibitem{R01} K.-H. Rehren, {\it Locality and modular invariance
  in 2D conformal QFT}, Fields Inst.\ Commun.\ 30, 341--354 (2001),
  arXiv:math-ph/0009004. 
\bibitem{ST78} B. Schroer, T.T. Truong, {\it The order/disorder
  quantum field operators associated with the two-dimensional Ising
  model in the continuum limit}, Nucl.\ Phys.\ B 144, 80--122 (1978).
\bibitem{X} F. Xu, {\it Mirror extensions of local nets}, Commun.\
  Math.\ Phys.\ 270, 835--847 (2007).
\bibitem{Y03} S. Yamagami, {\it C*-tensor categories and free
  product bimodules}, J. Funct.\ An.\ 197, 323--346 (2003). 


\end{thebibliography}
\end{document}